%% file: sphwhitarxiv.tex
\documentclass[12pt]{amsart}
\usepackage{amsfonts,amssymb,amscd,amsmath}
\usepackage{latexsym,amsbsy,bbold,stmaryrd}
\usepackage{tocvsec2}
\usepackage{mathrsfs}
\usepackage[all]{xy}

\usepackage[
hyperindex=true,pagebackref=true,bookmarks=true,
colorlinks=true,linkcolor=blue,citecolor=red]
{hyperref}
\def\~{{\rm --}} 

\settocdepth{subsubsection}

\input{macro}

\begin{document}
\newcommand{\comment}[1]{}
\comment{
Article-id: 0904.4324, Article password: v538r 
Spherical and Whittaker Functions via DAHA
Ivan Cherednik, Xiaoguang Ma

This work grew out of the  lectures given by the first author at 
Harvard in February and March, 2009. 
A draft of the lecture notes was prepared by the second author, 
and then expanded and made into their final form by the first author.
It begins with an introduction to the classical p-adic theory of
the Macdonald, Matsumoto and Whittaker functions. 
Its major directions are as follows:

1) extending the theory of DAHA to arbitrary levels;
2) the affine Satake map and Hall functions via DAHA;
3) the spinor Dunkl operators for the Q-Toda operators;
4) applications to nil-DAHA and Q-Whittaker functions; 
5) the technique of spinors in the differential theory. 

The latest variant is essentially mathematically equivalent
to variant 5, though there the were some improvements, extensive
editing and adding some justifications.

 
}


\title[Spherical and Whittaker functions via DAHA]
{\footnotesize A new take on spherical, Whittaker and Bessel functions\\
{\tiny (Spherical and Whittaker functions via DAHA I,II)}}


\author[Ivan Cherednik]{Ivan Cherednik $^\dag$}
\author[Xiaoguang Ma]{Xiaoguang Ma}

\begin{abstract}
This paper begins with an exposition of the classical p-adic theory of
the Macdonald, Matsumoto and Whittaker functions aimed at the
affine generalizations. The major directions are
the theory of DAHA for arbitrary levels and the affine 
Satake map and Hall functions via DAHA. The key result is the
proportionality of the two different formulas for the affine
symmetrizer, the Satake-type formula and that based on the
polynomial representation of DAHA. The latter approach results 
in two important formulas for the affine symmetrizer
generalizing the relations between the Kac-Moody characters and
Demazure characters. 

The second part of this paper is focused on the spinor 
(nonsymmetric) Whittaker functions in the rank one, related 
q-Toda-Dunkl operators, and other aspects of the spinor construction, 
including one-dimensional Bessel functions, 
and the isomorphism between the affine Knizhnik-Zamolodchikov 
equation and the Quantum Many-Body problem (the Heckman-Opdam 
system).
\end{abstract}

\date{October 24, 2012}

\address 
{\vskip -0.5cm (I.Cherednik) Department of Mathematics, UNC
Chapel Hill, North Carolina 27599, USA\\
chered@email.unc.edu}
\address
{\vskip -0.7cm (X. Ma) 632 CAB,
University of Alberta,
Edmonton, Alberta  T6G 2G1, CANADA\\
xiaoma@math.ualberta.ca}

\def\bysame{{\bf --- }}
\def\~{{\bf --}}

\def\~{-}
\settocdepth{subsubsection}

\renewcommand{\tilde}{\widetilde}
\renewcommand{\hat}{\widehat}
\newcommand{\dagx}{\hbox{\tiny\mathversion{bold}$\dag$}}
\newcommand{\ddagx}{\hbox{\tiny\mathversion{bold}$\ddag$}}

\thanks{$^\dag$ 
Partially supported by NSF grant DMS--0800642}

\maketitle



\renewcommand{\baselinestretch}{1.2}
{
{
\tableofcontents}
} 
\renewcommand{\baselinestretch}{1.0}
\vfill\eject

\renewcommand{\natural}{\wr}

\setcounter{section}{-1}
\setcounter{equation}{0}
\section{\sc{Introduction}}

\noindent
This work grew out of the 
lectures given by the first author at Harvard in 
February and March, 2009. 
A draft of the lecture notes was prepared by the second author, 
and then expanded and made into their final form by the first author.
It will be published by Selecta Mathematica in two parts 
``Spherical and Whittaker functions via DAHA, I,II," 
essentially corresponding to Sections 0,1,2,3 and Sections 4,5,6,7
of this preprint, which is a somewhat extended version of its
previous variant (posted in 2009). These two parts 
are related but relatively independent; the theory of spherical
functions is the main unifying theme.


\subsection{{\bf Objectives and main results}}
The first aim of the first part of this work is to connect DAHA with 
the theory of affine Hall functions using the approach to the 
classical Hall polynomials ($=\, p$\~adic spherical functions) 
via the
Matsumoto $p$\~adic functions, an important 
special case of the theory of nonsymmetric Macdonald polynomials. 

It is closely connected to the second major direction of this 
work, which is the nonsymmetric Whittaker theory.  
Classical Whittaker functions are already nonsymmetric, so we need a 
new theory of spinors (generally, $W$\~spinors) to achieve this; some 
instances already appeared in the related harmonic analysis. 
\smallskip

Dunkl-$q$-Toda operators and their eigenfunctions,
the spinor $q$\~Whit\-taker functions, are introduced and studied 
for $A_1$ in the second part of this paper.  
The $p$\~adic limits of these function are well defined and 
result in {\em new} Matsumoto-type (``nonsymmetric") 
$p$\~adic Whittaker functions; see Sections \ref{sect:Whitf}.
The definition can be given for any root systems.

Another (actually related) possible output of this project
could be the theory of nonsymmetric counterparts of 
the affine Hall functions and the corresponding Satake map,
including their connections with the DAHA 
elliptic-type representations from \cite{C5}; 
cf. Section \ref{sec: p-adic-e}. This work is in progress. 
\smallskip

More specifically, the results of this work (both parts)
can be grouped as follows.
 
(1) The theory of DAHA modules of arbitrary levels $l$
(not only $l=0,1$ as in \cite{C101}),
which technically means that its polynomial representation
can be multiplied by any powers of the Gaussian.

(2) The affine Satake isomorphism and affine Hall functions 
via DAHA; the latter functions attract growing attention, 
though not much is known so far for arbitrary $q$ and $t$, the
DAHA parameters.

(3) Establishing connections with the theory of
Kac-Moody characters, the $t\to\infty$ 
limits of the affine Hall functions, and the level
one Demazure characters.

(4) The theory of coinvariants of DAHA, their 
relations to the bilinear symmetric invariant
forms on DAHA of higher levels, and the 
corresponding spaces of Looijenga functions.

(5) Revisiting classical $p$\~adic theory
of the Satake-Macdonald, Matsumoto and Whittaker functions
with the focus on the Matsumoto functions and aiming at
the DAHA generalizations. 

(6) The study of new spinor Dunkl operators serving the
$q$\~Toda operators and the $q$\~Whittaker functions,
the related theory of the nil-DAHA and the spinor
Whittaker functions.

(7) Developing the technique of $W$\~spinors, including the
differential theory and its application to the Bessel functions,
symmetric and nonsymmetric, and the AKZ$\leftrightarrow$QMBP 
isomorphism theorem.

\subsubsection{\sf Affine Satake isomorphisms}
Among the main topics we consider, are the  
{\em DAHA-Satake map}, which is the infinite symmetrizer
on the affine Hecke subalgebra, and its relation
to the {\em affine Satake map} (and related constructions)
defined by the formulas used in \cite{Ka,FGT,BK}. The 
latter map is directly connected with the 
theory of Jackson integration developed in \cite{C1,C5,Sto}, which
provides exact formulas at levels $0,1$;
see also \cite{FGT}, Section 12.7 ``Lattice-hypergeometric sums."
Interestingly, the DAHA-Satake map
and the affine Satake map have different convergence ranges. 
The latter 
is well defined for any nonzero $t$, the former only as 
$\Re k < -1/h$ for $t=q^k$ and the Coxeter number $h$;
$|q|<1$ in the paper.
When both converge, they are proportional to each other.

The affine Satake series becomes essentially the 
{\em Weyl-Kac character formula} in the limit $t\to \infty$.
On the other hand, 
the DAHA-Satake map appeared to be related to the 
{\em Demazure characters}, due to the main proportionality
theorem and the $Y$\~formulas from 
Theorems \ref{YLEFT},\ref{YLEFTNEW}. We note that $t$\~counterparts
of the Kac-Moody string functions (and related matters)
are not discussed in this paper; see \cite{FGT,Vi}. 
Also, what seems promising to us is the study 
of the monodromy of the affine Hall functions (generalizing the
classical theorem due to Kac and Peterson); 
we hope to consider this problem in other works.

Concerning the algebraic theory of DAHA, the Satake
map and affine Hall functions are closely related to 
{\em DAHA coinvariants\,}, which, in turn,
are directly connected with the symmetric invariant
bilinear forms on DAHA of levels $l\ge 0$.
The bilinear forms of level $0$ and $1$ are exactly the key 
inner products from \cite{C101} and other works of the first author.
For arbitrary levels $l>0$, the space of DAHA coinvariants is 
isomorphic to the corresponding Looijenga space. Various applications
of the DAHA coinvariants are expected in mathematics and physics.

\subsubsection{\sf Spinor Whittaker functions}
The focus of the second part of this work 
is on the nonsymmetric Whittaker theory 
for $A_1$. The classical Whittaker functions are already nonsymmetric, 
so we need a new theory of spinors (generally $W$\~spinors) 
to achieve this; some its instances already appeared in the 
related harmonic analysis (we will discuss this).
\smallskip

The construction of the {\em spinor-Dunkl operators}
for the $q$\~Toda operators (also called chains) 
is an important and 
unexpected development in this classical field. It
can be presented as an isomorphism between the standard 
polynomial representation of {\em nil-DAHA} and the 
spinor-polynomial representation of its dual. The reproducing
kernel of this isomorphism is the {\em spinor nonsymmetric
Whittaker function}, which was mentioned in \cite{ChW}
as a possible major continuation of the theory
of $q$\~Whittaker functions. We note that the definition
of the difference (relativistic) Toda chain in the case 
of $A_n$ 
in the classical and quantum variants is essentially due 
to Ruijsenaars; see \cite{Rui} for a review.

In this paper the formula for the
{\em nonsymmetric} Whittaker function is discussed
for $A_1$ only. See \cite{ChW} for the 
theory of {\em global symmetric
$q$\~Whittaker functions}, which are closely connected with 
the theory of affine flag varieties and Givental-Lee 
theory. They may have other applications too; see \cite{GLO}.
Technically, the introduction of
{\em nonsymmetric} Whittaker functions is an important step for
using DAHA methods at their full potential.

It is important that the same
limit $t\to \infty$ serves the $q$\~Whittaker functions and
the passage to Kac-Moody theory. However, this limit must 
be calibrated in a very special way in the Whittaker case 
following the Ruijsenaars procedure (see \cite{Et1} and
\cite{ChW}). As a matter of fact, obtaining the 
Kac-Moody characters is also not immediate from DAHA; the affine 
Satake map is needed here, the major theme of the first part of
this work.  
The $q$\~Hermite polynomials emerge in the limit $t\to\infty$ 
for both, $q$\~Whittaker and
Kac-Moody theories. They play an important role in our analysis. 
The resulting connection between Kac-Moody 
theory and $q$\~Whittaker theory is expected
to be related to the geometric quantum Langlands program. 
\smallskip

\subsection{{\bf Dunkl operators via DAHA}}
To put this paper into perspective, let us briefly
outline the (current) status of DAHA theory from the viewpoint
of the constructions of the Dunkl operators. The families of the Dunkl
operators are essentially in one-to-one correspondence with  
the constructions of DAHA ``polynomial representations". The latter
are generally those induced from the affine Hecke subalgebras of 
DAHA, their variants and degenerations. Not all of them 
are really polynomial;  
{\em Fock representations} may be a better name.

Such approach to reviewing applications of DAHA is of course 
simplified, but maybe not too much. 
For instance, if the polynomial representation is known and 
well studied, then we know a lot about the corresponding DAHA. 
It gives the PBW theorem, the zeros of the corresponding Bernstein-Sato 
polynomial, the definition of the localization functor, the
construction of the corresponding spherical function and
more of these. 

The {\em spinor-polynomial} representation
needed for the $q$\~Toda-Dunkl operators appeared of a new type
(not exactly induced from AHA), which reflects interesting
new features of nil-DAHA. To explain it, let us begin with the
list of major families of Dunkl operators.

\subsubsection{\sf Main families of Dunkl operators}
We will stick to the crystallographic case; there
are important developments for the groups
generated by complex reflections and those generated by
symplectic ones (though the latter generally 
do not result in Dunkl-type operators).
With this reservation, the list of major known
families of Dunkl operators and corresponding 
polynomial representations is as follows.

(a) The rational-differential operators 
due to Charles Dunkl; {\em rational
DAHA\,} is self-dual and its theory (including the
polynomial representation) is the most developed now.

(b) Differential-trigonometric
and difference-rational polynomial representations
of {\em degenerate DAHA\,}; they are connected by
the generalized Harish-Chandra transform. 
 
(c) Macdonald theory and $q,t$\~DAHA, corresponding 
to the difference-trigonometric polynomial representation
and the corresponding Dunkl operators; it
is self-dual as in the rational case.

(d) Differential-elliptic representation of 
degenerate DAHA and the difference-elliptic 
representation of $q,t$\~DAHA \cite{Ch12,Ch13}; 
their dual counterparts have not been studied so far.

(e) The specializations of the representations
from $(b)$ in the theory of Yang-type systems of 
spin-particles. The references are \cite{Ug}
and \cite{EOS}; degenerate DAHA governs their theory.
\smallskip

Let us mention that the families from $(d)$
were introduced in  \cite{Ch12} and \cite{Ch13}, 
but there is no reasonably complete theory of these
representations so far. They are connected with the
{\em affine Hall functions\,}, the major
theme of the first part of the paper.
\smallskip

\subsubsection{\sf The Toda-Whittaker case}
The nonsymmetric $q$\~Whittaker functions are eigenfunctions
of new {\em spinor\,} Dunkl operators defined using {\em nil-DAHA},
which adds a new dimension to the list above.
The $q$\~Toda-Dunkl operators do require the 
spinors; they are different from those of the induced type defined in
\cite{Ch8} (and their degenerations).
 
The usual (``symmetric") $q$\~Whittaker functions have various 
applications, exceeding those of the difference spherical functions.
One of the reasons is that the coefficients of the 
$q$\~Whittaker functions are $q$\~integers.
\smallskip

There is a limiting procedure due to Ruijsenaars
that connects the $q$\~Toda operators and the difference QMBP;
see \cite{Rui, Et1}. 
It must be significantly modified in the nonsymmetric
case using the spinor setting and 
eventually leads to the {\em spinor polynomial representation},
an irreducible module of {\em nil-DAHA\,} of a new kind. 

To be more exact, the latter representation is a counterpart
of the polynomial representation 
multiplied by the Gaussian. Its nil-Fourier-dual equals 
the Gaussian times the standard polynomial representation of 
nil-DAHA. The map intertwining these two
representations is given in terms of the {\em nonsymmetric spinor 
global $q$\~Whittaker function\,}. The construction is a general
one, but we will stick to the $A_1$\~case in this work.
\smallskip

\subsection{{\bf The technique of spinors}}
It is an important tool in the QMBP 
(the Heckman-Opdam eigenvalue problem) and DAHA theory.
The main objective of the spinors is to address the problem 
that the Dunkl operators are not local; 
they become local in the space of spinors. Another (related)
purpose of this technique is to incorporate into DAHA theory 
all solution, not only $W$\~invariant, of the QMBP,
its generalizations and variants.
Solving QMBP in the class of all functions has interesting 
algebraic and analytic aspects. We will not try to review 
them here. 

As far as we know, this technique was used explicitly for the 
first time in \cite{C13}, when proving the so-called 
Matsuo- Cherednik isomorphism theorem. This theorem establishes 
an equivalence of the affine Knizhnik-Zamolodchikov equation, 
AKZ,  in the modules
of the degenerate Hecke algebra induced from 
(dominant) characters and the corresponding Heckman-Opdam
system (QMBP). See Chapter 1 of \cite{C101}, \cite{O2}
and Section \ref{sect:AKZ} below. 

Using the technique of spinors systematically (see
Section \ref{sect:AKZ}) makes the proof from \cite{C13}
entirely algebraic and establishes its direct connection
with the proof suggested (several years later) in \cite{O2};
compare Lemma 3.2 there with Theorem \ref{SPINDUN} below.
The approach from \cite{O2} is actually very close to the
justification of this theorem in our paper. Mathematically, 
Opdam's proof was essentially equivalent to the 
one from \cite{C13}, but this was done in \cite{O2} entirely
algebraically; the spinors in their algebraic
variant were certainly present there.

We note the technique of spinors (combined
with the explicit calculation of the AKZ-monodromy)
was actually used in \cite{C13} to obtain the {\em nonsymmetric
spherical function\,}, called the $G$\~function
in \cite{O2}. 

Generally speaking, there is nothing very new about
the definition of spinors, $W$\~spinors
to be more exact. They are simply sets of functions 
$\{f_w\}$ numbered by the elements from the Weyl group $W$
with the action of $W$ on the indices. 
The principle spinors are in the form $\{w^{-1}(f),\,w\in W\}$ for 
a global function $f$; generally $f_w$ are absolutely
independent functions. For instance, the {\em real spinors} 
are functions on the disjoint union of all Weyl chambers, 
collected (using $W$) in the fundamental Weyl chamber.
It is not surprising that they appeared in various contexts before. 


\subsubsection{\sf Connections to AKZ}
The Matsuo proof of the relation between AKZ and
QMBP from paper \cite{Mats}
was a direct algebraic verification. The Grothendieck-type 
notion of the monodromy {\em without
a fixed point} used in \cite{C13} made the proof very short 
and entirely conceptual; also, this paper was written for the
vector-valued solutions and included the rational QMBP.
Using this approach, such an equivalence was extended to the 
difference and elliptic cases. In the difference theory, this
map can be an embedding of the spaces of solutions (not an isomorphism); 
see \cite{Ch8}, which was finalized in \cite{Sto2}. 

The definition of the elliptic QMBP requires the trivial central 
charge condition, which is  $l=-kh$ for the Coxeter
number $h$ (where $t=q^k$); then the equivalence will hold too.
Apart from the elliptic case, the isomorphism theorems
from \cite{C13} and \cite{C101} (Chapter 1) can be
stated as follows.

\newtheorem*{keytheorem}{Theorem}

\begin{keytheorem}[{\sf AKZ$\to$Dunkl$\to$QMBP}]
Given an arbitrary weight $\la$, the space of AKZ-solutions
in the induced module $I_\la$ of the (degenerate) affine
Hecke algebra can be identified with the $\la$\~eigenspace of Dunkl 
operators in the corresponding DAHA {\sf spinor} representation. 
Then the latter eigenspace can be mapped to the space of all, not 
necessarily symmetric, solutions of the corresponding QMBP. 
For generic $\la$, this
map is an isomorphism (an embedding in the difference setting). 
\end{keytheorem}    
 
The spinors needed here are {\em complex}, defined in the domain
$U=\{z\}$ such that $\Im(z)$ belongs to the corresponding
fundamental Weyl chamber. They can be interpreted as
functions in the disjoint union $\cup_{w\in W} w(U)$; then
the principle spinors are global analytic functions. Using $W$,
we can gather these functions in $U$. Only functions in $U$ emerge 
in the spinor theory of the Dunkl-type eigenvalue problem, 
including the spinor integration and related inner products.
\smallskip

\subsubsection{\sf On the localization functor}
This construction is connected with the {\em localization 
functor}, one of the most powerful tools in the theory of DAHA. 
See \cite{GGOR} and \cite{VV1}. 
The localization construction assigns a local
system to a module of DAHA (from a proper category);
the case of induced representations is related to AKZ
and paper \cite{C13} as follows.

The starting point of the latter paper
was the AKZ with values in an arbitrary
finite-dimensional module $V$ of AHA (or degenerate AHA).
Then the spinor Dunkl operators were defined for these AKZ
via the monodromy representation. Combining
these Dunkl operators with the operators of 
multiplication by functions supplies the space 
of $V$\~valued analytic functions with the DAHA action. 

The relation of the spinor Dunkl operators to
the monodromy of AKZ is of independent interest.
The {\em monodromy cocycle} on $W$ from \cite{C13} 
(see also \cite{C101}, Chapter 1)
can be expressed in terms of the (usual) monodromy 
homomorphism of the braid group. This establishes a link 
to the localization functor.

We note that the construction AKZ$\to$Dunkl$\to$QMBP was aimed 
at applications to the corresponding eigenvalue 
problems and was done only within the class
of induced modules; the projective modules are of 
key importance for the theory of the localization functor.    
\medskip

\subsubsection{\sf The setting of the work}
We mainly use the standard affine root systems in contrast
to the {\em twisted} affine root systems considered in
\cite{C101} and many papers on DAHA.
The standard (untwisted) ``affinization" is (presumably) 
exactly the one 
compatible with the quantum  Langlands duality. For instance,
the {\em untwisted} affine exponents from \cite{C103},
describing the reducibility of the
polynomial representation, obey the quantum 
Langlands-type duality for the modular transformation
$q\mapsto \hat{q}$. This kind of
duality does not hold in the {\em twisted case} (at least,
we do not know how to formulate it). On the other hand,
the twisted affinization
has obvious merits (versus the standard setting)
for the theory of Gaussians. This is 
parallel to the advantages of the
twisted case for level-one character formulas 
in Kac-Moody theory.

Due to the standard (untwisted) setting, we need to state some of 
the results of this paper, especially where the Gaussians are 
involved, only for the simply-laced root systems.  We hope
to consider the corresponding {\em twisted} case in other 
publications. Using $t$  in this
paper is relaxed as well; we simply treat it as a single parameter.
Generally, $t$ (or $k$) are supposed to depend on the length
of the corresponding root.
In the second part of this work, we present some constructions 
only in the $A_1$\~case, where practically everything can be 
calculated explicitly.
However, the major results of this paper can be transferred to 
(or expected to hold for) arbitrary root systems.

The readers familiar with AHA and 
classical $p$\~adic theory can go directly to
the double affine generalizations, though the introduction
of the Macdonald's $p$\~adic spherical functions as symmetrizations
of Matsumoto functions, which are essentially delta functions, 
is not quite standard (even for specialists).
\medskip

\subsection{{\bf Acknowledgements}}
\subsubsection{\sf Harvard lectures}
The paper is based on a series of lectures 
delivered by the first author at Harvard (February-March 2009);
he is responsible for the scientific contents of this paper.

It was a somewhat unusual series, a sort
of reporting the current research activities on weekly basis.
The output of these lectures appeared better than the
lecturer expected (hopefully, for the listeners too). 
The initial TeX files of the lectures were prepared by 
Xiaoguang Ma.

Extensive usage of examples and  
exposition of the classical topics are an organic part of the 
design of this work. However, the focus is on  
general approaches and new results. Almost all examples
and direct verifications are needed to prepare 
affine and spinor generalizations, the main purpose of this work.

\subsubsection{\sf Special thanks}
My special thanks go to Dennis Gaitsgory and Pavel
Etingof for participating in these lectures, 
shaping their direction and contents,
and for various important discussions. 

Alexander Braverman and David Kazhdan significantly 
influenced the key topics of this series of lectures
and the papers. David Kazhdan greatly helped in 
improving the theory of the affine symmetrizers presented
in this work. I am grateful to Alexander Braverman, 
Ian Grojnowski and Manish Patnaik for discussions 
concerning the affine Hall functions. Talks to Victor Kac and 
Boris Feigin were helpful in establishing the connections with
Kac-Moody theory. I am thankful to Roman Bezrukavnikov, Michael
Finkelberg and Victor Ostrik for discussions
of the affine flag varieties and quantum groups. I thank 
Pavel Etingof, Eric Opdam, Simon Ruijsenaars, Jasper Stokman, 
the referee, and Ann Kostant for valuable comments.

The work is partially based on my notes on spinors (reported
at the University Paris VI in 2004 and at RIMS in 2005) and
on the DAHA approach to the decomposition of
the regular representation of AHA
(see \cite{ChL},\cite{HO2}) 
reported at CIRM (2006), MIT (2007) and at the University of 
Amsterdam (2008). Working on these papers continued at RIMS 
(Kyoto University, 2009) and completed at
the Hebrew University (2012). I am very grateful for 
the invitations.

Quite a few topics were stimulated by my talks to 
physicists; special thanks 
to Anton Gerasimov who introduced me to 
the brave new world of 
$q$\~Whittaker functions.   

I am grateful to the many people, mathematicians and physicists,
I talked to on these and related matters at Harvard, MIT, RIMS,
the Hebrew University, and in many other places.

\hfill {--Ivan Cherednik}

\setcounter{equation}{0} 
\section{\sc{P-adic theory revisited}}
\noindent
The area of affine Hecke
algebras, AHA, and spherical functions is vast. 
The classical $\mathfrak{p}$\~adic spherical
functions were subject to various generalizations. 
It is most important to note that 
they are limits as $q\to 0$ of the {\em symmetric Macdonald
polynomials}, due to Ian Macdonald. Similarly,
the limits of the {\em nonsymmetric Macdonald polynomials} are
the Matsumoto spherical functions, key to our approach.
The DAHA methods help a lot in clarifying the algebraic aspects 
of their theory.
See Section 2.11 from Chapter 2 in \cite{C101} (and references
therein) and \cite{O4}; see also \cite{Ion2,O5}. 

The purpose of this section is  
revisiting the $\mathfrak{p}$\~adic theory from the 
viewpoint of DAHA, which aims at establishing connections with the
affine Hall functions and $q$\~Whittaker functions.
  
\subsection{{\bf Affine Weyl group}} \label{sect:affineweyl} 
\subsubsection{\sf Root systems}
Concerning the classical theory of root systems and Weyl
groups, the standard references are \cite{Bo,Hu}; if these
sources are insufficient, then see \cite{C101}.

In this paper
$R=\{\al\}\subset \mathbb{R}^{n}$ is a simple reduced root system 
with respect to a nondegenerate symmetric bilinear form $(,)$ on 
$\mathbb{R}^{n}$. Let $\{\alpha_{i}\}_{i=1}^{n}\subset R$ be the set 
of simple roots and let $R_{+}$ (or $R_{-}$) be the set of positive (or 
negative) roots. The coroots are denoted by 
$\alpha^{\vee}=2\alpha/(\alpha,\alpha)$; $W$ is the Weyl group 
generated by $s_\al$. 

Let $Q=\bigoplus_{i=1}^{n}\mathbb{Z}\alpha_{i}$, 
$P=\bigoplus_{i=1}^{n}\mathbb{Z}\omega_i$, correspondingly, 
let $Q^{\vee}=\bigoplus_{i=1}^{n}\mathbb{Z}\alpha_{i}^{\vee}$ be the 
coroot lattice and $P^{\vee}=\bigoplus_{i=1}^{n}\mathbb{Z} 
\omega^{\vee}_{i}$ the coweight lattice, where 
$\{\omega_{i}^{\vee}\}$ are the fundamental coweights, i.e., 
$(\omega^{\vee}_{i}, \alpha_{j})=\delta_{ij}$. Replacing 
$\mathbb{Z}$ by $\mathbb{Z}_{+}=\Z_{\ge 0}$, 
we obtain $Q_+,Q^{\vee}_{+}$ and $P_+,P^{\vee}_{+}$. 

The maximal positive root will be denoted by $\theta$,
and the bilinear form will be normalized by the condition 
$(\theta,\theta)=2$; also
$\rho\equal\dfrac{1}{2}\sum_{\alpha\in R_{+}}\alpha$. 
Due to this normalization, 
\begin{align*}
&\ \ \ \ \ \ \ \ \ \ \ \ \ \ \ \ \ \ \ \ \ \ \ \ &&Q&      
&\subset&      &P&  
&\ \ \ \ \ \ \ \ \ \ \ \ \ \ \ \ \ \ \ \ \ \ \ \ & \\ 
&\ \ \ \ \ \ \ \ \ \ \ \ \ \ \ \ \ \ \ \ \ \ \ \ &&\cup&   
&&          &\cup&  
&\ \ \ \ \ \ \ \ \ \ \ \ \ \ \ \ \ \ \ \ \ \ \ \ & \\
&\ \ \ \ \ \ \ \ \ \ \ \ \ \ \ \ \ \ \ \ \ \ \ \ &&Q^\vee& 
&\subset& &P^\vee.& 
&\ \ \ \ \ \ \ \ \ \ \ \ \ \ \ \ \ \ \ \ \ \ \ \ &
\end{align*}

We stick to reduced root systems in this paper,
sometimes even to the $A-D-E$ systems.
Almost all results in the theory of DAHA 
and related Macdonald polynomials for reduced root systems
were transferred to the case of $C^\vee C$, the ultimate
nonreduced system, and to the corresponding theory of
Koornwinder polynomials.

\subsubsection{\sf Affine root systems} 
The vectors $\tilde{\alpha}= [\alpha,j]\in 
\mathbb{R}^{n}\times\mathbb{R}$, where $\alpha\in R$
and  $j\in \mathbb{Z}$, 
form the {\em standard affine root system} $\widetilde{R}$. 
The set of positive 
affine roots is 
 $\widetilde{R}_{+}=\{[\alpha,j]\,|\,j\in \mathbb{Z}_{>0}\}\cup 
 \{[\alpha,0]\,|\,\alpha\in R_{+}\}$.
Define $\alpha_{0}=[-\theta, 1]$, 
where $\theta$ is the maximal 
positive root in $R$.  
We will identify $\alpha\in R$ with 
$\tilde{\alpha}=[\alpha,0]\in \widetilde{R}$. The affine simple
roots $\{\al_i, 0\le i\le n\}$ form the extended (also called
affine) Dynkin diagram $Dyn^{\tt aff}
\supset Dyn=\{\al_i,1\le i\le n\}$.

For an arbitrary affine root $\tilde{\alpha}=[\alpha,j]$ and 
$\tilde{z}=[z, \zeta]\in \R^{n+1}$, the corresponding 
reflection is defined as follows: 
$$
s_{\tilde{\alpha}}(\tz)
=\tz-2\frac{(z, \alpha)}{(\alpha,\alpha)}\,\tilde{\alpha}
=\tz-(z,\alpha^{\vee})\,\tilde{\alpha}.
$$
We set $s_i=s_{\al_i}$ for $i=0,\ldots,n$.
The affine Weyl group $\widetilde{W}$ is generated by 
$\{s_{\tilde{\alpha}}\,|\,\tilde{\alpha}\in \widetilde{R}_{+}\}$;\,
$\{s_i\}$ for $i\ge 0$ are sufficient.

\begin{theorem}
We have an isomorphism
$$\widetilde{W}\cong W\ltimes Q^{\vee},$$ 
where the translation $\al^\vee\in Q^\vee$ is naturally
identified with the composition $s_{[-\al,1]}s_{\al}\in \tW$.
In terms of the action in $\R^{n+1}\ni \tz$, one has 
$b(\tz)=[z,\zeta-(b,z)]$ for $\tz=[z,\zeta],\, b\in Q^\vee;$
notice the sign of $(b,z)$.\sq
\end{theorem}

Define the {\em extended affine Weyl group} to be 
$\widehat{W}=W\ltimes 
P^{\vee}$ acting on $\R^{n+1}$ via the last formula from
the theorem with
$b\in P^\vee$. Then $\widetilde{W}\subset \widehat{W}$. 
Moreover, we have the following theorem. 
Let $\Aut=\Aut(Dyn^{\tt aff})$\, and
 $O\equal \{r\}$ for $\Aut(\al_0)=\{\al_r\}$, i.e., $O$
is formed by the indices of the simple roots from the
Aut\~orbit $\Aut(\al_0)$ of $\al_0$.

\begin{theorem} (i) The group
$\widetilde{W}$ is a normal subgroup of $\widehat{W}$ and 
$\widehat{W}/\widetilde{W}=P^{\vee}/Q^{\vee}$. 
The latter group can be identified with the group 
$\Pi=\{\pi_r\}$ of the elements of $\hW$ permuting simple 
affine roots under their action in $\R^{n+1}$. 
It is a normal commutative subgroup of 
$\Aut$; the quotient $\Aut/\Pi$
is isomorphic to the group $A_0=\Aut(Dyn)$ 
of the automorphisms preserving 
$\al_0$.

(ii) The
indices $r\in O^*\equal O\setminus \{0\}$ are exactly those for 
the minuscule coweights $\om_r^\vee$ satisfying the inequalities 
$(\al,\om_r^\vee)\le 1$ for all $\al\in R_+$. 
The elements $\pi_r\in \Pi$ 
are uniquely determined by the relations $\pi_r(\al_0)=\al_r$ 
($\pi_0=$id). An arbitrary element $\hw\in \hW$ can
be uniquely represented as $\hw=\pi_r\tw$ for $\tw\in \tW$.\sq
\end{theorem} 

It is not difficult to calculate $\pi_r$ explicitly
(see \cite{C101}):
\begin{align}\label{piromr}
&\pi_r=\om_r^\vee\, u_r^{-1} \hbox{\ \,for minuscule\ \,} 
\om_r^\vee\in P_+^\vee\subset \hW,\ u_r=w_0w_0^{(r)},
\end{align}
where $w_0^{(r)}$ is the element of maximal length in the
centralizer of $\om_r^\vee$ in $W$ for $r\in O^*$,
$w_0$ is the element of maximal length in $W$. Equivalently,
$u_r$ is of minimal possible length 
such that  $u_r(\om_r)
\in P_-=-P_+$ (see the next section). Note that 
$\pi_{r}s_{i}\pi_{r}^{-1}=s_{j}$ if 
$\pi_{r}(\alpha_{i})=\alpha_{j}$, $0\leq i\leq n$. 

\subsubsection{\sf The length function}\label{sect:length}
Any element $\hat{w}\in \widehat{W}$ can be written as 
$\hat{w}=\pi_{r}\tilde{w}$ for $\pi_{r}\in \Pi$ and $\tilde{w}\in 
\widetilde{W}$. The length $l(\hat{w})$ is defined to be the length 
of the {\em reduced decomposition} 
$\tilde{w}=s_{i_{l}}\cdots s_{i_{1}}$ 
(i.e., with minimal possible $l$) in terms of the simple reflections 
$s_{i}$. Thus, by definition, $l(\pi_{r})=0$. 

This is the standard {\em group-theoretical} definition. 
There are two other (equivalent)
definitions of the length for
the crystallographic groups, {\em combinatorial} and {\em
geometric}. Namely, the length 
$l(\hw)$ is the cardinality 
$|\tR_+\cap \hw^{-1}(\tR_-)|$ and can also be
interpreted as the ``distance" from the standard affine Weyl 
chamber to its image under $w$. Both definitions
readily give that  $l(\pi_{r})=0$; indeed,  $\pi_r$ sends 
positive roots $\tal$ to positive roots
and (therefore) leaves the standard affine Weyl chamber 
invariant.

Either the combinatorial or the geometric definition can be used 
to check that $l(w(b))=2(\rho,b)$ for arbitrary $b\in P^\vee_+$ and 
$w\in W$.

All three approaches to the length function are important 
in the combinatorial theory of affine Weyl groups, which is far 
from being simple and completed.

\subsubsection{\sf Twisted affinization}\label{sect:thinaff}
There is another affine extension $R^\nu$ of $R$,  
convenient in quite a few constructions (especially, when 
the DAHA Fourier transform and the Gaussians are studied).
This is the setting in \cite{C101} and in quite a few of
author's papers. This extension is defined
for the maximal {\em short} root $\vth$ instead of 
the maximal root $\th$. Accordingly,
$(\alpha,\alpha)=2\,$ for short roots and affine roots
are introduced as $\tilde{\alpha}=[\alpha,\nu_\alpha\, j]$
for $\nu_\al\equal\frac{(\al,\al)}{2}$ ($=1,2,3$). 
Adding $\al_0=[-\vth,1]$ for such $\vth$ to $\{\al_i,i>0\}$,
the resulting diagram is the extended Dynkin diagram 
$(Dyn^\vee)^{\tt aff}$ for $R^\vee$
where all the arrows are reversed. On can simply set
$\tR^\nu\equal((R^\vee)^{\tt aff})^\vee$, 
where
the form in $R^\vee$ is normalized by the (usual)
condition $(\al^\vee,\al^\vee)=2$ for long $\al^\vee$, 
which makes $\vth$ the maximal root in $R^\vee$. The second
check in $((R^\vee)^{\tt aff})^\vee$ is applied to the {\em affine} roots.
The formula $s_{[-\al,\nu_\al]}s_{\al}=\al\,$ naturally 
results in unchecked $Q,P$ in the {\em twisted affine Weyl group}:
$$\for R^\nu\,:\, \
\widetilde{W}\cong W\ltimes Q,\ \hW\cong W\ltimes P.
$$
In $\mathfrak{p}$-adic theory, the 
twisted Chevalley group is a {\em form} of the split group
for a proper Galois extension of the starting field.  
 
The appearance of $Q,P$ in $\tW,\hW$ results in the
invariance of the corresponding DAHA with respect
to the Fourier transform and other basic automorphisms.
This is the main reason why the book \cite{C101} is mainly 
written in such a ``self-dual" setting. Due to the special
choice of the normalization, $Q\subset Q^\vee$ in this
case; recall that $(\vth,\vth)=2$. The term ``twisted" matches 
similar terminology in Kac-Moody theory.

\subsection{{\bf AHA and spherical functions}}
\subsubsection{\sf Affine Hecke algebras}
The affine Hecke algebra $\mathcal{H}$ is generated by 
$T_{0}, T_{1}, \ldots, T_{n}$ and the group $\Pi=\{\pi_{r}\}$ 
with the relations: 
\begin{align}
&\underbrace{T_{i}T_{j}T_{i}\ldots}_{m_{ij} \text{ times}}=
\underbrace{T_{j}T_{i}T_{j}\ldots}_{m_{ij} \text{ times}}\,,\notag\\
&(T_{i}-t^{1/2})(T_{i}+t^{-1/2})=0,\label{heckerel}\\
&\pi_{r}T_{i}\pi_{r}^{-1}=T_{\pi_{r}(i)}.\notag
\end{align}
where $\pi_r(i)$ is the suffix of the simple root $\pi_r(\al_i)$;
$m_{ij}$ is the number of edges between vertex $i$ 
and vertex $j$ in the affine Dynkin diagram $Dyn^{\tt aff}$ and 
$t$ is a formal parameter
(later, mainly a nonzero number). 
\smallskip

\rmk
The above definition gives the affine Hecke algebra with 
{\em equal parameters}. More systematically, we can introduce a 
family of formal 
parameters $\{t_{\al}\}$ depending only on $|\al|$, setting
$t_{i}=t_{\al_i}$ for $0\le i\le n$.
Replacing relations (\ref{heckerel}) by the relations 
$(T_{i}-t_{i}^{1/2})(T_{i}+t_{i}^{-1/2})=0$, we come to the 
definition of the affine Hecke algebra standard in (modern)
geometric and/or algebraic theory (in the case
of {\em unequal parameters}).

The formulas below can be readily adjusted to this setting,
namely, $t_i$ must be used for $T_i$ 
and the subscript $\al$ must be added to $t$ in the formulas
involving  $Y_{\al^\vee}$. In DAHA theory, the same must be 
done for  $X_\al$; also, the relation $t=q^k$ below will become  
$t_\al=q^{k_\al}$. 
If $\tR^\nu$ is used instead of $\tR$, with $Y_\al$ instead of
$Y_{\al^\vee}$, then $q$ must be also replaced by 
$q_\al=q^{\nu_\al}$ in the formulas; accordingly, 
$t_\al=q_\al^{k_\al}$.   
\sq
\smallskip

For any element $\hat{w}\in \widehat{W}$, define 
$T_{\hat{w}}=\pi_{r}T_{i_{l}}\cdots T_{i_{1}}$, where 
$\hat{w}=\pi_{r}s_{i_{l}}\cdots s_{i_{1}}$ is a reduced 
representation of $\hat{w}$. The definition of $T_{\hat{w}}$ does 
not depend on the choice of the reduced decomposition. 

Setting $Y_b=T_b$ for  $b\in P^{\vee}_{+}\subset \widehat{W}$,  
one has $Y_{b}Y_{c}=Y_{c}Y_{b}$ for such (dominant) $b,c$;
use that $l(b)=2(\rho,b)$ for dominant $b$.
For any $a\in P^{\vee}$, we set $Y_{a}\equal Y_{b}Y_{c}^{-1}$,
where $a=b-c$ with some $b,c\in P^{\vee}_{+}$;
the commutativity guarantees that $Y_a$ depends only on $a$.
This definition is due to Bernstein, Zelevinsky, and Lusztig,
see, e.g., \cite{L}. 

Let $\mathscr{Y}\equal
\mathbb{C}[Y_{\omega_{i}^{\vee}}^{\pm}] \subset 
\mathcal{H}$. Then  
$$\mathcal{H}=\langle\mathscr{Y}, T_{1}, \ldots, T_{n}\rangle.$$
Indeed, $T_0=Y_\th T_{s_\th}^{-1}$ and $\pi_r=
Y_{\om_r^\vee}T_{u_r}^{-1}$ (see (\ref{piromr})).

\begin{theorem}\label{PBW}
(i) An arbitrary element $H\in \h$ can be uniquely
represented as $H=\sum c_{b,w}\,Y_b\,T_w$ (a finite sum)
for $b\in P^\vee, w\in W$, which is called the {\sf PBW Theorem}.

(ii) The subalgebra $\mathscr{Y}^W$ of $W$\~invariant
$Y$\~polynomials is the center of $\h$ (the {\sf Bernstein Lemma});
here $w(Y_b)=Y_{w(b)}$, see also Lemma \ref{Berns}.
\end{theorem}

\subsubsection{\sf Matsumoto functions}
Let $\H=\mathcal{H}_{\tt nonaff}$ be the Hecke algebra associated 
with the nonaffine root system $R$, i.e., generated by
$T_i$ for $1\le i\le n$. We define the 
$t$\~symmetrizer by 
the formula
$$
\mathscr{P}_{+}=\frac{\sum_{w\in W}t^{l(w)/2}T_{w}}
{\sum_{w\in W}t^{l(w)}}\in \H\,.
$$
One checks directly or using (\ref{tidelta}) below that
$$
\frac{(1+t^{1/2}T_{i})\mathscr{P}_{+}}{1+t}=
\mathscr{P}_{+},\ 1\le i\le n.
$$

The following renormalization 
$\delta_{\hat{w}}=t^{-l(\hat{w})/2}T_{\hat{w}}$ of $T_{\hw}$
(any $\hw\in \hW$) is convenient to
establish the connection with $\mathfrak{p}$-adic theory.
Then\,
\begin{equation}\label{tidelta}
T_{i}\delta_{\hat{w}}=\left\{\begin{array}{ccc}t^{1/2}
\delta_{s_i\hat{w}}, &  & \text{ if }l(s_{i}\hat{w})=l(\hat{w})+1; \\
t^{-1/2}\delta_{s_i\hat{w}}+(t^{1/2}-t^{-1/2})\delta_{\hat{w}}, 
& & \text{ otherwise.}\end{array}\right.
\end{equation}

Let 
$\Delta=\bigoplus_{\hat{w}\in\widehat{W}}\mathbb{C}\delta_{\hat{w}}$ 
be the (left) regular representation of $\mathcal{H}$. Its 
{\em  spherical submodule} is defined as follows: 
$$\Delta^{\sharp}=\Delta\mathscr{P}_{+}\cong 
\mathscr{Y}\mathscr{P}_{+}.$$
Identification with the Laurent $Y$\~polynomials
is based on claim $(i)$ (PBW) of Theorem  \ref{PBW}. 

From now on $\Delta^{\sharp}$ will be identified with 
$\mathscr{Y}$, i.e., $1\in \mathscr{Y}$ will be actually  
$\mathscr{P}_+$. By
$\delta^{\sharp}_{\hat{w}}$, we denote the image of 
$\delta_{\hat{w}}$ in $\Delta^{\sharp}$; explicitly,\,
$\delta_{\hat{w}}^{\sharp}\equal \delta_{\hat{w}}\mathscr{P}_{+}$. 

The {\em Matsumoto functions} \cite{Mat}, also called nonsymmetric
$\mathfrak{p}$-adic spherical functions, are defined 
(in this approach) to be 
$$\varepsilon_{b}=\delta_{b}^{\sharp}, \quad\forall\, 
b\in P^{\vee},$$
i.e., we simply restrict $\de^{\#}$ to $P^\vee$ here.
From this definition, 
$\varepsilon_{b}=t^{-(b,\rho)}Y_{b}$
for any $b\in P_{+}^{\vee}$. Representing (calculating)
$\varepsilon_{b}$ as a Laurent polynomial in terms of $Y$
for any $b\in P^{\vee}$ is of fundamental importance.

\subsubsection{\sf The rank-one case}\label{sec:rank1}
In the $A_{1}$ case, we can set 
$\om=\omega_1^{\vee}=\om^\vee$;
then  $\alpha=\al_1=2\omega$ and 
$\rho=\om$. The extended affine Weyl group $\widehat{W}$ 
is generated by $\pi=\pi_1$ and the reflection $s=s_{\alpha}$. As an 
element of $\widehat{W}$, $\om=\pi s$. Let 
$T=T_{1}\in \mathcal{H}$; then $Y=Y_{\om}=\pi T $. 

The affine Hecke algebra can be written as $\mathcal{H}=\lan Y, 
T\ran$ subject to
$T^{-1}YT^{-1}=Y^{-1}$ and $(T-t^{1/2})(T+t^{-1/2})=0$. 
The first of these relations is equivalent to $\pi^2=1$
for $\pi$ introduced as $YT^{-1}$. 
 
The symmetrizer is 
$$\mathscr{P}_{+}=\frac{1+t^{1/2}T}{1+t}.$$

For any $m\in \mathbb{Z}$, let $\delta_{m}=\delta_{m\om}$ 
and $\varepsilon_{m}=\delta_{m\om}^{\sharp} 
=t^{-|m|/2}T_{m\om}\mathscr{P}_{+}$. 
Then we have for $m\geq 0$, 
\begin{align}
&T\varepsilon_{m}\ =\ t^{1/2}\varepsilon_{-m},\label{Tvepm}\\
&T\varepsilon_{-m}\ =\ t^{-1/2}
\varepsilon_{-m}+(t^{1/2}-t^{-1/2})\varepsilon_{m}.\label{T-vepm}
\end{align}
Similarly, for $m\geq 0$, 
\begin{align*}
&T^{-1}\varepsilon_{-m}\ =\ t^{-1/2}\varepsilon_{m},\\
&T^{-1}\varepsilon_{m}\ =\ 
(T-(t^{1/2}-t^{-1/2}))\varepsilon_{m}=t^{1/2}
\varepsilon_{-m}-(t^{1/2}-t^{-1/2})\varepsilon_{m}.
\end{align*}

\begin{lemma}\label{LEMPI}
For any $m\in \Z$, $\pi\varepsilon_{m}=\varepsilon_{1-m}$. 
\end{lemma}
{\em Proof.}
Since $\pi^{2}=1$, it suffices to calculate
$\pi \varepsilon_{-m}$ for $m\ge 0$. 
Using that $Y\vep_m=t^{1/2}\vep_{m+1}$ (it results from the
definition of $\vep$ for such $m$),
$$\pi \varepsilon_{-m}=YT^{-1}\varepsilon_{-m}=
t^{-1/2}Y\varepsilon_{m}= 
\varepsilon_{1+m}.$$
\vskip -0.7cm
\sq

Let us apply the lemma to write down the action of $Y^{\pm1}$ on 
$\varepsilon_{m},\vep_{-m}$ for $m\ge 0$: 
 
\begin{eqnarray}
Y\varepsilon_{m}
&=&t^{1/2}\varepsilon_{m+1}\label{yepm},\\
Y\varepsilon_{-m}
&=&t^{-1/2}\varepsilon_{-m+1}+(t^{1/2}-t^{-1/2})\varepsilon_{m+1}
\label{yep-m},\\
Y^{-1}\varepsilon_{m+1}
&=&t^{-1/2}\varepsilon_{m}\label{y-epm},\\
Y^{-1}\varepsilon_{-m}
&=&t^{1/2}\varepsilon_{-m-1}-(t^{1/2}-t^{-1/2})\varepsilon_{m+1}.
\label{y-ep-m}
\end{eqnarray}
Note that (\ref{yepm})and (\ref{yep-m}) overlap at $m=0$,
as well as (\ref{Tvepm})and (\ref{T-vepm}).

The formulas for 
the action of $Y$ and $Y^{-1}$
are called {\em nonsymmetric Pieri rules}; they are {\em obviously}
sufficient to calculate the $\vep$\~functions (which holds 
in any ranks). However, the {\em technique of intertwiners} is 
generally more 
efficient for calculating the $\vep$\~polynomials and their 
variants than direct usage of the Pieri formulas (see, e.g., \cite{C101}).
In this particular example, formula (\ref{Tvepm}) 
is such an intertwiner. It is sufficient indeed:
\begin{eqnarray}\label{vepformulas}
\vep_m & = & t^{-\frac{m}{2}}Y^m \for m\ge 0 
\hbox{\ \, implies \ that}\\ 
\vep_{-m}
&=& t^{-\frac{1}{2}}T\vep_m=t^{-\frac{m+1}{2}}T(Y^m)\notag\\
&=&t^{-\frac{m+1}{2}}(t^{\frac{1}{2}}Y^{-m}+
(t^{\frac{1}{2}}-t^{-\frac{1}{2}})\frac{Y^{-m}-Y^{m}}
{Y^{-2}-1}).\notag
\end{eqnarray}

We are now ready to introduce the {\em $\mathfrak{p}$-adic spherical
functions}. In this (algebraic) approach, they are
\begin{equation*}
\varphi_{m}\equal\frac{1+t^{1/2}T}{1+t}\,\varepsilon_{m},\ m\ge 0.
\end{equation*}
Using formulas 
(\ref{yepm}), (\ref{y-epm}) and
the commutativity of $Y+Y^{-1}$ with $T$ (check it directly
or see below), we establish the {\em symmetric Pieri rules}:
\begin{align}
&(Y+Y^{-1})\varphi_m\ =\ t^{1/2}\varphi_{m+1}+t^{-1/2}\varphi_{m-1}
\hbox{\ as\ } m>0,\notag\\
&(Y+Y^{-1})\varphi_0\ =\ (t^{1/2}+t^{-1/2})\varphi_{1}. \label{pieri}
\end{align}

Note that the latter relation follows from the former if
one formally imposes the periodicity condition
$\varphi_{-1}=\varphi_1$. By construction, $\vph_0=1$; 
all other functions
can be calculated using the Pieri rules. 
All $\varphi_{i}$'s are 
invariant under $s:\,Y\mapsto Y^{-1}$ due to the commutativity
$[Y+Y^{-1},T]=0$.

The first three $\varphi_{m}$'s are as follows: 
\begin{eqnarray*}
\varphi_{0}=1,\quad \varphi_{1}=\frac{Y+Y^{-1}}{t^{1/2}+t^{-1/2}},
\quad
\varphi_{2}=\frac{(Y+Y^{-1})^{2}}{1+t}-t^{-1}.
\end{eqnarray*}

For the system $A_1$, the symmetric Pieri rules look simpler 
than their $\vep$\~counterparts, but this is exactly the other way 
around in higher ranks. Generally, there are no good formulas for the
action of $W$\~orbitsums in the form $\sum_w Y_{w(b)}$
on the spherical functions
(see (\ref{pieri})) except for the minuscule $b=\om_r^\vee$
and $b=\th$. Theoretically, the Pieri formulas are
sufficient to calculate all $\vph$\~polynomials,
but this can be used mainly for $A_n$
and in some cases of small ranks.
The nonsymmetric formulas of type (\ref{yepm}--\ref{y-ep-m}) 
exist (and are reasonably convenient to deal with) 
for arbitrary root systems.

\subsection{{\bf Spherical functions as Hall polynomials}}
\label{sec:Sph via Hall}
\subsubsection{\sf Macdonald's formula}
In general (for any root system $R$ as above), 
we can define the {\em spherical
function} as follows:
\begin{eqnarray*}
\varphi_{b}\equal\mathscr{P}_{+}\varepsilon_{b} =
t^{-(\rho,b)}\mathscr{P}_{+}Y_{b}\mathscr{P}_{+}\in \mathscr{Y},\ 
b\in P^{\vee}_+. 
\end{eqnarray*}
They become $W$\~invariant $Y$\~polynomials upon the identification
of $\De^{\#}$ and $\mathscr{Y}$  (the Bernstein Lemma), where
$w(Y_b)\equal Y_{w(b)}$  for $w\in W$.
Their ($\mathfrak{p}$-adic) theory was developed by Satake, 
Macdonald and others; 
we will mainly call them the {\em Macdonald spherical functions}. 
Macdonald established the following fundamental fact.

\begin{theorem}
Let $P(t)$ be the Poincar\'e polynomial, namely,
$P(t)=\sum_{w\in W}t^{l(w)}$. Then
\begin{eqnarray}\label{macdvarphi}
\varphi_{b}(Y)=\frac{t^{-(\rho,b)}}{P(t^{-1})}
\sum_{w\in W}Y_{w(b)}\prod_{\alpha\in R_{+}}
\frac{1-t^{-1}Y_{w(\alpha^{\vee})}^{-1}}
{1-Y_{w(\alpha^{\vee})}^{-1}}.
\end{eqnarray}\sq
\end{theorem}

The summation on the right-hand side is proportional to
the {\em Hall-Littlewood polynomial} associated with $b\in P_+^\vee$.
The potential poles (due to the denominators) will cancel each other,
so it is really a Laurent $Y$\~polynomial. It can be readily deduced 
from the fact that all anti-symmetric polynomials in $\mathscr{Y}$ are 
divisible by the {\em discriminant}, the common denominator on the 
right-hand side. The proof of this theorem will be given in the 
next section.

In the case of $A_{1}$, we obtain
\begin{align}
\varphi_{m}=&\frac{t^{-m/2}}{1+t^{-1}}
\left(\frac{Y^{m}-t^{-1}Y^{m-2}-Y^{-m-2}+t^{-1}Y^{-m}}
{1-Y^{-2}}\right)\notag\\
=&\frac{t^{-m/2}}{1+t^{-1}}
\left(\frac{(Y^{m+1}-Y^{-m-1})-t^{-1}(Y^{m-1}-Y^{1-m})}
{Y-Y^{-1}}\right),\label{vphiformula}
\end{align}
which matches our calculations above based on the Pieri rules. 
Compare with the ``nonsymmetric" formulas (\ref{vepformulas}). 
Macdonald established his formula by calculating the Satake 
$\mathfrak{p}$\~adic integral 
representing the spherical function (see below).

One can try to use the Pieri rules to justify the theorem, but as we 
noted above, reasonably simple explicit formulas exist only for 
$A_n$ and in some cases of small ranks. There is another, much more 
direct approach (any root systems), which can be 
generalized to DAHA theory. We will switch to it
after the following remarks clarifying the $\mathfrak{p}$\~adic
origins of the Pieri rules, to be continued in Section 
\ref{sect:S-Mtheory} on the classical $\mathfrak{p}$\~adic theory 
of spherical functions.

\subsubsection{\sf Comments on Pieri rules}
Formulas (\ref{pieri}) match the classical arithmetical definition 
of the (one-dimensional) Hecke operator. Let $t$ be the cardinality 
of the residue field of a $\mathfrak{p}$\~adic field $K$ ($t=p$ for 
$\Q_p$). The {\em Bruhat-Tits building} of type $A_1$ is a 
{\em tree} with $t+1$ edges from each vertex; 
the {\em vertices} $\{v\}$ 
correspond to the maximal parahoric subgroups of $G=PGL_2(K)$, 
which are (all) 
conjugated to $U=PGL_2(\o)\subset G=PGL_2(K)$ for the ring 
of integers $\o\subset K$. 
Two vertices are connected by an {\em edge} if their 
intersection is an Iwahori subgroup, i.e., is conjugated to
$B=\{g\in U\,\mid\, g_{21}\in \mathfrak{p}\}$ for
the maximal ideal $\mathfrak{p}\subset \o$. 
The group $G$ naturally acts on 
this tree by conjugation. Identifying the 
vertices with the cosets of 
$G/U$, the action of $G$ becomes left regular. 

Let $d(v)$ be the 
distance (in the tree) of the vertex $v$ from the origin $o$, 
which corresponds to $U$. 
The functions $f(m)$ on this tree depending only on the 
distance $m=d(v)\ge 0$ are exactly the functions on 
$G/\!/\,U=U\backslash G/U$. The figure is as follows ($t=p=3$): 
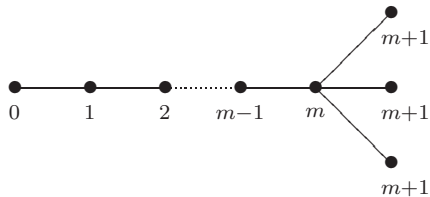
\begin{equation*}
\xy <1cm,0cm>: 
(0,1)*=0{\bullet}="-" ; 
(1,1)*=0{\bullet}="*", **@{-};  
(2,1)*=0{\bullet}="*", **@{-}; 
(3,1)*=0{\bullet}="*", **@{.}; 
(5,1)*=0{\bullet}="*", **@{-}; 
(4,1)*=0{\bullet}="*", **@{-}; 
(5,0)*=0{\bullet}="*", **@{-}, (4,1); 
(5,2)*=0{\bullet}="*", **@{-}, (2,1); 
(0,0.6)*+{^0};
(1,0.6)*+{^1};
(2,0.6)*+{^2};
(3,0.6)*+{^{m-1}};
(4,0.6)*+{^m};
(5.2,1.6)*+{^{m+1}};
(5.2,-0.4)*+{^{m+1}};
(5.2,0.6)*+{^{m+1}};
\endxy 
\end{equation*}
  
The classical {\em Hecke operator} is the (radial) Laplace 
operator $\De$ on this tree, the averaging over 
the neighbors. Explicitly, 
$$
\De f(m)=\frac{tf(m+1)+f(m-1)}{t+1} \for m>0, \ \,
\De f (0)=\frac{f(1)}{t+1}. 
$$
Thus (\ref{pieri}) is exactly the eigenvalue problem for
$(t^{1/2}+t^{-1/2})\De$ with the eigenvalue $Y+Y^{-1}$,
where $Y$ is treated as a free parameter.

For arbitrary Chevalley groups,
a combinatorial definition of the 
Laplace-type operator and its higher analogs in terms of
the Bruhat-Tits buildings becomes involved. 
The case of $A_n$ 
was considered by Drinfeld. 

The Bruhat-Tits building is equally useful in the
theory of {\em Whittaker functions}.
There is a unique infinite 
path from the origin such that the elements of the 
{\em unipotent} subgroup $N\subset G$ preserve its 
{\em direction to infinity}; only the direction, any finite
number of vertices can be ignored. Let us extend this path to
a {\em road}, infinite in both directions. Then any vertex
can be mapped onto this road (identified with 
$N\backslash G/U$) using $N$; its image is unique.
The Whittaker function can be interpreted as a 
function on this road, nonzero only on the original (positive)
path; see Section \ref{sect:Whitf} below for more detail. 

\subsubsection{\sf The major limits}\label{sect:thelimits}
Let us switch from the normalization we used 
(compatible with the $\mathfrak{p}$\~adic Hecke
operators), to the one more convenient algebraically.
Namely, we set $\tilde{\varphi}_{m}\equal t^{m/2}\varphi_{m}$,
which readily simplifies the (symmetric) Pieri rules:
$$
(Y+Y^{-1})\tilde{\varphi}_{m}=
\tilde{\varphi}_{m+1}+\tilde{\varphi}_{m-1}.
$$
This recurrence has the following elementary solutions
for $m\ge 0$.

1) The monomial symmetric functions (divided by $2$):
$$\m_{m}=(Y^{m}+Y^{-m})/2.
$$ 

2) The classical Schur functions $\chi_{m}$: 
$$
\chi_{m}=\frac{Y^{m+1}-Y^{-m-1}}{Y-Y^{-1}}.
$$

3) The renormalized Macdonald spherical functions:
$$
\tilde{\varphi}_{m}=\frac{1}{1+t^{-1}}\cdot
\frac{Y^{m+1}-Y^{-m-1}-t^{-1}(Y^{m-1}-Y^{1-m})}{Y-Y^{-1}}.
$$

All three sequences begin with $1$ at $m=0$. They are different
due to the {\em boundary conditions\,} 
at $m=-1$\,:
$$
1)\,\m_{-1}=\m_1,\ \ \,2)\,\chi_{-1}=0,\ \ \, 
3)\,\tilde{\varphi}_{-1}=
\tilde{\varphi}_{1}t^{-1}.
$$
The first two cases are limits of the third one:
\begin{equation*}
\xymatrix{ 
-\chi_{m-2} 
& \tilde{\varphi}_{m}\ar[l]_{\quad t\to 0} 
\ar[d]_{t\to 1}\ar[r]^{t\to \infty} &\chi_{m}\\ 
& \m_{m}}.
\end{equation*}
The limit $t\to \infty$ is actually the degeneration of
the Macdonald spherical functions to the Whittaker
functions; see Section \ref{sect:Whitf}.

\subsubsection{\sf The nonsymmetric case}
The Matsumoto spherical functions are right $U$\~invariant
and left Iwahori\~invariant, so they can be naturally identified 
with the functions depending on the distances from the
origin $\,o\,$ in the following two halves of the Bruhat-Tits 
building\,:\\ 
($+$) the paths from $\,o\,$ through the
{\em nonaffine} neighbors of $\,o\,$ ($t$ of them),\\  
($-$) the paths from $\,o\,$
through the {\em affine} neighbor $\,\hat{o}\,$ of $\,o\,$
(only one). 

The elements of $B\subset G$ are exactly those
preserving $\,o\,$ and the edge between $\,o\,$ and $\,\hat{o}\,$.
We will measure the distance using negative numbers 
in the second half ($-$).  Then the functions on 
$B\backslash G/U$ become 
$f(m)$ for  $m\in \Z$,
where $m=d'(v)\in Z$ for the new distance (may be
negative). 

Check that $d'(v)$ is the only invariant of the vertex under 
the action of the Iwahori subgroup and
interpret combinatorially
formulas  (\ref{yepm},\ref{yep-m}) in terms of $m=d'(v)$.
 
Let us switch in (\ref{vepformulas}) to
$\tilde{\vep}_m=t^{|m|/2}\vep_m$. Then
\begin{align}\label{veptilde}
&\tilde{\vep}_m=t^{m/2}\vep_m=Y^m,\ 
\tilde{\vep}_{-m}=Y^{-m}+
(1-t^{-1})\frac{Y^{-m}-Y^{m}}
{Y^{-2}-1},
\end{align}
where $m\ge 0$. 
There is no dependence on $t$ for nonnegative indices
(so the corresponding limits are trivial). 
The graph of the limits for $-m\,(m>0)\,$ 
reads as follows:
\begin{equation*}
\xymatrix{ 
\infty 
& \tilde{\vep}_{-m}\ar[l]_{\quad t\to 0\quad} 
\ar[d]_{t\to 1}\ar[r]^{t\to \infty} &\chi_m\\ 
& Y^{-m}}.
\end{equation*}

\subsubsection{\sf Proof of Macdonald's formula}
Recall that the affine Hecke algebra $\mathcal{H}$
in the $T$-$Y$\~presentation
is generated by the elements $T_{1}, \ldots, T_{n}$ and 
$Y_{b}$ for $b\in P^{\vee}$. The defining relations between
$T_{i}$'s and $Y_{b}$'s are:
\begin{align}\label{tyhecke}
&T_{i}^{-1}Y_{b}T_{i}^{-1}=Y_{b}Y^{-1}_{\al_{i}},\ 
\text{ if }(b, \al_{i})=1,\\
&T_{i}Y_{b}=Y_{b}T_{i}, \text{ if }(b, \al_{i})=0,\ i>0.
\label{t-yhecke}
\end{align}
The connection with the original definition is as follows: 
\begin{align*}
T_{0}=Y_{\theta}T_{s_{\theta}}^{-1},
\  \pi_{r}=Y_{\omega^{\vee}_{r}}T_{u_{r}}^{-1},
\end{align*}
where $u_r$ are from (\ref{piromr}). 

Formulas  (\ref{tyhecke}),(\ref{t-yhecke}) are actually 
the relations of the orbifold braid group
of $\C^*/W$. Using the quadratic 
relations,
\begin{equation}{\label{lus}}
T_{i}Y_{b}-Y_{s_{i}(b)}T_{i}=(t^{1/2}-t^{-1/2})
\frac{Y_{s_{i}(b)}-Y_{b}}{Y^{-1}_{\al_{i}^{\vee}}-1},\, i>0.
\end{equation}
These formulas are due to Lusztig (see e.g.,
\cite{L}).

\begin{lemma}\label{Berns}
The center of the affine Hecke algebra is 
$$Z(\mathcal{H})=\mathscr{Y}^{W}=\mathbb{C}[Y_{b}]^{W}.$$
\end{lemma}
{\em Proof.}
By regarding both sides of \eqref{lus} as operators on 
$\mathscr{Y}\ni f(Y)$, we have
\begin{equation}\label{lusoper}
T_{i}(f)=t^{1/2}s_i(f)+(t^{1/2}-t^{-1/2})
\frac{s_{i}(f)-f}{Y^{-1}_{\al_{i}^{\vee}}-1}.
\end{equation}
Thus $T_{i}(f)=t^{1/2}f$ for all $i>0$ are equivalent to
the relations $s_{i}(f)=f$ 
for all $i>0$, which means that
$f\in \mathscr{Y}^{W}$.
\sq

\begin{theorem}[{\sf Operator Macdonald 
Formula}]{\label{MAC}}
Let $$\tilde{M}\equal
\prod_{\al\in R_+}\frac{1-t^{-1}Y_{\al^{\vee}}^{-1}}
{1-Y_{\al^{\vee}}^{-1}}.$$
Then we have the following identity of operators acting 
in $\mathscr{Y}$ 
\begin{equation}{\label{mac1}}
P(t^{-1})\mathscr{P}_{+}=(\sum_{w\in W}w)\circ \tilde{M}, 
\end{equation}
Using the definition of $\mathscr{P}_{+}$,
\begin{equation}{\label{mac2}}
\sum_{w\in W}T_{w}^{-1}t^{-l(w)/2}=(\sum_{w\in W}w)\circ 
\tilde{M}.
\end{equation}
Equivalently, (\ref{mac1}) holds in the (abstract)
algebra \,$\b$\, of operators  generated by $W\ni w$ and 
rational functions in terms of \,$\{Y_b\}$ subject to
the relations \ $w Y_b w^{-1}=Y_{w(b)}$\, ($w\in W,b\in P$).
\end{theorem}

{\em Proof.}
The equivalence of 
\eqref{mac1} and \eqref{mac2} is due to 
\begin{equation*}
\mathscr{P}_{+}=\frac{\sum_{w\in W}t^{l(w)/2}T_{w}}
{\sum_{w\in W}t^{l(w)}}
=\frac{\sum_{w\in W}t^{-l(w)/2}T^{-1}_{w}}
{\sum_{w\in W}t^{-l(w)}}.
\end{equation*}
Indeed, both operators are divisible by $1+t^{1/2}T_i$
on the right
and on the left for any $i>0$, 
and act identically on $1\in \mathscr{Y}$ (which provides
the exact normalization factors).

Following \cite{Ch10} (upon the affine degeneration), 
let us introduce the following involution acting
on the operators from the algebra $\b$, 
\begin{align}
\iota: Y_{b}\mapsto Y_{b},\ t^{1/2}\mapsto -t^{-1/2},
\ s_{i}\mapsto -s_{i}.
\end{align}
Applying this involution to the operator from (\ref{lusoper}), 
\begin{equation*}
T_{i}=t^{1/2}s_{i}+\frac{t^{1/2}-t^{-1/2}}
{Y_{\al_{i}^{\vee}}^{-1}-1}(s_{i}-1),
\end{equation*}
one readily obtains
\begin{equation*}
T_{i}^\iota\ =\ t^{-1/2}s_{i}-
\frac{t^{1/2}-t^{-1/2}}{Y_{\al_{i}^{\vee}}^{-1}-1}
(s_{i}+1).
\end{equation*}
 
The $q\to 0$ limit of the $\mu$\~function\, from the 
DAHA theory is 
$$
M\equal
\prod_{\al\in R_+}\frac{1-Y_{\al^{\vee}}^{-1}}
{1-tY_{\al^{\vee}}^{-1}}.
$$
This function is {\em equivalent}
($\leftrightarrows$)
to $\tilde{M}$ in the following sense:
they coincide up to a $W$\~invariant factor. Indeed,
$$
\tilde{M}\leftrightarrows
\tilde{M}'\equal\prod_{\al\in R_+}\frac{1-Y_{\al^{\vee}}}
{1-t^{-1}Y_{\al^{\vee}}}\leftrightarrows M .$$

\begin{lemma}\label{IOTAM}
$MT_{i}M^{-1}=T_{i}^{\iota}$ for $i=1,\ldots,n$
(see \cite{Ch10}). 
\sq\end{lemma}

\begin{lemma}\label{TLEFTRIGHT} For $i\ge 1$,

$T_{i}\,+t^{-1/2}=(s_{i}+1)\cdot F_i$ for a rational function 
$F_i(Y)$,

$T_{i}^{\iota}+t^{-1/2}=G_i\cdot(s_{i}+1)$ for 
a rational function $G_i(Y)$.
\sq\end{lemma}
 
Returning to the proof of the theorem, 
$\mathscr{P}_{+}\circ \tilde{M}^{-1}\leftrightarrows
\mathscr{P}_{+}\circ M^{-1}$, and
these operators are divisible by $(1+t^{1/2}T_i)$ on the
left and by $(1+t^{1/2}T_i^\iota)$ on
the right. The left divisibility is
straight from that of $\mathscr{P}_{+}$; the right
divisibility results from Lemma \ref{IOTAM}.

Using Lemma \ref{TLEFTRIGHT}, we obtain that 
$\mathscr{P}_{+}\circ \tilde{M}^{-1}$
is divisible on the right and on the left by
$(s_i+1)$. Thus it commutes with the operators of multiplication
by functions from $\mathscr{Y}^W$ and must be in the form 
$G(Y)\circ\sum_{w\in W}w$ for a $W$\~invariant (rational)
function $G(Y)$. Hence, $G=P(t^{-1})^{-1}$ due to
$\sum_{w\in W}w(\tilde{M})=P(t^{-1})$. The latter is an 
immediate corollary of the divisibility of 
antisymmetric Laurent polynomials by the discriminant;
see \cite{Bo} and \cite{Hu}, formula
(35), Section 3.20.   
\sq
\medskip

The operator Macdonald formula is actually from
\cite{Ma5}, formula (5.5.14).  We deduced this from
Lemma \ref{IOTAM}; Macdonald checks the divisibility
of the operator $(\sum_{w\in W}w)\circ \tilde{M}$ by $1+t^{1/2}T_i$
on the left and on the right directly. Then he 
equates the leading terms in (\ref{mac1}), the coefficients of 
the longest element $w_0\in W$. Note that his last step cannot 
be used in DAHA theory (the longest element does not exist in
$\hW$). We think that the interpretation of $M$ and 
$\mu$ from \cite{Ch10} as intertwiners 
between the symmetric and antisymmetric 
polynomial representations 
clarifies well their appearance in this context. 

\subsection{{\bf Satake-Macdonald theory}}\label{sect:S-Mtheory}
\subsubsection{\sf Chevalley groups}\label{sect:Cheva}
Let $K$ be a $\mathfrak{p}$\~adic field and $\o\subset K$
the valuation ring in $K$ with the (unique) prime ideal
$\mathfrak{p}=(\varpi)$ for the uniformizing element $\varpi$. 
We set $t=|k|$, where $k$ is the residue field $\o/(\varpi)$.

For an irreducible reduced root system $R$ 
as above and the coweight lattice $P^{\vee}$,
the Lie algebra $\mathfrak{g}_{K}$ is
defined as the $\mathfrak{g}\otimes K$
for the Lie algebra $\mathfrak{g}$ defined
over $\Z$ as the span of
$\{x_{\al}, h_{b}\}$ for $\al\in R,
b\in P^{\vee}$ subject to the relations 
\begin{align*}
&[h_a,h_b]=0,\ [h_{b},x_{\al}]=(b, \al)x_{\al},
\ [x_{\al}, x_{-\al}]=h_{\al^\vee}, \\
&[x_{\al},x_{\be}]=N_{\al,\be}x_{\al+\be}\
\text{ if }\al+\be\in R,\, \text{ otherwise } 0.
\end{align*}
Accordingly, $\mathfrak{g}_{\o}=\mathfrak{g}\otimes \o$.
The integers $N_{\al,\be}$ can be chosen here uniquely up to
signs; we will omit their discussion.

The unipotent groups 
$X_{\al}$ are defined for $\al\in R$ as
``exponents" of $Kx_{\al}$; 
$H$ is the $K$\~torus corresponding to $P^\vee$.
By construction, these groups act on $\mathfrak{g}_K$.
We will also need the group lattice formed by the elements 
$\varpi^b\in H$ for $b\in P^\vee$ defined 
as follows:
$$
\varpi^{b}(x_{\al})=\varpi^{(b,\al)}x_{\al},\, \forall\, \al\in R.
$$
Finally, the (split) {\em Chevalley group} 
$G$ is the span of $X_\al$ 
for all $\al\in R$ and $H$. The standard {\em unipotent
subgroup} $N$ is the group span of $X_\al$ for
$\al\in R_+$. The {\em maximal parahoric subgroup}
$U$ is the centralizer of $\mathfrak{g}_{\o}$ in $G$.
Note that $P^\vee$ is used here; if it is replaced by
$Q^\vee$, then the corresponding group is the group of 
$K$\~points of the connected
simply connected split algebraic group associated with $R$. 

We have the Cartan decomposition of $G$
\begin{equation}
G=UH_{+}U=\bigcup_{b\in P^{\vee}_{+}} U\varpi^{b}U,
\end{equation}
and the Iwasawa decomposition
\begin{equation}{\label{Iwasawa}}
G=UHN=\bigcup_{b\in P^{\vee}} U\varpi^{b}N\,;
\end{equation}
the unions are disjoint.

As an exercise, introduce the Chevalley group corresponding
to the {\em twisted affinization} $\tR^\nu$ of $R$ considered
in Section \ref{sect:thinaff}. Using algebraic groups,
it will be a group of $K$\~points of a 
nonsplit group over $K$, which splits
over certain ramified extension of $K$.  

\subsubsection{\sf The Satake integral}
Let $L(G,U)$ be the space of complex valued
functions $f$ on $G$, compactly supported,
satisfying the bi-$U$-invariance condition:
\begin{equation*}
f(u_{1}x u_{2})=f(x) \text{ for all } x\in G, \text{ and any }
u_{1}, u_{2}\in U.
\end{equation*}

This is a ring; the product of two functions 
$f,g\in L(G,U)$ is 
defined by the {\em convolution\,}
\begin{equation*}
f\ast g(x)=\int_{G}f(xy^{-1})g(y)dy, 
\end{equation*}
where $dy$ is the Haar measure on $G$ normalized 
by $\int_{U}dy=1$. Moreover, it is a commutative
ring (use the ``$-1$"\~automorphism of $R$ and $R^\vee$ extended
to $G$). 

The {\em zonal spherical function} on $G$ relative to $U$ is a
continuous bi-$U$-invariant 
complex-valued function $\Phi$ on $G$ 
satisfying the following condition:
\begin{align} \label{sphdef}
&\Phi\ast f =c_{f}\Phi \hbox{\ \, for\ any\ \,} f\in L(G,U), 
\end{align}  
and for the constants $c_{f}$ depending on $f$. In other words, 
$\Phi$ is a common eigenfunction of all the convolution operators
with the elements $f\in L(G,U)$; then $c_{f}$ are the corresponding
eigenvalues. The normalization is $\Phi(1)=1$.

Satake (following Harish-Chandra)
found that an {\em arbitrary} zonal spherical function
can be uniformly described in terms of the vector 
$\la\in \C\otimes_{\Z}P\cong \C^n$.
Using the Iwasawa decomposition \eqref{Iwasawa},
let us define the projection map onto $P^\vee$  
\begin{align}
&\mathrm{pr}:G\ \to\  P^{\vee}, \ 
x\in U\varpi^{b}N\ \mapsto \ b.
\end{align}
Using this map, the zonal spherical functions are given as follows:
\begin{equation}\label{satake}
\Phi_{\la}(x)=\int_{U}t^{(\mathrm{pr} (x^{-1}u),\rho-\la)}du
\end{equation}
for the Haar measure restricted to $U$.

Macdonald calculated this integral in \cite{Ma1} using the 
combinatorics of $U$. This was not too simple; see his
Madras lectures \cite{Ma2}
(the lectures also include relations to the real theory,
positivity matters and other issues).
It suffices to evaluate $\Phi_{\la}$ at 
$\varpi^b$. His formula reads as 
\begin{align}\label{macdformula}
&\Phi_\la (\varpi^b)=\frac{1}{P(t^{-1})}
\sum_{w\in W}t^{(b,w(\la)-\rho)}\,\prod_{\al\in R_+}
\frac{1-t^{-1-(\al^\vee,w(\la))}}{1-t^{-(\al^\vee,w(\la))}}.
\end{align}
Connecting $\mathfrak{p}$\~adic theory
and our algebraic approach can be achieved by 
replacing $Y_{b}$ by $t^{(b,\la)}$, namely,
\begin{align*}
&\vph_b(Y)\ =\ \Phi_{\la}(\varpi^b)\,\bigl
[t^{(b,\la)}\mapsto Y_b
\bigr].
\end{align*}
Recall that in (\ref{macdvarphi}), 
\begin{eqnarray*}
\varphi_{b}(Y)=\frac{t^{-(\rho,b)}}{P(t^{-1})}
\sum_{w\in W}Y_{w(b)}\prod_{\alpha\in R_+}
\frac{1-t^{-1}Y_{w(\alpha^{\vee})}^{-1}}
{1-Y_{w(\alpha^{\vee})}^{-1}}.
\end{eqnarray*}

\subsubsection{\sf The universality principle}
The approach via the Matsumoto spherical
functions establishes a bridge between
the algebraic theory above and the $\mathfrak{p}$\~adic theory,
and {\em proves} 
(\ref{macdformula}) without taking a single 
$\mathfrak{p}$\~adic integral.

The coincidence of these two theories, algebraic and
$\mathfrak{p}$\~adic,  can be also seen
by observing that the defining relations 
from (\ref{sphdef}) are nothing but the Pieri rules in the 
algebraic theory. However this is with the reservation that the 
(symmetric) Pieri rules are generally not explicit.  

One can also use the following {\em universality principle}.

Formula (\ref{sphdef}) ensures that there
exists a family of {\em pairwise commutative} 
{\em difference} operators in terms of $b$; they are 
convolutions with different $f\in L(G,U)$. 
It is not necessary to know exactly how the convolution is defined;
it can be of any origin, say, from geometric theories.
Provided there exist such operators (differential or
difference) and certain natural {\em symmetries}, such
a family is essentially unique for a given root system. 
 This claim can be made  
rigorous if more information on the structure of
difference or differential operators under consideration
is available. 

The key point is that we have very few such families
in mathematics  
(subject to certain symmetries and boundary conditions). 
Cf. the discussion in Section \ref{sect:thelimits}.
Major examples come from
the theory of Macdonald polynomials and DAHA, from
their counterparts, generalizations and degenerations.
In physics, the same phenomenon is the universality of the 
quantum many body problem.

Thus, one can expect {\em a priori}
(or even conclude rigorously) that
$\mathfrak{p}$\~adic spherical functions must be proper
specializations of the Macdonald polynomials. In our case, 
specialization of the general $q,t$\~theory is by letting 
$q\to 0$ and under minor renormalizations. The justification of
this connection is straightforward if the algebraic approach 
via the Matsumoto functions is used. However, it is not obvious
at all if the spherical functions and the operators are defined
$\mathfrak{p}$\~adically, via the convolution on $G$.

\subsubsection{\sf Whittaker functions}\label{sect:Whitf}
The universality  principle discussed above  
works equally well for the Whittaker functions. 
We introduce them following
\cite{CS} with some simplifications;
see also \cite{Shi} for the $GL_n$\~case. The notation
is from Section \ref{sect:Cheva}. 

The theory of $q$\~Whittaker functions will be discussed
in the second part of this work, including the nonsymmetric
(spinor) functions. A natural challenge is to define
the Matsumoto-type (``nonsymmetric") p-adic Whittaker functions
(which can be only spinor ones); their 
definition is outlined below.

The unramified $\mathfrak{p}$\~adic Whittaker function 
$\w$ is introduced for an additive character $\psi$, the
product of the ($K$\~additive) characters 
$\psi_i: K\to K/\o \to \C^*\ (i=1,\ldots,n)$; 
each $\psi_i$ must be nontrivial on
$\varpi^{-1}\o/\o$. This can be naturally extended to a 
character of the group $N$ (vanishing on $X_\al$ for nonsimple
roots $\al>0$). 

For an algebra homomorphism 
$\chi:\,L(G,U)\to \C$, there is a unique function
$\w_{\chi}$ on $G$ such that $\w_{\chi}(1)=1$,
\begin{align} \label{whit}
&\w_\chi(ngu)=\psi(n)\,\w_\chi(g) 
\for n\in N,\, u\in U,\, g\in G, \\
&\and \w_\chi\ast f =\chi(f)\,\w \hbox{\ \, for\ any\ \,} f\in L(G,U).
\notag 
\end{align}  

Similar to the spherical function $\Phi$,
it suffices to know the values $\w_{\chi}(\varpi^b)$
for $b\in P^\vee$. However, 
$\w_{\chi}(\varpi^b)$ is {\em not} a $W$\~invariant function
of $b$. Moreover,  $\w_{\chi}(\varpi^b)=0$ unless
$b\in P^\vee_+$ (anti-dominant in Lemma 5.1 from \cite{CS}). 

The universality principle 
is actually sufficient to
conclude/expect that, up to a certain renormalization, 
$\w_{\chi}(\varpi^b)$ {\em does not depend on $t$} (a surprising
fact!) and that it 
is a classical finite-dimensional character
of the Langlands dual group of $G$. Here the corresponding
dominant weight is $b$ and $\chi$ must be treated as the argument. 
See Theorem 5.4 from \cite{CS} and \cite{Shi} for the 
precise statements. 

The fact that $\w_{\chi}(\varpi^b)$ vanishes for 
$b\not\in P^\vee_+$ is the key here. 
 It provides the boundary condition sufficient to identify 
the Whittaker functions with the characters (practically
without calculations).
Cf. Section \ref{sect:thelimits}, case (2).
A counterpart of this property in the theory of real and 
complex Whittaker functions is a certain decay condition; see
\cite{ChW} for the $q$\~Whittaker functions.

Let us demonstrate the mechanism of this vanishing property
in the case of $GL_2(K)$. 
Using the first relation from the definition of $\w=\w_\chi$,

\begin{align*}
&\psi(\varpi^{-1})
\w\Biggl(\left(\begin{array}{cc}
\varpi^n & 0 \\0 & \varpi^{n+1} \\ \end{array}\right)\Biggr)
=\w\Biggl(
\left(\begin{array}{cc}
1& \varpi^{-1} \\0 & 1 \\ \end{array}\right)
\left(\begin{array}{cc}
\varpi^n & 0 \\0 & \varpi^{n+1} \\ \end{array}\right)
\Biggr)\\
&=\w\Biggl(
\left(\begin{array}{cc}
\varpi^n& \varpi^n \\0 & \varpi^{n+1} \\ \end{array}\right)
\Biggr)
=\w\Biggl(
\left(\begin{array}{cc}
\varpi^n& \varpi^n \\0 & \varpi^{n+1} \\ \end{array}\right)
\left(\begin{array}{cc}
1& -1 \\0 & 1  \\ \end{array}\right)
\Biggr)\\
&=\w\Biggl(\left(\begin{array}{cc}
\varpi^n & 0 \\0 & \varpi^{n+1} \\ \end{array}\right)\Biggr)
=0 \hbox{\ \, due to \ \,} \psi(\varpi^{-1})\neq 1.
\end{align*}
At the level of formulas,  
$\w_{\chi}(\varpi^b)$ for $\la=\chi$ is the
limit $t\to \infty$ of the $\mathfrak{p}$\~adic spherical 
function from (\ref{macdformula}); 
see Section \ref{sect:thelimits}
for the demonstration in the $A_1$\~case.
\smallskip

Generalizing (any root systems),
we claim that the $\mathfrak{p}$\~adic
Whittaker functions can be obtained as limits
of the properly normalized spherical functions 
when the cardinality of the residue field $k$ tends 
to  $\infty$. I.e., we replace the starting 
$\mathfrak{p}$\~adic field by (the completion of) 
its maximal unramified extension; the limiting
procedure can be correctly defined. It results in 
the switch from the affine Hecke algebra to the
affine nil-Hecke algebra. The Matsumoto functions
go to new {\em spinor-Whittaker functions} in this limit.
\smallskip

Let us make this explicit for $A_1$.
The quadratic relation becomes $T(T-1)=0$
in such a limit. Correspondingly,
$T^{-1}$ in the formulas must be replaced by 
$T'\equal T-1$. For instance, the relation $TYT=Y^{-1}$
now becomes 
$T' Y=Y^{-1} T$; more generally,
\begin{equation*}
T Y^n-Y^{-n}T=
\frac{Y^{-n}-Y^n}{Y^{-2}-1} \for n\in \Z.
\end{equation*}
Cf. Section \ref{sec:rank1} above.

The definition
of the Matsu\-mo\-to- Whittaker function
remains $T_{\hw}\mathscr{P}_+$ for $\hw\in \hW$ 
and for the symmetrizer $\mathscr{P}_+$, which 
is now simply $T$ (for $A_1$). Following (\ref{vepformulas}),
they must be expressed in terms of $Y^{\pm 1}$.
Setting $\psi_{-n}=Y^n T$ for $n\ge 0$,
the nil-counterpart of  $(T\pi)^n \mathscr{P}_+\,$  is
\begin{align*}
&\psi_n\, =\, T Y^n T= Y^{-n} T+ \frac{Y^{-n}-Y^n}{Y^{-2}-1}T
=(\sum_{m=0}^n Y^{n-2m})T \for n\ge 0.
\end{align*}
Thus $\psi_n=\{\,Y^{|n|}\hbox{ for } n\le 0,\ 
(Y^{n+1}-Y^{-n-1})/(Y-Y^{-1})\hbox{ for } n> 0\,\}.$

The identities 
$T\psi_{-n}=\psi_n=T\psi_{n}\, (n\ge 0)$ are directly connected
with the theory of the second part of this work. The symmetrization 
of the spinor Whittaker function (here applying $T$)
must be the {\em diagonal spinor} (under the symmetry
$n\to -n$) constructed from the ``symmetric" Whittaker function.

The connection to the {\em spinor $q$\~Whittaker function}
from the second part of this work is direct. Namely, it is the 
limit $t\to 0$ where $\La$ is replaced by $Y$. Recall
that $t$, the cardinality of the residue field,
changes to  $t^{-1}$ in the $q,t$\~theory.  
The theory of the spinor $q$\~Whittaker functions
for arbitrary root systems is in progress, including
the $\mathfrak{p}$\~adic applications.
 
\setcounter{equation}{0}
\section{\sc{Double affine generalizations}}{\label{sect:daha}}
\subsection{{\bf Double affine Hecke algebra}}
We continue to use the notations from Section \ref{sect:affineweyl}.
Let $\hat{P}=\{\hat{a}=[a,j]\,|\,a\in P, j\in \mathbb{Z}\}\subset 
\mathbb{R}^{n}\times \mathbb{R}$ be the {\em affine weight
lattice}. Correspondingly, let $X_{[a,j]}\equal X_a q^j$ for 
pairwise commutative $X_a$ ($X_{a+b}=X_aX_b$) and 
a parameter $q$ 
(later, a nonzero number). Setting $X_j=X_{\om_j}$ for  
$j=1,\ldots,n$ (they are algebraically independent)\,:
$$X_{a}=\prod_{j=1}^{n}X_{j}^{l_j},
\hbox{\ \, where\,\ } l_j=(a,\al_j^\vee) \text{\ \,due to \ \,}  
a=\sum_{j=1}^{n}l_{j}\omega_{j}.
$$

Recall the definition of the action of the extended affine 
Weyl group $\hat{W}=W\ltimes P^{\vee}$  in
$\R^{n+1}$:
$$
b[z,\xi]=[z, \xi-(b,z)]\, (b\in P^{\vee}),\  
w[z,\xi]=[w(z),\xi]\, (w\in W).
$$ 
Accordingly, we set $\hat{w}(X_{\hat{a}})\equal
X_{\hat{w}(\hat{a})}.$

This action is dual to  the {\em standard affine action}
of $\hat{W}\ni \hw$ in $\R^n\ni x$ via the 
translations defined
as $wb(x)=w(x+b)$ for $w\in W,\, b\in P^\vee$. In the 
space of functions of $x$, this reads as $\hw(f)(x)=f(\hw^{-1}(x))$ 
(notice the sign). Applying $\hw=wb\in \hW$ to
$\, X_a\equal q^{x_a}\,$ for $x_a\equal (x,a)$, one has
\begin{align}\label{xbaction}
&\hw(X_a)=q^{(w^{-1}x-b,a)}=q^{(x,w(a)-(b,a))}=X_{[w(a),-(b,a)]}=
X_{\hw(a)}.
\end{align}

The {\em double affine Hecke algebra} (DAHA), denoted by
$\HH$, is 
defined over the ring of constants 
$\Z[q^{\pm 1/m},t^{\pm 1/2}]$ 
for $m\in Z_+$ such that $(P,P^\vee)= \frac{1}{m}\Z$.
In this paper we will mainly consider DAHA over the field
$\C_{q,t}\equal \C(q^{1/m},t^{1/2}).$  This algebra is 
generated by the affine Hecke algebra $\h=
\lan T_{i}, i=0,\ldots,n, \Pi\ran$ defined above and
pairwise commutative elements $\{X_{a}, a\in P\}$ 
subject to the following {\em cross-relations}:
\begin{align}
&T_{i}X_{a}T_{i}=X_{a}X_{\al_{i}}^{-1} \text{\ if\ } 
(a,\al_{i}^{\vee})=1,
\notag\\
&T_{i}X_{a}=X_{a}T_{i} \text{\ if\ } 
(a, \al_{i}^{\vee})=0,\label{dahadef}
\\
&\pi_{r}X_{b}\pi_{r}^{-1}=X_{\pi_{r}(b)},\notag
\end{align}
where $r\in O$ is from the orbit $O$ of $\al_{0}$ in  
$Dyn^{\tt aff}$; see (\ref{heckerel}).

Recall that the $Y_{b}$ for $b\in P^\vee$ from (\ref{tyhecke})  
satisfy the {\em dual} cross-relations:
\begin{align*}
&T_{i}Y_{b}T_{i}=Y_{b}Y_{\al_{i}^\vee}^{-1},\  
\text{ if }(b,\al_{i})=1,\\
&T_{i}Y_{b}=Y_{b}T_{i}, 
\text{ if } (b, \al_{i})=0. 
\end{align*}
Using $Y_b$ instead of $\{\pi_r,T_0\}$,
$\HH=\lan X_{a}\,(a\in P),\, Y_{b}\,(b\in P^\vee), 
T_{1}, \ldots, T_{n}\ran$.

\subsubsection{\sf The PBW Theorem}
An important fact is the PBW Theorem (actually,
there are $6$ of them depending on the ordering of $X,T,Y$):
\begin{theorem}[{\sf PBW for DAHA}]\label{PBWDAHA}
Every element in $\HH$ can be uniquely written in the form
\begin{equation}\label{pbwcdaha}
\sum_{a,w,b}C_{a,w,b}\,X_{a}T_{w}Y_b \text{\ \, for \,} 
C_{a,w,b}\in \C_{q,t},\ a\in P,\ w\in W,\ b\in P^\vee. 
\end{equation}
\vskip -1cm\sq
\end{theorem} 

The theorem readily results in the definition of the {\em polynomial 
representation} of $\HH$ in  
$\mathscr{X}\equal\C_{q, t}[X_b]=\C_{q,t}[X_{\om_i}]$; 
the ring $\Z[q^{\pm 1/m},t^{\pm 1/2}]$  is sufficient
in its definition. 
Using Theorem \ref{PBWDAHA}, we can  
identify $\mathscr{X}$ with the induced representation
Ind$_\h^{\HH}\,\C_+$, where $\C_+$ is the one-dimensional
module of $\h$ such that $T_{\hw}\mapsto t^{l(\hw)/2}$.
 
The generators $X_{b}$ act by multiplication; 
$T_i (i\ge 0)$ and $\pi_r (r\in O^*)$ 
act in $\mathscr{X}$ as follows:
\begin{eqnarray}\label{pitpolyn}
&\pi_{r}\mapsto \pi_{r},\ \, 
T_{i}\mapsto t^{1/2}s_{i}+\dfrac{t^{1/2}-t^{-1/2}}
{X_{\al_{i}}-1}(s_{i}-1).
\end{eqnarray}
Here 
$s_{0}(X_{b})=X_{b}X_{\theta}^{-(b, \theta)}q^{(b, \theta)}$.

\rmk
If one begins with formulas (\ref{pitpolyn}),
then the DAHA relations 
for these operators are not difficult to check directly. 
This approach gives the PBW Theorem for $\HH$
(the polynomial representation is faithful if $q$ is not a 
root of unity). In the affine case the
deduction of the PBW Theorem from the (nonaffine) formulas 
(\ref{pitpolyn}), checked directly, is actually due to Lusztig 
(in one of his first papers on AHA). 
Kato interpreted these formulas as those in
Ind$_{\H}^\h\C_+$ for nonaffine $\H$ and the 
plus-representation $\C_+$ (but then you need to
use the PBW Theorem). 

In the DAHA case the best way to obtain the PBW Theorem
is by defining the representation $\mathscr{X}$ 
via the formulas from (\ref{pitpolyn})and checking that
it is faithful for generic $q$.
There is no problem to order 
$X,Y,T$ as in (\ref{pbwcdaha}) for any $q,t\in \C^*$
using the DAHA relations, so the polynomial representation 
for generic $q$ (when this representation is faithful) 
provides the uniqueness of such 
expansions (which is the key) for all $q$.
\smallskip

\subsubsection{\sf The mu-functions}
We set
\begin{equation}\label{mutildemu}
\mu(X;q,t)=\prod_{\tilde{\al}>0}\frac{1-X_{\tilde{\al}}}
{1-tX_{\tilde{\al}}},\ \ \tilde{\mu}(X;q,t)=
\prod_{\tilde{\al}>0}\frac{1-t^{-1}X_{\tilde{\al}}}
{1-X_{\tilde{\al}}}.
\end{equation}
Following Section \ref{sect:length},
\begin{align}\label{Lahw}
&\La(\hw)\,\equal\,\tilde{R}_{+}\cap \hat{w}^{-1}(\tilde{R}_{-})
=\{\tal>0\,|\,\hw(\tal)<0\} \for  \hat{w}\in \hat{W}
\end{align}
consists of $l(\hat{w})$ positive roots.
The following are the key relations
for the functions $\mu,\tilde{\mu}$:
\begin{align}\label{murelations}
&\frac{\hw^{-1}(\mu)}{\mu}=
\frac{\hw^{-1}(\tilde{\mu})}{\tilde{\mu}}=
\prod_{\tilde{\al}\in \Lambda(\hat{w})}
\frac{1-t^{-1}X_{\tilde{\al}}^{-1}}{1-X_{\tilde{\al}}^{-1}}\cdot
\frac{1-X_{\tilde{\al}}}{1-t^{-1}X_{\tilde{\al}}}\\
=&\prod_{\tilde{\al}\in \Lambda(\hat{w})}
\frac{1-t^{-1}X_{\tilde{\al}}^{-1}}{1-t^{-1}X_{\tilde{\al}}}\cdot
\frac{1-X_{\tilde{\al}}}{\ \,1-X_{\tilde{\al}}^{-1}}= 
\prod_{\tilde{\al}\in \Lambda(\hat{w})}
\frac{t^{-1}-X_{\tilde{\al}}}{1-t^{-1}X_{\tilde{\al}}}.\notag
\end{align}
We see that $\mu/\tilde{\mu}$ is (formally)
a $\hW$\~invariant function.
Note that both functions, $\mu$ and $\tilde{\mu}$, are
invariant under the action of $\Pi=\{\pi_r,r\in O\}$.

We will need the formula for the constant term ct$(t)$ 
of $\mu$ (the coefficient of $X^0$):
\begin{align}\label{consterm}
&\hbox{ct}(t)=\prod_{\al\in R_{+}}\prod_{i=1}^{\infty}
\frac{(1-t^{(\al,\rho^{\vee})}q^{i})^{2}}
{(1-t^{(\al,\rho^{\vee})+1}q^{i})(1-t^{(\al,\rho^{\vee})-1}q^{i})}.
\end{align}
It will be treated as an element in $\mathbb{C}[t][[q]]$;
we will use this formula mainly for $t^{-1}$ instead of $t$.

\subsection{{\bf Affine symmetrizers}}
\subsubsection{\sf The hat-symmetrizers}
Let us introduce formally the infinite counterpart of
the $P$\~symmetrizer as follows:
\begin{equation}\label{defhatP}
\hat{\mathscr{P}}_{+}
=\sum_{\hat{w}\in \hat{W}}t^{-l(\hat{w})/2}
T^{-1}_{\hat{w}}/\hat{P}(t^{-1}) \for 
\hat{P}(t)=\sum_{\hat{w}\in \hat{W}}
t^{l(\hat{w})},
\end{equation}
the affine Poincar\'e series, which is a rational function of $t$.
Here and below $\hat{P}(t^{-1})^{-1}$ is expanded 
with respect to $t^{-1}$. We also set

$$
\hat{\mathscr{P}}\,'_+\equal
\sum_{\hat{w}\in \hat{W}}t^{-l(\hat{w})/2}
T^{-1}_{\hat{w}},\ \
\hat{\mathscr{S}}\,'_+\equal\sum_{\hat{w}\in \hat{W}}\hat{w},\ \
\hat{\mathscr{I}}\equal \hat{\mathscr{S}}\,'_+
\circ \tilde{\mu}\,.
$$
All constructions below can be extended to the
{\em minus-symmetrizers} (generally, to arbitrary
characters of the affine Hecke algebra), but we will
stick to the plus-case in this paper.

We understand these operators in this paper mainly
(but not always) as follows. 
Let us move all $\hw\in \hW$ in the series for
$\hat{\mathscr{S}}\,'_+
\circ \tilde{\mu}$ to the right and expand the coefficients
in terms of $X_{\al_i}$ for $i=0,\ldots,n$. 
Such expansions will contain only nonnegative powers of $q$.
Similarly, $t^{-l(\hu)/2}T_{\hu}^{-1}$ are understood 
as operators in the  polynomial representation,
where we move all $\hw$ to the right. The resulting coefficients 
will be {\em infinite} sums in terms of $X_{\al_i}\,(i\ge 0)$
by construction, to be analyzed in the next theorem,
which extends Theorem \ref{MAC} to the affine case.

\subsubsection{\sf The kernel and the image}
\begin{theorem}\label{KERIMAGE}
(i) The coefficients of $\hw$ in the above representations
of $\hat{\mathscr{S}}\,'_+
\circ \tilde{\mu}$ and  $\hat{\mathscr{P}}\,'_+$ 
will contain only nonpositive powers of $t$. These coefficients
are well defined as formal series in terms of $X_{\al_i}$ for
$i\ge 0$ and $t^{-1}$.
Moreover, provided that $|q|<1$ and $|t|>1$, the coefficients
of individual $X_a \hw\,(a\in Q\subset P, \hw\in \hW)\,$ 
will converge as series in terms of $q,t^{-1}$. 
 
(ii) Letting $\a=\hat{\mathscr{P}}\,'_+$ or 
$\a=\hat{\mathscr{S}}\,'_+\circ \tilde{\mu}$,
the following annihilation properties hold:
\begin{align}\label{annihil}
(\hw-1)\a\ =\ 0\ =\ &(t^{-\frac{l(\hw)}{2}}T_{\hw}-1)\a \notag\\
\ =\ 0\ =\ &\a(t^{-\frac{l(\hw)}{2}}T_{\hw}-1).
\end{align}
The products in (\ref{annihil}) must be 
transformed in the same way as $\hat{\mathscr{S}}\,'_+
\circ \tilde{\mu}$ and  $\hat{\mathscr{P}}\,'_+$. Namely,
all $\{T_{\hw}^{\pm 1}\}$ must be expressed
via $\{\hw\}$ using (\ref{pitpolyn})\,; then all $\hw$ 
must be moved to the right and, finally, the resulting 
coefficients of $\hw$  must be expanded as series from 
$\Z[[t^{-1/2},X_{\al_i},i\ge 0]]$. 

(iii) The right multiplication by 
\, $(t^{-\frac{l(\hw)}{2}}T_{\hw}-1)$\,
is well defined for any series $\c=\sum_{\hu}\, C_{\hu}\, \hu$
with the coefficients in $\Z[[t^{-1/2},X_{\al_i},i\ge 0]]$
or its localization by $t$. Namely, given $\hw\in \hW$,
\begin{align*}
&\c\,(t^{-\frac{l(\hw)}{2}}T_{\hw}-1)\ =\ 
\sum_{\hu,\hv} \,C_{\hu} B^{\hu}_{\hv}\,\hu \hv \for \\
&\hu \,(t^{-\frac{l(\hw)}{2}}T_{\hw}-1)\ =\ 
\sum_{\hv}\, B^{\hu}_{\hv}\,\hu \hv,\ \,
B^{\hu}_{\hv}\in \Z[[t^{-1/2},X_{\al_i},i\ge 0]],
\end{align*}
where $\hv$ are taken from 
the (finite) Bruhat set of the element $\hw$.
\end{theorem}
{\em Proof.}
To check $(i)$ for 
$\hat{\mathscr{S}}\,'_+\circ \tilde{\mu}$, let us divide it
by $\tilde{\mu}$ on the left. 
Then, using (\ref{murelations}),
\begin{eqnarray}\label{tmudivided}
\tilde{\mu}^{-1}\circ\hat{\mathscr{S}}\,'_{+}\circ\tilde{\mu}
=\sum_{\hat{w}\in \hat{W}}\prod_{\tilde{\al}\in 
\Lambda(\hat{w})}
\frac{t^{-1}-X_{\tilde{\al}}}{1-t^{-1}X_{\tilde{\al}}}\circ
\hat{w}^{-1},
\end{eqnarray}
which can be readily expanded in terms of $t^{-1}$.
Multiplying (\ref{tmudivided}) by the expansion of
$\tilde{\mu}$ in terms of $X_{\al_i}$ for $i\ge 0$, 
we obtain the required.
 
Only the nonnegative powers of $t^{-1}$ appear in the 
expressions of $t^{-l(\hw)/2}T_{\hw}^{-1}$ and
$\hat{\mathscr{P}}\,'_+$. Indeed,
using (\ref{pitpolyn}),
\begin{align*}
t^{-1/2}T_{i}^{-1}\ &=\ 
t^{-1/2}(t^{-1/2}s_{i}+\dfrac{t^{-1/2}-t^{1/2}}
{X_{\al_{i}}^{-1}-1}(s_{i}-1))\\
&=\ t^{-1}s_{i}+\dfrac{(t^{-1}-1)X_{\al_i}}
{1-X_{\al_{i}}}(s_{i}-1).
\end{align*}
The $\hw$\~coefficients
of $\hat{\mathscr{P}}\,'_+$ 
are infinite sums, well defined 
due to part $(e)$ of Lemma \ref{LEMThatw} below. 

We note that the operators $\hat{\mathscr{S}}\,'_+
\circ \tilde{\mu}$ and  $\hat{\mathscr{P}}\,'_+$
will be used later in concrete spaces; then
their coefficients will be treated as 
(meromorphic) functions of $X,q,t$. 

The convergence  of the coefficients of 
$\hat{\mathscr{P}}\,'_+$ subject to $|q|<1<|t|$ 
is part of Theorem \ref{PSBULLET}. 
It can be also obtained from Theorem \ref{LEVZERO}; 
see an outline of its proof in 
Section \ref{sect: Coefficient}.
The sharp estimate 
is actually $|t|>q^{1/h}$ (see below). 
\smallskip

Let $\iota$ be the involution, not an anti-involution, in 
$\mathscr{X}$ or acting in a proper localization of $\HH$ given by
\begin{equation*}
\iota: s_{i}\mapsto -s_{i}\, (i\ge 0),\, \pi_r\to\pi_r,\,
\ X_{a}\mapsto X_{a},\, q\mapsto q,
\, t^{1/2}\mapsto -t^{-1/2}. 
\end{equation*}

We have the following two lemmas extending 
the corresponding nonaffine Lemmas \ref{IOTAM}
and \ref{TLEFTRIGHT} (used for verifying the 
Macdonald formula).

\begin{lemma}\label{IOTLEM}
$\mu T_{i}\mu^{-1}=T_{i}^{\iota}$,\, for\, $i=0, \ldots, n$\,
(see \cite{Ch10}). \sq 
\end{lemma}
\begin{lemma} For $i\geq 0$,

$t^{1/2}T_{i}+1=(s_{i}+1)\cdot F_i$ for a rational function 
$F_i$,

$t^{1/2}T_{i}^{\iota}+1=G_i\cdot(s_{i}+1)$ for 
a rational function $G_i$. \sq
\end{lemma}
Note that
the automorphism $\,H\mapsto \mu^\iota H^\iota (\mu^\iota)^{-1}\,$
acts trivially on the element 
$T_i(i\ge 0),\,X_a,\,Y_a,\,q$, changing only $t$.

These lemmas are sufficient to establish $(ii)$. Claim $(iii)$
is straightforward.
\sq

\subsubsection{\sf Employing the E-polynomials}
From now on we will frequently represent $t$ in the 
form $t=q^k$. Given $a\in P$,
$u_{a}$ will be the element of
minimal possible length in $W$ such that 
$u_{a}(a)\in P_{-}$. We set
\begin{align}\label{pibdef}
&a_-\equal u_a(a)\in P_-,\ \, \pi_a\equal au_a^{-1}.
\end{align}
Here $l(\pi_a w)=l(\pi_a)+l(w)$ for an arbitrary
$w\in W$, which is the defining property of $\{\pi_a\}$.

The Macdonald polynomials $E_{a}$, $a\in P$ are $Y$\~eigenvectors:
\begin{align}\label{epolyno}
&Y_{b}^{-1}(E_{a})\,=\,q^{(b, a_{\#})}E_{a}, \ b\in P^\vee,\ 
a_{\#}=a-k u_{a}^{-1}(\rho), 
\end{align}
which fix them uniquely up to proportionality 
for generic $k$. The standard normalization condition
is $E_a=X_a+$(lower terms); see books \cite{Ma4,C101}.  
Note that $u_0=$id and $0_\#=-k\rho$. More generally, $u_{a}=\,$id for
$a\in P_{-}$ and 
$Y_{b}^{-1}(E_{a})=q^{(b, a-k \rho)}E_{a}$ for such $a$ and
any $b\in P^\vee$. 
\smallskip

These polynomials
were introduced by Heckman and Opdam in the differential
setting, then by Macdonald for  $t=q^k$ for integers $k$ 
and then in \cite{C4} in complete generality (in the reduced
case).  They are orthogonal Laurent polynomials with respect to 
the inner product 
$$
\hbox{Constant Term\,}(fg^*\mu) \for f,g\in \mathscr{X},\ 
q^*=q^{-1},t^*=t^{-1}, X_b^*=X_b^{-1}.
$$
See  \cite{Ma4,C101} and also \cite{OS}; the latter contains 
historic remarks and references including the important 
$C^\vee C$\~case, which we do not discuss here.
The symmetric Macdonald polynomials 
for the classical root systems were defined (and used)
for the first time by Kevin Kadell.
\smallskip

Among quite a few properties of the $E$\~polynomials, let 
us mention the nonsymmetric
Macdonald conjectures, namely, the norm-formula, the
duality-evaluation formula and the Pieri rules. 
They are now established in an entirely
conceptual way (see \cite{C101}
and \cite{C103}); these properties can be deduced from
the self-duality of DAHA practically without calculations. 

In a sense the duality claim
is the starting (and the simplest) in this 
chain of properties and the constant term formula is 
the endpoint. The nonsymmetric Pieri rules
do not belong to the standard list of Macdonald's conjectures, 
but they are the key to connect the duality with the evaluation 
and norm formulas. We note that their proof in \cite{C4} goes 
via the reduction to the roots of unity. 
\smallskip

The symmetric (usual) Macdonald conjectures can be deduced
from the nonsymmetric ones or can be obtained directly from the
DAHA theory upon symmetrization. The key feature 
of the nonsymmetric theory, which has no symmetric counterpart,
is the technique of intertwiners. It simplifies dealing with the 
$E$\~polynomials significantly vs. the symmetric theory (the
$P$\~polynomials).

We note that \cite{C101} and other works
of the first author  are mainly written for
the {\em twisted affinization} $\tR^\nu$
(in the reduced case). A natural notation
is $\HH(\tR^\nu;\tR^\nu)$, which 
means that the $X$\~generators and $Y$\~generators
are labeled by same lattice $P$. Then the $\HH$ from this 
paper must be denoted by $\HH(\tR;\tilde{R^\vee})$.

The technique of intertwiners can be 
transferred to $\HH(\tR;\tilde{R^\vee})$
(which is the setting of this paper). The norm and evaluation
formulas for $\tR^\vee$ hold for 
$\HH(\tR;\tilde{R^\vee})$ upon natural modifications 
of the formulas. For instance, the evaluation
formula for $E_{a}(t^{-\rho^\vee})$ can be obtained from
the one in \cite{C101} or from the Main Theorem of 
\cite{C4} (formula (5.4)) by the following transformations:

(a) adding check to $\rho$,\ 
(b) replacing $q_\al$ by $q$,\  and\ 
(c) setting $t_\al=q^{k_\al}$.

Explicitly, for $b\in P$,
\begin{align}
&E_{b}(t^{-\rho^\vee}) \ =\ t^{(\rho^\vee,\,b_-)}
\prod_{[\al,j]\in \La'(\pi_b)}
\Bigl(
\frac{
1- q^{j}t^{1+(\rho^\vee,\,\al)}
 }{
1- q^{j}t^{(\rho^\vee,\,\al)}
}
\Bigr),\where
\label{evaluform}\\
&\La'(\pi_b)\ =\
\{[\al,j]\ |\  [-\al,\nu_\al j]\in \La(\pi_b)\} \for 
\pi_b\equal bu_b^{-1},
\notag
\end{align}
and we use the elements $u_b,\pi_b$ from
(\ref{epolyno}),(\ref{pibdef}).
The same transformation must be performed with the 
norm-formula (5.5) from \cite{C4}.  

\rmk
We note that the DAHA of {\em untwisted} type
$\HH(\tR;\tilde{R^\vee})$ are expected
to satisfy the quantum Langlands duality (see \cite{C103}).
Trying to help the readers interested in this setting,
let us discuss briefly the changes with the key
DAHA\~automorphisms from \cite{C101} needed in the untwisted
case. The $\si$ from \cite{C101} (coinciding with $\om^{-1}$ 
from \cite{C4}) maps now  
$\HH(\tR;\tilde{R^\vee})$ to $\HH(\tilde{R^\vee},\tR)$.
The automorphism $\tau_+$ acts in the former,
$\tau_-$ in the latter. One has
$$\si \tau_+^{-1}=\tau_-\si,\ \si \tau_+=\tau_-^{-1}\si.$$
There are unsettled questions with the difference Mehta-Macdonald 
formulas from \cite{C5} in the untwisted case; they will be 
partially addressed when discussing the 
affine Hall functions of level one. \sq

\subsubsection{\sf Convergence at level zero}
Let us begin with the remark that the summation formula 
for ct$(t)$ from \cite{Ma3} was interpreted in \cite{C1} 
as the Jackson integration version
of the {\em constant term conjecture}. It was generalized there 
to the Jackson-type norm formulas for arbitrary $E$\~polynomials.
The relation of \cite{C1} to the present paper 
is direct; the definition of the Jackson integral of $f(X)$ from
\cite{C1} is nothing but 
$$
\hat{\mathscr{S}}\,'_+(\tilde{\mu} f(X))\,[X\mapsto q^\xi] \for
\xi\in \C^n;
$$ 
the vector $\xi$ (arbitrary) is called the origin of the
Jackson integral, which is a summation.
The following theorem is a particular case of 
the Jackson norm-formulas from \cite{C1}, Proposition 5.7.

\begin{theorem}\label{THMSE}
For $|q|<1$, $t=q^k$ and $a\in P$ such that $E_{a'}$ are 
well defined for all $a'\in W(a)$, the sums
$\hat{\mathscr{S}}\,'_{+} (\tilde{\mu}E_{a'})$ 
absolutely converge if and only if\,
$\Re(2k\rho+a_{+},\om_i)<0$ for all $i=1,\ldots,n$.
Here $\{a_-\}=W(a)\cap P_-$, 
$a_+=w_0(a_-)$ for the element $w_0$ of maximal length
in $W$, $\Re$ denotes the real part. Under this condition,
$\hat{\mathscr{S}}\,'_{+} (\tilde{\mu}E_{a'})=0$ for $a\neq 0$
and all $a'\in W(a)$. \sq
\end{theorem}
To give some examples, 
the (absolute) convergence range for $a=\rho=
\al_1+\al_2$ in the case of $A_2$ is $\{\Re k>-1/2\}$; 
it becomes $\{\Re k>-1/3\}$ for $a=\om_1=\om_1^\vee
=(2\al_1+\al_2)/3$.  

We continue to assume that $k$ is generic
(we will need this to employ the $E$\~polynomials). Considering 
generic $k$ in Theorem \ref{THMSE} and in a similar convergence
statement is sufficient for us. Indeed, the 
inequalities for $\Re k$ that provide
the convergence (in a given finite-dimensional 
subspace of $\mathscr{X}$) for all but finitely many special
$k$ hold automatically for such special values. The convergence
can be better at such special values, but no worse than at
generic $k$, which is sufficient in what will follow.

\begin{theorem}\label{LEVZERO}
The sum $\hat{\mathscr{P}}'_{+}(E_{a'})
=\sum_{\hat{w}\in \hat{W}}t^{-l(\hat{w})/2}
T_{\hat{w}}^{-1}(E_{a'})$
absolutely converges for any $a'\in W(a)$ if and only the following
its sub-sum converges absolutely:
$\sum_{b\in P^\vee_{+}}t^{-(\rho, b)}Y_{b}^{-1}(E_{a_-})$.
Using (\ref{epolyno}), this readily results in the same condition as 
from the previous theorem, namely,
$\Re(2k\rho+a_{+},\om_i)<0$ for all 
$i=1,\ldots,n$. Provided the convergence,
\begin{align}\label{propor}
&\hat{\mathscr{P}}\,'_+
=ct(t^{-1})\hat{\mathscr{S}}\,'_+
\circ \tilde{\mu} \hbox{\ \, as operators acting in\ \,} \mathscr{X},
\end{align} 
where ct$(t^{-1})$ is the constant term 
of $\mu(X;q,t^{-1})$:
$$\hbox{ct}(t^{-1})=\prod_{\al\in R_{+}}\prod_{i=1}^{\infty}
\frac{(1-t^{-(\al,\rho^{\vee})}q^{i})^{2}}
{(1-t^{-(\al,\rho^{\vee})-1}q^{i})
(1-t^{-(\al,\rho^{\vee})+1}q^{i})}
\in \mathbb{C}[t^{-1}][[q]].$$
\end{theorem}
{\em Proof.}
Let us begin with establishing the proportionality claim from
(\ref{propor}) assuming the convergence. 
Copying the affine case,
$\hat{\mathscr{P}}\,'_{+}\circ\,\tilde{\mu}^{-1}$ is divisible by 
$(t^{1/2}T_{i}+1)$ on the left and by $(t^{1/2}T_{i}^{\iota}+1)$
on the right. Hence it is divisible by $(s_{i}+1)$ on the left and
on the right. Therefore 
$$
\hat{\mathscr{P}}\,'_{+}\circ\,\tilde{\mu}^{-1}\ =\
G(X)\sum_{\hat{w}\in \hat{W}}\hat{w}\ =\ 
G(X)\cdot \hat{\mathscr{S}}\,'_+
$$
for a certain $\hW$\~invariant function $G(X)$.
Using \cite{Ma3},  $G=ct(t^{-1})$. 

More directly, we can check that $\hat{\mathscr{P}}\,'_{+}(E_a)=0$ for 
any $a\in P\setminus \{0\}$; combining this with
Theorem \ref{THMSE} we readily establish the required proportionality.

The operator $\hat{\mathscr{S}}\,'_+$ of course 
diverges in (the whole) $\mathscr{X}$, so we must apply the 
argument above as follows. Given $N\in \N$, 
formulas (\ref{annihil}) guarantee that
the images and the kernels of 
$\hat{\mathscr{P}}\,'_+$ and 
$\hat{\mathscr{S}}\,'_+\circ \tilde{\mu}$
coincide upon acting in the linear spaces 
$V_N=\oplus_{(\rho,a_+)<N} \C X_a$, provided that
$\Re k<0$ and $|\Re k|$ is sufficiently large
(depending on $N$). 
Thus these operators are proportional in every $V_N$ and the 
coefficient of proportionality (a constant) does not depend
on $N$.

The convergence analysis for
$\hat{\mathscr{P}}\,'_+$ in $\mathscr{X}$ is different from
that for $\hat{\mathscr{S}}\,'_+\circ \tilde{\mu}$.
First, it suffices to assume that $a\in P_-$, using the
standard relations between the polynomials $E_{a\,'}$ for 
$a'\in W(a)$. Second, we observe that
the convergence is the worst for  terms $Y_b^{-1}(E_a)$ 
with $b\in P_+^\vee$ and $a\in P_-$. Thus, we need to analyze
\begin{eqnarray*}
\sum_{b\in P^\vee_{+}}t^{-(\rho, b)}Y_{b}^{-1}(E_{a})
&=&\sum_{b\in P_{+}}q^{(b, a-2k\rho)}E_{a};
\end{eqnarray*}
this sum converges absolutely if and only if \,
$\Re(2k\rho+a_{+})\in \R_{>0}Q_+$.
The completion of this argument is based
on the following theorem. 

\subsubsection{\sf Y-formulas for P-hat}
\label{sect:RatPhat}
Recall that 
$\hat{\mathscr{P}}\,'_+$ is the plus-symmetrizer
without the exact projector normalization, i.e., without
the division by $\hat{P}(t^{-1})$. By
$P(t)$, we denote the {\em nonaffine} Poincar\'e 
polynomial. For a subset
$\I\subset \{1,2,\ldots,n\}$, the Poincar\'e
polynomial of the root subsystem $R_{\I}\subset R$ generated
by the simple roots $\{\al_i\,\mid\, i\in \I\}$ will be
denoted by $P_{\I}(t)$. It is $1$ if $\I=\emptyset$.

\begin{theorem}\label{P+FORMULA}
The symmetrizer $\hat{\mathscr{P}}\,'_+$
can be presented as the following summation
over all subsets $\I\subset \{1,2,\ldots,n\}$
including the empty set and $\I=\{1,\ldots,n\}$\,:
\begin{align}\label{hatPrat}
&\hat{\mathscr{P}}\,'_+ = 
P(t^{-1})\mathscr{P}_+\Bigl(\,
\sum_{\I}\frac{P(t)}{P_{\I}(t)}
\prod_{i\not\in \I}
\frac{\,t^{-(\om_i^\vee,\,\rho)}\,Y_{\om_i^\vee}^{-1}}
{1-t^{-(\om_i^\vee,\,\rho)}\,Y_{\om_i^\vee}^{-1}}\,\Bigr) 
\mathscr{P}_+,
\end{align}
which is understood coefficient-wise upon the expansion 
of the rational expressions in the products 
in terms of $t^{-1}$ (a set of identities in
$\h_Y$).
\end{theorem}
{\em Proof.} We employ the key property of
the elements $\pi_b$ from (\ref{pibdef}), namely,
the equality $l(\pi_b w)=l(\pi_b)+l(w)$
for any $w\in W$. Since $\pi_b=b u_b^{-1}$, one has
$\pi_b w=u_b^{-1}\,b_-\,w$.
The element $u=u_b$ can be arbitrary such that 
its length is minimal possible for a given $b=u^{-1}(b_-)$,
i.e., minimal in the coset $Z(b_-)u$ for the centralizer
$Z(b_-)$ of $b_-$ in $W$.
It results in (\ref{hatPrat}).\sq
\smallskip

Note that formula (\ref{hatPrat}) gives a {\em rational}
expression for the affine Poincar\'e series 
$\hat{P}(t^{-1})$ $=\hat{\mathscr{P}}\,'_+(1)$.
Provided that $\hat{P}(t^{-1})\neq 0$,
the theorem gives a universal map 
{\em onto} the space of $Y$\~{\em spherical vectors}
$$
\{\, v\,\mid\,
T_{\hw}(v)=t^{l(\hw)/2)}v \for \hw\in \hW\,\}, 
$$
which is applicable to $\HH\,$\~modules that are 
unions of {\em finite-dimensional}
$Y$\~invariant subspaces, including $\mathscr{X}$.
Theorem \ref{SYMRANK1}
can be readily extended to arbitrary one-dimensional
characters of $\h_Y$; the case of the 
{\em affine minus-symmetrizer}, corresponding to
$\{T_{\hw}\mapsto (-t^{-1/2})^{l(\hw)}\}$, is of importance.

The right-hand side of formula (\ref{hatPrat}) 
is a rational function and can be used as such without 
the $t^{-1}$\~expansion. However, one has to ensure 
that the denominators in (\ref{hatPrat}) are nonzero.
For instance, this formula can be used in the (whole)
polynomial representation $\mathscr{X}$ 
for $A_1$ with any $q,t$ unless $t^2\in q^{-1-\Z_+}$
and for $A_2$ unless $t^6\in q^{-1-\Z_+}$ or 
$t^3\in q^{1+\Z_+}$. It is under the assumption that 
$q$ is not a root of unity and $\hat{P}(t^{-1})\neq 0$.
At roots of unity, this formula can be applied only in certain 
quotients of $\mathscr{X}$.

Formula (\ref{hatPrat}) is the subject of 
Theorem \ref{SYMRANK1} in the case of $A_1$. 
For $A_2$, it reads as follows:
\begin{align*}
&\hat{\mathscr{P}}\,'_+\ =\  
P(t)P(t^{-1})\mathscr{P}_+\Bigl(\,
\frac{\,t^{-2}\,Y_{\om_1+\om_2}^{-1}}
{(1-t^{-1}\,Y_{\om_1}^{-1})(1-t^{-1}\,Y_{\om_2}^{-1})}\\
&+\,
\frac{1}{1+t}\,
\bigl(\frac{\,t^{-1}\,Y_{\om_1}^{-1}}
{1-t^{-1}\,Y_{\om_1}^{-1}}\,
+\frac{\,t^{-1}\,Y_{\om_2}^{-1}}
{1-t^{-1}\,Y_{\om_2}^{-1}}\bigr)\,
+\,\frac{1}{(1+t)(1+t+t^2)}\,
\Bigr)\mathscr{P}_+.
\end{align*}
Here $\rho=\al_1+\al_2$ and $(\rho,\om_i)=1$ for $i=1,2$;
$P(t)= (1+t)(1+t+t^2)$. Recall that $\om_i=\om_i^\vee.$
Applying this formula to 
$1\in \mathscr{X}$ and using that
$t^{-1}\,Y_{\om_i}^{-1}(1)$$=t^{-2}$, the resulting series
is the $t^{-1}$\~expansion of $\hat{P}(t^{-1})$; we arrive
at the formula $\hat{P}(t^{-1})=3(1-t^{-3})/(1-t^{-1})^3$. 
\smallskip

The expression on the right-hand side of
(\ref{hatPrat}) treated as an
element in the localization of affine Hecke subalgebra
$\h_Y=\lan T_{\hw},\hw\in \hW\ran$ must be {\em identically}
zero. Indeed, no affine symmetrizer exists in
$\h_Y$ or its localizations unless completions are
allowed. Similarly, this expression becomes identically zero
when applied in $\HH$\~modules that are unions of 
finite-dimensional $\h_Y$\~modules containing {\em no}
$Y$\~spherical vectors. 
This is the key point of the following theorem;
we mention that the $A_1$\~case is considered in full detail in 
Theorem \ref{SYMRANK1Y} below.

\begin{theorem}\label{YLEFT}
Given a set of representatives 
$\mathbf b=\{b^1,\ldots,b^p\}\subset P_+^\vee$ for 
the group $\Pi=P^\vee/Q^\vee$
(of cardinality $p$), let
\begin{align}\label{SiYgen}
&\tilde{\Sigma}_{\mathbf b}=
\prod_{\al\in R_+}\frac{(1-tY_{\al^\vee}^{-1})}
{(1-Y_{\al^\vee}^{-1})}\,
\frac{\sum_{j=1}^p
t^{-(b^j,\rho)}Y_{b^j}}
{\prod_{i=1}^n (1-tY_{\al_i^\vee}^{-1})}\,,\\
&\overline{\Sigma}_{\mathbf b}=\frac
{\prod_{\al\in R_+\setminus\{\al_1,\,
\ldots,\,\al_n\}}\,(1-t^{1-(\al^\vee,\rho)})}
{\prod_{\al>0}(1-t^{-(\al^\vee,\rho)})}
\sum_{j=1}^p\, t^{-(b^j,\rho)}Y_{b^j}.
\label{SiYgenbar}
\end{align}
We consider $\hat{\mathscr{P}}\,'_+$  as a 
standard formal series $\sum_{\hw}C_{\hw}\hw$ provided 
the convergence of the coefficients as formal series or
point-wise or as an operator acting
in any representations of $\h_Y$ where it is well defined. 
If $t$ is treated as a number, $\hat{P}(t^{-1})$ is 
supposed to be invertible.

Let $b^j\to\infty$, which means that $(b^j,\al_i)\to \infty$ for all 
$1\le j\le p,$\, $i>0$. We also assume that 
\begin{align}\label{tYvanish} 
&\lim_{b^j\to\infty}
t^{-(b^j,\rho)}Y_{w(b^j)}\mathscr{P}\,'_+=
\left\{\begin{array}{c}
 \hbox{\, exists\, for\, all \,\,} w\in W\\
 \hbox{\, equals zero for\ } w\neq \hbox{\sf id} 
\end{array}\right\}
\end{align} 
coefficient-wise in the standard $\hw$\~expansions 
(provided then that $|q|$ is sufficiently small if the 
coefficients are treated as meromorphic functions) 
or element-wise in a given $\h_Y$\~module. Then
\begin{align}\label{PSiY}
&\hat{\mathscr{P}}\,'_+\, =\, 
\lim_{\mathbf b\to\infty}
\tilde{\Sigma}_{b}\,\mathscr{P}\,'_+\,=\,
\lim_{\mathbf b\to\infty}
\overline{\Sigma}_{b}\,\mathscr{P}\,'_+\for
\mathscr{P}\,'_+\, =\, P(t^{-1})\,\mathscr{P}_{\,+}.
\end{align}
In the one-dimensional representation
of $\h_Y$ corresponding to ``$+$", (\ref{PSiY}) results  
in formula (5.9) from \cite{Ma3} for the affine 
Poincar\'e series 
$\hat{P}(t)$ in terms of the degrees $d_i$:
\begin{align}\label{rataffpoin}
\hat{P}(t) =\frac{|\Pi|}{(1-t)^{n}}
\,\prod_{i=1}^n\frac{1-t^{d_i}}{1-t^{d_i-1}},\hbox{\ where\ }
P(t)=\frac{\prod_{i=1}^n
(1-t^{d_i})}{(1-t)^n}.
\end{align}
\end{theorem} 
{\em Sketch of the proof.} 
\comment{
In the limit ${\mathbf b}\to\infty$,
the convergence of $\hat{\mathscr{P}}\,'_+$ implies that
\begin{align}\label{tYvanish} 
&t^{-l(a)/2}Y_a\to 0\for a\not\in P^\vee_+,
\where
l(a)=l(a_+)=2(a_+,\rho); 
\end{align} 
we omit the justification. 
}
Relation (\ref{tYvanish}) implies  
that $\overline{\Sigma}_{b}\,\mathscr{P}\,'_+\ $
from (\ref{PSiY}) converges to the affine symmetrizer
up to proportionality, i.e., satisfies the  
invariance properties upon multiplication by $T_{\hw}$ 
($\hw\in \hW$) on the right and on the left.
It is obvious when $T_{\hw}=Y_a(a\in P_+^\vee)$, which
is sufficient. 
Cf. Theorem \ref{SYMRANK1Y} below for $A_1$.

A straightforward
calculation of the coefficient of proportionality 
results in the first equality in (\ref{PSiY}). 
It readily gives that
$\tilde{\Si}_{\mathbf b}\mathscr{P}\,'_+$ and 
$\overline{\Si}_{\mathbf b}\mathscr{P}\,'_+$
must coincide in the limit provided the convergence
of the latter expression. 
Indeed, the multiplication or division by the ratio 
$(1-CY_{\al^\vee})/(1-Ct^{(\al^\vee,\rho)})$ 
will not change  $\tilde{\Si}_{\mathbf b}$  in the limit for
a sufficiently general constant $C$.
\smallskip

As noted above,
the first equality in (\ref{PSiY})
can be deduced directly from relation
(\ref{hatPrat}); let us outline the main steps.

We introduce the truncation
$\Upsilon_{\mathbf b}$ of the $Y$\~expression between the two
$\mathscr{P}_+$ in formula (\ref{hatPrat}) as follows.
Upon the $Y^{-1}$\~expansion, only the monomials $Y_a^{-1}$
subject to $b^j- Q_+\ni a\in P_+$ for $b^j=a\mod Q$ 
will be kept. Let $b^j=\sum_{i=1}^n r_i^j \al_i^\vee$;
recall that $(b^j,\al_i)\to \infty$, so 
the whole $\hat{\mathscr{P}}\,'_+\,$ will be obtained
in this limit. 
The {\em finite} sum  
\begin{align}\label{Upsilonb}
\Upsilon_{\mathbf b}^\flat\,\equal\,P(t)\sum_{j=1}^p
\prod_{i=1}^n
\frac{(1-t^{-1-r_i^j(\al_i^\vee,\,\rho)}\,
Y_{\al_i^\vee}^{-r_i^j-1})}
{(1-t^{-(\al_i^\vee,\,\rho)}\,Y_{\om_i^\vee}^{-1})}
\end{align}
contains all such $Y_a^{-1}$, i.e., contains $\Upsilon_{\mathbf b}$, 
but there will be extra 
(nondominant) terms there with $a\not\in P_+^\vee$. 

We are going now to use nonaffine formulas (\ref{mac1}) 
and (\ref{mac2})\,: 
\begin{equation}{\label{mac12}}
P(t^{-1})\mathscr{P}_{+}=(\sum_{w\in W}w)\circ \tilde{M}\for
\tilde{M}\equal
\prod_{\al\in R_+}\frac{1-t^{-1}Y_{\al^{\vee}}^{-1}}
{1-Y_{\al^{\vee}}^{-1}}.
\end{equation}

Due to these formulas combined with the vanishing
property from (\ref{tYvanish}), the contributions
of $Y_a^{-1}$ in (\ref{Upsilonb})
with $(a,\rho)\ll (b^j,\rho)$ for $a\in b^j-Q_+$ 
tend to zero in the limit. Thus the nondominant
terms can be disregarded in $\Upsilon_{\mathbf b}^\flat$.
Moreover, it suffices to consider only $\I=\emptyset$ in
Theorem \ref{P+FORMULA} in the limit
upon applying the operator from (\ref{mac12}).


Similarly, the numerator in formula (\ref{Upsilonb}) can be 
actually reduced to  $P(t)\,\bigl(1+(-1)^n\prod_{i=1}^n
t^{-1-r_i^j(\al_i^\vee,\,\rho)}\,Y_{\al_i^\vee}^{-r_i^j-1}\bigr)$.
Using that (\ref{hatPrat}) is zero in localizations of $\HH\,$,
$
\mathscr{P}_{+}\,\Upsilon^\flat_{\infty}\,
\mathscr{P}_+=0 \for \Upsilon^\flat_{\infty}
$
for $\Upsilon_{\infty}^{\flat}$ given by (\ref{Upsilonb}) upon
making the numerators $1$, i.e., by deleting the terms
that contain any $r_i^j$.
This identity can be obtained
directly from (\ref{mac12}); use the divisibility of 
the anti-invariant Laurent polynomials by the discriminant. 

This makes it possible to switch to  $Y_{b^j+a}^{-1}$ with 
$a\in Q_+^\vee$ in the limit; the terms here apart from
the initial truncation will not contribute to the limit.
Therefore $\Upsilon_{\mathbf b}^{\flat}$
can be replaced by
\begin{align*}
\Upsilon^\sharp_b\,\equal\,
P(t)&\,(-1)^n\,\sum_{j=1}^p\,t^{-(b^j,\rho)}\,Y_{b^j}^{-1}
\,\prod_{i=1}^n \frac
{t^{-(\al_i^\vee,\rho)}Y_{\al_i^\vee}^{-1}}
{(1-t^{-1}Y_{\al_i^\vee}^{-1})},\\
\hbox{and\ } & \hat{\mathscr{P}}\,'_+\ =\ 
\lim_{\mathbf b\to\infty}
\sum_{w\in W}\,
w\bigl(\Upsilon^\sharp_b\tilde{M}\bigr)\mathscr{P}_+.
\end{align*}

Using the vanishing condition from (\ref{tYvanish})
once again, we see that
only $w=w_0$ here really contributes to  
$\hat{\mathscr{P}}\,'_+$ in the limit.
Let us substitute $b\mapsto -w_0(b)\,$ in the resulting 
expression. Then
$\hat{\mathscr{P}}\,'_+$ becomes the limit of  
\begin{align*}
(-1)^n\,\sum_{j=1}^p\frac 
{t^{-(b^j+\rho^\vee,\rho)}Y_{b^j+\rho^\vee}}
{\prod_{i=1}^n(1-t^{-(\al_i^\vee,\rho)}Y_{\al_i^\vee})}
\prod_{\al\in R_+}\frac
{(1-t^{-1}Y_{\al^\vee})}
{(1-Y_{\al^\vee})}\,
P(t)\mathscr{P}_+,
\end{align*}
where we use that $\rho^\vee=\sum_{i=1}^n \om_i^\vee$.
Rewriting the latter formula in terms of $Y_{\al^\vee}^{-1}$,
we finalize (\ref{PSiY}).

Applying (\ref{PSiY}) to $1$ in the standard one-dimensional
representation of $\h_Y$, one arrives at
(\ref{rataffpoin}). Indeed, $Y_a$ become
$t^{(a,\rho)}$ upon this evaluation and $\mathscr{P}(1)=1$.
This formula is due to Matsumoto and Macdonald; see formula
(5.9) from \cite{Ma3}.
\sq
\smallskip

The conditions from (\ref{tYvanish}) 
hold coefficient-wise via the action of $Y_b$ in the 
polynomial representation followed by the standard 
expansion $Y_b= \sum _{\hw\in \hW} C_{\hw}\hw$ and in
the representations  $\mathscr{X}q^{\,lx^{2}/2}$
for $l>0$. See Theorem \ref{YLEFTNEW} below; 
the standard expansions of $Y_{w(b)}$ are discussed
there in detail.

Formula (\ref{SiYgen}) 
for $\overline{\Si}_b$ coincides with formula 
(\ref{YSiabsbar})
for $\overline{\Sigma}_M$ below in the case of $A_1$.
One needs to set 
$b^1=M\om, b^2=(M-1)\om$ for $\om=\om_1$.
We note that $b^1$ and
$b^2$ can be taken arbitrary (approaching infinity);
the $\hw$\~expansions of $Y_{b^j}$ in (\ref{YSiabsbar})
are for two disjoint sets of $\hw$, for $j=1$ and $j=2$.  

The vanishing condition from (\ref{tYvanish}) becomes
(\ref{YMzero}) for $A_1$ and always holds
provided the existence of $\hat{\Si}_\infty^+$ in
Theorem \ref{SYMRANK1Y}.

\rmk
In the Kac-Moody limit $t\to\infty$, (\ref{SiYgenbar})
combined with the proportionality claim from
(\ref{propor}) give a presentation of the Kac-Moody
characters as limits of the (affine)
Demazure characters. The latter are directly related to 
the operators 
$T^\infty_{\hw}=\lim_{t\to\infty}t^{-l(\hw)/2}T_{\hw}$. 
Namely, the corresponding
Demazure characters are proportional to  
$q^{-l\frac{x^2}{2}}\,
T^\infty_{\hw}(X_{-a}\,q^{l\frac{x^2}{2}})$ upon
the substitution $X_b\mapsto e^{-b}$.
Here $a$ are affine $l$\~dominant weights, i.e.,
$a\in P_+$ and $(a,\th)\le l$. For $\hw=b\in P_+$
as $b\to\infty$, they
approach $\prod_{i=1}^\infty
\frac{1}{(1-q^i)^n}\hat{\chi}_a^{(l)}$;
see (\ref{Kac-Moody}) below. 

Here it is not necessary to stick to the affine 
dominant weights $a$ of level $l$. One can define
the Kac-Moody characters formally for arbitrary $a\in P$
using the Kac-Weyl formula. 
The proportionality claim (\ref{propor})
itself provides that the Kac-Moody characters are {\em sums}
of properly normalized Demazure characters, which is connected
with the (infinite-dimensional) Demazure modules associated with
the opposite Borel subalgebra (to that used for the highest
vectors).

For arbitrary $t$, (\ref{SiYgenbar})
states that the corresponding affine Hall functions 
from (\ref{affineHall}) are limits of the 
{\em Demazure $t$\~characters} for $a\in P_+$
defined (formally) as $q^{-l\frac{x^2}{2}}\,
\overline{\Si} \mathscr{P}_+'(X_{-a}\,q^{l\frac{x^2}{2}})$,
where actually we do not need $\mathscr{P}_+'$ (see below).
The summation formula also holds and is equally important.

\smallskip

\subsubsection{\sf Coefficient-wise proportionality}
\label{sect: Coefficient}
Theorem \ref{LEVZERO} is sufficient to claim 
the existence of the coefficients of the operator 
$\hat{\mathscr{P}}\,'_+ $ as meromorphic functions 
and the coefficient-wise proportionality from
(\ref{propor}). We will outline
here an analytic version of this approach based on a natural
analytic extension of the polynomial representation.

\begin{theorem} \label{TQ1H}
Let $|t|>q^{1/h}$ for the Coxeter
number $h=(\th,\rho)+1$. Expanding 
$\hat{\mathscr{P}}\,'_+=\sum_{\hw\in \hW} F_{\hw}(X)\,\hw$,
the coefficients $F_{\hw}$ converge absolutely and to an
analytic function on any given compact subsets in 
$\{0\neq X_\al\not\in q^{\Z}, \al\in R\}$
for sufficiently small $|q|$ depending on this subset.
Moreover, $F_{\hw}$ coincide with the corresponding coefficients
of\, $ct(t^{-1})\,\hat{\mathscr{S}}\,'_+\circ \tilde{\mu}$
in this range; for instance, 
$F_{\hbox{\tiny id}}=ct(t^{-1})\,\tilde{\mu}(X;q,t)$.
\end{theorem}

The proof of this theorem, including the proportionality
claim and the sharp estimate of the radius of convergence
with respect to $t$ of the coefficients of $\hw$,
results from Theorem \ref{PSBULLET} below, 
based on the representations of $\HH$ in the space of 
delta functions. Also,
the existence of $\{F_{\hw}\}$ as meromorphic functions can
be obtained using direct estimates for the coefficients of operators 
$Y_b$\,; see Lemma \ref{LEMThatw} below and 
Theorem \ref{YEXPA} in the case of $A_1$.
Nevertheless, it is quite natural to try to deduce the convergence 
and proportionality directly from the properties of  
$\hat{\mathscr{P}}\,'_+$, considered as an operator
acting in the polynomial representation and its extensions. 
\smallskip

Let us outline here an  
approach to the coefficient-wise existence and the 
proportionality utilizing the following analytic modification of 
Theorem \ref{LEVZERO}. As a matter fact, the approach from
Theorem \ref{PSBULLET} (entirely algebraic) is very similar
to the following considerations.

We will assume in the sketch below that  $|t|>1$.
When dealing with the affine symmetrizers analytically,
it is convenient to replace $\mathscr{X}$ by the 
union of  Paley\~Wiener type
spaces $\mathscr{P\!W\!}_{M}(\u)$ of analytic functions 
in a given $\hW$\~invariant
domain $\R^n\subset \u\subset\C^n$. Here $M\in \Z_+$ and  
the growth condition is as follows:
$$
f(x)\in  \mathscr{P\!W\!}_{M}(\u)\,\Rightarrow\,
{}^{bw}f(x)<C_{x}(M)\, q^{-M(b_+,\rho)},\ b\in P^\vee,\, w\in W, 
$$
for a constant $C_{x}(M)$ continuously depending on $x\in \u$.
For $M=0$, this space includes $1$ and all  
$\hW$\~invariant functions analytic in $\u$, for 
instance, the images of $\hat{\mathscr{P}}\,'_+$ and 
$\hat{\mathscr{S}}\,'_+\circ \tilde{\mu}$. These two operators
act in $\mathscr{P\!W\!}_{M}(\u)$ for sufficiently large negative 
$\,\Re k\,$, depending on $M$, and for sufficiently small
$\u$ containing $\R^n$.

The kernels and images of these operators in $\HH$\~invariant
subspaces of $\cup_{M\ge 0}\mathscr{P\!W\!}_{M}(\u)$ coincide and
Theorem \ref{KERIMAGE} (in an analytic variant) implies the
proportionality 
\begin{align}\label{proporan}
&\hat{\mathscr{P}}\,'_+=
ct(t^{-1})\hat{\mathscr{S}}\,'_+\circ \tilde{\mu}
\hbox{\ \, provided the convergence.\,\ }
\end{align} 
\smallskip

To extract and then equate the coefficients 
of the operators under consideration, we need certain 
modifications of delta functions in 
the space $\mathscr{P\!W\!}_{0}(\u)$ in a sufficiently small
neighborhood $\u$ of $0\in\R^n$. Let 
$\a=\sum_{\hw\in\hW}F_{\hw}(X)\hw$, assuming that this operator
is convergent
with the coefficients analytic in $\u$ and satisfies
the conditions from (\ref{annihil}).
It suffices to know $\tilde{F}_b\equal \sum_{w\in W}F_{bw}(X)$
for $b\in P^\vee$;
expand $\a$ in terms of $bT_{w}$ for $\hw=bw$
to see it (use that $q,t$ are generic).

Let us extract from $\a$ 
the value of the coefficient $\tilde{F}_0$
at $x=0$. Recall the notation $X=q^x$, $x_\al=(\al,x)$. 
The following {\em probe function} from 
$\mathscr{P\!W\!}_{0}(\u)$ 
can be used, a substitute for the delta function at zero:\,  
$$
\zeta_N(x)=-\prod_{\al\in R_+}
\frac{(\exp({N\pi \imath x_\al)-\exp(-N\pi \imath x_\al}))^2}
{(\exp(N\pi x_\al)-\exp(-N\pi x_\al))^2},
$$ 
where $N\in \N$, $\imath^2=-1$.
This function is of order $1+O(|x|^2/N)$ near $x=0$ and of order
$O\bigl(\,|x-b|^2\cdot\frac{\exp(-CN)}{N}\,\bigr)$ for 
$x\approx b\in P^\vee\setminus 0$
for some constant $C>0$. Obviously,
$\a(\zeta_N)(x=0)=\tilde{F}_0(x=0)$, and we recover
the value of $\tilde{F}$ at $x=0$.

Using the function $\sum_w \zeta_N(w(x)-x_0)$ in the same manner,
we can find the values $\tilde{F}_0(x=x_0)$ for any given
$x_0$ in a sufficiently small neighborhood of $x=0$. 
This gives the function 
$\tilde{F}_0$ in $\u$ pointwise in terms of the action
of $\a$ in $\mathscr{P\!W\!}_{0}(\u)$.
Alternatively, 
recovering $\tilde{F}_0(x)$ for small $x$  can be achieved by
tending $N$ to $\infty$ (we will omit details).
 
The same approach can be used for extracting any $\tilde{F}_b$
from $\a$
upon applying the translations by $b\in P^\vee$ to the argument 
$x$ in the probe function (which fix its numerator).

This is of course based on the existence of $\a$ when applied
to $\zeta_N$ in a neighborhood of $x=0$.
The numerator of $\zeta_N$ is a 
{\em pseudo-constant}, a $\hW$\~invariant
function. Thus, the rate of convergence depends only on
the denominator and the convergence of the
operators $\hat{\mathscr{P}}\,'_+$
and $\hat{\mathscr{S}}\,'_+\circ \tilde{\mu}$\, 
applied to $\zeta_N$ is no worse than that 
for constants (or pseudo-constants). Actually, it 
is better than this; it holds for small {\em positive} $\Re k$ 
too (presumably, the inequality $\Re k<1/h$ is sufficient
here). 

As a matter of fact, we need to know here the convergence 
only for  
large negative $\Re k$ (for recovering the coefficients), 
a weaker fact. Indeed, the coefficients of
$\hat{\mathscr{P}}\,'_+$
and $\hat{\mathscr{S}}\,'_+\circ \tilde{\mu}$\, 
are meromorphic functions in $k$ (provided the convergence).
If the proportionality of these operators is known
for $\Re k \ll 0$, then it holds coefficient-wise.
Thus, the coefficient-wise existence and 
proportionality require only Theorem \ref{propor} extended 
analytically to the functions similar to $\zeta_N$; the 
proportionality factor will be automatically $ct(t^{-1})$. 
Theorem \ref{PSBULLET} below is an algebraic variant of this approach. 

\subsection{{\bf Affine Hall functions}}
\subsubsection{\sf Main definition}
The above considerations were for the 
$0$\~level case of the general theory of 
{\em affine Hall functions of arbitrary levels},
which will be the subject of this section. 
We continue to assume that $|q|<1$. 

Expressing  $X_{a}=q^{x_{a}}=q^{(x,a)}$,
let us introduce the {\em $l$\~Gaussian} as $q^{\,l\,x^2/2}$ for
$x^{2}\equal\sum_{i=1}^{n}x_{\omega_{i}}x_{\alpha_{i}^{\vee}}.$
In the case of $A_{2}$, for example, we have
$\al_{1}=\al_1^\vee=2\omega_{1}-\omega_{2}$, 
$\al_{2}=\al_2^\vee=2\omega_{2}-\omega_{1}$ and
$$\frac{x^{2}}{2}=
\frac{x_{1}(2x_{1}-x_{2})}{2}+\frac{x_{2}(2x_{2}-x_{1})}{2}=
x_{1}^{2}-x_{1}x_{2}+x_{2}^{2}.$$ 
One readily checks that
$$
\hw(q^{\,lx^2/2})=
q^{\,lb^2/2} X_{lb}^{-1} q^{\,lx^2/2} \for
\hw=bw,\, b\in P^\vee, w\in W. 
$$ 
These formulas
are actually the defining relations of the Gaussian
in what will follow. Recall that 
$\,bw(X_a)=q^{-(b,w(a))}X_a\,$ for $\,a\in P$.

To simplify notations, we set
\begin{align}\label{Jproj}
\hat{\mathscr{I}}\equal
\hat{\mathscr{S}}\,'_{+}\circ\tilde{\mu},\ 
\ \mathscr{I} \hbox{\ stays here for ``integration"}.
\end{align}

The {\em Hall functions of level $l>0$} are
defined as
\begin{align}\label{affineHall}
&H_{a}^{(l)}\equal\hat{\mathscr{I}}(X_{a}q^{\,lx^{2}/2}),
\ a\in P,\  \,\mathscr{H}_{l}\equal
\hat{\mathscr{I}}(\mathscr{X}q^{\,lx^{2}/2}).
\end{align}
Thanks to the presence of the Gaussian, 
$q^{\,-lx^2/2}\,H_a^{(l)}$
are absolutely convergent series
in terms of $X_b$ ($b\in P$) 
for all $x$ and $t$ (no poles due to the denominator
of $\tilde{\mu}$ will occur).
This is known and can be readily checked using 
$\hat{\mathscr{P}}_+^{\,'\,}$, which preserves the
Laurent polynomials. Indeed, 
the residues at (potential) poles of $H_a^{(l)}$
are meromorphic functions in terms of $q,t$; however they 
must vanish for sufficiently general $t$ due to (\ref{propor})
or (\ref{proporan}), the proportionality.
 
The absolute convergence actually holds here for any $l\in \C$ such
that $\Re l>0$, but then we will not be able to represent the
funcuions 
$q^{-\,lx^{2}/2}\,H_{a}^{(l)}$ as Laurent series. Also,
singularities in $x$ can appear for nonintegral 
$\,l\,$ at {\em nonreal 
poles} of $\tilde{\mu}(q^x)$, which are as follows:
\begin{align}\label{polesmu}
\bigl\{x\,\mid\, (x,\al)+j\in 2\pi \imath\, \log(q)\,
\{P^\vee\setminus 0\},\
[\al,j]\in \tR_+\bigr\},
\end{align}
where $\imath$ is the imaginary unit.
There will be no singularities in a sufficiently small
neighborhood of $\R^n\subset \C^n$ for nonintegral levels.


Note that for any $\hW$\~invariant function $f$, called a
{\em pseudo-constant},
\begin{equation}\label{PSonf}
\hat{\mathscr{P}}\,'_{+}(f)\,=\,\hat{P}(t^{-1})f\,=\,
\hbox{ct}(t^{-1})
\hat{\mathscr{I}}(f),
\end{equation}
where we need to assume that $\Re k<0$ to ensure 
convergence. Here  $\hat{P}(t)$ is the affine Poincar\'e
series. 

The coefficient of proportionality
is the same as in (\ref{propor}) because the 
action of our operators on any pseudo-constants $f$ is 
no different from the action on $1\in \mathscr{X}$.
For instance, (\ref{PSonf}) holds for functions
from $\mathscr{H}_l$ provided that $\Re k<0$.
\smallskip

\rmk
The proportionality from (\ref{propor}) cannot hold
for all $k$; otherwise 
$H_a^{(l)}$ would vanishes {\em identically} for all $a\in P$ at 
the poles of $\hbox{ct}(t^{-1})$, which is not the case.
For instance, $\mathscr{H}_l$ must be $\{0\}$ 
as $t=q^{1/h}$ for the Coxeter number $h$ if $\mathscr{P}\,'_+$
is well defined at this point, which happens only for $l=1$.
Indeed, the proportionality always holds when both operators are
well defined. 

We claim that for any (integral) $l>0$, the space
$\mathscr{H}_{l}$ is always smaller than
the corresponding Looijenga space (see the definition below)
at $t=q^{1/h}$ and at other zeros of 
$H_{a}^{(l=1)}$ from part $(ii)$ of the next 
Theorem \ref{HALLONE} (the simply-laced case).
However it is generally nonzero.
The justification of this and similar facts is based
on diminishing the level due to formula (\ref{PSonf}).

Numerical calculations of the space 
${\mathscr H}_{l}$ 
$=\hat{\mathscr{I}}(\mathscr{X})$
for $A_1,A_2,B_2$ show that this space is
really nonzero at $t=q^{1/h}$, i.e., that, generally,
$\hat{\mathscr{P}}\,'_+$  cannot be continued analytically
to $\Re k\ge 1/h$.
The latter inequality seems sharp for $l>1$, namely,
the convergence of $\hat{\mathscr{P}}\,'_+$
and (its corollary) the vanishing property 
$\mathscr{H}_l(k=1/h)=\{0\}$ 
are not expected to hold for $l=1\pm \vep$ for arbitrarily
small $\vep>0$. Only integral $l$ are considered
in this paper, but the definition of the corresponding 
spaces for any complex $l$ with $\Re l>0$ is straightforward.
\sq
\smallskip

\comment{
{\em DAHA-Matsumoto polynomials}.
Let us connect this definition with the 
approach from Section \ref{sec:Sph via Hall} based on
the Matsumoto (nonsymmetric) spherical
functions, which were introduced there in an 
entirely algebraic way. 
We define the {\em DAHA-Matsumoto polynomials}
$H\!\!\!E_{a}^{(l)}\in \mathscr{X} q^{\,lx^2/2}\,$
of level $l\ge 0$ by the relation:
\begin{align}\label{affineHallns}
&H\!\!\!E_{a}^{(l)}\hat{\mathscr{P}}_+
\ =\ \si^{-1}(T_{a})\,q^{\,lx^{2}/2}\,
\hat{\mathscr{P}}_+ \,\in\, 
\HH\, q^{\,lx^2/2}\hat{\mathscr{P}}_+ .
\end{align}
Here $\tau_{\pm}$ and $\si=\tau_+\tau_-^{-1}\tau_+$
are the standard automorphisms of DAHA; see \cite{C101}.
We will use that $\si^{-1}(Y_a)=X_a^{-1}$ for $a\in P$. 

Because of this definition,
we are supposed to add formally all Gaussians 
$q^{\,lx^2/2}\, (l\in \Z_+)$ to $\HH$ with the 
relations  $q^{\,lx^2/2} H=\tau_+^l(H) q^{\,lx^2/2}$ 
for $H\in \HH$, which allow moving  $q^{\,lx^2/2}$ 
to the right. However, such an extension is not really
necessary; the strict meaning of (\ref{affineHallns}) is
as follows.

First, we apply the DAHA-PBW theorem in terms of the 
basis $X_b\, (b\in P)$, $T_{w}\,(w\in W)$ 
and  $\tau_+^l(Y_b)\, (b\in P)$
to the element $\si^{-1}(T_{a})$. Second, we replace the elements  
$\tau_+^l(T_{\hw})$ (placed on the right) by $t^{l(\hw)/2}$;
the result will be a {\em polynomial} in terms of $X_b$. 
Then $H\!\!\!E_{a}^{(l)}$ will be this polynomial upon its
multiplication by $q^{\,lx^2/2}$. 

For example,
$H\!\!\!E_{a}^{(l)}=X_a^{-1} q^{\,lx^2/2}$ for $a\in P_+$, since
$T_a=Y_a$ for such $a$. It readily results in 
\begin{align}\label{hnonsym-sym}
\hat{\mathscr{P}}\,'_+ (H\!\!\!E_{a}^{(l)})\ =
\ ct(t^{-1}) H_a^{(l)}
\for a\in P_+
\end{align}
due to 
$\hat{\mathscr{P}}\,'=ct(t^{-1})\,\hat{\mathscr{I}}$.
It of course includes the analysis of the coefficient-wise
convergence of the operator $\hat{\mathscr{P}}\,'$ and when 
it acts in $\mathscr{X} q^{\,lx^2/2}$ from \cite{ChNT}
and this work.

We think that such an interpretation of our proportionality
theorem (it was not stated in \cite{ChNT}, v5) can be of
importance. 
A statement similar to (\ref{hnonsym-sym}) was known to 
Alexander Braverman (unpublished), with a reservation
that  only dominant weights $a$ subject
to $(a,\th)\le l$ appear in his approach. 
He did not define the nonsymmetric affine Hall polynomials and 
we do not know how he justified the proportionality theorem
(presumably, using the theory from \cite{ChNT}). The convergence of 
$\hat{\mathscr{P}}\,'$ is of course the key here.

We are very thankful to him for triggering (\ref{hnonsym-sym}).
As it was explained above, the extension of $\HH$ by the powers 
of the Gaussians is not necessary here and in (\ref{affineHallns}); 
the PBW theorem twisted by $\tau_+^l$ is the main algebraic 
ingredient of the construction.
}
 
\subsubsection{\sf Discussion, some references}
The formula for the affine Satake-type operator
$\hat{\mathscr{S}}\,'_+\circ \tilde{\mu}$ was considered
by several specialists as a ``natural" extension of
the Macdonald $\mathfrak{p}$\~formula, including 
certain geometric aspects and applications.

The main reference is \cite{Ka}; see also \cite{FGT,BK}.
Equivalent definitions of the affine Hall-Littlewood functions 
were suggested by several authors (not always published),
for instance, by Feigin and Grojnowski; let us also mention 
Garland's works.
\smallskip

Independently, the affine Hall functions of level one were 
explicitly calculated in \cite{C5} in the context of 
Jackson integrals (see also \cite{Sto}). The paper \cite{FGT}
contains an important interpretation of the affine Hall
functions via the Dolbeault cohomology of the 
{\em affine Grassmannian} and related flag
varieties. The appearance of the $ct(t^{-1})$   
in the formulas is interpreted there as the 
``failure of the Hodge decomposition." See also Section 
12.7 in \cite{FGT} concerning the level-one formulas. 

The definition of $\hat{\mathscr{P}}\,'_+$ is straight;
it belongs to a completion of the corresponding affine 
Hecke algebra. It becomes really interesting when acting
in DAHA modules; this theory is new.

Both operators,
$\hat{\mathscr{I}}\,=\,\hat{\mathscr{S}}\,'_+\circ \tilde{\mu}$\, 
and $\hat{\mathscr{P}}\,'_+$, are proportional whenever 
the operator $\hat{\mathscr{P}}\,'_+$ exists (see 
Theorem \ref{HALLONE}). They complement each other
in the following sense.

The convergence of the $\hat{\mathscr{I}}\,$ for $l>0$  
is better and much simpler to manage than that of 
$\hat{\mathscr{P}}\,'_+$. 
However, the latter operator acts naturally in 
DAHA modules and, importantly, does not 
require a priori knowledge of the $\mu$\~function;
for instance, this provides an alternative 
way to supply the polynomial and similar representations
with inner products. Accordingly, this operator
 has no singularities
(at the denominator of $\tilde{\mu}$).
Also, $\hat{\mathscr{P}}\,'_+\,$ is an exact DAHA-version of the 
classical Satake isomorphism in the AHA theory and it is 
closely connected
with the theory of Demazure characters. Let us comment on the latter. 

Under the limit $t\to \infty$, the operator
$\hat{\mathscr{I}}\,$ is directly connected with 
the {\em Weyl-Kac formula} for Kac-Moody characters;
the functions $ct(t^{-1})H_{-b}^{(l)}$ tend to the corresponding 
characters for the affine dominant weights $b$.  Theorem \ref{YLEFT}
generalizes the presentation
of the corresponding Kac-Moody character as an inductive limit
of the {\em Demazure characters}\,. 
The proportionality itself is an operator 
$t$\~variant of the presentation of the Kac-Moody characters as 
{\em sums} of properly normalized
Demazure characters associated with the Demazure modules
the Borel subalgebra opposite to the one used for the highest
vectors.
\smallskip


\subsubsection{\sf Proportionality for
\texorpdfstring{$l>0$}{level>0}}
Let us begin with the level-one case. 
Then we have a reasonably complete
theory from \cite{C5} (see also \cite{C101}) and 
\cite{Sto} devoted to the $C^\vee C$\~case.
Let us mention \cite{Vi}, where the level-one case is
addressed in the simply-laced case. Theorem 2 there
is a special case of Theorem 7.1 from
\cite{C5} (for simply-laced root systems).
The relation of Theorem 2 to the difference 
Mehta-Macdonald formulas from \cite{C5} 
{\em in the compact case} is discussed in \cite{Vi}.
The compact case is that based on the constant term 
inner product
(more generally, on the imaginary integration).
However, it is the {\em noncompact case}, 
namely the Jackson integration
formula from \cite{C5} (not mentioned in
\cite{Vi}), that is {\em directly} connected with the
affine Hall functions of level one. 

Works \cite{C5,C101} were written in the self-dual setting, 
i.e., for the {\em twisted} affine root system $\tR^\nu$, 
where the same lattice $P$ is used in $\hW$ and
for $X_a$ (and $E_a$).
Accordingly, the operator $T_0$ changes to the one with 
$\al_0=[-\vth,1]$ for the maximal {\em short} root $\vth$. 
Restricting ourselves to the simply-laced
case, the results from \cite{C5} on the Mehta-Macdonald
formulas in the context of Jackson integration can be
formulated as follows. Recall that $\al^\vee=\al$,
$\om_i^\vee=\om_i$ in this  case due to the normalization
$(\al,\al)=2$ for $\al\in R$.

\begin{theorem}\label{HALLONE}
Let $R$ be a simply-laced root system. We set
$\gamma(x)\equal
|W|^{-1}\,\sum_{\hat{w}\in \hat{W}}\hat{w}(q^{x^{2}/2})$
$=q^{x^2/2}\,\sum_{b\in P}\,X_{b}\,q^{b^{2}/2}$ for the
order $|W|$ of the nonaffine Weyl group $W$.
Let $X_b(q^a)\equal q^{(b,a)}$, $\hat{P}(t^{-1})$ 
is from (\ref{rataffpoin}). The level will be $l=1$.

(i) The series $\hat{\mathscr{P}}\,'_+$ considered
as an operator in $\mathscr{X}q^{x^2/2}$ converges  
element-wise for all $t\in \C^*$. The proportionality
relation 
$$
\hat{\mathscr{I}}\equal \hat{\mathscr{S}}\,'_+
\circ \tilde{\mu}\ =\ \ct(t^{-1})^{-1}\hat{\mathscr{P}}\,'_+
$$ 
holds for any $t\neq 0$ as well; cf. (\ref{propor}).

(ii) Assuming that $E_a$ is well defined, 
\begin{align}\label{levonehall}
&\hat{\mathscr{S}}\,'_{+}(\tilde{\mu}\,E_{a}q^{x^{2}/2})=
\frac{\hat{P}(t^{-1})}{\hbox{ct}(t^{-1})}
\hat{\mathscr{P}}_{+}(E_{a}q^{x^{2}/2})\\
&=E_{a}(q^{-k\rho})\,
q^{-a^{2}/2-k(a_{+},\rho)}
\cdot\prod_{\al\in R_{+}}
\prod_{j=0}^{\infty}\frac{1-t^{-1-(\rho,\al)}q^{j}}
{1-t^{-(\rho,\al)}q^{j}}\cdot \gamma(x).\notag
\end{align}

(iii) If $t$ is not a root of 
unity, then the linear map $\hat{\mathscr{P}}\,'_{+}$ is identically 
zero in $\mathscr{X}q^{x^{2}/2}$ if and only if $t^{m_i}=q^j$ for 
$j\in \N$ (for instance, for $t=q$). Here $\{m_1,m_2,\ldots,m_n\}$ 
are the exponents of $R$; $m_i=d_i-1$ for the degrees $\{d_i\}$. 
The map $\hat{\mathscr{I}}$ is identically 
zero on  $\mathscr{X}q^{x^{2}/2}$ if and only if $t^{d_i}=q^j$ 
for $j\in \N$ and
$j/d_{i}\not\in \N$ (for instance, this map vanishes identically at
$t=q^{1/h}$, where $h=(\th,\rho)+1$ is the 
Coxeter number). 
\end{theorem}

{\em Sketch of the proof.}
The existence of $\hat{\mathscr{P}}\,'_+$ for all $k\in \C$
and the corresponding extension of the proportionality from
$(\ref{propor})$  is due to the fact
that the image of this operator is one-dimensional
for generic $k$ and therefore proportional to $\ga(x)$.
The best way to proceed here is via the level-one variant
of Theorem \ref{PSBULLET}, namely, by considering
the inner product 
\begin{align*}
&\langle f,g\rangle_1\ =\ 
\bigl(\hat{\mathscr{P}}\,'_+(fg q^{x^2/2})\bigr)(id)\,.
\end{align*}

Paper \cite{C5} contains the formula for 
$\hat{\mathscr{S}}\,'_{+}(\tilde{\mu}\,E_{a}q^{x^{2}/2})$
from (\ref{levonehall}). 
To check $(iii)$, use the explicit formula for 
ct$(t^{-1})$ and the fact that all $E_a$ are well defined
with nonzero $E_a(q^{-k\rho})$ for positive $\Re k$. 
\sq

\rmk  
The levels $0$ and $1$ are exceptional
from the viewpoint of convergence.
For $l=0$, the convergence of both,  $\hat{\mathscr{I}}$ and 
$\hat{\mathscr{P}}\,'_+$, is (naturally) significantly
worse than the convergence in the presence of the Gaussian. For $l=1$,
$\hat{\mathscr{P}}\,'_+$ converges much better than for (integral)
$l>1$ due to
the fact that its image is one-dimensional. 
Recall that $\hat{\mathscr{I}}$ always converges for $l>0$.
\sq
 
\begin{theorem}\label{GENPROPR}
We continue to assume that
$R$ is simply-laced, but $l$ can be an arbitrary complex
number now such that $\Re l>0$. If $l\not\in \Z$, then we need
to avoid the {\sf\, nonreal
singularities\,} of the function $\tilde{\mu}(X;q,t)$;
see (\ref{mutildemu}) and (\ref{polesmu}). Restricting
the functions to a sufficiently small neighborhood of $x=0$
is sufficient. Considering $\hat{\mathscr{I}}$ and  
$\hat{\mathscr{P}}\,'_+$
as operators acting in the space $\mathscr{X}q^{\,lx^{2}/2}$,
the former operator converges absolutely element-wise 
for any $k$ and
the latter converges absolutely as $\Re k <1/h$ 
for the Coxeter number $h$. 
Under the condition $\Re k <1/h$, the proportionality
holds: $\ct(t^{-1})\hat{\mathscr{I}}\ =\ \hat{\mathscr{P}}\,'_+.$
\end{theorem}
{\em Proof.}
The convergence and proportionality here can be deduced 
from the corresponding coefficient-wise claims from Theorem 
\ref{TQ1H} in  Section \ref{sect: Coefficient}.
See also Lemma \ref{LEMThatw} concerning the convergence.
The estimates in $(e)$ there and the fact that 
the growth of the coefficients of $\hat{\mathscr{P}}\,'_+$
is no greater than exponential
are sufficient for the convergence due to the presence 
of the Gaussian.
Theorem \ref{YEXPA} below provides sharp estimates 
for the coefficients of $Y$\~operators in the case of $A_1$.

For $\Re k<0$, the absolute convergence of 
$\hat{\mathscr{P}}\,'_+$ 
and, therefore, the proportionality follow from  
the convergence of this operator in the space 
$\mathscr{PW}_0(\u)$ there.
The estimates from Theorem \ref{LEVZERO} (the level
zero case) can be almost directly used for such $k$ as well; 
the convergence will be no worse than it was for $a=0$ 
in this theorem. 
\sq

\rmk
Let us mention the 
symmetrizer $ \sum_{\hw}t^{l(\hw)/2}T_{\hw}$,
with $t,T$ instead of $t^{-1},T^{-1}$. 
Its convergence range in the space $\mathscr{X}q^{-lx^2/2}$
is $\Re k>-1/h$ (unless $l=0,1$), i.e., negating 
the range for $\hat{\mathscr{P}}\,'_+$ (acting in
$\mathscr{X}q^{+lx^2/2}$). We continue to assume that $|q|<1$.

This symmetrizer corresponds
to the theory of {\em imaginary integration}. Applying it to 
$\mathscr{X}q^{+lx^2/2}$ with positive $\Re l$ is possible too
provided that $\Re k>0$, however the result will be zero identically. 
\sq

\subsubsection{\sf The Looijenga spaces}{\label{sect:loo}}
For positive integral levels $l>0$, let us introduce the 
{\em Looijenga space} 
\begin{equation*}
\mathcal{L}_{l}=
\{\sum_{\hat{w}\in \hat{W}}\hat{w}(X_{a}q^{lx^{2}/2}),\ 
a\in P\}.
\end{equation*}

It can be identified with the space
Funct$\,(P/lP^\vee)^{\Pi W}$ formed by
the $\Pi W$\~invariant functions on the
set $P/lP^\vee$.
Recall that $P^\vee\subset P$ due to the normalization 
$(\th,\th)=2$. 
The action of $W$ is natural. The action of the group 
$\Pi=\{\pi_r=\om_r u_r^{-1}\,\mid\, r\in O\}$
is as follows.

Let us identify the space 
Funct$\,(P/lP^\vee)^{W}$  
with the space Funct$\,(\c_l)$, defined for the set
$\c_l\equal\{b\in P_+\,|\, (b,\th)\le l\}$.
The group $\Pi$ naturally acts on the 
{\,\em set\,} $\c_l$ through
its action on the {\em closed  fundamental affine Weyl 
chamber} $\{x\in \R_+\cdot P_+\,|\, (x,\th)\le 1\}$ ``multiplied" 
by $l$.  More algebraically,
we can identifying $\Pi$ with the group
$\{(l\om_r) u_r\,\mid\, r\in O\}$ and consider
the affine action of the latter on the points of
the set $\c_l$. Then $\mathcal{L}_{l}$ becomes isomorphic to
Funct$\,(\c_l)^\Pi$.

For instance,
 the permutation induced by
$\pi_1\in \Pi$ on $\c_2$ in the case of $A_2$
reads as follows:
\begin{align*}
\c_2=&\{0,\om_1,\om_2,\om_1+\om_2, 2\om_1,2\om_2\}\\
\pi_1(\c_2)=&\{2\om_1,\om_1+\om_2,\om_1,\om_2,2\om_2,0\}.
\end{align*}
Thus dim\,$\l_2=6/|\Pi|=2$ in this example. Only the
sets $\c_{3p}$ contain a (unique) $\Pi$\~invariant point, which is
$p(\om_1+\om_2)$. The general dimension formula for $A_2$ $(l>0)$ is
$$
\hbox{dim}\,\l_{l}=(\frac{(l+2)(l+1)}{2} + \de_l)/3
\for \de_{3p}=2, \de_{3p\pm 1}=0.
$$ 

For $A_1$, $\hbox{dim}\,\l_{l}=1+[l/2]$,
where $[\cdot ]$ is the integer part. Indeed, $\pi_1$
transposes $0$ and $l\om_1$ in this case and has a fixed
point if and only if $l$ is even.

\begin{theorem}
The space  $\mathscr{H}_{l}=
\hat{\mathscr{I}}(\mathscr{X}q^{lx^{2}/2})$\, 
belongs to $\mathcal{L}_{l}.$
For generic $k$, for instance, provided that
$\Re k<0$, this space coincides with $\l_l$. 
\end{theorem}
{\em Proof}. The surjectivity of the map
$\hat{\mathscr{I}}: \mathscr{X}q^{lx^{2}/2}\to$
$\mathcal{L}_{l}$ for generic $k$ is straightforward;
adding $\tilde{\mu}$ does not change the image. One can
also use that this map is zero on  
$\mathcal{J}_{l}(\mathscr{X})q^{lx^{2}/2}$ (see below)
and apply Theorem \ref{THMDIMCOINV}.\sq

Note that the group of the automorphisms of 
the {\em nonaffine} Dynkin diagram acts
in Funct$\,(\c_l)^\Pi$. This action commutes 
with the action of this groups on $\mathscr{X}$ under the
map $\hat{\mathscr{I}}$,
since the Gaussian is invariant with respect 
to these automorphisms.
For instance,   
\begin{align}\label{varsigmaf}
(H^{(l)}_{a})^{\varsigma}=H^{(l)}_{\varsigma(a)}
\for \varsigma(a)=-w_{{}_0}(a),\  X^\varsigma_a=X_{\varsigma(a)}.
\end{align}

\subsection{{\bf DAHA coinvariants}}\label{sec:coinvariants}
\subsubsection{\sf Polynomial coinvariants}\label{sec:pcoinv}
We will introduce the coinvariants only in the context of 
the polynomial representation. 
The space of {\em coinvariants of level $l$}
is $\mathscr{X}/\mathcal{J}_{l}(\mathscr{X})$ 
for the subspace  
$$
\mathcal{J}_{l}(\mathscr{X})\equal\lan\, q^{-lx^{2}/2}\,
T_{\hw}\,q^{\,lx^{2}/2}(X_a)-
t^{l(\hw)/2}X_a\,|\, 
\hw\in \hW, \, a\in P\,\ran\subset \mathscr{X}.
$$
We note that taking only finitely many  $X_a$ is sufficient
in this definition (and all $\hw$). 
For instance, it suffices to make $a=0$ 
if the quotient is one-dimensional
(say, when $l=1$ in the simply-laced
case). 

By construction, $\mathcal{J}_{l}(\mathscr{X})q^{\,lx^{2}/2}$ 
belongs to the kernel of the map $\hat{\mathscr{I}}$. Denoting
the map $\HH\ni A\mapsto q^{x^{2}/2}A\,q^{-x^{2}/2}$ by
$\tau$ (it is an automorphism of $\HH$), 
$\mathcal{J}_{l}(\mathscr{X})=
\tau^{-l}(\mathcal{J}_{0}(\mathscr{X}))$.

We claim that the dimension
of $\mathscr{X}/\mathcal{J}_{l}(\mathscr{X})$ {\em always}
coincides with that of the Looijenga space (defined above). 
The dimension of the 
space of coinvariants can be calculated without any reference
to the Looijenga space.   

\begin{theorem}\label{THMDIMCOINV}
For any $q,t\in \C^*$ and $l>0$, 
\begin{align*}
\dim(\mathscr{X}/\mathcal{J}_{l}(\mathscr{X}))
\,=\,\dim( \text{Funct\,}(\c_l)^{\Pi}).
\end{align*}
\end{theorem}
{\em Sketch of the proof.}
We use the PBW theorem to 
establish the inequality  
\begin{align}\label{dimineqality}
\dim(\mathscr{X}/\mathcal{J}_{l}(\mathscr{X}))
\,\leq\,\dim( \text{Funct\,}
(\c_l)^{\Pi}).
\end{align}
Let $k\to 0$ $(t=q^k\to 1)$. Then $T_{\hw}\to \hw$ and
$\HH(t=1)$ becomes the classical Weyl algebra generated by
$X_a$ and $Y_b$ extended by $W$. The dimension can be
readily calculated at $k=0$; it equals  
$\dim( \text{Funct\,}(\{b\in P_{+}, (b,\theta)\leq l\})^{\Pi})$.
Due to (\ref{dimineqality}),
this dimension must remain the same for all $q,t$.\sq

\subsubsection{\sf The B-case}\label{sect:THEBCASE} 
Avoiding the non-simply-laced root systems
in Theorem \ref{HALLONE} is not only a technicality.
The dimension of $\l_{1}$ is greater than one if $P\neq P^\vee$,
so it is not true (generally) that all level-one Hall functions 
are proportional to $\ga(x)$, as stated in this theorem. 
However for $B_n$, there
is the following possibility to make the image really 
one-dimensional (for $l=1$).

We use that $Q=P^\vee$ in this case and consider 
$\mathscr{X}\,'=\C_{q,t}[X_a,a\in Q]$ 
instead of the complete polynomial representation
$\mathscr{X}$. The space $\mathscr{X}\,'$
is a module over the {\em little DAHA} (in the terminology
from \cite{C101}), which is generated by 
$\mathscr{X}\,'$ and the same $\{T_{\hw},\hw\in \hW\}$;
all the considerations above hold under this restriction. 
The corresponding level-one Looijenga space will be isomorphic to 
Funct$\,(Q/lQ^\vee)^{\,W}$, i.e., will be of dimension one as $l=1$.
The formula (\ref{levonehall}) holds if $\rho$ is replaced by
$\rho^\vee$ and $a\in Q$.

Generally, if there is any DAHA-submodule $\mathscr{X}\,'$,
then, automatically, 
$$\hat{\mathscr{I}}(\mathscr{X}\,'\,q^{lx^{2}/2})
\subset
\{\sum_{\hat{w}\in \hat{W}}\hat{w}(G(X)\,q^{lx^{2}/2}),\, 
G(X)\in \mathscr{X}\,'\} \hbox{\ \, for\ any\ } l>0.
$$

\subsubsection{\sf Levels 0 and 1}
Let us consider the (simplest) 
cases when the space of coinvariants is one-dimensional.
\begin{theorem}\label{thm:coinv}
In the level-zero case,
provided that the space of $Y$\~eigenvectors with the
eigenvalue $t^{\rho}$ (i.e., containing
$E_{0}=1$) is one dimensional in $\mathscr{X}$,
$$
\dim(\mathscr{X}/\mathcal{J}_0(\mathscr{X}))=1\ \, \hbox{and\ \,} 
\oplus_{\,q^\la\,\neq\, t^{\rho}\,}\mathbb{C}\mathscr{X}_{\la}=
\mathcal{J}_0(\mathscr{X}),
$$
where $\mathscr{X}_\la=\{f\in \mathscr{X}\,|\,
(Y_a-q^{(\la,a)})^N (f)=0\}$ for sufficiently large $N$;
we identify $q^\la$ if they give coinciding $Y$\~eigenvalues.   
This dimension is one for $l=1$ as well in the simply-laced case;
then $q,t$ can be {\sf arbitrary} nonzero.
\end{theorem}
{\em Proof.}
If the nonsymmetric Macdonald polynomials
$E_a$ are well defined, then they form a basis 
for $\mathscr{X}$. Otherwise, use the generalized
$Y$\~eigenvectors in the following reasoning.
Recall that the action of $Y_{b}$ is given by 
$Y_{b}^{-1}(E_{a})=q^{(a_{\#},b)}E_{a}$ for
$a\in P,\, b\in P^\vee$.
So for any $a\in P$ such that $q^{(a_{\#},b)}\neq q^{k(\rho,b)}$,
we have $E_{a}\in \mathcal{J}_0(\mathscr{X})$.
Then $E_{0}=1$ is of multiplicity one in $\mathscr{X}$ and
$\dim(\mathscr{X}/\mathcal{J}(\mathscr{X}))=1$. 

In the case $l=1$, we use that  $\tau_-\tau_+^{-1}(Y_b)=$
$\tau_-\tau_+^{-1}\tau_-(Y_b)=
\si^{-1}(Y_b)=X_b^{-1}$ and apply $\tau_-^{-1}$ to the
triple 
$\{\,\{X_a\}, \{T_w\}, \{Y_b\}\,\}$, 
satisfying the PBW theorem. See \cite{C101} for the
definitions of $\tau_{\pm}, \si$ and also see
Lemma \ref{lem:basis} below for the case of $A_1$. \sq

\subsection{{\bf Kac-Moody limit}}
The limiting case $t\to\infty$ ($k\to -\infty$) is important. Then
the Hall function $\tilde{H}_{a}^{(l)}$ for a
weight $a\in P_+$ subject to $(a,\th)\le l$ becomes
proportional to the character
of the corresponding integrable Kac-Moody module. 
The level $\,l\in \N\,$ equals the action of the central element 
$c$ in the standard normalization; we consider here
only the case of standard (split) Kac-Moody algebras.  

Notice that we use the extended affine 
Weyl group $\hW$ with $P^\vee$ instead of $Q^\vee$ (usual
in Kac-Moody theory) and that, in our approach, 
the weights $a\in P$ are not supposed to be $l$\~dominant. 
The Hall functions can be defined for any $a$, but
their interpretation as characters of integrable modules 
of level $\,l\,$ in the limit does require $a\in P_+$ and 
the inequality $(a,\th)\le l$. This connection with
the Kac-Moody characters is known; see
e.g., \cite{Vi}. Let us discuss this in detail. 

\subsubsection{\sf Explicit formulas} 
From (\ref{mutildemu}) and (\ref{consterm}),
\begin{align}\label{tildemulim}
&\tilde{\mu}(t\to\infty) =
\prod_{\tilde{\al}>0}\frac{1}
{1-X_{\tilde{\al}}},\ 
\lim_{t\to\infty}\hbox{ct}(t^{-1})=\prod_{i=1}^\infty\frac{1}
{(1-q^{i})^n}. 
\end{align} 
Also, $\hat{P}(t^{-1})\to |\Pi|$
as $t\to\infty$\, for\, $\Pi=P^\vee/Q^\vee$. Setting 
$$
\hat{\chi}_a ^{(l)}
\equal q^{-l\frac{x^2}{2}}\,
\lim_{t\to \infty}\tilde{H}_{-a}^{(l)} \for
a\in P \hbox{\ \,(notice\ $-a$)},
$$
\begin{align}\label{Kac-Moody}
\hat{\chi}_a^{(l)}\ &=\ q^{-l\frac{x^2}{2}}
\sum_{\hw\in \hW}\, 
\hw(X_{a}^{-1}\,\tilde{\mu}(t\to\infty)\,q^{l\frac{x^2}{2}})\\
&=\ \bigl(\sum_{\hw=bw} (-1)^{l(\hw)} 
X_{\hw(\hat{\rho}+a)-\hat{\rho}+lb}^{-1}\ q^{\,lb^2/2}\bigr)/
\prod_{\tal\in \tR_+}(1-X_{\tal}).
\notag
\end{align}
Here the summation is over all $b\in P^\vee, w\in W$ and 
we set (symbolically) $\hat{\rho}=\frac{1}{2}
\sum_{\tal\in \tR_+}\tal$ (as for Kac-Moody algebras).
What we really need is the relation
$$\sum_{\tal\in \La(\hw^{-1})}\tal\ =\ \hat{\rho}-\hw(\hat{\rho})
$$
for the sets $\La(\hw^{-1})$ defined in (\ref{Lahw});
note $\hw^{-1}$ here.
Using the level-zero and level-one formulas 
for $\hat{\mathscr{S}}\,'_+\circ \tilde{\mu}$,
\begin{align}\label{thetaprod}
\prod_{\tal\in \tR_+}(1-X_{\tal})&=
\frac{\sum_{\hw=bw} (-1)^{l(\hw)} 
X_{\hat{\rho}-\hw(\hat{\rho})}}
{|\Pi|\,\prod_{j=1}^\infty (1-q^j)^n}\\
&=\frac{\sum_{\hw=bw} (-1)^{l(\hw)} 
X_{\hat{\rho}-\hw(\hat{\rho})-b}\ 
q^{\,b^2/2}}{\sum_{b\in P}\,X_{b}\,q^{\,b^2/2}}.
\label{thetaprodone}
\end{align}
Formula (\ref{thetaprodone}) is stated here in the simply-laced case
as in (\ref{levonehall}). One can readily adjust this formula to 
the setting of \cite{C101}, i.e., to the case of {\em twisted}
$\tR^\nu$\~affinization (then an arbitrary reduced  
nonaffine $R$ can be used). 

These two formulas are the denominator identity and
the level-one formula due to Kac. See 
Theorem 10.4, Lemma 12.7 and (12.13.6) from \cite{Kac}.
We conclude that $ct(t^{-1})|_{t\to \infty}\,\hat{\chi}_a^{(l)}$ is
the character of the corresponding Kac-Moody integrable
module of level $l$ provided that
$a\in P_+$ and $(a,\th)\le l$; this is upon the substitution
$X_{b}\mapsto e^{-b}$.

Let us provide the
first few terms of the numerators of these formulas 
in the case of $A_1$:
\begin{align*}
&
\sum_{\hw=bw} (-1)^{l(\hw)} X_{\hat{\rho}-\hw(\hat{\rho})-lb}\ 
q^{\,lb^2/2}\mod (q^2)\\
=&\left\{\begin{array}{ccc} 
1-X^2+q^{1/4}(X^{-1}-X^3-qX^{-3}+qX^5) 
&\for l=1,&  \\
2(1-X^2+qX^{4}-qX^{-2})
&\for l=0,&
\end{array}\right.
\end{align*}
where $X=X_{\om_1}$.
Compare this with the left-hand side of
(\ref{thetaprod}) and (\ref{thetaprodone}) multiplied by the
corresponding denominators (here the calculations are 
direct).

In our approach, there are no clear reasons to stick here
to {\em affine $l$\~dominant} weights,
i.e., to $a\in P_+$ subject to $(a,\th)\le l$.
Apart from the weights of integrable modules, i.e.,
for arbitrary $a\in P$, the following level-one
formulas in terms of the polynomials 
$\tilde{E}_a\equal E_a(t\to \infty)$
are worth mentioning:
\begin{align}\label{Kac-M-Herm}
\hat{\mathscr{S}}_+'\bigl(
\tilde{\mu}(t\to\infty)\, \tilde{E}_a\,q^{x^2/2}\bigr)=
\left\{\begin{array}{ccc}q^{-a^2/2}\ga(x), &  & 
\hbox{if\  } a\in P_-, \\0, &  & \hbox{otherwise}.
\end{array}\right.
\end{align}

We use formula (\ref{levonehall}). The polynomials 
 $\tilde{E}_a$ are closely connected with  
the $q$\~Hermite polynomials
$E_a(t\to 0)$ studied in \cite{ChW} (and which play the
key role in the theory of $q$\~Whittaker functions).
\smallskip

\rmk
Let us consider briefly the limit 
$t\to 0$ ($k\to \infty$). Then the series $\tilde{\mu}^{-1}\circ
\hat{\mathscr{S}}\,'_+\circ \tilde{\mu}$ can also be
interpreted via the Kac-Moody characters. 
Due to (\ref{murelations}),
\begin{align*}
&q^{-l\frac{x^2}{2}}\,
\lim_{t\to 0}\tilde{H}_{a}^{(l)} \ =\ 
\frac{\sum_{\hw=bw} (-1)^{l(\hw)} 
X_{\hw(\hat{\rho}+a)-\hat{\rho}-lb}\ q^{\,lb^2/2}}{
\prod_{\tal\in \tR_+}(1-X_{\tal})}.
\end{align*}

\subsubsection{\sf Match at level one} 
We note that (12.13.6) from Kac's book is stated 
in the simply-laced case, which matches the setting we
use for formulas (\ref{levonehall}) and (\ref{thetaprodone}).
Calculating the level-one characters in the cases 
$B_n, F_4, G_2$ is due to Kac and Peterson. As for the $k$\~case
(i.e., when $t$ is added),
we explained in Section \ref{sect:THEBCASE} how to proceed in the
$B$\~case for the lattice $Q^\vee$. 
The root systems  $F_4$ and $G_2$ with $k$ also seem doable.  

The most difficult case in the theory of level-one
Kac-Moody characters is $C_n$ (managed by Kac and Wakimoto); 
it seems exactly
parallel to the problem with explicit formulas for the 
affine Hall functions of type $C_n$ for $l=1$ (untwisted). 
The paper \cite{Sto} devoted to the $C^\vee C$
may contain the methods and results sufficient to manage
this case.  

The above discussion and considerations of this section
are in the {\em untwisted} case.  
The formulas for the {\em twisted Kac-Moody characters} 
are known for any root systems. The twisted KM\~characters 
correspond (with some 
reservations) to our using $\tR^\nu$, the twisted
affinization from Section \ref{sect:thinaff}. Similar to
Kac-Moody theory, the level-one formulas with $k$ 
were obtained (uniformly) in \cite{C5}
for {\em any} reduced root systems.

It is worth mentioning that
the classification of Kac-Moody algebras is {\em not} the
same as that for DAHA (which continues the classical classification
of symmetric spaces). However, when they intersect,
it seems that there is almost an exact match between the problems 
arising in the theory of Kac-Moody characters
and those for the affine Hall functions (with $k$).
At least, this is so in the level-one case. 
Recall that the affine Hall functions belong
to the same Looijenga space as the Kac-Moody characters. 
We do not discuss explicit formulas for $l>1$,
where not much is actually known; see \cite{Vi}. 

Let us mention that in the level-one case,
the affine Demazure characters are directly connected with the
nonsymmetric $q$\~Hermite polynomials $E_a(t\to 0)$ (see
above). They become $W$\~invariant for $a\in P_-$ and
their coefficients in this case are given in terms of 
the $q$\~Kostka numbers (see \cite{San},\cite{Ion1}).
\smallskip

We are grateful to Victor Kac who helped us
establish the correspondence between the two theories,
the classical KM theory and the one for arbitrary $k$.
We thank Boris Feigin for a helpful discussion.
As a matter of fact, we introduce in this paper certain
$t$\~deformations of the Demazure characters,
but our definition is of a technical nature and we do not now
how far this can go.
\smallskip

\subsection{{\bf Shapovalov forms}}
We will begin with a very general approach to constructing
inner products (in functional analysis, known as $GNS$ 
construction). Let $\f$ be a cyclic $\HH$\~module, 
i.e., $\f=\HH(vac)$ 
for some $vac\in \f$. Actually $\f$ can be absolutely
arbitrary in the following (formal) considerations, but
we prefer to restrict ourselves to cyclic modules here. 
We assume that $\HH$ and $\f$
are defined over a field $\tilde{\C}$. It can be 
$\C_{q,t}$, the definition field
for the polynomial representation of $\HH$,
or its extension by the parameters of $\f$
(treated as independent variables).
If $q,t$ and the parameters of $\f$
are considered as nonzero complex numbers, then 
$\tilde{\C}=\C$. 

\subsubsection{\sf Symmetric J-coinvariants}
We set $\mathcal{J}=\{A\in \HH\,|\,A(vac)=0\}$ (a left ideal).
Then $\f\cong \HH/\mathcal{J}$. 
{\em Any} form on $\f$ which is symmetric and $\HH$\~invariant
with respect to a given anti-involution $\star$ can be obtained
as follows.

Let $\star$ be an anti-involution $\star$ on $\HH$\,
($\star^2=1$ is required because the form must be symmetric)
and let $\varrho: \HH\to \tilde{\C}$ be a functional on $\HH$ 
such that $\varrho(A^{\star})=\varrho(A)$ and 
$\varrho(\mathcal{J})=0$. 
Automatically, we have that 
$\varrho(\mathcal{J}^{\star})=0$ ($\mathcal{J}^{\star}$ is
a right ideal in $\HH$).
Since $\varrho(\mathcal{J}+\mathcal{J}^{\star})=0$, it comes 
from a functional 
$\varrho': \f\to \f/\mathcal{J}^{\star}(\f)\to \tilde{\C}$. 

Then the form on $\f$ is introduced as follows: 
$$\lan f, g\ran\equal\varrho(\bar{f}^{\star}\,\bar{g})=
\varrho'(f^{*}\,g), \, f,g\in \f,$$
where we lift $f,g$ to 
$\bar{f},\bar{g}\in \HH$ and set $f^{*}=\bar{f}^{\star}(vac)$.

This form $\lan\,,\,\ran$ is obviously symmetric and 
$\star$\~invariant:
$$\lan A(f), g\ran=\lan f,  A^{\star}(g)\ran, \text{ where }
f,g\in \f, A\in \HH.
$$

To describe all such forms, let us introduce the space  
\begin{align}\label{starcoinv}
&\HH/(\mathcal{J}+\mathcal{J}^{\star})=\f/\mathcal{J}^{\star}(\f).
\end{align}
and its dual 
$\Hom_{\tilde{\C}}(\f/\mathcal{J}^{\star}(\f),\tilde{\C})$. 
Both have a natural action of $\star$ and are direct
sums of $\pm 1$\~eigenspaces.

The subspace of $\star$\,\~invariant elements of 
$\Hom_{\tilde{\C}}(\f/\mathcal{J}^{\star}(\f),\tilde{\C})$
will be called the {\em space of  $\star$\,\~symmetric
$\j$\~coinvariants}. We will always assume that $1^\star=1$, 
correspondingly, $vac^\star=vac$.

The $\pm 1$\~eigenvectors of $\star$ from
$\Hom_{\tilde{\C}}(\f/\mathcal{J}^{\star}(\f),\tilde{\C})$ 
lead to either $\star$\,\~invariant forms or 
to $\star$\,\~anti-invariant ones, respectively. In the examples 
we consider, the action of $\star$ is trivial in the whole
space from (\ref{starcoinv}) and its dual, but 
the minus-sign (equally interesting) may occur as well.

Let us discuss basic examples.

\subsubsection{\sf Shapovalov pairs} 
We call the nonzero form $\lan\,,\,\ran$ 
a {\em  Shapovalov  form} if
$$\hbox{dim}_{\tilde{\C}} 
\bigl(\HH/(\mathcal{J}+\mathcal{J}^{\star})\bigr)=1=
\hbox{dim}_{\tilde{\C}}\bigl(\f/\mathcal{J}^{\star}(\f)\bigr),
$$
and therefore this form is a unique symmetric 
$\star$\,\~invariant form in $\f$ up to proportionality.
Accordingly, $\{\mathcal{J},\star\}$ is called a 
{\em Shapovalov pair}. 

This terminology may be somewhat misleading.
The anti-involutions we are going to consider generally have 
little to do with those
in Lie theory; the connection with the
Heisenberg and Weyl algebras is significantly more direct.
However, our usage of the PBW theorem is really 
similar to the original Shapovalov construction.

Given a Shapovalov pair $\{\mathcal{J},\star\}$,
finding $\lan f,g\ran$ is
purely algebraic problem directly related to 
the PBW theorem. For
instance,  $\lan f,g\ran$ always depends rationally 
on the parameters $t,q$ of $\HH$. 
It is valuable, since the forms given by
integrals (or similar) are generally well defined only
for some $t,q$. Their meromorphic
continuation to other values of $q,t$ can be involved.
\smallskip

\rmk
We follow in this section unpublished notes by the
first author devoted to the  Arthur-Heckman-Opdam 
formulas \cite{HO2} in the theory of the spectral 
decomposition of AHA (due to Lusztig and many others).
This approach is  based on a relatively direct 
(without geometry) meromorphic continuation of the 
corresponding Plancherel formula and 
``picking the residues".

The DAHA version of this decomposition is completed (by now) only 
for $A_n$ (unpublished). The best reference we can give 
so far is \cite{ChL}. The main theorem is that 
the Shapovalov form coincides with the analytic continuation 
of the corresponding inner product defined in terms of the 
standard integration over $i\R^n$ subject to $\Re k>0$.
A direct analytic continuation of the latter to negative
$\Re k$ appeared a certain generalization of the 
``picking the residues" 
in AHA theory. In contrast to the Arthur-Heckman-Opdam 
method \cite{HO2}, the result of this procedure is
known {\em a priori}. It is the Shapovalov form, which
is defined entirely algebraically, and is rational
or even regular in terms of $t$; see
Theorem \ref{Shapind} below.
\sq
\smallskip
 
The case of the standard form associated with the
anti-involution $\ast$ of the polynomial representation,
sending $t,q,X_a,Y_b,T_i$ to their inverses,
was considered
in \cite{C101},Proposition 3.3.2. The rational dependence
of the corresponding inner products in terms of $q,t$ 
was deduced there from the uniqueness of such a form up to
proportionality. 
A similar approach was applied to the anti-involution
$\phi$ (governing the duality) in \cite{C101} and to the 
bilinear invariant forms involving the $q$\~Gaussians
(generalizations of the Mehta-Macdonald integrals).
See (\ref{diamondef}) below.

\subsubsection{\sf Y-induced modules}
Let us discuss the Shapovalov forms for
the $Y$\~induced modules $\f=\i_\la$, where $\la\in \tilde{\C}^n$.
By definition, $\i_\la$ is a free $\HH$\~module over $\tilde{\C}$
generated by $vac$ with the defining relations 
$Y_b(vac)=q^{(\la,b)}\,vac$. It belongs to the category $\o$
with respect to the action of $Y$\~elements, i.e., it can be 
represented as a direct sum of the {\em finite-dimensional\,}
spaces of generalized $Y$\~eigenvectors.
For the sake of definiteness,
let us assume that $T_i^{\star}=T_i$ for $i=1,\ldots,n$. Then
the corresponding  $\varrho$ satisfies the following:
\begin{align}\label{Shapo}
\varrho(Y_a^{\star}T_w Y_b)=q^{(\la,a+b)}\varrho(T_w),\ 
\varrho(T_w)=\varrho(T_{w^{-1}}) \for w\in W.
\end{align}
The latter relation simply means that $\varrho$ is a trace functional
on the nonaffine Hecke algebra $\H$. 

We call the anti-involution 
$\,\star\, $ of {\em strong Shapovalov type} with respect to
$\mathscr{Y}$ if $\HH$ satisfies the PBW condition for 
$\mathscr{Y}$, 
$\H$ and $\mathscr{Y}^{\star}$ (replacing $\mathscr{X}$).
Namely, if an arbitrary $A\in \HH$ can be uniquely represented
as $c_{awb}\,Y_a^{\star}T_w Y_b$ for $a,b\in P^\vee$ and 
$w\in W$. Then the conditions from (\ref{Shapo}) determine $\rho$
completely. We see that the simply-laced root systems are
generally needed here, unless in the twisted (self-dual) 
setting for the affine root system $\tR^\nu$, as in \cite{C101}. 
Note that the definition
of strong Shapovalov anti-involutions
depends only on $\star$ and $\mathscr{Y}$, 
not on the module $\i_\la$ ($\la$ can be arbitrary).
 
An important example of the {\em weak\,} (not strong)
Shapovalov anti-involution in $\i_{\la}$ is when
$\mathscr{Y}^{\star}=\mathscr{Y}$, i.e., $\mathscr{Y}$ is
a {\em normal subalgebra} with respect to $\star$. Then
the Shapovalov condition holds for $\i_\la$ provided that 
the {\em generalized\,} $Y$\~eigenspace containing 
$vac$ is one-dimensional in $\i_\la$. Indeed, the 
linear span of the spaces 
$(Y_a-q^{(a,\la)})\i_\la\subset$Ker$(\varrho)$ 
is of codimension one in $\i_\la$ in this case. 
Here $\star$ can be arbitrary, provided 
$\mathscr{Y}$ is normal.

There are actually only a few {\em strong} Shapovalov anti-involutions
in DAHA theory, 
essentially the examples (1) and (3) considered below
(for the subalgebra $\mathscr{Y}$). 
They play a significant role. The corresponding PBW property 
holds for any (nonzero) $q$ and $t$ for these anti-involutions.

The following rationality theorem clarifies the importance of the 
Shapovalov property in both, the weak and strong variants.
The first generally guarantees rational dependence of 
the inner products on the parameters (including $q,t$);
the second provides regular dependence.

We follow Proposition 3.3.2 from \cite{C101}. Let the algebra
$\HH$, the representation $\f$, and the functional $\rho$
be defined over the same field $\tilde{\C}$. For instance,
the field of rationals $\C(q^{1/m},t^{1/2})$ can be taken
for the polynomial representation (generally this field 
is supposed to contain the parameters of the module $\f$).

\begin{theorem}\label{Shapind}
(i) A form $\lan\,,\,\ran$ on $\f$ corresponding
to a Shapovalov pair $\{\j,\star\}$ is a unique symmetric 
$\star$\,\~invariant form in $\f$ up to 
proportionality; let us normalize it by the condition
$\lan 1,1\ran=1$. Then given $f,g\in \f$, their inner 
product $\lan f,g\ran$ belongs to the field $\tilde{\C}$
(which may include the parameters of $\f$).

(ii) Assuming that $\star$ satisfies the strong Shapovalov 
property for any nonzero $q$ and $t$, let $f,g$ be taken from
$\HH_{int}(vac)$, where  
\begin{align}\label{HHinteg}
\HH_{int}=\C[\,q^{\pm 1/m},t^{\pm 1/2}\,]\,[X_a,Y_b,T_w]_
{\hbox{\tiny nc}}\ \subset\ \HH.
\end{align}
The ring of coefficients here is the standard $\C$\~algebra
necessary for the defining DAHA relations and by 
$[\ ]_{\hbox{\tiny nc}}$ we mean the noncommutative
algebraic span.
Then the inner product
$\lan f,g\ran$ is well defined for any nonzero $q,t$. 
In other words, if the PBW property holds for $\mathscr{Y}$,\, 
$\H$ and $\mathscr{Y}^{\star}$, then the
corresponding form is regular in terms of 
$q^{\pm 1/m},t^{\pm 1/2}$.   
\end{theorem}

\subsubsection{\sf The polynomial case}{\label{sect:poly}}
Let us discuss the Shapovalov condition for an arbitrary
anti-involution $\star$, fixing $T_i$ for $i>0$, combined with
the polynomial representation 
$\mathscr{X}$. This representation is a quotient of $\i_{\la}$ for
$\la=k\rho$; the vacuum element (the cyclic generator of 
$\i_{k\rho}$) becomes $1\in \mathscr{X}$. One has
$$\HH/(\HH\mathcal{J}+\mathcal{J}^{\star}\HH)
\cong \mathscr{X}/\mathcal{J}^{\star}(\mathscr{X})$$
for the left ideal $\mathcal{J}$ linearly generated
by the spaces $\HH(T_{\hw}-t^{l(w)/2})$. 
This results in\  
$
\varrho(Y_a^{\star}T_w Y_b)=t^{(\rho,a+b)+l(w)/2}.\ 
$ 

Chapter 3 of \cite{C101} is actually the theory of 
the following three anti-involutions and the corresponding
symmetric forms:
\begin{align}\label{diamondef}
(1)&\ \  \varphi: X_a\leftrightarrow Y_a^{-1}, 
T_w\mapsto T_{w^{-1}},\notag \\
(2)&\ \  \Diamond: X_a\mapsto T_{w_0}^{-1}X_{-w_0(a)}T_{w_0}, 
Y_b\mapsto Y_b, T_w\mapsto T_{w^{-1}},\\
(3)&\ \  \Diamond_1= q^{-x^2/2}\circ\Diamond\,\circ q^{x^2/2}:
Y_a\mapsto q^{-x^2/2}Y_a q^{x^2/2}.\notag
\end{align}
We assume that $R$ is simply-laced in (1) (it is arbitrary
in \cite{C101} because $\tR^\nu$ is considered there). 
Let us provide some details. 

(1) This anti-involution controls the duality and evaluation 
conjectures and is related to the Fourier transform. 
The Shapovalov property for $\varphi$ is 
{\em exactly} the PBW Theorem (any $q,t$). The corresponding form 
is well defined for any $q,t$ and the study of its radical is an 
important tool in the theory of the polynomial representation 
of DAHA.

(2) The second anti-involution governs the inner product 
in $\mathscr{X}$ (without conjugating $q,t$);
$\Diamond$ is of Shapovalov type only for generic $k$ (and
there is no immediate relation to the PBW theorem). So it is 
{\em weak\,}. The corresponding bilinear form is the key in the 
DAHA harmonic analysis, including 
the Plancherel formula for $\mathscr{X}$ and its Fourier
image, the representation 
of $\HH$ in delta functions.  

(3) The third appears in the 
difference Mehta-Macdonald formulas and is used to prove that 
the Fourier transform of the DAHA module $\mathscr{X}q^{-x^2/2}$ is 
$\mathscr{X}q^{+x^2/2}$. The strong Shapovalov property holds here, 
so the form is well defined for any $q,t$. 
The radical of the corresponding pairing is closely related to that 
in (1) (they coincide in the rational theory).

\subsection{{\bf Using induced modules}}
\subsubsection{\sf Level-zero forms}{\label{sect:form}}
Let us consider the coinvariants in the case $l=0$ via 
the affine symmetrizer $\hat{\mathscr{P}}_+$. 
The $\hat{\mathscr{P}}$\~symmetrizer 
is more convenient here than $\hat{\mathscr{I}}$. 
The definition is in (\ref{defhatP}); we will also 
use the rational formula of Theorem \ref{P+FORMULA},
which gives a $t$\~meromorphic continuation of this
operator when acting in $\mathscr{X}$. 
   
Recall that  
$
\,\hat{\mathscr{P}}_{+}(f)\,=\,\hat{\mathscr{P}}\,'_{+}(f)/
\hat{P}(t^{-1}),\,
$ 
where $\hat{P}(t)$ is the affine
Poincar\'e series; see (\ref{defhatP}).
We continue using the notation $\mathcal{J}\subset \HH$
for the ideal such that $\mathscr{X}=\HH/\mathcal{J}$;
it is the linear span of subspaces
\begin{equation*}
\HH(T_{\hat{w}}-t^{l(\hat{w})/2})
\for \hat{w}\in \hat{W}.
\end{equation*}

For the anti-involution  $\Diamond$ in
(\ref{diamondef}), the 
functional 
\begin{align*}
\varrho_+:\, \HH \, \to\, \C_{q,t}\hbox{\ \,sending\ \,} 
\  A\, \mapsto\, 
\hat{\mathscr{P}}_+A(1) 
\end{align*}
satisfies the $\Diamond$\~invariance property
$
\varrho_+(\mathcal{J}^{\Diamond}+\mathcal{J})=0.
$
Indeed, 
$$
\varrho_+(f)=\hat{\mathscr{P}}_+(f),\  
\varrho_+((T_{\hat{w}}^{\Diamond}-t^{l(\hat{w})/2})f)=0
\for f\in\mathscr{X},
$$
since $\Diamond$ preserves $\h=\lan T_{\hw}\ran$. Thus,
$\varrho_+$  can be used to construct a symmetric 
form on $\mathscr{X}$ corresponding
to the anti-involution $\Diamond$. 

This argument is of course {\em formal}; 
one needs to address the existence of 
$\hat{\mathscr{P}}_+(f)$.   
Theorem \ref{P+FORMULA} provides the existence
of  $\hat{\mathscr{P}}_+$ if
there are no $Y_{\om_i^\vee}$\~eigenvectors in 
$\mathscr{X}$ with the eigenvalue $t^{-(\rho,\om_i^\vee)}$
for $i=1,2,\ldots,n$.  

\rmk
The rational formula for  $\hat{\mathscr{P}}_+(f)$ 
from Theorem \ref{P+FORMULA} cannot be used
in (the whole) $\mathscr{X}$ if $q$ is a root of unity even if
$t$ is sufficiently general. Indeed, recall that
the $Y$\~eigenvalue of $1\in \mathscr{X}$ is $t^{\rho}$.  
For generic $q$, the
parameter $t$ can be an $N$-th root of unity for 
sufficiently large $N$. The latter is needed to
avoid the zeros of $\hat{P}(t^{-1})$.
\sq
\smallskip 

Under these conditions, the space of $\{\varrho_+\}$ is
one-dimensional and 
$\hat{\mathscr{P}}_{+}$ becomes a {\em universal}\,
$\Diamond$\~coinvariant, which leads to the following
construction. Recall that 
$\varsigma(a)=-w_{{}_0}(a),\  
X^\varsigma_a=X_{\varsigma(a)}$; see (\ref{varsigmaf}).

\begin{theorem}\label{DIAZERO}
(i) Let us assume that $\mathscr{X}$ has a nonzero symmetric
form $\lan f,g \ran$ with the anti-involution $\Diamond$
normalized by $\lan 1,1\ran=1$. Given any
$f,g\in \mathscr{X}$, $\lan f,g\ran$ is a rational
function in terms of $q,t$. Provided that  $\Re k<0$
and $|\Re k|$ is sufficiently large (depending on $f,g$), 
\begin{align}\label{psymformula}
\lan f,g\ran=
t^{-l(w_0)/2}\,\hat{\mathscr{P}}_{+}(f T_{w_0}(g^\varsigma)).
\end{align}

(ii) Let $\hat{P}(t^{-1})\neq 0$ for the affine
Poincar\'e series expressed as in (\ref{rataffpoin}),
$\f$ be a $\HH$\~quotient of $\mathscr{X}$  such that 
it has no  $Y_{\om_i^\vee}$\~eigenvectors  
with the eigenvalue $t^{-(\rho,\om_i^\vee)}$.
for any  $i=1,2,\ldots,n.$
Using the rational presentation 
for $\hat{\mathscr{P}}\,'_{+}$ from Theorem \ref{P+FORMULA}, 
formula (\ref{psymformula}) supplies $\f$ with  
a bilinear symmetric form associated with the 
anti-involution $\Diamond$
and satisfying $\lan1,1\ran=1$.\sq 
\end{theorem} 

Compare with Proposition 3.3.2 from \cite{C101} and
with Theorem \ref{Shapind} above. 

\subsubsection{\sf X-induced modules}
A modification of formula (\ref{psymformula}) 
can be used in $X$\~induced $\HH$\~ modules. They are defined
as universal $\HH$\~modules $\i^X_\xi$ generated
by $v$ subject to $X_a(v)=q^{(\xi,a)}v$ for
$\xi\in \C^n$, $a\in P$. If
$\xi$ is generic, then the module 
$\i^X_\xi$ is $X$\~semisimple and can be identified 
with the {\em delta-representation} of $\HH$ in the space
$$
\De_\xi\equal \oplus_{\hw\in \hW}\,\C_{q,t}\,\chi_{\hw}
$$ 
in terms of the {\em characteristic functions} 
$\chi_{\hw}$ defined as follows:
$$
\chi_{\hw}(\hu)= \de_{\hw,\hu}\,,\ \,
\chi_{\hw}\chi_{\hu} = \de_{\hw,\hu}\,\chi_{\hw}
\hbox{\ \ for the Kronecker delta }.
$$
The action of the $X$\~operators
is via their evaluations at $\{q^{\hw(\xi)}\}$:
\begin{align*}
&X_a(\chi_{\hw})\equal 
X_a(\hw)\chi_{\hw} \for a\in P, \hw\in \hW,\\
&X_a(bw)\equal X_a(q^{b+ w(\xi)})\, =\,
q^{(a,b)}X_{w^{-1}(a)}(q^{\xi}).\notag
\end{align*}

The group $\hW$ acts on the characteristic functions
through their indices:\,  $\hu(\chi_{\hw})=\chi_{\hu\hw}$
\, for $\hu,\hw\in \hW$.
Accordingly,
\begin{align*}
T_i(\chi_{\hw})\,&=\,
\frac{t^{1/2}X_{\al_i}^{-1}(q^{w(\xi)})
q^{-(\al_i,b)} - t^{-1/2}}{
{X_{\al_i}^{-1}(q^{w(\xi)})q^{-(\al_i,b)} - 1}  }\,
\chi_{s_i\hw}\notag\\
&-\,\frac{t^{1/2}-t^{-1/2}}{
{X_{\al_i}(q^{w(\xi)})q^{(\al_i,b)} - 1}  }\,
\chi_{\hw} \for \hw=bw\in \hW, \notag\\
\pi_r(\chi_{\hw})\,=\, &\chi_{\pi_r\hw}, \where
\pi_r\in \Pi,\, 0\le i\le n,\, X_{\al_0}=qX_{\th}^{-1}.
\end{align*}

The $X$\~weight $q^\xi$ is assumed generic in this formula
and below.
We follow Section 3.4.2, ``Discretization", from \cite{C101}.

The {\em delta functions\,} are defined as 
$
\de_{\hw}(\hu)\ =\ \mu_\bullet(\hw)^{-1}\chi_{\hw}
$
for $\mu_\bullet(\hw)\equal\mu(\hw)/\mu(\hbox{id})$,
the measure function in the following inner product:
\begin{align}
&\langle f,g\rangle_\bullet = 
\sum_{\hw\in \hW}
\mu_\bullet(\hw) f(\hw)\ g(\hw)\ =\
\langle g,f\rangle_\bullet\,.
\label{innerdelnew}
\end{align}
Here $f,g$ are finite or infinite (provided the convergence)
linear combinations of the characteristic functions considered
as functions on $\hW$. By construction, 
$\lan \chi_{\hu},\de_{\hw}\ran_{\bullet}=
\de_{\hu,\hw}$ for $\hu,\hw\in \hW$ and Kronecker's 
$\de_{\hu,\hw}$. 

The values
$\mu_\bullet(\hw)$ are given by formulas in 
(\ref{murelations}); replace in this formula $X$ by $q^\xi$ and
$\hw$ by $\hw^{-1}$. We see that (\ref{innerdelnew}) is
directly connected
with the affine symmetrizer 
$\hat{\mathscr{S}}\,'_+\circ\tilde{\mu}$:
\begin{align}\label{bullettildes}
&\langle f,g\rangle_\bullet\ =\ 
\tilde{\mu}^{-1}
\hat{\mathscr{S}}\,'_+(\tilde{\mu}fg)(\hbox{id});
\end{align}
recall that $F(X)(\hbox{id})=F(q^\xi)$ for functions $F$ of $X$
and $\chi_{\hw}(\hbox{id})=\delta_{\hw,\hbox{\tiny id}}$.

The anti-involution of $\HH$ associated with 
$\lan\,,\,\ran_\bullet$ is
\begin{align}\label{diambul}
&\Diamond_\bullet:\,T_i\mapsto T_i(i\ge 0),\ \,X_a\mapsto X_a
(a\in P),\ \,
\Pi\ni\pi_r\mapsto \pi_r^{-1}.
\end{align}
See Section 3.2.2 from \cite{C101} and formula (3.9.4) from
Section 3.9.1; 
compare with the definition of $\Diamond$ from (\ref{diamondef}). 
The (ideal of the) module $\De_\xi$
and the anti-involution $\Diamond_\bullet$ satisfy the 
nonstrong Shapovalov property (for generic $q^\xi$).

\smallskip
\subsubsection{\sf Theorems \ref{TQ1H}, \ref{GENPROPR}
revisited} \label{sec:appl-to-conv}
The Shapovalov property of 
$\De_\xi$ and $\Diamond_\bullet$ 
guarantees that this module has a unique up to proportionality 
bilinear form associated with $\Diamond_\bullet$ 
(for sufficiently general $\xi$). 
Using $\hat{\mathscr{P}}\,'_+$ instead of 
$\hat{\mathscr{S}}\,'_+\circ\tilde{\mu}$ in (\ref{bullettildes}), 
one can establish the coefficient-wise proportionality of these 
operators. A direct usage of the 
divisibility argument as in Theorem \ref{KERIMAGE} can be now 
avoided;\, though it is of course present in this approach.
The justification goes as follows. 

\begin{theorem}\label{PSBULLET}
Let $\hat{\mathscr{P}}\,'_+=\sum_{\hw\in \hW}C_{\hw}\hw$
be the expansion from Theorem \ref{KERIMAGE},$(i)$
(see also Lemma \ref{LEMThatw} below). We set
$$\hat{\mathscr{P}}_+^\circledast\,=\,
\sum_{\hw\in \hW}\,C_{\hw}^\circledast\hw
\for C_{\hw}^\circledast\equal C_{\hw}/C_{\hbox{\tiny id}},
$$ 
where $C_{\hw},C_{\hw}^\circledast\,\in\, 
\Z[[t^{-1/2},\,X_{\al_i},\,i\ge 0]]$.
Then for $f,g\in \De_\xi$ (with the coefficient-wise multiplication), 
$$
\bigl(\hat{\mathscr{P}}_+^\circledast(fg)\bigr)(\hbox{id})\,
=\,\lan f,g\ran_\bullet\,=\,
\bigl(\tilde{\mu}^{-1}\hat{\mathscr{S}}\,'_+(\tilde{\mu}fg)\bigr)
(\hbox{id})\,,
$$
where the values are in the algebra
$\C[[t^{-1/2},\,q^{(\xi,\al_i)}\, ,i\ge 0]]$.
In particular, $\,C_{\hw}^\circledast(q^\xi)=\mu_\bullet(\hw^{-1})\,$
for any $\hw\in \hW$ (when $f=\chi_{\hw^{-1}}=g$ are taken). 
Thus, $\hat{\mathscr{P}}_+\,'$ and 
$\hat{\mathscr{S}}\,'_+\circ\tilde{\mu}$
are proportional to each other, which readily results in the 
proportionality claim from (\ref{propor}) or (\ref{proporan}).
\end{theorem}

{\em Proof.} We define the inner product $\lan f,g\ran'\,$ 
for $f,g\in \De_\xi$ using a direct counterpart of
(\ref{bullettildes}):
\begin{align}\label{bullettildep}
&\langle f,g\rangle'\ =\ 
\bigl(\hat{\mathscr{P}}^\circledast_+(fg)\bigr)(id)\,;
\end{align}
cf. formula (\ref{psymformula}).
Here 
$\hat{\mathscr{P}}\,'_+$ is considered 
as in its original definition from Theorem 
\ref{KERIMAGE}, i.e., with the coefficients 
$C_{\hw}\in \Z[[t^{-1/2},\,X_{\al_i}\, ,i\ge 0]]$ in
its decomposition $\sum_{\hw}C_{\hw}\hw$.
Then this series is applied to $fg\in \De_\xi$ and
finally the coefficient of $\de_{id}=\chi_{id}$ has
to be considered. The output will be a finite 
linear combination of proper $C_{\hw}$
evaluated at $q^\xi$, i.e., an element of
$\Z[[t^{-1/2},\,q^{(\xi,\al_i)}\, ,i\ge 0]]$.
We assume here that the coefficients of the expansion
of $fg$ in terms of $\chi_{\hw}$ are from $\Z$ or from
this algebra.

Due to formula (\ref{annihil}), 
\begin{align}\label{Pa-right}
\hat{\mathscr{P}}\,'_+\,T_{\hw}\ =\ 
t^{\frac{l(\hw)}{2}}\, \hat{\mathscr{P}}\,'_+\,.
\end{align}
Note that $T_{\hw}$ are placed here on the right, which is 
covered by part $(iii)$ of Theorem \ref{KERIMAGE}. 

Relation (\ref{Pa-right}) provides
that all images $\hat{\mathscr{P}}\,'_+(\de_{\hw})$ are
proportional to each other with certain constant coefficients of
proportionality. Thus, taking the evaluation at
any $\hu$ instead of  $id$ in
(\ref{bullettildep}) will not change this bilinear form 
up to proportionality.

These relations are sufficient to conclude that
the form $\lan f,g\ran'$ satisfies all properties of the form
from (\ref{bullettildes}). It can be checked directly, but
this can be avoided since we already know that the vanishing
conditions from (\ref{annihil}) are the same for 
$\hat{\mathscr{P}}\,'_+$ and for 
$\hat{\mathscr{S}}\,'_+\circ\tilde{\mu}$.

We see that $\lan f,g\ran'$ for
$f,g\in \De_\xi$
is $\lan\,,\,\ran_\bullet$ times a constant, which may, generally
speaking, depend on $\xi$. This constant is easy to find by 
taking $f=1=g$.
Here $1=\sum_{\hw}\chi_{\hw}$ is an infinite sum in $\De_\xi$, 
but the limits $\lan 1,1 \ran_\bullet$ and $\lan 1,1 \ran'$
are well defined.
\sq
\smallskip

Theorem \ref {PSBULLET} establishes that the coefficients
of the expansion of  $\hat{\mathscr{P}}\,'_+\,$
(in its initial definition from Theorem \ref{KERIMAGE}) 
are actually those of $\hat{\mathscr{S}}\,'_+\circ \tilde{\mu}\,$
up to a general (functional) coefficient of proportionality.
The coefficient of proportionality is immediate;
it is $ct(t^{-1})$ due to Macdonald.  

Moreover, there is a common radius of convergence 
of the coefficients of $\hat{\mathscr{P}}\,'_+$ 
with respect to $t^{-1}$, which depends only on the 
``first appearance" of the singularities in $ct(t^{-1})$ 
and readily results in the estimate $\Re k <1/h$, 
equivalently, $|t|\,>q^{1/h}\,.$ For such $t$ and $|q|<1$,
$$ 
\hat{\mathscr{P}}\,'_+=
ct(t^{-1})\,\hat{\mathscr{S}}\,'_+\circ \tilde{\mu}\,,
$$
which finalizes Theorem \ref{TQ1H}. 
We use that the coefficients of  
$\hat{\mathscr{S}}\,'_+\circ \tilde{\mu}$ are well defined 
for any $t$. 
\smallskip

{\em Theorem \ref{GENPROPR}.}
Similarly, the operator $\hat{\mathscr{S}}\,'_+\circ \tilde{\mu}$ 
converges for $|q|<1$ (any $t$) when 
applied to the functions from the spaces  
$\mathscr{X}q^{\,lx^{2}/2}$ for the levels $0<l\in \Z_+$. 
It is apart from potential $X$\~singularities, which are actually
not present due to cancelations of residues (see below).

The operator $\hat{\mathscr{P}}\,'_+$ acts there too; 
by construction, its images are
certain series multiplied by $q^{\,lx^{2}/2}$. The 
coefficient-wise
proportionality provides that these images are actually 
expansions of meromorphic $X$\~functions when  $|t|\,>q^{1/h}\,.$
This finalizes Theorem \ref{GENPROPR}.
\smallskip

As a byproduct, we obtain that 
$\hat{\mathscr{S}}\,'_+\circ \tilde{\mu}$ 
has no singularities when acting in
$\mathscr{X}q^{\,lx^{2}/2} \, (0<l\in \Z_+)$. A direct
justification of this (known) fact is by establishing
the cancelation of singularities, which is not needed now
due to the proportionality theorem. Indeed, it suffices 
to assume here that $|q|$ is small. 
Then $\hat{\mathscr{P}}\,'_+$ converges and has no singularities 
because it is defined in terms of the divided differences, which
preserve Laurent polynomials.

We note that given $t\in \C^*$ and $ b\in P$, 
it is not too difficult to check directly
that $\hat{\mathscr{P}}\,'_+(X_b q^{\,lx^{2}/2}) $
is an analytic function for $|q|<1$ and sufficiently
large $|t|$.

\smallskip
To conclude the convergence and proportionality matters, 
let us emphasize that there are two major approaches to the 
analysis of  $\hat{\mathscr{P}}\,'_+$. The first is based
on its $Y$\~{\em rational} presentation from formula
(\ref{hatPrat}), which, for instance,
results in Theorem \ref{YLEFT}. The second is
Theorem \ref {PSBULLET} (and its predecessors), which equates 
this operator with  
$ct(t^{-1})\hat{\mathscr{S}}\,'_+\circ \tilde{\mu}$ and then
the theory of the latter can be used.
\smallskip

\subsubsection{\sf Application to Theorem \ref{YLEFT}}
A similar approach can be used to finalize
Theorem \ref{YLEFT}. We will prove here that 
the convergence and vanishing assumptions in this theorem hold
when $\, \hat{\mathscr{P}}\,'_+$,\,
$\hat{\mathscr{S}}\,'_+\circ \tilde{\mu}\,,$
$\Si_{\infty}$\, 
and other operators involved
are understood coefficient-wise.  

The coefficients will be treated as the elements of the algebra
\begin{align}\label{ZX-prime}
&\Z\!X_+\equal 
\Z[[t^{-1},\,X_{\al_i}\, ,i\ge 0]]\,;
\end{align}
recall that it contains all positive powers of $q$.
We will use  $\hat{P}(t)$, the affine Poincar\'e series 
from (\ref{rataffpoin}). By $b\to \infty$, we mean that
$b\in P_+$ and $(b,\al_i)\to \infty$ for all $i>0$.
Similarly,\,  
$\mathbf b\to\infty$ for a set $\mathbf b=\{b^j\}$
if and only if $b^j\to \infty$ for all $j$.

\begin{theorem}\label{YLEFTNEW}
Given a system of representatives 
$\mathbf b=\{b^1,\ldots,b^p\}\subset P_+^\vee$ for 
the group $\Pi=P^\vee/Q^\vee$
(of cardinality $p$), we set
\begin{align}\label{SiYgennew}
&\Sigma_{\mathbf b}\, =\, 
\frac{1}{|\Pi|}\sum_{j=1}^p\, t^{-(b^j,\rho)}Y_{b^j},
\ \,\Sigma_\infty\,  =\, \frac{1}{|\Pi|}\lim_{\mathbf b\to\infty}\,
\sum_{j=1}^p t^{-(b^j,\rho)}\,Y_{b^j}.
\end{align}

(i) Given $W\ni u\neq \hbox{id}\,$
and $\hw\in \hW$, there 
exists a constant $\de=\de_{\hw}>0$ such that
for all  $\mathbf b$ 
sufficiently close to $\infty$, 
\begin{align}\label{tYlimdefw} 
&C_{\hw}^u\,\in \, q^v \,\Z\!X_+ \for v>
\de\,(\sum_{j=1}^p(b^j,\rho)) \hbox{,\ \, where }
\notag \\
&\Sigma_{\infty}^u\, \equal\, \frac{1}{|\Pi|}
\sum_{j=1}^p \lim_{\mathbf b\to\infty}
t^{-(b^j,\rho)}Y_{u(b^j)}\,=\sum_{\hw\in \hW} C^u_{\hw} \hw.
\end{align}

(ii) The limit $\Sigma_\infty\,$
exists as a series $\sum_{\hw\in \hW} C_{\hw}\,\hw$
with the coefficients $C_{\hw}=C_{\hw}^{\hbox{\tiny id}}$ 
in the algebra $\Z\!X_+$.
Recall that we replace $Y_b$ by the corresponding operators
acting in the polynomial representation, move $\hw$ to the
right and then expand the resulting $X$\~rational coefficients 
of $\hw$ in terms of $X_{\al_i} (i\ge 0)$.

Given $\hw\in \hW$ and a compact subset belonging 
to $\{0\neq X_\al\not\in q^{\Z}, \al\in R\}$, 
the coefficients $C_{\hw}^u$ converges uniformly
in this subset provided that $|t|>1$ and $|q|$ is
sufficiently small (depending on $|t|$ and this subset);
moreover, $C_{\hw}^u\to 0$ for $u\neq 0.$

(iii) Treating the $C$\~coefficients as
elements from $\Z\!X_+$,
\begin{align}\label{tYliminv} 
&\Sigma_\infty Y_a\  =\ \ 
t^{(a,\rho)}\,\Sigma_\infty \ \for\, a\in P,\notag\\ 
&\Sigma_\infty T_{\hw}\ =\ 
t^{l(\hw)/2}\,\Sigma_\infty \for \hw\in \hW.
\end{align}
These identities formally result in
\begin{align}\label{tYliminvS} 
\Sigma_\infty\ =\ \hat{\mathscr{P}}_+\ =\ 
(ct(t^{-1})/\hat{P}(t^{-1}))\,
\hat{\mathscr{S}}\,'_+\circ \tilde{\mu}\,. 
\end{align}

(iv) We continue (ii) and (iii). For $0\le |q|<1$ and
$X$ from a given compact subset of 
$\{0\neq X_\al\not\in q^{\Z}, \al\in R\}$, the condition
$|t|>1$ is sufficient for the uniform point-wise 
convergence of the coefficients $C_{\hw}^u$ of 
$\Sigma_\infty^u$; the convergence is to $0$ for $u\neq$id.

Correspondingly, the coefficients of the 
$\hw$\~expansions in the identities from (\ref{tYliminvS}) coincide 
point-wise provided that $0\le |q|, |t|^{-1}<1$ subject to
$\{0\neq X_\al\not\in q^{\Z}, \al\in R\}$.
\end{theorem}
{\em Proof.}
We will use the following presentation  
of $T_{\hu}$ ($\hu\in \hW$) acting in $\mathscr{X}$,
which is especially convenient  
for $Y_b=T_b\, (b\in P_+)$. Let  
\begin{align}\label{Gtaldef}
&G_{\tal} \ \equal\ 1+\frac{1-t^{-1}}
{X_{\tal}^{-1}-1}\,(1-s_{\tal})\ =\ 
\frac{X_{\tal}^{-1}-t^{-1}}{X_{\tal}^{-1}-1}-
\frac{1-t^{-1}}{X_{\tal}^{-1}-1}\,s_{\tal},\\ 
&G'_{\tal}\ \equal\ G_{\tal}(X\mapsto X^{-1})\ =\
G_{-\tal}\ =\ \frac{X_{\tal}-t^{-1}}{X_{\tal}-1}-
\frac{1-t^{-1}}{X_{\tal}-1}\,s_{\tal}\notag
\end{align}
for $\tal\in \tilde{R}$;\ 
recall that $X_{\tal}=X_{\al}q^j$\, for $\tal=[\al,j]$.

Given  a reduced decomposition 
$b= \pi_r s_{j_l}\cdots s_{j_1}$ for $\ l=l(b)$
and $r\in O$,
one has 
\begin{align}\label{prodforgb}
&t^{-(\rho,b)}Y_b\ =\ b\, G_{\tal^{\,l}}\cdots G_{\tal^1}=
\ G'_{\tbe^{\,l}}\cdots G'_{\tbe^1}\,b,\\
&\hbox{for\ \ }  
\tal^1=\al_{j_1},
\tal^2=s_{j_1}(\al_{j_2}),\,\ldots\,,\, 
\tbe^r=-b(\tal^{\,r})\in \tR_+.
\notag 
\end{align} 
The set $\{\tal^1,\tal^2,\, \ldots\,\}=\La(b)\subset \tR_+$
is the $\La$\~set defined in (\ref{Lahw}). Note that 
$b s_{\tal^l}\,\cdots\, s_{\tal^1}=$ $\pi_r.$
See e.g., \cite{Ch13}, where the main construction is
actually close to the limits we consider here.

Here and below we need some basic properties of the
$\La$\~sets and the Bruhat ordering in $\hW$. 

Only $X_{\tbe}$ with positive $\tbe$ are present in the
terms of the $G'$\~product from (\ref{prodforgb}), but
negative roots do appear when calculating
its (right) $\hw$\~expansion when the terms with $s_{\tal}$
are taken from the corresponding binomials. The resulting
$\{\hw\}$ will form the Bruhat set for $b$ (pure $b$
is obtained if no single $s_{\tal}$ is taken).

As always in this paper, we 
expand the resulting $X$\~rational coefficients in terms of
$X_{\al_i}$ $(i\ge 0)$ and can readily check that 
they are actually from  $\Z\!X_+$ for any given
$\hu\in \hW$. The problem is to justify the existence of
these coefficiients when $l(\hu)\to \infty$.

Claim $(ii)$ is the key in this theorem;
it will be deduced from Theorem \ref{YLEFT}, where
$\Si_{\infty\,}^+$ was obtained from $\hat{\mathscr{P}}\,'_+$ 
assuming $(i)$. The fact that the $C$\~coefficients 
of $\Si_{\infty\,}$ converge cannot be justified at the moment 
directly from (\ref{prodforgb}), at least
for arbitrary root systems. Potentially, there can be terms
in the resulting summation destroying the convergence; 
they cancel each other, which follows 
from Theorem \ref{YLEFT}.

\smallskip

{\em Claim $(i)$.\,} We need the following modification
of (\ref{prodforgb}).

Given a reduced decomposition $\hu=\pi_r s_{j_l}\cdots s_{j_1}$,
let
\begin{align}\label{prodforg}
&t^{-l(\hu)/2}T_{\hu}^{-1}\ =
\tilde{G}_{\tal^1}\cdots \tilde{G}_{\tal^l}\, \hu^{-1},\\
&\tilde{G}_{\tal}=
\frac{1-t^{-1}X_{\tal}^{-1}}{1-X_{\tal}^{-1}}+
\frac{1-t^{-1}}{1-X_{\tal}^{-1}}s_{\tal}.\notag
\end{align}
Note that $\tal^1=\al_{j_1}$, $\tal^2=s_{j_1}(\al_{j_2})$,
$\tal^3=s_{j_1}s_{j_2}(\al_{j_3})$ and so on constitute
the set $\La(\hu)$.

First, it is simple
to calculate the ``greatest" $C$\~coefficient in this
expression, which is that of $\hu^{-1}$. It can be obtained only
by picking the terms without $\,s_{\tal}\,$ from all
binomials in (\ref{prodforg}). Thus, 
\begin{align}\label{cminbt}
C_{\hu^{-1}}\  =\ \prod_{[\al,j]\in \La(\hu)}
\frac{t^{-1}-q^j X_{\al}}{1-q^j X_{\al}}\ =\ 
t^{-l(\hu)}\!\!\prod_{[\tal]\in \La(\hu)}
\frac{1-t X_{\tal}}{1-X_{\tal}}\,.
\end{align}
It readily converges to zero in the sense of $(i)$, i.e.,
it will become divisible in $\Z\!X_+$ by growing positive powers 
either of $t^{-1}$ or of $q$ as $l(\hu)\to \infty$.
Recall that we expand the
denominators here and in any other products in terms of 
nonnegative powers of $X_{\al_i}$ for $i\ge 0$.

The case of fixed (bounded) $\hw$ as $l(\hu)\to \infty$ is, 
in a sense, opposite to this example. The following lemma
addresses it. 

\subsubsection{\sf Combinatorics of C-coefficients}
Let us examine the individual 
products contributing to the coefficients $C_{\hw}$ in the 
standard decomposition
$$
t^{-l(\hu)/2}T_{\hu}^{-1}=\sum_{\hw\in \hW} C_{\hw} \hw.
$$

\begin{lemma}\label{LEMThatw}
(a) There exists a constant $\de_{\hbox{\tiny total}}\,>\,0$ such that 
for any $\hu\in \hW \ni\hw $ all such individual products belong to
$$
(q^v+t^{-v})\,\Z\!X_+ \for v\,>\,\de_{\hbox{\tiny total}}\ l(\hu).
$$

(b) Given $\,\hw\in \hW$, there exists a constant $\de_{\hw}>0$
such that for any $\hu\in P_+$, the corresponding individual
products from (a) belong to 
$$
q^v\,\Z\!X_+ \for v\,>\,\de_{\hw}\, l(\hu).
$$

(c) Given \ id\ $\neq u\in W$,
the same holds for the standard decomposition of 
\ $t^{-(\rho,b)}Y_{u(b)}^{-1}\,$, where $b\in P_+$ and
$l(\hu)$ is replaced by $l(b)=2(\rho,b)$;\,
we assume that $b\to\infty$, i.e., $(b,\al_i)\to\infty$
for all $i>0$. 

(d) Claims from (b,c) hold when the algebra $\Z\!X_+$ from
(\ref{ZX-prime}) is changed to algebra
$\Z\!X'_+\equal 
\Z[[\,t'\, ,\ X_{\al_i}\, ,i\ge 0\,]] \for  t'\equal 1-t^{-1}$,
i.e., when we expand the coefficients at the point $\,t=1\,$ 
instead of $\,t=0$.

(e) The $C$\~coefficients of the decomposition
$\hat{\mathscr{P}}\,'_+\ =\ \sum_{\hw\in \hW}\,C_{\hw}\,\hw
$
are well-defined elements of $\Z\!X_+$ or $\Z\!X'$.
Moreover, they belong to
$$
\sum_{j=0}^\infty q^j\,
\Bigl(\,\Z[[\,X_{\al_i}\, ,i\ge 0\,]]\,[t^{-1}]\,\Bigr)\ \subset\
\Z\!X_+\,\cap\,\Z\!X_+'.
$$
\end{lemma}
{\em Proof.}
The smallest possible $\hw$ that can be obtained
from $\hu$ is when we always
pick the terms with $\,s_{\tal}\,$ from the binomials
in the product (\ref{prodforg}) for $\,t^{-l(\hu)/2}T_{\hu}^{-1}\,.$
It will contribute to the $C$\~coefficient of $\pi_r^{-1}$ 
(maximally distant from $\hw^{-1}$).  There can be of course
other products that contribute to $C_{\pi_r^{-1}}\,$; 
their number grows exponentially in terms of $l(\hu)$.
The corresponding product is as follows:
\begin{align}\label{prodforpi}
\prod_{r=1}^l
\frac{1-t^{-1}}{1-X_{\al_{j_r}}^{-1}}\ &=\ 
(-q/X_\th)^{l_0}\frac{(1-t^{-1})^l}{(1-q/X_{\th})^{l_0}}
\,\prod_{j_r\neq 0}\frac{1}{1-X_{\al_{j_r}}^{-1}}\notag\\
&=\ (-1)^l\,\frac{(1-t^{-1})^l(q/X_\th)^{l_0}}{(1-q/X_{\th})^{l_0}}
\,\prod_{j_r\neq 0}\frac{X_{\al_{j_r}}}{1-X_{\al_{j_r}}},
\end{align}
where $l_0$ is the number of indices $j_r=0\ (1\le r\le l)$
in the reduced decomposition of $\hu$; recall that $\al_j$ are
simple roots. 

Obviously, this product satisfies $(a)$. 
Moreover, its minimal power of $q$ will grow linearly 
with respect to $l(\hu)$. We use that for any given $i\ge 0$, 
the number of $j_r$ such that $j_r=i$ must grow linearly 
with $l(\hw)$. Indeed, if in a certain connected portion 
of the reduced decomposition of $\hu$, the simple reflection
$s_i$ is missing, then this portion comes from a finite Weyl 
subgroup of $\hW$ (for instance, $W$ for $i=0$). Therefore, 
the maximal possible length of such a segment is bounded by the 
length $l(w_0)$ of the element $w_0$ of the maximal length in $W$. 

We conclude that the number of $X_{\th}^{-1}$ from 
$(q/X_{\al_0})^{l_0}$ that we can terminate using $X_{\al_{j_r}}$ 
will grow linearly
together with $l(\hu)$. This will release the power of $q$
growing linearly with respect to $l(\hu)$.

\smallskip

{\em Omitting some $s_{\tal}$.\,} 
Let us take now one $\tilde{G}_{\tal^p}$
for some $p$  in (\ref{prodforg}) and pick the term there 
without 
$s_{\tal}$ for $\tal=\tal^p$; anything else remains unchanged.
The corresponding contribution will be
to the coefficient  $C_{\hw}$ for $\hw= s_{j_p}\pi_r^{-1}$. 
It is a pure product equal to 
\begin{align}\label{prodforpi1}
&\Bigl(\,\prod_{r=1}^{p-1}
\frac{1-t^{-1}}{1-X_{\al_{j_r}}^{-1}}\,\Bigr)\,
\frac{1-t^{-1}X_{\al_{j_p}}^{-1}}{1-X_{\al_{j_p}}^{-1}}\,
\Bigl(\,\prod_{r=p+1}^{l}
\frac{1-t^{-1}}{1-X_{\be_{r}}^{-1}}\,\Bigr),\notag\\
&\ \ \ \ \ \ \ \hbox{\ where\ } \be_{r}\ =\ 
s_{j_p}(\al_{j_r})\, \for r>p\,. 
\end{align}
Unless $j_p=0$, the estimate of the $q$\~power is completely 
parallel to that for (\ref{prodforpi}). 

If $j_p=0$, then the indices $\{j_r=0\}$
{\em after} $j_p$ will not contribute any longer to the
total power of $q$, since $s_0(\al_0)=-\al_0$. However,
$s_0(\al_i)=\al_0+\al_i$ for simple $\,\al_i\,(i>0)\,$ 
neighboring to $\al_0$ in the completed Dynkin diagram. 
The number of such indices 
$\,i\,$ in the reduced decomposition of $\hu$ will tend to 
infinity together with $l(\hu)\to\infty$. The corresponding
$X_{\al_i-\th}$ in the numerator 
can be terminated by using nonaffine $X_{\al}$ with $\al>0$ 
exactly in the same way as it was done in (\ref{prodforpi})
for $X_{-\th}$.
The released $q$ will provide the required growth of the total
power of $q$. 
\smallskip

If there are two places $p<p'$ where the terms without
$s_{\tal}$ are taken, then the resulting product will 
contribute to $C_{\hw}$ for 
$\hw=s_{j_p}s_{j_{p'}}\pi_r^{-1}\,;\,$ it reads
\begin{align}\label{prodforpi2}
&\Bigl(\,\prod_{r=1}^{p-1}
\frac{1-t^{-1}}{1-X_{\al_{j_r}}^{-1}}\,\Bigr)\,
\frac{1-t^{-1}X_{\al_{j_p}}^{-1}}{1-X_{\al_{j_p}}^{-1}}\,
\Bigl(\,\prod_{r=p+1}^{p'-1}
\frac{1-t^{-1}}{1-X_{\be_{r}}^{-1}}\,\Bigr)\\
&\ \ \times\,\frac{1-t^{-1}X_{\be_{p'}}^{-1}}{1-X_{\be_{p'}}^{-1}}\,
\Bigl(\,\prod_{r=p'+1}^{l}
\frac{1-t^{-1}}{1-X_{\be_r}^{-1}}\,\Bigr)\,, \where
\notag\\
\be_{r}=&s_{j_p}(\al_{j_r}) \hbox{\, if \,} p'\ge r>p\,,\
\be_{r}=s_{j_p}s_{j_{p'}}(\al_{j_r}) \hbox{\, if \,} r>p'\,.\notag
\end{align}
For the sake of uniformity, we will replace here 
the remaining $\al_{j_r}$ by $\be_r$ as well, setting 
$\beta_r=\al_{i_r}$ for $p\ge r >0$.
\smallskip

{\em Minimal $q$\~powers\,.}
The analysis of the minimal power of $\,q\,$
remains essentially the same in the case of two $p$.
The number of the indices $r>p\,$ with {\em affine negative\,} 
$\be_r\,$ (they do not
contribute to the minimal power of $q$) {\em approximately}, i.e.,
in the limit $l(\hu)\to\infty$\,, is no smaller than the number 
of positive affine $\be_r$ (which do contribute). It can be readily 
generalized to any number of indices $\,p\,$ such that
the corresponding terms without $s_{\tal}$ are taken.
Let $\hw$ will be the corresponding 
index of the $C$\~coefficient; 
$\hw=s_{j_p}s_{j_{p'}}\pi_r^{-1}$ for two $p$.
\smallskip

We can assume that $\hw(\al_0)<0$ for $\hw\in \hW$.
Indeed, if the reduced decomposition of $\hu\,$ grows to infinity 
after $\{p\}$, then we can assume that $\hw(\al_0)<0$ for 
$\hw\in \hW$; otherwise such $\hw$ will not change the positivity
of $X_{\al_0}$ after $\{p\}$. If such a growing interval in
the reduced decomposition of $\hu$ occurs between some $p$, 
we can diminish $\{p\}$ to end this sequence before this interval. 
The positivity of the terms before $\{p\}$ remains
unchanged, which provides the required power of $q$ if such
growing interval occurs before $\{p\}$.  
 
Representing $\hw=v a$ for $v\in W,\,a\in P$, the condition
$\hw(\al_0)<0$ can happen only if $(a,\th)<0$. Indeed,
$\hw(\al_0)=[-v(\th),1+(a,\th)]$. However,
$\th$ is a sum of simple roots with positive coefficients. 
Thus $(a,\al_j)=-d<0$ for at least one $j>0$ and
$\hw(\al_j)=[v(\al_j),d]>0$.

Finally, $d$ here will tend to $\infty$
together with $l(\hu)$ because, as we already used, 
the number of $s_j$  (for any given $j>0$) in the reduced 
decomposition of $\hu$ grows linearly with respect to $l(\hu)$.  
\smallskip

{\em Omitting many $p$\,.} Let 
$\mathbf p=\{\cdots > p''>p'>p\,\}\,$ be the sequence  
of the terms (binomials) where we omit the corresponding  $s_{j_p}\,$.
If $\mathbf p=\La(\hu)$, then we arrive at 
(\ref{cminbt}); however, now we are interested in the case when
the corresponding $\hw\,$ remains bounded.

One has
$$ 
\be_p=\al_{j_p},\ 
\be_{p'}=s_{j_p}(\al_{j_{p'}}),\
\be_{p''}=s_{j_p}s_{j_{p'}}(\al_{j_{p''}})
\hbox{\,,\ and so on}.
$$ 

The corresponding product will contribute to 
the coefficient $C_{\hw}$ for 
$\hw=s_{j_{p}}s_{j_{p'}}s_{j_{p''}}\,\cdots\,\pi_r^{-1}$.
The set $\La(\hw^{-1})$ is very explicit;
it is obtained from the set 
$\mathbf{\be}=\{\,\ldots \be_{p''},\be_{p'},\be_{p}\,\}\,$ by 
removing all pairs in this set in the form $\,\{\tal,-\tal\}\,.$ 

The main problem we have to address is that $\hw$ can be small 
for arbitrarily large $\hu$. This is actually the key point
of the justification of existence of $\hat{\mathscr{P}}^{\,'}_+$
and other operators under consideration. Given $\hw$, the number of 
contributions to $C_{\hw}$ of this kind will go to $\infty$
together with $l(\hu)$. 

\smallskip
For instance, one can take $\,\hu=(-b_+)w_0 b_+\,$ for any
$b_+\in P_+$ such that $-b_+=w_0(b_+)$. Then the corresponding length 
will be $l(\hu)=2\,l(b_+)+l(w_0)$, but if we delete $w_0$ it 
will drop to zero.  Actually, this is a typical example.
The corresponding product will be that from (\ref{prodforpi2}) 
with only one group in the
parentheses, corresponding to $w_0$, and with the products
before and after it (without the parentheses)
corresponding to $b_+$ and $(-b_+).$
\smallskip

Assuming that $|\mathbf p|$ is large and $l(\hw)$ is bounded by
a certain constant, almost all elements
of $\Lambda(\hw^{-1})$  will appear in the pairs $\{\tal,-\tal\}$ 
for $\tal=[\al,j]>0.$ Any such a pair will
contribute either $t^{-1}$ or at least $q^{j-1}$ to the 
resulting product. Indeed, the product of the corresponding 
quantities will be (before the expansion in terms 
$X_{\tbe}$ with $\tbe>0$)
$$
\frac{(1-t^{-1}X_{\tal})(t^{-1}-X_{\tal})}{(1-X_{\tal})^2}\,.
$$ 
See (\ref{prodforpi2}), the terms there without the parentheses. 
This concludes $(a)$, but the resulting powers of $q$ 
can be estimated better than needed in $(a)$, which is part 
$(b)$ of the lemma.
Note that this argument is not applicable to $(d)$, though the 
estimates for the $q$\~powers below hold in this case.
\smallskip

{\em Part (b).}
Since the original decomposition of $\hu$ was reduced,
there will be $\al_{j_r}$ for $r$ with positive $\be_r$ 
between some of $p$ with their affine components
approaching $\infty$. Indeed,
if all these affine components
remain bounded, then $l(\hw)\to \infty$
together with $l(\hu)$.
\smallskip

For instance, in the example of $\,\hu=(-b_+)w_0 b_+\,$ with 
$b_+\in P_+\,$, the element $(-b_+)$ must be on the left and
$b_+$ on the right to ensure that the corresponding combined
decomposition is reduced.
Thus $b_+^{-1}(\al_i)=[\al_i, (b_+,\al_i)]$ and there must 
be growing positive affine components at least for some $i>0$. 

To demonstrate the essence of our estimates in the case of large
$|\mathbf b|$ with bounded $l(\hw)$, let us insert here any 
$v\in W$ instead of $w_0$. 
For $c=v(b), \,b\in P_+\,$, the length of $\hu=c\,v\, b=v\cdot(2b)$ 
is $l(c)+l(v)+l(b)$. We pick the terms with $s_{\tal}$ only from
$v$ here (the corresponding portion of the reduced decomposition
of $\hu$) assuming that $v(b)+b$ is bounded. 
The set $\La(v)$ is formed by certain positive linear combinations
of $\al_{i_r}$ with {\em nonnegative} integral coefficients 
for the indices $I_v=\{i_r\}$ from a given reduced decomposition 
$v=s_{i_m}\cdots s_{i_1}$. 

One has $b^{-1}(\La(v))\in \tR_+$.
If all scalar products $(b,\al_i)$ here are no greater than
a certain constant for all $i\in I_v$, then there must
exist an index $0<k\not\in I_v$  such that $M=(b,\al_k)$ is 
positive and large compared to $(b,\rho)$.  
Representing $b=b'+b''$ for $b'=M\om_k$, we can assume that 
$v(b')=b'$; otherwise $M$ will contribute to the power of
$q$ (the terms inside the parentheses) and we will obtain the 
required growth. However the relation $v(b')=b'$ readily
contradicts the assumption that $v(b)+b$ is bounded. 
\smallskip

{\em The general case.} We need to examine the products in the 
form $\hu=\hu^*\,\hv\,\hu'$, where $\hv\in \hW \ni \hu',\hu^*$, 
such that 
$l(\hu^*\,\hv\,\hu')=l(\hu^*)+l(\hv)+l(\hu')$ and 
the set $\La(\hu^*\,\hu')$ remains
bounded as $l(\hu^*\,\hv\,\hu')\to \infty.$

It is a generalization of the example considered
above, where $\hv$ substitutes for $v$ and 
$\hu'$ replaces $b\in P_+$.
 
Let $I_{\hv}=\{i_r\}$ for a given reduced decomposition 
$\hv=s_{i_m}\cdots s_{i_1}\,$. We set 
$\,\hu=a u\, (a\in P, u\in W)$, $a=\sum_{i=1}^nM_i \om_i$
and $a=a'+a''$ for $a'=\sum_k M_k\om_k$ for the set
$K=\{k\}$ of all indices $i$ such that $|M_i|\to \infty$. 
Using that 
$a^{-1}(\La(\hv))\in \tR_+$, we can check by induction that
$K\cap I_{\hv}=\emptyset$ unless the power of $q$ tends
to infinity together with  $l(\hu^*\,\hv\,\hu')$ (the fact we need
to establish). Let us demonstrate it.

Indeed, for the first appearance
of $i$ in the sequence $I_{\hv}$, the root
$(a')^{-1}(\al_i)$ and the inner product $(a',\al_i)$ 
must be nonnegative to ensure the positivity of
$(\hu')^{-1}(\al_i)\in \La(\hv\,\hu')$. This holds
assuming that we already know that 
$(a',\al_{i'})=0$ for all previous $i'$ in
$I_{\hv}$. The corresponding contribution
to the resulting power of $q$ will be
$(a',\al_i)\to \infty$ 
as $l(\hu^*\,\hv\,\hu')\to \infty$ 
if $(a',\al_{i})>0$. Therefore,
$(a',\al_i)=0$. Check that here $\al_i$
can be allowed to be $\al_0$. Thus we conclude 
that $\hv(a')=a'$.
\smallskip

To finalize this reasoning (and part $(b)$), one
can assume that $\hu^*$ in $\hu^*\,\hv\,\hu'$ can be 
represented as $\hu^*=\hw^*\, b^*$
for $b^*\in P, \hw^*\in\hW$ such that $l(\hu^*)=l(\hw^*)+l(b^*)$
and the sum $b^*+a'$ remains bounded in the limit. However then 
$b^*\hv a'$ cannot be reduced. This contradiction shows that
the powers of $q$ in the products  under consideration
(contributing to $C_{\hu^*\hu'}$ for bounded $\hu^*\hu'$)
go to infinity as  $l(\hu^*\,\hv\,\hu')\to \infty$.
\smallskip

\rmk
We note that continuing the (combinatorial) analysis of the
products contributing to the coefficients 
$C_{\hw}$ in the decomposition of $t^{-l(\hu)/2}T_{\hu}^{-1}$
for fixed $\hw$ and growing $\hu\in \hW$, one can eventually  
arrive at the sharp convergence range $|t|>q^{1/h}$ for
the $C$\~coefficients of the operator  
$\hat{\mathscr{P}}^{\,'}_+$; see $(e)$.  We will not demonstrate this
here, since it formally results from 
the proportionality of this operator with $\hat{\mathscr{I}}$,
where sufficiently small $|q|$ and $|t|^{-1}$ are sufficient
(the coefficients of $\hat{\mathscr{I}}$ are explicit).
Let us mention that the Coxeter number $h$ appears
in our considerations due to counting the ``density" of $s_{0}$
in the reduced decomposition of $\hu$ and in similar estimates.
\sq
\smallskip


{\em Part (c).\,}
Let us apply the $q$\~estimates from $(a,b)$ to $Y_{u(b)}$ 
for id$\neq u\in W$, $b\in P_+$. 
The above analysis of the expansion
$t^{-(\rho,b)}Y_b^{-1}=\sum_{\hw\in \hW} C_{\hw}\,\hw$\
for $b\in P_+$\, corresponds to $u=w_0$ and 
can be readily extended to arbitrary such $u$. 
The main step is as follows.

Let $c=u(b)$ for $b\in P_+$. Then
$l(c)=l(b)$ and the representation
$c= b'-b''$ for $b',b''\in P_+$ 
results in the following.
Setting
$c=u(b)=\pi_r s_{j_l}\,\cdots\, s_{j_1}\ (l=l(c))\,$, 
\begin{align}\label{prodforgc}
&Y_c\ =\ \pi_r T_{j_l}^{\ep_l}\,\cdots T_{j_1}^{\ep_1}\,,
\where \,\ep_j=\pm 1;\hbox{\ \, correspondingly\,,}\notag\\ 
&t^{-(\rho,b)}Y_c = c\,\hat{G}_{\tal^l}\,\cdots\, 
\hat{G}_{\tal^1}
\for \tal^1=\al_{j_1}, \tal^2=s_{j_1}(\al_{j_2}),\,\ldots\,,\\
&\hbox{where\, }\hat{G}_{\tal^j}=G_{\tal^j}\, \hbox{ for }\ep_j=+1\, 
\hbox{ and } 
\hat{G}_{\tal^j}=\tilde{G}_{\tal^j} \hbox{\, otherwise}.\notag
\end{align}

Then we move $\,c\,$ to the right and focus on
the terms $\tilde{G}_{\tal^j}$ for $\ep_j=-1$;
the condition $b\to \infty$ guarantees complete analogy
with the considerations for $c=-b\, (b\in P_+)$. 
\smallskip

{\em Parts (d,e).} The estimates for the power of $q$,
and therefore
the claims from $(b)$ and $(c)$ hold
when we replace $t^{-1}$ by $1-t'$ for $t'=1-t^{-1}$
and analyze the resulting expressions. This is 
claimed in $(d)$. It provides the existence
of the coefficients of  $\hat{\mathscr{P}}^{\,'}_+$ and
those of $\Sigma_{\infty}$ as formal series in terms 
of $q$ an $t'$.

Part $(e)$ follows from  $(b)$ and the observation (obvious
from its justification) that the power of $t^{-1}$ is bounded in 
the products contributing to $C_{\hw}$ unless they are 
divisible in $\Z\!X_+$ by powers of $q$ approaching $\infty\,$.
\sq
\smallskip

\subsubsection{\sf Back to Theorem \ref{YLEFTNEW}} 
{\em Claim $(ii)$.} Combining $(i)$ and
Theorem \ref{YLEFT}, we obtain that 
$\Si_\infty^+=\Si_\infty\,\mathscr{P}_+$
exists as a series with coefficients
in $\Z\!X_+$. A justification of this
fact based directly on (\ref{prodforgb}) is not known
at the moment (at least for arbitrary root systems).

Recall that the connection of $\hat{\mathscr{P}}\,'_+$ 
and $\Si_{\infty\,}^+$ is a sequence of algebraic 
manipulations based on the fact that 
$\hat{\mathscr{P}}\,'_+$ treated as a rational
function is identically zero. See (\ref{hatPrat})
and also Theorem \ref{SYMRANK1} below (the case
of $A_1$). 

Since $\Si_{\infty\,}$ is given in terms of $b^j$
representing all elements in $\Pi$, this formally
results in the existence of 
$(\Si_{\infty\,}^{\pi})^+=\Si_{\infty\,}^{\pi}\mathscr{P}_+$
for any $\pi=\pi_r\in \Pi$, where   
\begin{align*}
\Si_{\infty\,}^{\pi}\equal
\lim_{b\to \infty}\, t^{-(\rho,b)}Y_{b} \for 
b \hbox{\ \ such that\ \ } b-\om_r\in Q.
\end{align*}

We will use the nonaffine $(i>0)$ intertwining
operators; cf. Section \ref{sect:intertw}. One has
\begin{align}\label{y-intert}
&\Phi_i Y_b=Y_{s_i(b)}\Phi_i\, \for \,
\Phi_i\equal \frac{T_i+(t^{1/2}-t^{-1/2})/(Y_{\al_i}^{-1}-1)}
{t^{1/2}+(t^{1/2}-t^{-1/2})/(t^{-(\rho,\al_i)}-1)},\\
&\mathscr{P}_+\ =\ \sum_{u\in W} \Phi_{u}, \where\,
\Phi_{uv}=\Phi_{u}\Phi_{v}\, (u,v\in W), \ \,
\Phi_i=\Phi_{s_i}. \label{phivias}
\end{align}
Also, $\mathscr{P}_+\,$ is divisible by
$\sum_{u\in W} u$ on the left; see (\ref{mac1}).

Applying $(i)$, 
\begin{align*}
(\Si_{\infty\,}^{\pi})^+=(\tilde{\Si}_{\infty\,}^{\pi})^+\
\for \tilde{\Si}_{\infty\,}^{\pi}\equal
\lim_{b\to \infty}\, t^{-(\rho,b)}(\sum_{u\in W}Y_{u(b)}),
\ \,b\in \om_r+ Q.
\end{align*}
Moreover, representing $\mathscr{P}_+$ via the
$Y$\~intertwiners, we can place it on the left: 
\begin{align}\label{sipinftyx}
(\Si_{\infty\,}^{\pi})^+\ =\ {}^+\tilde{\Si}_{\infty\,}^{\pi}\,\equal\,
\mathscr{P}_+\,
\lim_{b\to \infty}\, t^{-(\rho,b)}(\sum_{u\in W}Y_{u(b)})
\end{align}
for $b$ from $\om_r+Q$.

Since $\mathscr{P}_+\,$ is divisible by $\,\sum_{u\in W} u\,$
on the left, the $C$\~coefficients of 
${}^+\tilde{\Si}_{\infty\,}^{\pi}$ must satisfy the $W$\~invariance
relations $C_{u\hw}=C_{\hw}$ for $u\in W$ and $\hw\in \hW$.
The $C$\~coefficients of sums $\sum_{u\in W}Y_{u(b)}$ have
the same invariance condition up to the terms from  $q^N\Z\!X_+$ for  
$N$ growing together with $(\rho,b)$. Use the $Y$\~intertwiners
and part $(i)$ (see also part $(iii)$ below). 

However,  
$\mathscr{P}_+\,$ for generic $t$ has no kernel when 
acting in the space of $W$\~invariant delta function
defined as follows: 
$\de_a(X_{w(c)})=q^{(a,c)}$ for $a,c\in P_+, \, w\in W$. Therefore,
the $C$\~coefficients of $\sum_{u\in W}Y_{u(b)}$ can be 
uniquely recovered from those of  
$\mathscr{P}_+\,\sum{u\in W}Y_{u(b)}$ modulo
$q^N\Z\!X_+$. This provides the existence of
$\Si_{\infty\,}^{\pi}$ and justifies the first part of $(ii)$.

\smallskip

{\em Claim $(iii)$}. Actually, we have already used 
the main arguments needed here in $(ii)$ above. Nevertheless, let 
us see how the proportionality claims can be obtained
directly from the existence of $\Si_{\infty\,}$.

The first of the formulas from
(\ref{tYliminv}) results from $(ii)$. Let us 
demonstrate that
the second follows directly from $(i)$.
Using (\ref{y-intert}),
$\Si_\infty \Phi_i =
\Phi_i \Si^{s_i}_\infty=0$ upon the $\hw$\~expansions
with $C_{\hw}$ treated
as a formal series with the coefficients in $\Z\!X_+$ (or 
point-wise for sufficiently small $|q|$).
Therefore 
$$
\Si_\infty T_i=-\frac{t^{1/2}-t^{-1/2}}{t^{-1}-1}
\Si_\infty =t^{1/2}\Si_\infty.
$$

These formulas show that we can use
$\Si_\infty$ exactly in the same way as
$\hat{\mathscr{P}}\,'_+\,$ 
in Theorem \ref{PSBULLET}, i.e., it can be used to define
the corresponding form $\lan f, g\ran$ in $\De_\xi$. 
The uniqueness of such a bilinear form in $\De_\xi$ up to
proportionality results in the coefficient-wise
proportionality of $\Si_\infty$,
$\hat{\mathscr{P}}\,'_+\,$ and
$\hat{\mathscr{S}}\,'_+\circ \tilde{\mu}\,$.
Upon evaluation at $1$, we see that 
$\Si_\infty=\hat{\mathscr{P}}_+\,.$ 
\smallskip

{\em Point-wise convergence; (ii) and (iv).}
The considerations from $(i)$ can be equally used for 
the point-wise convergence to zero of the 
$C$\~coefficients of $t^{-(\rho,b)}Y_b^{-1}$ in the 
limit $b\to\infty$ and, more generally, the
coefficients of $\Si^u_{\infty}$ for id$\neq u\in W.$

If $\hw$ is fixed, then the minimal common power
of $q$ in the expansion of  $C_{\hw}$ 
will grow linearly together with $(\rho,b)$.  
Therefore the sum of absolute 
values of all coefficients of $C_{\hw}$ expanded in terms of 
the powers of $t'$ and $X_{\al_i}\,(i\ge 0)$
can grow  no greater than exponentially. Thus the functional 
convergence of the series for $C_{\hw}$ to zero can be achieved 
by making $\,|q|\,$ sufficiently small, depending on $t$ and the
compact set were $X$ is taken, naturally apart from the
singularities. This can be readily extended to
any $u\neq$id. 
\smallskip

\rmk
We note that a direct 
justification of the functional (point-wise) convergence 
to $0$ in (\ref{prodforpi}) for the whole $C_{\pi_r^{-1}}$ and for any 
$C_{\hw}$ is doable as well (without the $X,q,t^{-1}$\~expansions), 
though it follows essentially
the same lines.
Let us also mention that the statements from  $(ii)$ are discussed 
in detail in the case of $A_1$  in (\ref{YexpB});
see the first formula there and Theorem \ref{SIMINV}. \sq
\smallskip

Similar estimates show that the coefficients
of $\Si_{\infty}\,$ exist as analytic functions 
for sufficiently small $|q|,|t'|$. Thus
(\ref{tYliminvS}) holds for such $q,t$, where
the $C$\~coefficients are treated analytically
(in this range), which is the first part of $(iv)$.
\smallskip

{\em Sharp estimates.\,}
The exact estimates from $(iv)$ including the coincidence
statements formally follow 
from its first part, which is for sufficiently 
small $|q|$ and $|t'|$. The following actually repeats the
argument that have been already used for the sharp estimates 
of the convergence of the coefficients of 
$\hat{\mathscr{P}}\,'_+\,$.

Let $|q|<1$.  The coefficients 
of $(ct(t^{-1})/\hat{P}(t^{-1}))\,
\hat{\mathscr{S}}\,'_+\circ \tilde{\mu}$
are very explicit and the first singularity with respect to
$t$ is at a certain root of unity due to the zeros
of $\hat{P}(t^{-1}))$. Due to (\ref{tYliminvS}),
this gives that the radius of convergence is
$|t|^{-1}<1$ for all coefficients of $\Sigma_\infty\ $.
Note that this is different from the answer obtained for 
$\hat{\mathscr{P}}\,'_+\,$ due to the presence of the
Poincar\'e series in the coefficient of proportionality.

This concludes the justification of Theorem \ref{YLEFTNEW}.
\sq
\smallskip

\rmk
Let us discuss very briefly
the case $|t|<1$. We set
$\Sigma^+_{\mathbf b}\,=\, 
(1/|\Pi|)\sum_{j=1}^p\, (-t)^{(b^j,\rho)}Y_{b^j}\,$
instead of that in (\ref{SiYgennew}). Then 
$\Sigma^+_{\infty}$ defined as $
\lim_{{\mathbf b} \to \infty}\Sigma^+_{\mathbf b}\,$
will be proportional
to \,$\sum_{\hw\in \hW} (-t)^{l(\hw)/2}T_{\hw}^{-1}$\, and to
$\mu(X; q,t)^{-1}\circ \sum_{\hw\in \hW}
(-1)^{l(\hw)}\,\hw$ for
$\mu$ from (\ref{mutildemu}), the $|t|<1$ 
counterpart of $\hat{\mathscr{S}}\,'_+\circ \tilde{\mu}$.
This 
makes the standard {\em right} decomposition 
$\Sigma^+_{\infty}=\sum C_{\hw}^+\hw$ naturally much 
simpler than that for $\Sigma_{\infty}$. We note that the exact 
coefficients of proportionality can be obtained from 
the theory $|t|>1$ using the DAHA involution
(not an anti-involution) sending $t^{1/2}\mapsto -t^{-1/2}$ 
and fixing $q$ and the generators of $\HH$. It coincides with  
$\,H\mapsto \mu^{-1}\, H^\iota \,\mu\,$
for the involution $\iota\,$ used in Lemma \ref{IOTLEM}, when
acting on operators in a proper completion of the polynomial 
representation. It results in
$$
\Sigma^+_\infty\ \,=\ \, 
(ct(t)/\hat{P}(t))\,
\mu(X; q,t)^{-1}\circ \sum_{\hw\in \hW}
(-1)^{l(\hw)}\,\hw.
$$ 
Compare with (\ref{YexpA}) in the case of $A_1$.
\sq

\subsubsection{\sf Higher levels}
Conjugating $\Diamond$ from (\ref{diamondef})
by $q^{\,lx^2/2}$
for an integer $l\ge 0$, one obtains the following 
anti-involution: 
\begin{eqnarray*}
&\Diamond_l:\, T_{i}\, \mapsto\, T_{i},\  (i>0), 
Y_{b}\, \mapsto\, q^{-lx^2/2}\,Y_{b}\,q^{lx^2/2}, 
\  (b\in P^{\vee}), \\
&X_{a}\, \mapsto\,  
T_{w_{0}}^{-1}X_{a^{\varsigma}}T_{w_{0}},\ 
X_{a^{\varsigma}}=\varsigma(X_{a})=X_{-w_{0}(a)}, \ a\in P.
\end{eqnarray*}
The formulas for $T_0$ and $\pi_r$ can also be calculated but
they are not that direct. 
Let us discuss the invariant forms corresponding to
$\Diamond_l$ for $l>0$.
The $\HH$\~module will be
the polynomial representation $\mathscr{X}$.

We use that $\hat{\mathscr{I}}$ identifies 
the space of coinvariants 
$\mathscr{X}/\mathcal{J}_{l}(\mathscr{X})$, 
from Section \ref{sec:coinvariants} with the  Looijenga
space $\l_l\, (l\in \N)$ for generic $k$. Recall that
$\mathcal{J}_{l}(\mathscr{X})$ is the span of linear spaces
$$q^{-lx^{2}/2}\,(T_{\hw}-
t^{l(\hw)/2})(\mathscr{X}q^{\,lx^{2}/2}) \for
\hw\in \hW.
$$
Thus this is exactly the space of $\Diamond_l$\~coinvariants
from (\ref{starcoinv}):
\begin{align*}
&\HH/(\mathcal{J}+\mathcal{J}^{\Diamond_l})=
\mathscr{X}/\mathcal{J}^{\Diamond_l}(\mathscr{X}),\ 
\mathcal{J}=\,\hbox{Ker\,}(\HH\ni A\mapsto 
A(1)\in \mathscr{X});
\end{align*}
the subspaces $\mathcal{J}_l\subset \mathscr{X}$ 
from Section \ref{sec:pcoinv} and $\mathcal{J}^{\Diamond_l}$ 
coincide. 

The action of $\Diamond_l$ is trivial in this quotient; 
use the limit $t\to 1$ to see this. Therefore every functional 
on this space can be used to construct a form associated 
with $\Diamond_l$, and every such a form can be obtained in this
way. Using $\hat{\mathscr{I}}$, we come to the following 
extension of Theorem \ref{DIAZERO} from $l=0$ to $l>0$.  

\begin{theorem}\label{DIAGAU}
Let us assume that $\mathscr{X}$ has a nonzero symmetric
form $\lan f,g \ran$ corresponding to the anti-involution 
$\Diamond_l$ and normalized by $\lan 1,1\ran=1$. 
Provided that $\hat{\mathscr{I}}(\mathscr{X})=\l_l$,
this form can be represented as follows:
$$\lan f,g\ran=
\psi(\hat{\mathscr{I}}(f T_{w_0}(g^\varsigma)\,q^{l x^2/2}))
$$
for a proper linear functional $\psi: \l_l\to \C$.
When $l=1$, the resulting symmetric form satisfies the
Shapovalov property; the corresponding anti-involution
is of strong type with respect to $\mathscr{Y}$.
\sq
\end{theorem} 
\smallskip

{\bf Analytic theories.}
 \ Let us take a function $\phi(x)$ such that
$\phi\, q^{-l x^2/2}$ is $\hW$\~invariant, for
instance $\phi=q^{-l x^2/2}$. Then the form   
\begin{equation}\label{diaform}
\lan f,g\ran_\phi=
t^{-l(w_{0})/2}\int f T_{w_{0}}(g^{\varsigma})\phi\mu',\, 
\mu'=\mu (l<0),\, \mu'=\tilde{\mu} (l>0).
\end{equation}
is symmetric and is served by $\Diamond_{l}$ for the
following major choices of the integration (``theories"):

(a) imaginary integration $\int_{e+\imath \R^n}$ for
$e\in \R^n$ subject to $l< 0$,

(b) real integration $\sum_{w\in W}\int_{w(e)+\R^n}$ for
$\,e\not\in \R^n\,$, where $l>0$,

(c) Jackson integration $\int_\xi f\,=\,
\sum_{\hw\in \hW}\,f(q^{\hw(\xi)})$,  where $l>0$.

The function $\phi$ must be analytic in a neighborhood of
the integration contour for $(a)$ and everywhere 
for $(b)$ to ensure that the integral does not depend
on the choice of $e\in \R^n$. Finding the kernels
of the linear map $\phi\mapsto
\lan\, \cdot\,,\,\cdot\,\ran_\phi$ is an interesting problem;
the dimension of its image equals dim\,$\l_l$.
 
Establishing connections between these theories 
is fundamental in harmonic analysis. Relating
them to {\em algebraic} Shapovalov-type 
inner products is equally important. The latter
inner products do not involve integrations and are well defined 
for all or almost all $q,t$. This problem is directly linked to 
the DAHA-generalization of the  
Arthur-Heckman-Opdam approach from \cite{HO2}, which can be 
stated as the problem of
{\em finding presentations of algebraically defined inner
products in DAHA\~modules} (Shapovalov-type ones) 
{\em in terms of the integrations} 
(with respect to the affine residual subtori). 

\setcounter{equation}{0}
\section{\sc{The rank-one case}}
\subsection{{\bf Polynomial representation}}
\subsubsection{\sf Basic definitions}
Let us consider the root system $A_{1}$.
Following Section \ref{sec:rank1},
$\HH$ is generated by $Y=Y_{\om_1},T=T_1,X=X_{\om_1}$  
subject to the quadratic relation $(T-t^{1/2})(T+t^{-1/2})=0$
and the cross-relations:
\begin{align}\label{dahaone}
&TXT=X^{-1},\ T^{-1}YT^{-1}=Y^{-1},\ Y^{-1}X^{-1}YXT^2q^{1/2}=1.
\end{align}
Using $\pi\equal YT^{-1}$, the second relation becomes $\pi^2=1$.
The field of definition will be $\C(q^{1/4},t^{1/2})$, though
$\Z[q^{\pm 1/2},t^{\pm 1/2}]$ is sufficient for many constructions;
actually $q^{\pm 1/4}$ will be needed only in the automorphisms 
$\tau_{\pm}$ below. We will frequently treat $q,t$ as numbers;
then the field of definition will be $\C$.

The following map can be extended to an anti-involution on $\HH$\,
$\varphi: X \leftrightarrow Y^{-1}, T\to T$. The first two
relations in (\ref{dahaone}) are obviously fixed by $\varphi$;
as for the third, check that 
$\varphi(Y^{-1}X^{-1}YX)=Y^{-1}X^{-1}YX$. 

The following DAHA automorphism is of key importance in this paper:
$$
\tau_+(X)=X,\ \tau_+(T)=T,\ \tau_+(Y)=q^{-1/4}XY,\ 
\tau_+(\pi)=q^{-1/4}X\pi,
$$
which can be interpreted as conjugation by the Gaussian $q^{x^2}$ 
for $X=q^x$.
Check that $T^{-1}YT^{-1}=Y^{-1}$ is transformed to
$Y^{-1}X^{-1}YXT^2q^{1/2}=1$ under $\tau_+$. 
Applying $\varphi$ we obtain an automorphism
$\tau_-=\varphi\tau_+\varphi$\,:
$$
\tau_-(Y)=Y,\ \tau_-(T)=T,\ \tau_-(X)=q^{1/4}YX.
$$

The Fourier transform corresponds to the
following automorphism of $\HH$ (it is not an
involution)\,:
\begin{align}
\si(X)= Y^{-1},\ \si(T)=T,\ &\si(Y)=q^{-1/2}Y^{-1}XY=XT^2,\ 
\si(\pi)=XT,
\notag\\ 
\si\ &=\  \tau_+\tau_-^{-1}\tau_+\ =\ \tau_-^{-1}\tau_+\tau_-^{-1}.
\label{tautautau}
\end{align}
Check that $\,\si\tau_+=\tau_-^{-1}\si,\ \, 
\si\tau_+^{-1}=\tau_-\si.$ 

The polynomial representation is defined as 
$\mathscr{X}=\C_{q,t}[X^{\pm 1}]$ over the field 
$\C_{q,t}=\C(q^{1/4},t^{1/2})$ with $X$ acting by
the multiplication. The formulas for the 
other generators are 
\begin{eqnarray*}
T=t^{1/2}s+\frac{t^{1/2}-t^{-1/2}}{X^{2}-1}
\circ (s-1),\ \, Y=\pi T
\end{eqnarray*}
in terms of the multiplicative reflection $s(X^n)=X^{-n}$
and $\pi(X^n)=q^{n/2}X^{-n}$ for $n\in\Z$.

The Gaussian $q^{x^2}$ is an element of
a completion of $\mathscr{X}$. 
However the conjugation $A\mapsto
q^{x^2}\,A\,q^{-x^2}$ for $A\in \HH$
preserves $\HH$ and
coincides with $\tau_+$. To see this use that
$$
Y=\om\circ(t^{1/2}+\frac{t^{1/2}-t^{-1/2}}{X^{-2}-1}
\circ (1-s)).
$$
Recall that $X=q^x$ and
\begin{align*}
&s(x)=-x,\ \om(f(x))=f(x-1/2),\ 
\pi=\om s,\ \pi(x)=1/2-x,\\ 
&\om(q^{x^2})\,=\,q^{1/4}X^{-1}q^{x^2},\ \, 
Y(q^{-x^2})\,=\,\om(q^{-x^2})\,=\,q^{-1/4}X q^{-x^2}.
\end{align*}

It is important that $\HH$ at $t=1$ becomes the 
Weyl algebra defined as the span
$\lan X,Y\ran /(Y^{-1}X^{-1}YXq^{1/2}=1)$ extended
by the inversion $s=T(t=1)$ sending 
$X\mapsto X^{-1}$ and $Y\mapsto Y^{-1}$. 

\subsubsection{\sf The E-polynomials}{\label{sect:macpoly}}
Let us assume that
$k$ is generic; we set $t=q^k$.
The definition is as follows: 
\begin{align}\label{nonsymp}
YE_{n}=q^{-n_{\#}}E_{n}\for n\in Z,\ \,E_n\in \mathscr{X}, &&&&\\
n_{\#}=\left\{\begin{array}{ccc}\frac{n+k}{2} 
&  & n>0, \\\frac{n-k}{2} &  & n\le 0,\end{array}\right\}, 
\text{\, note that }\,  0_{\#}=-\frac{k}{2}.\label{nonsymp1}
\end{align}
The normalization is
$E_{n}=X^{n}+\text{ ``lower terms'' },$
where by ``lower terms'', we mean
polynomials in terms of $X^{\pm m}$ as $|m|<n$
and, additionally,  $X^{|n|}$ for negative $n$.
It gives a filtration in $\mathscr{X}$ with the
consecutive quotients of dimension $1$.
Check that $Y$ preserves it, which justifies that $Y$
is diagonalizable in $\mathscr{X}$ and readily
provides the formulas for the eigenvalues 
from (\ref{nonsymp}),(\ref{nonsymp1}).
 
The $E_{n} (n\in \Z)$ are called
{\em nonsymmetric Macdonald polynomials} or simply 
$E$\~polynomials.
Obviously,  $E_{0}=1, E_{1}=X$. 

\subsubsection{\sf The intertwiners}\label{sect:intertw}
The first intertwiner comes from AHA theory:
$$
\Phi\equal T+\frac{t^{1/2}-t^{-1/2}}{Y^{-2}-1}\,:\ 
\Phi Y=Y^{-1}\Phi.
$$
The second is $\Pi\equal q^{1/4}\tau_+(\pi)$;
obviously, $\Pi^2=q^{1/2}$.
Explicitly,
$$
\Pi=X\pi=q^{1/2}\pi X^{-1}\,:\ 
\Pi Y=q^{-1/2}Y^{-1}\Pi.
$$
Use that $\phi(\Pi)=\Pi$ to deduce the latter
relation from $\Pi X\Pi^{-1}=q^{1/2}X^{-1}$.
The $\Pi$\~type intertwiner is due to Knop and Sahi for $A_n$ 
(the case of arbitrary reduced systems was considered in \cite{C1}).
Since $\Phi,\Pi$ ``intertwine" $\mathscr{Y}$, they can be
used for generating the $E$\~polynomials. Namely,
\begin{align}{\label{signE}}
&E_{n+1}=q^{n/2}\Pi (E_{-n}) \for n\ge 0,\\
&E_{-n}=t^{1/2}(T+\frac{t^{1/2}-t^{-1/2}}{q^{2n_{\#}}-1})E_{n}.
\label{interphi}
\end{align}
Beginning with $E_0=1$, one can readily construct the
whole family of $E$\~polynomials.
For instance,
\begin{eqnarray*}
T(X)&=&t^{1/2}X^{-1}+\frac{(t^{1/2}-t^{-1/2})(X^{-1}-X)}{X^{2}-1}\\
&=&t^{1/2}X^{-1}-(t^{1/2}-t^{-1/2})X^{-1}\ =\ t^{-1/2}X^{-1}, \\
E_{-1}&=&t^{1/2}(T+\frac{t^{1/2}-t^{-1/2}}{qt-1})E_{1}
=X^{-1}+\frac{1-t}{1-t q}X.
\end{eqnarray*}
Using $\Pi$, 
\begin{align*}
&E_{2} = q^{1/2}\Pi E_{-1}=  X^{2}+q\frac{1-t}{1-tq}.
\end{align*}
Applying $\Phi$ and then $\Pi$,
\begin{align*}
&E_{-2} = X^{-2}+\frac{1-t}{1-tq^{2}}X^{2}
+\frac{(1-t)(1-q^{2})}{(1-tq^{2})(1-q)},\\
&E_{3}=X^{3}+q^2\frac{1-t}{1-tq^{2}}X^{-1}+q\frac{(1-t)(1-q^{2})}
{(1-tq)(1-q)}X.
\end{align*}
It is not difficult to find the general formula. See e.g.,
(6.2.7) from \cite{Ma4} for integral $k$. However, recalculating
these formulas from integral $k$ to generic $k$ is not
too simple; we will provide the exact formulas for 
the $E$\~polynomials
below (in the form we need).

\subsubsection{\sf The E-Pieri rules}
For any $n\in \Z$, we have the {\em evaluation formula}
\begin{eqnarray*}
E_{n}(t^{-1/2})=t^{-|n|/2}\prod_{0<j<|\tilde{n}|}
\frac{1-q^{j}t^{2}}{1-q^{j}t},
\end{eqnarray*}
where $|\tilde{n}|=|n|+1$ if $n\leq 0$ and 
$|\tilde{n}|=|n|$ if $n>0$.

It is used to introduce the {\em nonsymmetric spherical polynomials}
$$\mathcal{E}_{n}=\frac{E_{n}}{E_{n}(t^{-1/2})}.$$
This normalization is important in many constructions
due to the 
{\em duality formula} $\e_m(q^{n_\#})=\e_n(q^{m_\#})$.
The Pieri rules are the simplest for the
$E$\~spherical polynomials:
\begin{eqnarray}\label{pierie}
X\mathcal{E}_{n}
=\frac{t^{-1/2\pm 1}q^{-n}-t^{1/2}}
{t^{\pm1}q^{-n}-1}\mathcal{E}_{n+1}
+\frac{t^{1/2}-t^{-1/2}}{t^{\pm1}q^{-n}-1}\mathcal{E}_{1-n}.
\end{eqnarray}
Here the sign is $\pm=+$\, if\, $n\leq 0$\  and $\pm=-$ if $n>0$.
These formulas give an alternative approach to
constructing the $E$\~polynomials and establishing
their connections with other theories, for instance, with
$\mathfrak{p}$\~adic Matsumoto functions.

\subsubsection{\sf Rogers' polynomials}{\label{sect:Rogers}}
Let us introduce the {\em Rogers polynomials} for $n\ge 0$:  
$$
P_{n}=(1+t^{1/2}T)\bigl(E_{n}\bigr)=(1+s)
\bigl(\frac{t-X^2}{1-X^2}E_n\bigr)=E_{-n}+\frac{t-tq^n}{1-tq^n}E_n,
$$ 
$P_n=X^n+X^{-n}+$``lower terms", where the latter are $X^m+X^{-m}$
for $0\le m<n$.
They are eigenfunctions
of the following well-known operator
\begin{align}\label{Lopertaor}
\l=\frac{t^{1/2}X-t^{-1/2}X^{-1}}{X-X^{-1}}\Ga+
\frac{t^{1/2}X^{-1}-t^{-1/2}X}{X^{-1}-X^{1}}\Ga^{-1},
\end{align}
where we set $\Ga(f(x))=f(x+1/2),\, \Ga(X)=q^{1/2}X$,
i.e., $\Ga$ acts as $-\om$ in $\mathscr{X}$. This operator
is the restriction of the operator $Y+Y^{-1}$ to symmetric
polynomials, which is the key point of the DAHA approach
to the theory of the Macdonald polynomials.

The exact eigenvalues are as follows:
\begin{equation}{\label{eqn:LP}}
\l(P_n)\ = \ (q^{n/2}t^{1/2}+q^{-n/2}t^{-1/2})\,P_n,\ n\ge 0.
\end{equation}
The evaluation formula reads
$$
P_n(t^{\pm1/2})=t^{-n/2}\prod_{0\le j\le n-1}
\frac{1-q^{j}t^{2}}{1-q^{j}t}.
$$
The spherical $P$\~polynomials $\p_n\equal P_n/P_n(t^{1/2})$
satisfy the duality $\p_n(t^{1/2}q^{m/2})=\p_m(t^{1/2}q^{n/2})$.

\subsubsection{\sf Explicit formulas}{\label{sect:formulas}}
Let us begin with the well-known formulas for the 
Rogers polynomials ($n\ge 0$): 
\begin{align}\label{exactp}
&P_{n}=X^{n}+X^{-n}+\sum_{j=1}^{[n/2]}M_{n-2j}\prod_{i=0}^{j-1}
\frac{(1-q^{n-i})}{(1-q^{1+i})}\,\frac{(1-tq^{i})\  }
{(1-tq^{n-i-1})},
\end{align}
where $M_{n}=X^{n}+X^{-n}\, (n>0)$ and $M_{0}=1$. 

The formulas for the $E$\~polynomials are as follows ($n>0$):
\begin{align}\label{exacte-}
E_{-n}&=X^{-n}+X^n\frac{1-t}{1-tq^n}+
\sum_{j=1}^{[n/2]}X^{2j-n}\,\prod_{i=0}^{j-1}
\frac{(1-q^{n-i})}{(1-q^{1+i})}\,\frac{(1-tq^{i})}{(1-tq^{n-i})}
\notag\\
&+\sum_{j=1}^{[(n-1)/2]}X^{n-2j}\,
\frac{(1-tq^{j})}{(1-tq^{n-j})}\prod_{i=0}^{j-1}
\frac{(1-q^{n-i})}{(1-q^{1+i})}\,\frac{(1-tq^{i})}{(1-tq^{n-i})},
\end{align}
\begin{align}\label{exacte+}
E_{n}=X^{n}+&
\sum_{j=1}^{[n/2]}X^{2j-n}\,q^{n-j}\,
\frac{(1-q^{j})}{(1-q^{n-j})}\prod_{i=0}^{j-1}
\frac{(1-q^{n-i-1})}{(1-q^{1+i})}\,
\frac{(1-tq^{i}\ ) }{(1-tq^{n-i-1})}
\notag\\
&+\sum_{j=1}^{[(n-1)/2]}X^{n-2j}\,q^j\,
\prod_{i=0}^{j-1}\frac{(1-q^{n-i-1})}{(1-q^{1+i})}\,
\frac{(1-tq^{i}\ ) }{(1-tq^{n-i-1})}.
\end{align}

\subsection{{\bf The p-adic limit}}
Let us ``separate" $t$ and $q$; they will not be connected
any longer by the relation $t=q^k$ in this section.

\subsubsection{\sf The limits of P-polynomials}
We will begin with the symmetric case.
Formula (\ref{exactp}) readily gives that 
\begin{align*}
P^{0}_{n}\equal\lim_{q\to 0}P_{n}&=X^{n}+X^{-n}+
\sum_{j=1}^{[n/2]}M_{n-2j}
\prod_{i=0}^{j-1}(X^{n-2j}+X^{2j-n})(1-t)\\
&=X^{n}+X^{-n}+(1-t)\chi_{n-2}
=\chi_{n}-t\chi_{n-2}
\end{align*}
for the monomial symmetric functions $M_n$ and the
classical characters $\chi_n=(X^{n+1}-X^{-n-1})/
(X-X^{-1})$. In the spherical normalization, 
$\p_n=P_n/P_n(t^{1/2})$, where $P^0_n(t^{1/2})=t^{-n/2}(1+t)$.
One has
$$
\mathcal{P}_{n}^{0}=(\chi_{n}-t\chi_{n-2})\frac{t^{n/2}}{1+t}.
$$
By letting $t\to t^{-1}$ and $X\to Y$, we obtain that
$\mathcal{P}^{0}_{n}$ coincides with the spherical
function $\varphi_{n}$.

Let us obtain this fact directly from the 
definition of  the Rogers polynomials $P_{n}$ in terms of
the operator $\l$:
\begin{align*}
\left(\frac{t^{1/2}X-t^{-1/2}X^{-1}}{X-X^{-1}}\Ga+
\frac{t^{1/2}X^{-1}-t^{-1/2}X}{X^{-1}-X^{1}}\Ga^{-1}\right)
&P_{n}\\
=\ (q^{n/2}t^{1/2}+q^{-n/2}t^{-1/2})&P_{n}\,;
\end{align*}
see Section \ref{sect:Rogers}. 

Recall that $\Ga(X^{m})=q^{m/2}X^{m}\,$ for any $\,m\in \Z\,$.
It gives that 
\begin{align}\label{limgax}
&\lim_{q\to 0}q^{n/2}\Ga^{\pm 1}(X^{\pm m})=0
\for |m|\le n \hbox{\ \ unless\ }\notag\\ 
&\lim_{q\to 0}q^{n/2}\Ga(X^{-n})=X^{-n},\  
\lim_{q\to 0}q^{n/2}\Ga^{-1}(X^{n})=X^{n} \for n\ge 0.
\end{align}
Therefore
\begin{align*}
t^{-1/2}P_{n}^{0}=\frac{t^{1/2}X-t^{-1/2}X^{-1}}
{X-X^{-1}}X^{-n}
+\frac{t^{1/2}X^{-1}-t^{-1/2}X}{X^{-1}-X}X^{n}.
\end{align*}
Using that $P^0_n(t^{1/2})=t^{-n/2}(1+t)$, we obtain that
\begin{align*}
\mathcal{P}_{n}^{0}=\left(\frac{tX^{2}-1}{X^{2}-1}X^{-n}
+\frac{tX^{-2}-1}{X^{-2}-1}X^{n}\right)\frac{t^{n/2}}{1+t},
\end{align*}
which is exactly the Macdonald summation formula
(\ref{vphiformula})
under the substitution $X\mapsto Y,t\mapsto t^{-1}$:
\begin{align}\label{Macdsummat}
\varphi_{n}\ =\  
\frac{t^{-n/2}}{1+t^{-1}}
\left(\frac{1-t^{-1}Y^{-2}}{1-Y^{-2}}Y^{-n}
+\frac{1-t^{-1}Y^{2}}{1-Y^{2}}Y^{n}\right).
\end{align}
We see that the right-hand side of (\ref{Macdsummat}) is
actually the limit of the operator $\l$; this is a general
fact (true for any root systems). 
\smallskip

\subsubsection{\sf The limits of E-polynomials}\label{sec: p-adic-e}

We mainly follow \cite{C101}, however, with certain technical
modifications.
\begin{theorem}
The limit $\e^0_n(X)=
\lim_{q\to 0}\mathcal{E}_{n}\equal \mathcal{E}_{n}^{0}$ 
exists. The Matsumoto functions 
$\varepsilon_{n}$ from (\ref{vepformulas}) are connected
with $\e^0_n$ as follows:
$$
\varepsilon_{n}=
\mathcal{E}_{n}^{0}(t\to t^{-1}, X\to Y).
$$
\end{theorem}
{\em Proof.} 
First, $\lim_{q\to 0}E_{n}(t^{-1/2})=t^{-|n|/2}$.
For $n>0$, we have
\begin{eqnarray*}
X\mathcal{E}_{n}^{0} & = & t^{-1/2}\mathcal{E}_{n+1}^{0},\\
X\mathcal{E}_{-n}^{0} & = & t^{1/2}\mathcal{E}_{-n+1}^{0}
-(t^{1/2}-t^{-1/2})\mathcal{E}_{n+1}^{0}.
\end{eqnarray*}
These are exactly the Pieri relations for the Matsumoto
functions from (\ref{yepm}--\ref{y-epm}) upon the
substitution $Y\mapsto X,t\mapsto t^{-1}$.\sq

We know from (\ref{vepformulas}) that for $n\ge 0$,
\begin{align}\label{vepnfor}
&\vep_n =  t^{-\frac{n}{2}}Y^n,\  
\vep_{-n}
=t^{-\frac{n+1}{2}}(t^{\frac{1}{2}}Y^{-n}+
(t^{\frac{1}{2}}-t^{-\frac{1}{2}})\frac{Y^{-n}-Y^{n}}
{Y^{-2}-1}).
\end{align}
Obtaining these formulas directly from (\ref{exacte-}) and
(\ref{exacte+}) is of some interest. 
\smallskip

Let us show how to use here the $Y$\~operator, namely,
the formulas for the action of $Y$ and $Y^{-1}$, correspondingly,
on $E_n (n>0)$ and $E_{-n}(n\ge 0)$. One can present (\ref{nonsymp})
as follows: 
\begin{align*}
&\frac{1-qtX^{-2}}{1-qX^{-2}}\,q^{n/2}\Ga^{-1}(E_{n})+
\frac{1-t}{1-q X^{-2}}\,s(q^{n/2}\Ga(E_{n}))\,=\, E_n \ (n>0),\\
&(\frac{1-t^{-1}X^{-2}}{1-X^{-2}}-
\frac{1-t^{-1}}{1-X^{-2}}\,s)\,(q^{-n/2}\Ga(E_{-n}))\,=\, 
E_{-n} \for n\ge 0.
\end{align*}
Setting $E_m^0=E_m(X; q=0)$ and applying (\ref{limgax}),
\begin{align*}
&E_{n}^0=X^n  \and  E_{-n}^0=
(\frac{1-t^{-1}X^{-2}}{1-X^{-2}}-
\frac{1-t^{-1}}{1-X^{-2}}\,s)(X^{-n}) \for n\ge 0,
\end{align*}
where in the first formula we use that $q^{n/2}\Ga(E_{n})=0$
in the limit $n\to\infty\,$
because $E_{n>0}$ does not contain $X^{-n}$.
Switching from $E_m^0$ to $\e_m^0$ and then to $\vep_m$ 
(any $m\in \Z$), we arrive at (\ref{vepnfor}).
\medskip

\rmk
We expect that a
similar connection holds between the {\em difference-elliptic\,}
symmetric Macdonald-type Looijenga functions (which can be
called elliptic $P$\~functions) and the affine
Hall functions. There is no 
general theory of such Macdonald-Looijenga functions 
so far; the paper \cite{Ch13} dealt with
the difference-elliptic theory only at the level of
operators. Their diagonalization was not performed there.
The elliptic Ruijsenaars operators and their 
generalizations to arbitrary root systems are really  
connected with the affine symmetrizers
in the corresponding limit. 

Paper \cite{Ch13} indicates that the corresponding
{\em nonsymmetric Hall polynomials\,} can be obtained
from the nonsymmetric elliptic Macdonald-Looijenga 
functions, {\em elliptic $E$\~functions}, for the direct 
counterpart of the limit $q\to 0$. The nonsymmetric elliptic 
$E$\~functions are actually simpler to define than
the $P$\~functions. The limit $q\to 0$ 
is well defined for the properly 
normalized difference-elliptic $Y$ operators from
\cite{Ch13}, which provides certain nonsymmetric variant 
of the DAHA symmetrizer considered in this paper;
it is in progress. A connection is expected with
\cite{EFMV}.

Note that there are other theories of 
elliptic orthogonal polynomials. The most advanced theory
we know is \cite{Ra}; however, this seems not 
what is needed here. 

\subsection{{\bf Coinvariants and symmetrizers}}
\subsubsection{\sf DAHA coinvariants}
Let us prove Theorem \ref{thm:coinv} for the level $l=1$
in the case of $A_1$.
\begin{theorem}
For any $q,t=q^k$,  
$\dim(\mathscr{X}/\mathcal{J}_{1}(\mathscr{X}))=1$.
\end{theorem}
{\em Proof.}
Let $\varrho: \HH\to \mathbb{C}$ be a functional on $\HH$ such that 
\begin{align}\label{eqn:varrho1}
\varrho(\HH\cdot(T_{\hat{w}}-t^{l(\hat{w})/2}))\,=\,0
\text{\ \,and\ }\\
\label{eqn:varrho2}
\varrho(\tau_+^{-1}(T_{\hat{w}}-t^{l(\hat{w})/2})\cdot \HH)\,=\,0
\end{align}
for all $\hw\in \hW=W\lsmash P^\vee=\S_2\lsmash \Z\om$.
\begin{lemma} \label{lem:basis}
An arbitrary $A\in \HH$ can be uniquely represented as
\begin{align*}
A=\sum c_{n,\varepsilon,m}\tau_+^{-1}(Y^{n})T^{\varepsilon}Y^{m},
\end{align*}
where $\varepsilon=0$ or $1$,\  $m, n$ are integers and 
$c_{n,\varepsilon,m}$ are constants. 
\end{lemma}
{\em Proof of Lemma \ref{lem:basis}.}
One has $\tau_-\tau_+^{-1}(Y)=$
$\tau_-\tau_+^{-1}\tau_-(Y)=
\si^{-1}(Y)=X^{-1}$. Applying $\tau_-^{-1}$ to
$X^{-n}$, $T^{\varepsilon}$ and $Y^m$,
we obtain that the elements 
$\{\tau_+^{-1}(Y^{n})T^{\varepsilon}Y^{m}\}$ form a 
PBW basis for $\HH$.
\sq

Now, for $A\in \HH$,
relations \eqref{eqn:varrho1} and \eqref{eqn:varrho2} give
that
\begin{align*}
\varrho(\tau_+^{-1}(Y^{n})A)\,=\,
t^{n/2}\varrho(A) \text{\ and \ } 
\varrho(AT_{\hw})\,=\,t^{l(\hw)/2}\varrho(A). 
\end{align*}
Representing $A$ as in  Lemma \ref{lem:basis}, 
\begin{align*}
\varrho(A)=\sum c_{n,\varepsilon,m}
\varrho(\tau_+^{-1}(Y^{n})T^{\varepsilon}Y^{-m})\,=
\sum c_{n,\varepsilon,m}\,t^{n/2+\vep/2-m/2}.
\end{align*}
Thus
$\dim(\mathscr{X}/\mathcal{J}_{1}(\mathscr{X}))=1$.
\sq

\rmk
A similar argument can be employed for arbitrary simply-laced
root systems (or if the twisted setting is used). 
A counterpart of Lemma \ref{lem:basis} is the claim that
an arbitrary $A\in \HH$ can be uniquely represented as
\begin{align*}
A=\sum c_{b,w,a}\tau_+^{-1}(Y_{b})T_{w}Y_{a},
\end{align*}
where $w\in W$,\, $a,b\in P$\, and $c_{b,w,a}$
are constants. 

For any level $l>0$, $\tau_+^{-l}(Y)=q^{-l/4}X^{-l}Y$. 
Calculating the space of coinvariants generally requires 
knowing $\tau_+^{-l}(Y^m)$. The latter can be computed
using the relation 
$Y^{-1}X^{-1}YXT^2q^{1/2}=1$, but explicit formulas are
involved. Nevertheless, they are sufficient for finding the
dimension of the space of coinvariants (for arbitrary
simply-laced root systems as well).

\subsubsection{\sf P-hat in rank one}\label{sect:symmrank1}
Let us discuss the rank-one version of Theorem
\ref{P+FORMULA}. The explicit list of
the elements $\hw\in \hW$ (there
are four types) and the corresponding 
$T_{\hw}'\equal t^{-l(w)/2}T_{\hw}^{-1}$,
presented in terms of $Y,T$,
is as follows:
\begin{align*}
&1)\ \hw= &m\om\cdot s\, &(m>0),\ \ l(\hw)=&m-1&,\ 
\ T_{\hw}'= &t^{-\frac{m-1}{2}}TY^{-m},\\
&2)      &m\om\, &(m>0),\, &m&,\, &t^{-\frac{m}{2}}Y^{-m},\\
&3)     &-m\om\, &(m\ge 0),\, &m&,\, 
&t^{-\frac{m}{2}}TY^{-m}T^{-1},\\
&4)     &(-m\om)\cdot s\, &(m\ge 0),\, &m+1&,\, 
&t^{-\frac{m+1}{2}}Y^{-m}T^{-1}.
\end{align*}
Note that we use a presentation that is somewhat different from
the one used in the justification of Theorem \ref{P+FORMULA}.

\begin{theorem}\label{SYMRANK1}
The affine symmetrizer $\hat{\mathscr{P}}\,'_+$ (the prime
here indicates that there is no division by 
$\hat{P}(t^{-1})$) can be expressed as follows:
\begin{align}\label{hatprank1}
&\hat{\mathscr{P}}\,'_+ = 
(1+t^{\frac{1}{2}}T)\,\Bigl(\,
\frac{\,t^{-\frac{1}{2}}\,Y^{-1}}
{1-t^{-\frac{1}{2}}\,Y^{-1}}\, 
(1+t^{-\frac{1}{2}}T^{-1})+t^{-\frac{1}{2}}T^{-1}\,\Bigr). 
\end{align}
In particular, 
$\hat{\mathscr{P}}\,'_+(1)=2\,\frac{1+t^{-1}}{1-t^{-1}}=
\hat{P}(t^{-1})=2+\sum_{m=1}^{\infty}4t^{-m}$ for $|t|>1$.
\sq\end{theorem}

The formula for $\hat{\mathscr{P}}\,'_+$ from the
theorem in terms of
$t^{-1/2}$ is exactly the definition of the $P$\~hat symmetrizer
upon using (1,2,3,4) above, as well as its particular case, the sum 
$2+\sum_{m=1}^{\infty}4t^{-m}$.
Note that
$
(1+t^{1/2}T)t^{-1/2}T^{-1}=1+t^{-1/2}T^{-1}.
$

As remarked in Section \ref{sect:RatPhat}
concerning formula (\ref{hatPrat}),
the right-hand side of (\ref{hatprank1})
becomes {\em identically} zero when treated
as an element of a proper localization of the affine
Hecke algebra  
$\h_Y=\lan T,Y^{\pm 1}\ran$. Indeed, this expression can be only
zero because the localization is not sufficient to
construct such an affine symmetrizer in $\h_Y$ 
(a completion is needed). 
One can deduce that (\ref{hatprank1}) vanishes directly from
the relation
\begin{align}\label{aharat}
Tf(Y)-f(Y^{-1})T=\frac{t^{1/2}-t^{-1/2}}{Y^{-2}-1}
(f(Y^{-1})-f(Y))
\end{align}
extended to rational functions $f(Y)$; let us demonstrate it.
First, the extension of (\ref{aharat}) to rational functions
is straightforward since an arbitrary
rational function in terms of $Y$ can be represented as
a Laurent polynomial divided by a {\em $W$\~invariant} Laurent
polynomial, commuting with $T$. 

Let  $1^+\equal (1+t^{-1/2}T^{-1})$,
\begin{align}\label{UUplus}
U\equal \frac{t^{-1/2}\,Y^{-1}}{
1-t^{-1/2}\,Y^{-1}},\ U^+=U(1+t^{-1/2}T^{-1}).
\end{align}
Then
\begin{align*}
& TU^+=-\frac{t^{-1/2}}{1-t^{-1/2}Y^{-1}}(1+t^{-1/2}T^{-1})=
-t^{-1/2}U^+ - t^{-1/2}1^+\,;
\end{align*}
therefore, $(1+t^{1/2}T)U^+ +1^+=0$, which is exactly
vanishing the right-hand side of (\ref{hatprank1}).

This identity is the key point of the formula for
$\hat{\mathscr{P}}\,'_+$ as a limit of the powers of $Y$.
Let us discuss in detail the corresponding deduction
of Theorem \ref{YLEFT} from Theorem \ref{P+FORMULA}
in the case of $A_1$.
   
For integers $M>0$, we introduce the 
{\em truncated symmetrizers}
\begin{align}\label{hatprank1M}
&\hat{\mathscr{P}}\,'_M = 
(1+t^{\frac{1}{2}}T)\,\Bigl(\,
\,\sum_{j=1}^M\,t^{-\frac{j}{2}}\,Y^{-j}
(1+t^{-\frac{1}{2}}T^{-1})\Bigr)+ 
1+t^{-\frac{1}{2}}T^{-1}.
\end{align} 
\begin{theorem}\label{SYMRANK1Y}
(i) Moving $T$ in $\hat{\mathscr{P}}\,'$ via 
(\ref{aharat}), one arrives at identities in $\h_Y$:  
\begin{align}\label{hatsirank1}
&\hat{\mathscr{P}}\,'_M=
\hat{\Sigma}^+_M\equal\,\hat{\Sigma}_M\,(1+t^{-1/2}T^{-1}),
\for\\
&\hat{\Sigma}_M \equal t^{-[\frac{M}{2}]}+
\sum_{j=1}^M\, t^{-[\frac{M-j}{2}]-\frac{j}{2}}\,(Y^j+Y^{-j}),
\notag
\end{align}
where $[a/b]$ is the integer part.

(ii) The operator
$\hat{\mathscr{P}}\,'_+$ is well defined if and only if
the limit $\hat{\Si}_\infty^+=\lim_{M\to\infty}
\hat{\Sigma}^+_M$ exists; then these operators coincide.
The existence of  $\hat{\Si}_\infty^+$ formally results in
the following condition:
\begin{align}\label{YMzero}
&\lim_{M\to\infty} t^{-M/2}(Y^{-M})^+=0.
\end{align}
In its turn, (\ref{YMzero}) ensures that 
$\hat{\Si}_\infty^+$ is an {\sf affine 
symmetrizer} if it exists, i.e., 
satisfies the symmetries
\begin{align}\label{YSiabs}
&Y\,\hat{\Si}^+_\infty=
\hat{\Si}^+_\infty\, Y=
t^{\frac{1}{2}}\,\hat{\Si}^+_\infty
=T\,\hat{\Si}^+_\infty =
\hat{\Si}^+_\infty\, T.
\end{align}

(iii) Finally, we claim that the existence of
$\hat{\mathscr{P}}\,'_+$, for instance its
coefficient-wise convergence in the 
$\hw$\~decomposition, results in the identity 
\begin{align}
&\hat{\mathscr{P}}\,'_+=
\lim_{M\to\infty}\overline{\Sigma}^+_M
\,\for\, \overline{\Si}_M\equal\frac
{t^{-\frac{M}{2}}\,Y^M+t^{-\frac{M-1}{2}}\,Y^{M-1}}{1-t^{-1}}, 
\label{YSiabsbar}
\end{align}
which includes the existence (convergence) of 
$\overline{\Sigma}^+_{\infty\,}\equal
\lim_{M\to\infty}\overline{\Sigma}^+_M$
in the same sense as that for $\hat{\mathscr{P}}\,'_+$.
\end{theorem}

\subsubsection{\sf Proof of Theorem \ref{SYMRANK1Y}}
Let us prove this theorem (and Theorem \ref{YLEFT} for $A_1$);
(\ref{YSiabsbar}) is its main part, called {\em the sigma-formula}.
The following was outlined in
Theorem \ref{YLEFT} for arbitrary root systems.

Only the $t$\~powers 
$t^{-M/2}$ and $t^{(1-M)/2}$ appear in 
the formula for $\hat{\Sigma}_M$\,:
\begin{align*}
\hat{\Sigma}_M\ =\ 
&t^{-\frac{M}{2}}(Y^M+Y^{-M})+
t^{\frac{1-M}{2}}(Y^{M-1}+Y^{1-M})\\
+\ &t^{-\frac{M}{2}}(Y^{M-2}+Y^{2-M})+\,\ldots\, 
+t^{-[\frac{M}{2}]}.
\end{align*}
For instance, in the case of even $M$,
$$
\hat{\Sigma}_M(1)=\sum_{j=2l} t^{-M/2}(t^{j/2}+t^{-j/2})+
\sum_{j=2l-1}t^{-M/2+1/2}(t^{j/2}+t^{-j/2})
$$
for $l=1,2\ldots,M/2$. The resulting
$t^{-1}$\~series  is 
$2+2t^{-1}+2t^{-2}+\ldots$\,; we obtain that 
$$
\lim_{M\to\infty}\hat{\Sigma}_M\cdot(1+t^{-1/2}T^{-1})(1)=
2\,\frac{1+t^{-1}}{1-t^{-1}}=\hat{P}(t^{-1})
\for |t|>1.
$$

Let us check (\ref{hatsirank1}); we use the truncation
$U_M=\sum_{j=1}^M t^{-j/2}Y^{-j}$ of the series $U$ 
introduced in (\ref{UUplus}) and set 
$U_M^+=U_M\,(1+t^{-1/2}T^{-1})$ for $U_M$ and other
operators. Then 
\begin{align*}
\hat{\mathscr{P}}\,'_M-1^+\ =\ &(1+t^{\frac{1}{2}}T)\,U_M^+\\
=\ &\,U_M^+ + ts_{{}_Y}(U_M)^+ +
\frac{t-1}{Y^{-2}-1}(ts_{{}_Y}(U_M) -U_M)^+\\
=\ \sum_{j=1}^M 
t^{-\frac{j}{2}}&\Bigl( (Y^{-j}+tY^{j}) +
(1-t)\bigl(Y^j+Y^{j-2}+\ldots+Y^{2-j}\bigr)\Bigr)^+
\end{align*}
for $s_{{}_Y}(Y^j)=Y^{-j}$.
Collecting the terms with $Y^{\pm i}$, we obtain that
\begin{align*}
\hat{\mathscr{P}}\,'_M\ =\ 
\sum_{i=1}^M &\Bigl(\bigl(\frac{1-t}{1-t^{-1}}\, 
t^{-\frac{i}{2}}(1-t^{-1-[\frac{M-i}{2}]})+t^{1-\frac{i}{2}}
\bigr)\,Y^i\Bigr)^+\\
+\sum_{i=0}^{M-2} &\Bigl(\bigl(\frac{1-t}{1-t^{-1}}\, 
t^{-1-\frac{i}{2}}(1-t^{-[\frac{M-i}{2}]})+t^{-\frac{i}{2}}
\bigr)\,Y^{-i}\Bigr)^+ \\
+\ &\Bigl(
t^{-\frac{M}{2}}Y^{-M}\Bigr)^+
+\Bigl(t^{\frac{1}{2}-\frac{M}{2}}Y^{1-M}
\Bigr)^+,
\end{align*}
where the last term is present only for $M\ge 2$.
For $M=1$:
$$\hat{\mathscr{P}}\,'_M=1^+ 
+(1+t^{1/2}T)(t^{-1/2}Y^{-1})^+=
1^+ +t^{-1/2}(Y+Y^{-1})^+,$$ 
which immediately follows from (\ref{dahaone}).

As we have already checked, this sum becomes identically
zero as $M\to \infty$.
Therefore significant algebraic simplifications are granted;
only the terms containing $M$ will contribute.

Finally,
\begin{align*}
\hat{\mathscr{P}}\,'_M=
\Bigl(
\sum_{i=1}^M &t^{-\frac{i}{2}-[\frac{M-i}{2}]})Y^i+
\sum_{i=0}^{M-2}t^{-\frac{i}{2}-[\frac{M-i}{2}]})Y^{-i}\\
+\,&t^{-\frac{M}{2}}Y^{-M}+t^{\frac{1}{2}-\frac{M}{2}}
Y^{1-M}\Bigr)^+,
\end{align*}
which can be readily transformed to formula
(\ref{hatsirank1}). Claim $(i)$ is checked.
\smallskip

{\em Claim (ii).} Let us demonstrate that
\begin{align}\label{YSiplus}
t^{-\frac{1}{2}}Y\,\hat{\Si}^+_\infty\ =\ 
\hat{\Si}^+_\infty\ =\ 
t^{-\frac{1}{2}}T\,\hat{\Si}^+_\infty\,. 
\end{align}
The second of these formulas 
is an immediate corollary of the
$s_{{}_Y}$\~invariance of $\hat{\Si}^+_\infty$.

Provided the convergence of $\hat{\mathscr{P}}\,'_+$
or (equivalently) $\hat{\Si}^+_{\infty\,}$,  
the first relation from (\ref{YSiplus}) is
formally  equivalent to the
condition
\begin{align}\label{limYM}
&\lim_{M\to \infty} t^{-M/2}(Y^{-M})^+=0. 
\end{align}
Indeed, if 
$\hat{\Si}^+_M$ converges, then so does  
\begin{align*}
t^{-1/2}Y\,\hat{\Si}^+_M\ =\ 
\hat{\Si}^+_{M+1}- (t^{-(M+1)/2}Y^{-M-1}
+t^{-M/2}Y^{-M})^+.
\end{align*}
Thus the condition 
$(t^{-(M+1)/2}Y^{-M-1}-t^{-M/2}Y^{-M})^+\to 0$ as $M\to \infty$
is necessary for the existence of $\hat{\Si}^+_\infty$.
This condition holds if and only if
it is satisfied for each of the two terms separately, which 
is (\ref{limYM}). We conclude that the existence of
$\hat{\Si}^+_\infty$ results in (\ref{limYM}) and the latter,
in its turn, gives the $t^{-\frac{1}{2}}Y$\~invariance
condition from (\ref{YSiplus}).

Then
$$
(1+t^{-1})\hat{\Si}^+_\infty\ 
=\lim_{M\to\infty}\ (1+t^{-1/2}T^{-1})\hat{\Si}_M
(1+t^{-1/2}T^{-1}),
$$
and we see that $\hat{\Si}^+_\infty$ is invariant 
under the action of the {\em anti-involution} of $\h_Y$ 
sending $Y\mapsto Y$ and $T\mapsto T$ (and fixing $t,q$).
Applying this anti-involution to (\ref{YSiplus}), 
we arrive at the counterpart of these relations
with $\hat{\Si}^+_\infty\ $ placed on the left and 
$Y,T$ on the right. 
\smallskip

{\em Claim (iii).}
Finally, relation (\ref{limYM}) 
readily results in (\ref{YSiabsbar}).
\sq
\smallskip

\rmk
It is worth mentioning that the sigma formula for 
$\hat{\mathscr{P}}\,'_+$ makes it
possible to calculate its $C$\~coefficients
{\em directly}\, and establish 
the proportionality with
$\hat{\mathscr{S}}\,'_+\circ \tilde{\mu}$
in the most explicit way.

Theorem \ref{YLEFTNEW}, which is a
continuation of Theorem \ref {YLEFT}, establishes that
the right multiplication of $\overline{\Sigma}_M$ by 
$(1+t^{-1/2}T^{-1})$
(the notation was $\overline{\Sigma}_M^+$) 
is actually not necessary in (\ref{YSiabsbar}).
The following holds:
\begin{align}
&\hat{\mathscr{P}}\,'_+\ =\ 
\lim_{M\to\infty}\overline{\Sigma}_M. 
\label{YSiabsbar1}
\end{align}
The coefficients here
can be treated as formal series in terms of $X_{\al_1}=X^2,$
$X_{\al_0}=q X^{-2},$ $t^{-1}$
or as functions provided that $|t|>1>|q|$.
One needs to check that  
$\lim_{M\to \infty} t^{-M/2}\,Y^{-M}\,=0\,$ without ${}^+$
as in (\ref{limYM}); this formally results in
$$
\overline{\Sigma}^+_{\infty\,}=
(1+t^{-1})\overline{\Sigma}_{\infty\,}.
$$ 
The convergence in the algebraic variant means here that
$t^{-M/2}\,Y^{-M}$ is getting divisible by powers of
 $q$ growing together with $M.$ 
See Theorem \ref{YLEFTNEW} and below, the second formula 
in (\ref{YexpB}) and  Theorem \ref{SIMINV}.
\smallskip

\rmk
We note that under the Kac-Moody limit $t\to\infty$,
formula (\ref{YSiabsbar}) leads to a presentation
of the Kac-Moody characters introduced for affine
dominant weights as
inductive limits of the corresponding Demazure characters.
It can be used of course for arbitrary weights, not necessarily
dominant, or even for arbitrary functions provided the
convergence, which is an interesting development of this
classical direction. Actually 
$$
\overline{\Si}_M\equal\frac
{t^{-\frac{M}{2}}\,Y^M+t^{-\frac{M-1}{2}}\,Y^{M-1}}{1-t^{-1}},
$$
applied to $q^{l x^2/4}$,
and its generalization via $\overline{\Sigma}_{\mathbf b}$
from Theorem \ref{YLEFT} 
can be considered as certain $q,t$\~Demazure characters.
\sq

\subsubsection{\sf Stabilization of Y-powers}
Let us provide explicit analysis of the limits of the
powers $Y$\~operators, including the coefficient-wise convergence 
of $\overline{\Sigma}_\infty$ and $\hat{\mathscr{P}}\,'_+$;
see Theorem \ref{SYMRANK1Y}.
Recall that we expand operators in the form $\sum_{\hw}C_{\hw}\hw$,
where $C_{\hw}$ can be considered as formal series or 
functions of $X$. Let us treat them as (meromorphic) functions. 
  


Note that if we know that the $C$\~coefficients are meromorphic
functions, this does not  guarantees that
this operator converges in the corresponding
space. For instance, when acting in the polynomial
representation $\mathscr{X}$,
it is well defined at 
a given Laurent polynomials $P(X)$ only for sufficiently
large negative $\Re k$  (depending on $P$), which is
significantly worse than the condition $|qt^{-2}|<1$ 
(necessary and) sufficient for the coefficient-wise 
convergence of  $\hat{\mathscr{P}}\,'_+$.

In contrast to the case $l=0$, 
the convergence of $\hat{\mathscr{P}}\,'_+$
in the spaces $\mathscr{X}q^{lx^2}$ for $l>0$ 
is equivalent to the existence of the corresponding
$\{C_{\hw}\}$ (considered in the next theorem).  
It is with a reservation concerning $l=1$, where
the operator $\hat{\mathscr{P}}\,'_+$  
is well defined for any $t$. This fact is not very surprising due
to the presence of the Gaussians;
the growth of the $C_{\hw}$\~coefficients is no greater
than exponential in terms of $l(\hw)$. 

The following theorem is directly related to 
Theorems  \ref{TQ1H} and \ref{GENPROPR}.
 
\begin{theorem} \label{YEXPA} Continuing to
assume that $|q|<1$, we represent:
\begin{align}\label{YexpA}
\for |t|<1\,:\ \,t^{\frac{m}{2}}q^{-\frac{m}{2}}Y^{-m}&=
\sum_{\hw\in \hW}A^{(-m)}_{\hw}(X)\,\hw \\ 
\and \ \,t^{\frac{m}{2}}Y^{m}&=
\sum_{\hw\in \hW}A^{(m)}_{\hw}(X)\,\hw \,,\notag\\
\label{YexpB}
\for |t|>1\,:\ t^{-\frac{m}{2}}q^{-\frac{m}{2}}Y^{-m}&=
\sum_{\hw\in \hW}B^{(-m)}_{\hw}(X)\,\hw \\ 
\and \ t^{-\frac{m}{2}}Y^{m}&=
\sum_{\hw\in \hW}B^{(m)}_{\hw}(X)\,\hw \,,\notag
\end{align}
where $m\in \Z_+$. These are just algebraic expansions
in the polynomial representations;
the sums are finite.

The claim is that, given $\hw\in \hW$, the limits 
$A^{\pm\infty}_{\hw}=\lim_{m\to\infty} A^{(\pm m)}_{\hw}$
and
$B^{\pm\infty}_{\hw}=\lim_{m\to\infty} B^{(\pm m)}_{\hw}$
exist and are meromorphic functions in terms of $X^2$
analytic apart from $0\neq X^2 \not\in q^{\Z}$. 
\sq\end{theorem}

Using the second formula, we see that the operator
$t^{-m/2}Y^{-m}$ for $|t|>1$  has the coefficients tending to
zero as $m\to\infty$. Indeed, given $\hw\in \hW$, the
coefficient $B^{(-m)}_{\hw}(X)$ behaves as 
$q^{m/2}B^{-\infty}_{\hw}(X)$ in the limit of large $m>0$.
Similarly, $t^{-m/2}Y^{-m}$ has the $A$\~coefficients (for
$|t|<1$) convergent to zero as $m\to\infty$ if 
$|qt^{-2}|<1$. 

It readily results in formula (\ref{limYM}) needed
above. We obtain that the $C_{\hw}$\~coefficients of 
$\hat{\mathscr{P}}\,'_+$ are meromorphic functions 
when $|qt^{-2}|<1$. Here one can use (\ref{YSiabsbar})
or directly (\ref{hatprank1}).
\smallskip

\rmk
Note that the case $|t|=1$ is not covered by 
Theorem  \ref{YEXPA}. In this case, the $A,B$\~coefficients
remain bounded for large $-m$, which is sufficient
for the application to $\hat{\mathscr{P}}\,'_+$. 
\sq
\smallskip

The theorem is closely connected with the action of
$Y^{\pm m}$ in the polynomial representation. For instance,
the first line of (\ref{YexpB}) is related to
the fact that for any given $n\in \Z_+$,\,
$\lim_{m\to \infty}t^{-m/2}Y^{-m}(X^{\pm n})= 0$,
provided that $|q|<1$ and $|tq^{n/2}|>1$, i.e.,
for sufficiently large $t$, exactly, when
$|t|>|q|^{-n/2}$. This fact was actually used 
in Theorem \ref{LEVZERO} (for arbitrary root systems). 
It can be readily checked by expressing
$X^{\pm n}$ in terms of the $E$\~polynomials. 
For the latter,  
\begin{align*}
&t^{-m/2}Y^{-m}(E_{-n})=t^{-m}q^{-mn/2}E_{-n}=
(t q^{n/2})^{-m} E_{-n} \for n\ge 0, \\
&t^{-m/2}Y^{-m}(E_{n})\ =\ q^{mn/2}E_{n} \for n>0.
\end{align*} 
\smallskip

\subsubsection{\sf More on stabilization}
The expansions from (\ref{YexpB}) for the $B$\~coefficients 
and the relations from (\ref{YSiabs}) are of clear
algebraic nature. Let us demonstrate it. The expansion of
the operators
in the following theorem will be considered in the polynomial 
representation as above, however we will now treat their 
coefficients as formal $q$\~series.

\begin{theorem} \label{SIMINV}
(i) The $C$\~coefficients in the expansion
$t^{-m/2}Y^{-m}=\sum_{\hw\in \hW}C^{(-m)}_{\hw}\hw$ 
\, for $m\ge 0$\, are from the ring  
$$
\mathbb{X}\ 
=\ \Z[t^{-1},q^{1/2},X^{\pm 2},(1-q^l X^{\pm 2r})^{-1})],
$$
where $l,r\in \Z_+,\, r>0,\, l>0$ for $-2r$.
Moreover, the coefficient $C^{(-m)}_{\,w\cdot b\,}$, where 
$w=1,s$ and $b=\pm n$ for $m\ge n\ge 0$, belongs to
$q^{(m-n)/2}\,\mathbb{X}\subset\mathbb{X}$.

(ii) In particular, the coefficients of
$w\cdot(\pm n)$ in the $\hw$\~expansions of 
\begin{align}\label{YSiprime}
&t^{-\frac{1}{2}}Y\,\hat{\Si}^+_M-\hat{\Si}^+_M \and 
t^{-\frac{1}{2}}T\,\hat{\Si}^+_M-\hat{\Si}^+_M
\end{align}
belong to the ideal\,
$q^{(M-n)/2}\,\mathbb{X}$\,
for  $0\le n\le M$. If $n$ is fixed and $M\to\infty$,
these coefficients tend to zero with respect to the system 
of ideals $q^{m}\,\mathbb{X}$ for $m\to\infty$.\sq
\end{theorem}
\smallskip

This theorem is a refined  $A_1$\~version of the corresponding 
(algebraic) part $(ii)$ of Theorem \ref{YLEFTNEW}, which 
established that 
given $\hw\in \hW$, the coefficients $C_{\hw}$,
counterparts of the coefficients $C^{(-m)}_{\hw}$, are divisible
by powers of $q$ growing linearly in the limit 
$\mathbf b \to \infty$, corresponding to $m\to \infty$.

\smallskip

{\em Exact formulas as $t=0$ and $t\to \infty$}.
Let us provide the first several exact formulas for the coefficients
$A_{\hw}^{(\pm m)}$ and $B_{\hw}^{(\pm m)}$ from Theorem
\ref{YEXPA} in the corresponding limits $t\to 0$ and
$t\to \infty$. We will consider only the case of even
powers. Then $\hat{w}$ that appear in the formulas
can be represented as 
$$
[n;\ep]\,\equal\, \Ga^{2n}\, s^\ep\,=\,(-2n\om)\, s^\ep\, 
\for n\in \Z,\ \ep=0,1.
$$ 
When $m\ge 1$, the range of nonzero terms in the 
$A,B$\~coefficients is
\begin{align}\label{ab-range}
&-m\le n\le m-1, \where \ep=0,1
\for A_{[n;\,\ep]}^{(+2m)}, B_{[n;\,\ep]}^{(+2m)},\\
&-m\le n\le m,\ \ep=1 \hbox{\, if\, }\ n=-m,\,
\for A_{[n;\,\ep]}^{(-2m)}, B_{[n;\,\ep]}^{(-2m)}.\notag
\end{align}
The elements $\hw$ not in the form $[n;\ep]$ will not
contribute.

We set $\bar{A},\bar{B}$ for the limits of these coefficients
respectively for $t\to 0$ and $t\to \infty$.
The formulas below (they are known and are 
``pure" products for any coefficients) 
are of importance when analyzing the relations to
the Demazure characters in Kac-Moody theory and, hopefully,
for the study of the $t$\~deformations of the Demazure 
characters. 
\smallskip

First, 
\begin{align}\label{ab-zero}
&\bar{A}_{[0;0]}^{(\pm 2m)}\ =\ \bar{B}_{[0;0]}^{(\pm 2m)}\ =\ 
C_0^m\,\equal\, \prod_{i=2}^{2m-1}(1-q^i)\\
\times\bigl(&\,\prod_{i=2}^{2[\frac{m}{2}]}(1-q^i)\,
\prod_{i=1}^{2[\frac{m+1}{2}]-1}(1-q^i)\,
\prod_{i=0}^{m-1}
(1-q^{i+1}X^{-2})(1-q^i X^2)\,\bigr)^{-1}\notag
\end{align}
for $m>1$  and with $C_0^1\equal 1$ as $m=1$
(any signs of $\pm 2m$); here
$[m/2]$ is the integer part of $m/2$. 

\smallskip
Second, the case of the reflection,
\begin{align}\label{ab-min}
&\bar{A}_{[0;1]}^{(2m)}\ =\ -C_0^m,\ \ \ \ \,
\bar{A}_{[0;1]}^{(-2m)}\ =\ -X^2 C_0^m \,
\frac{1-q^m X^{-2}}{1-q^m X^2},\\
&\bar{B}_{[0;1]}^{(2m)}\ =\ -X^2 C_0^m ,\ 
\bar{B}_{[0;1]}^{(-2m)}\ =\ -X^4 C_0^m \,
\frac{1-q^m X^{-2}}{1-q^m X^2}.
\notag
\end{align}
\smallskip

Then for $\Ga^2$ (an element of length $2$)\,:
\begin{align}\label{abc10}
&\bar{A}_{[1;0]}^{(2m)}=C_0^m \frac{(1-q^{m-1})}{(1-q^{m+1})}
\frac{(1-q^m X^{-2})}{(1-q^m X^2)},\\ 
&\bar{A}_{[1;0]}^{(-2m)}=q^{-1} C_0^m \frac{1-q^m X^{-2}}{1-q^m X^2},
\notag\\
&\bar{B}_{[1;0]}^{(2m)}=q X^4 C_0^m \frac{(1-q^{m-1})}{(1-q^{m+1})}
\frac{(1-q^m X^{-2})}{(1-q^m X^2)},\notag\\ 
&\bar{B}_{[1;0]}^{(-2m)}=X^4 C_0^m \frac{1-q^m X^{-2}}{1-q^m X^2},
\notag
\end{align}
and for $s \Ga^2$ (its length is $1$):
\begin{align}\label{abc-10}
&\bar{A}_{[-1;1]}^{(2m)}=-C_0^m 
\frac{1-q^{m-1}X^2}{1-q^{m+1}X^{-2}},\ \
\bar{A}_{[-1;1]}^{(-2m)}=-q^{-1}X^2 C_0^m ,\\
&\bar{B}_{[-1;1]}^{(2m)}=-q X^{-2} C_0^m 
\frac{1-q^{m-1}X^2}{1-q^{m+1}X^{-2}},\ \
\bar{B}_{[-1;1]}^{(-2m)}=-C_0^m.
\notag
\end{align}
\smallskip

Finally, $\Ga^2 s$ (an element of length $3$),
\begin{align}\label{abc11}
&\bar{A}_{[1;1]}^{(2m)}=-C_0^m \frac{(1-q^{m-1})}{(1-q^{m+1})}
\frac{(1-q^{m}X^{-2})}{(1-q^{m}X^{2})},\\
&\bar{A}_{[1;1]}^{(-2m)}=-q X^2 C_0^m \,
\frac{(1-q^{m-1})}{(1-q^{m+1})}
\frac{(1-q^{m-1}X^{-2})(1-q^{m}X^{-2})}
{(1-q^{m}X^{2})(1-q^{m+1}X^{2})}\notag,\\
&\bar{B}_{[1;1]}^{(2m)}=-q^3 X^{6} C_0^m 
\frac{(1-q^{m-1})}{(1-q^{m+1})}
\frac{(1-q^{m-1}X^{-2})}{(1-q^{m+1}X^2)},\notag\\ 
&\bar{B}_{[1;1]}^{(-2m)}=-q^4 X^{8}C_0^m 
\frac{(1-q^{m-1})}{(1-q^{m+1})}
\frac{(1-q^{m-1}X^{-2})(1-q^{m}X^{-2})}
{(1-q^{m}X^{2})(1-q^{m+1}X^{2})}.
\notag
\end{align}
\medskip

\setcounter{equation}{0}
\section{\sc{Spinor Whittaker function}}
\subsection{{\bf Q-Hermite polynomials}}
We will begin with the limiting procedures
connecting $q$\~Toda theory with the difference
QMBP.

\subsubsection{\sf The Ruijsenaars limit} 
\label{sect:whitlim}
Recall the definition of the $L$\~operator
from (\ref{Lopertaor}) :
\begin{align}\label{Lopertaorx}
\l=\frac{t^{1/2}X-t^{-1/2}X^{-1}}{X-X^{-1}}\Ga+
\frac{t^{1/2}X^{-1}-t^{-1/2}X}{X^{-1}-X^{1}}\Ga^{-1},
\end{align}
where we set $\Ga(f(x))=f(x+1/2),\, \Ga(X^n)=q^{n/2}X$
for $X=q^x$. It is symmetric with respect to the action of
$s\,:\, X\mapsto X^{-1},\, \Ga\mapsto \Ga^{-1}$.

This operator preserves
the space of {\em symmetric} Laurent polynomials.
The space of all Laurent polynomials will be denoted by
$\mathscr{X}=\C_{q,t}[X^{\pm}]$, where the field
of definition is  $\C_{q,t}\equal\C(q^{1/2},t^{1/2})$.
 
The Rogers polynomials $P_n\in \mathscr{X}\, (n\ge 0)$
are the eigenfunctions of $\l$ normalized by the 
conditions $P_n=X^n+X^{-n}+$``lower terms".
The eigenvalues are as follows (see (\ref{eqn:LP}))\,:
\begin{equation}{\label{eqn:LPx}}
\l(P_n)\ = \ (q^{n/2}t^{1/2}+q^{-n/2}t^{-1/2})\,P_n,\ n\ge 0.
\end{equation}

In this section, $|q|<1$ and $t=q^k$ for $k\in \C$.
We will use the difference operator 
$\Ga_k(X^n)\equal t^{k/2}X^n$, 

Following Ruijsenaars, Etingof demonstrates in
\cite{Et1} that 
$$
\lim_{k\to -\infty}q^{-kx}\Ga_{k}\,\mathcal{L}\,\Ga_{-k}q^{kx}
$$
becomes the so-called 
$q$\~Toda (difference) operator. To be exact, they
considered the case of $A_n$. The {\em difference} Toda operators 
of type $A_n$ are due to Ruijsenaars too;
see e.g., \cite{Rui}. Inozemtsev extended Ruijsenaars'
limiting procedure to the case of differential {\em periodic} 
Toda lattice (which we do not consider here).

The $A_n$ is exceptional because 
all fundamental weights are minuscule and the formulas
for the Macdonald-Ruijsenaars difference QMBP operators are
explicit. The justification of this limiting procedure
in the case of arbitrary (reduced) root systems (conjectured
by Etingof) was obtained in \cite{ChW}; one can employ the 
Dunkl operators in Macdonald theory or use directly the 
formula for the {\em global} $q$\~Whittaker function from 
\cite{ChW}. It is worth mentioning that the classical 
integrability (at the level of the Poisson brackets) of QMBP 
and the classical Toda chain is significantly simpler
than that of its quantum (operator) generalization.

Following \cite{ChW}, we tend $k$ to $\infty$ 
($t\to 0$) in this section. Let
\begin{align*}
\hbox{\ae}(\mathcal{L})&\equal 
q^{kx}\Ga_{k}^{-1}\,\mathcal{L}\,\Ga_{k}q^{-kx},\ 
R\!E(\l)\equal \lim_{k\to \infty}
\hbox{\ae}(\l),
\end{align*}
where the second limit is
the {\em Ruijsenaars-Etingof procedure}.
At the level of functions $F(X)$:
$$R\!E(F)=\lim_{k\to \infty}q^{kx}\,F(q^{-k/2}X)
=\lim_{k\to \infty}q^{kx}\Ga_{k}^{-1}(F).
$$
Generally, the $R\!E$ procedure 
requires very specific functions $F$ to
be well defined. Formally, if
$\l(\Phi)=(\Lambda+\Lambda^{-1})\Phi$, then
$$
R\!E(\l)(\w)
=(\Lambda+\Lambda^{-1})\mathcal{W} \for \w=R\!E(\Phi)
\hbox{\, provided its existence}.
$$
At the level of operators,
\begin{align}\label{whitoper}
\hbox{\ae}(\mathcal{L})
&=\,\frac{X-X^{-1}}{t^{-1/2}X-t^{1/2}X^{-1}}t^{-1/2}\Ga
+\frac{tX^{-1}-t^{-1}X}{t^{1/2}X^{-1}-t^{-1/2}X}t^{1/2}\Ga
\notag\\
&=\frac{X-X^{-1}}{X-tX^{-1}}\Ga
+\frac{t^{2}X^{-1}-X}{tX^{-1}-X}\Ga^{-1}.
\end{align}
Therefore  
\begin{align}\label{qToda}
R\!E(\l)
=\frac{X-X^{-1}}{X}\Ga+\Ga^{-1}=
(1-X^{-2})\Ga+\Ga^{-1},
\end{align}
where the latter is the {\em $q$\~Toda operator}.

One of the main results of \cite{ChW} states that 
the $R\!E$\~image of the
{\em global 
$q,t$\~spherical function}  (arbitrary reduced root
systems; see the definition there)
is as follows:
\begin{align}\label{Whitsym}
\w_{q}(X,\La)\,=\,
\sum_{m=0}^{\infty}q^{m^{2}/4}\,X^m \overline{P}_{m}(\La)\,
\prod_{s=1}^{m}\frac{1}{1-q^{s}}\,q^{x^{2}}\,q^{\la^{2}},
\end{align}
where  $\prod_{s=1}^0=1$,  
$\La=q^\la$ as for $X$, $\overline{P}_m$ are
the symmetric $q$\~Hermite polynomials, defined as
the specializations of $P_m$ at $t=0$. 

The existence of  $\{\overline{P}_m\}$
can be readily deduced from the explicit formulas
from the previous part of the paper. It will be
discussed systematically (from scratch) below.

One of the key properties of $\w_{q}(X,\La)$
is the Shintani-type formula; see \cite{ChW}.
Setting
$
\widetilde{\w}_q(X,\La)
\equal \sum_{m=0}^\infty 
q^{\frac{m^2}{4}} 
\frac{X^m\ \overline{P}_m(\La) } 
{\prod_{s=1}^{m} (1-q^{s})
}\,$
one has:\, $\widetilde{\w}_q(q^{n/2},\La)=0$\, for\, $n>0$\, and
\begin{align}\label{shintqa1}
&q^{n^2/4}\,\widetilde{\w}_q(q^{-n/2},\La)\ =\
\th(\La)\overline{P}_n(\La) 
\prod_{j=1}^{\infty}\Bigl(
\frac{1}{1-q^j }\Bigr),
\end{align}
where $\,n\ge 0, \ \th(\La)\, =\, 
\sum_{j=-\infty}^\infty q^{j^2/4}\La^j.$

\subsubsection{\sf Nonsymmetric polynomials}
\label{sect:QHermite}
We will use the $E$\~polynomials $E_{a}\in \mathscr{X}$
from the previous part of the paper, which are the eigenfunctions
of the difference Dunkl operator
$$
Y\equal\Ga^{-1}\circ(t^{1/2}+\frac{t^{1/2}-t^{-1/2}}{X^{-2}-1}
\circ (1-s)).
$$
Namely, see (\ref{nonsymp}) above,
\begin{align}\label{nonsympx}
YE_{n}=q^{-n_{\#}}E_{n}\for n\in Z,&&&&\\
n_{\#}=\left\{\begin{array}{ccc}\frac{n+k}{2} 
&  & n>0, \\\frac{n-k}{2} &  & n\le 0,\end{array}\right\}, 
\text{\, note that }\,  0_{\#}=-\frac{k}{2}.
\end{align}
The normalization is
$E_{n}=X^{n}+\text{ ``lower terms'' },$
where by ``lower terms'', we mean
polynomials in terms of $X^{\pm m}$ as $|m|<n$
and, additionally,  $X^{|n|}$ for negative $n$.

 Let us define  their two limits:
$$
\tilde{E}_{a}=\lim_{t\to\infty}E_{a} \and
\overline{E}_{a}=\lim_{t\to 0}E_{a}.
$$
Both limits exist (use the explicit formulas or the intertwining
operators from the previous part of the
paper) and are closely connected to each other. The
following theorem provides the connection.

\begin{theorem}
For $n\geq 0$,
\begin{align}\label{tildee}
\tilde{E}_{-n}\,
=\,\bigl(q^{\frac{n}{2}}\overline{E}_{-n}(Xq^{\frac{1}{2}})\bigr)
\Big|_{q\to q^{-1}},\ 
\tilde{E}_{n}\,
=\,\bigl(q^{-\frac{n}{2}}\overline{E}_{n}(Xq^{\frac{1}{2}})\bigr)
\Big|_{q\to q^{-1}}.
\end{align}
\end{theorem}
\hfill\sq 

The polynomials $\overline{E}_{a}$ are called 
{\em nonsymmetric (continuous)
$q$\~Hermite polynomials} (see \cite{ChW} and references
therein). Upon the substitution $X\mapsto X^{-1}$, 
the polynomials $\overline{E}_a$ are directly 
connected with the Demazure characters 
of level-one Kac-Moody integrable modules; see
\cite{San} for the $GL_n$\~case. Generally it holds
only for the twisted affinization; see \cite{Ion1}. 
These polynomials also appear naturally in formulas 
$\hat{\chi}_{a}^{(l=1)}$ from 
(\ref{Kac-Moody}),
when the latter are used for arbitrary $a\in P$; 
see also 
(\ref{Kac-M-Herm}) there. 

More systematically, let us define 
\begin{align*}
\overline{T}\,\equal\, \lim_{t\to 0}t^{1/2}T=
\frac{1}{1-X^{2}}\circ (s-1),\ 
\overline{T}(\overline{T}+1)=0.
\end{align*}
Using intertwiners, $\overline{E_0}=1$,
\begin{align*}
\overline{E}_{1+n}\,=\,q^{n/2}\Pi \overline{E}_{-n},\\
\overline{E}_{-n}\,=\,(\overline{T}+1)\overline{E}_{n}
\end{align*}
for $n\ge 0$; the raising operator $\Pi\equal X\pi$ was discussed 
in 
Section \ref{sect:intertw}.

From the divisibility condition 
$\overline{T}+1=(s+1)\cdot\{\,\}$, we obtain that
$\overline{E}_{-n}$ is symmetric ($s$\~invariant) and
$\overline{P}_{n}=\overline{E}_{-n}$ for $n\ge 0.$

Explicitly,
\begin{align*}
\overline{E}_{-n-1}\,=\,((\overline{T}+1)\Pi q^{n/2})
\overline{E}_{-n},\\
(\overline{T}+1)\Pi
\,=\,\frac{X^{2}\Ga^{-1}-X^{-2}\Ga}{X-X^{-1}}.
\end{align*}

The bar-Pieri rules read as follows:
\begin{align}\label{pienilp+}
&X^{-1}\overline{E}_{-n}=\overline{E}_{-n-1}-\overline{E}_{n+1}
\ (n\ge 0),\\
&X^{-1}\overline{E}_n=(1-q^{n-1})\overline{E}_{n-1}+q^{n-1}
\overline{E}_{1-n}\ (n\ge 1),\notag\\
\label{pienilp-}
&X\overline{E}_{-n}=(1-q^n)\overline{E}_{1-n}+\overline{E}_{n+1}
\ (n\ge 0),\\
&X\overline{E}_n=\overline{E}_{n+1}-q^{n}
\overline{E}_{1-n}\ (n\ge 1).\notag
\end{align}

Let $\overline{Y}=\pi \overline{T}=\lim_{t\to 0}t^{1/2}Y$.
Recall that
\begin{eqnarray*}
YE_{n}=\left\{\begin{array}{ccc}t^{-1/2}q^{-n/2}E_n, &  & n>0, \\
t^{1/2}q^{n/2}E_{n}, &  & n\leq 0.\end{array}\right.
\end{eqnarray*}
In the limit,
\begin{eqnarray}\label{Ynilp+}
\overline{Y}\,\overline{E}_{n}
=\left\{\begin{array}{ccc}q^{-|n|/2}\overline{E}_{n}, &  & n>0, \\
0, &  & n\le 0.\end{array}\right.
\end{eqnarray}
Since $\overline{Y}$ is not invertible, we need to introduce
$$
\overline{Y}\,'=\lim_{t\to 0}t^{1/2}Y^{-1}=
\lim_{t\to 0}t^{1/2}T^{-1}\pi=\overline{T}\,'\pi
$$ 
for $\overline{T}\,'=\overline{T}+1$. Then 
$\overline{Y}\,\overline{Y}\,'=0=\overline{Y}\,'\,\overline{Y}$ and
\begin{eqnarray}\label{Ynilp-}
\overline{Y}\,'\,\overline{E}_{n}
=\left\{\begin{array}{ccc}q^{-|n|/2}\overline{E}_{n}, 
&  & n\le 0, \\
0, &  & n>0.\end{array}\right.
\end{eqnarray}
Finally, see 
\eqref{eqn:LP}, 
\begin{align*}
\overline{\l}=\lim_{t\to 0}t^{1/2}\l=\overline{Y}\,'
+\overline{Y}
\,=\,\frac{1}{1-X^{2}}\Ga+\frac{1}{1-X^{-2}}\Ga^{-1}
\end{align*}
and $\overline{\l}\,\overline{P}_{n}=
q^{-n/2}\overline{P}_{n},\, n\ge 0$; recall that 
$\overline{P}_n=\overline{E}_{-n}$.

\subsubsection{\sf Nil-DAHA}
We come to the following definition of {\em
nil-DAHA\,} (which can be 
readily adjusted to any reduced root systems).

\begin{theorem}\label{TWONIL}
(i) {\sf Nil-DAHA}\, $\overline{\HH}_+$ is 
generated by $T,\pi_+, X^{\pm 1}$ over the ring $\C[q^{\pm 1/4}]$
with the defining relations $T(T+1)=0$,
\begin{align}\label{nildahax}
\pi_+^2=1,\ \pi_+ X\pi_+=q^{1/2}X^{-1}, 
\ TX-X^{-1}T=X^{-1},
\end{align}
resulting in $X^{-1}=XT-TX^{-1}$.
Setting $\,Y\equal\pi_+ T\,$ and $\,Y'\equal T'\pi_+\,$
for $\,T'\equal (T+1)\,$, the relation $TT'=0$ gives that
(\ref{nildahax}) gives that 
$\,TY-Y'T=-Y,\  TY'=0=YT'$,\, which results in \, $TY'-YT=Y$.
 
(ii) Similarly, one can define
$\overline{\HH}_-=\C[q^{\pm 1/4}]\lan T,\pi_-, Y^{\pm 1}\ran$
subject to $T(T+1)=0$\, and
\begin{align}\label{nildahay}
&\pi_-^2=1,\ \pi_- Y\pi_-=q^{-1/2}Y^{-1},\\ 
&TY-Y^{-1}T=-Y\, \Rightarrow\, YT-TY^{-1}=-Y.\notag 
\end{align}
Setting $X\equal\pi_- T',\,
X'\equal T\pi_-,\, T'=T+1,$\ one has
$$
TX-X'T=X',\ T'X'=0=XT\, \Rightarrow\, TX'-XT=-X'.
$$ 

(iii) The algebra $\,\overline{\HH}_-\,$ is the image of
the algebra  
$\,\overline{\HH}_+\,$ under the 
anti-isomorphism
$$
\varphi: T\mapsto T,\, \pi_+\mapsto\pi_-,\, X\mapsto Y^{-1}.
$$
Correspondingly, $\varphi: Y\mapsto X', Y'\mapsto X$. 
There is also an isomorphism $\si:$
$\,\overline{\HH}_+\,\to\,\overline{\HH}_-\,$ sending
\begin{align*}
&\si: T\mapsto T,\ \,X\mapsto Y,\ \,\pi_+\mapsto \pi_-,\\
&\si:\, Y\mapsto \pi_-T,\ \,
Y'\mapsto T'\pi_-.
\end{align*}

(iv) The automorphism $\tau_+$ fixing $T,X$ and sending
$Y\mapsto q^{-1/4}XY$ acts in $\,\overline{\HH}_+\,$. 
Correspondingly,
$\tau_-\equal\varphi\tau_+\varphi^{-1}$ acts in $\overline{\HH}_-$
preserving $T,Y$ and sending \,
$X\mapsto q^{1/4}YX$. One has
the relations 
\begin{align}\label{sitausi}
&\si\tau_+\ =\ \tau_-^{-1}\si, \ \,
\si\tau_+^{-1}\ =\ \tau_-\si, 
\end{align}
matching the identity from 
(\ref{tautautau}) 
in the generic case. \sq
\end{theorem}
Both algebras $\,\overline{\HH}_{\pm}\,$ satisfy the PBW
Theorem, so $\HH$ is their {\em flat} deformation. 
The formulas above give
an explicit description of the {\em bar-polynomial\,} 
representation of $\overline{\HH}_+$ in $\mathscr{X}=
\C_{q}[X^{\pm 1}]$; recall that 
$T,\pi_+, X^{\pm 1},Y,Y'$ are mapped to the operators 
$\overline{T},\pi, X^{\pm 1}, \overline{Y},\overline{Y}\,'$.
It holds even if $q$ is a root of unity, including 
the construction of the $q$\~Hermite polynomials
(use the intertwiners).

A surprising fact is that the construction of
nonsymmetric Whittaker functions
naturally leads to a module over
$\overline{\HH}_-\,$ that is similar to $\mathscr{X}$
as a vector space but has a very different module structure.
We will call it later the {\em hat-polynomial\,} 
representation; this will require using the
{\em spinors}, to be discussed next.

\rmk
Let us mention the relation of 
nil-DAHA $\overline{\HH}_+\,$
to the $T$\~equivariant $K_T(\b)$ for affine flag 
varieties $\b$ from \cite{KK} and the Demazure-type
operators on this (commutative) ring considered
in this paper. Here
$T$ is the maximal torus in the Lie group $G$ constructed
by the root system $R$.  

The exact $K$\~theoretic 
interpretation of DAHA was obtained in \cite{GG} (see also 
\cite{GKV}). Namely,
$\HH$ is essentially $K^{T\times \C^*}(\La)$ for a certain 
canonical Lagrangian subspace $\La\subset \t^*(\b\times\b)$,
that is the Grothendieck group of the (derived) category of
$T\times \C^*$\~equivariant coherent sheaves on $\La$.

This interpretation is for arbitrary $q,t$. Switching from $\b$ in 
\cite{KK} to $\La\subset \t^*(\b\times\b)$ is important
because it gives the definition of convolution and, therefore,
supplies $K^{T\times C^*}(\La)$ with a structure of
algebra (isomorphic to $\HH$). We note that 
the Gaussians were added to the definition of DAHA in
\cite{GG}. We prefer not to consider the Gaussians as part 
of the definition of DAHA, 
treating them as outer automorphisms of $\HH\,$,
as in the classical theory of Heisenberg-Weyl algebras 
and metaplectic representations.

\subsection{{\bf Nonsymmetric Q-Toda theory}}
The problem of finding Dunkl operators for 
the {\em  $q$\~Toda operator} from
(\ref{qToda}) seems not well defined since the
Toda operators are not symmetric. Nevertheless, it
has a solution ( below). It provides a spinor variant
of the representation $\l=Y+Y^{-1}$ (upon the restriction
to the symmetric functions) for $\l$ from (\ref{Lopertaorx}).

The spinor-Dunkl operators make it possible to use
DAHA methods at their full potential algebraically
and in the theory of the $q$\~Whittaker functions.
The construction can be extended to arbitrary
root systems (in progress).
We will begin with the introduction of the spinors.

\subsubsection{\sf The spinors}\label{sect:spinors}
Generally, the $W$\~spinors are needed in DAHA theory 
as discussed in the introduction. In the $A_1$\~case,
we will call them simply {\em spinors}. In this case,
they are really connected with spinors
from the theory of the Dirac operator (and with
super-algebras). Under the rational degeneration,
the Dunkl operator for $A_1$ becomes
the square root of the (radial part of the) 
Laplace operator, i.e., the {\em Dirac operator}.
However, this relation (and using super-variables) 
is a special feature of the root system $A_1$. 

For practical calculations with spinors, the language of 
$\Z_2$\~graded algebras can be used in the $A_1$\~case
(see the differential
theory below). However, we prefer to proceed here
in a way that does not rely on the special 
symmetry of the $A_1$\~case and can be transferred
to $W$\~spinors for arbitrary root systems.
\smallskip

The {\em spinors\,} are simply pairs $\{f_1,f_2\}$ of
elements (functions) from a space $\f$ with an action
of  $s$; the addition or multiplication
(if applicable) of spinors is componentwise.
The space of spinors will be denoted by $\hat{\f}$.

The involution $s$ on spinors is defined as follows 
$s\{f_1,f_2\}=\{f_2,f_1\}$, so this does not involve the
action of $s$ in $\f$. There is a ``natural"
embedding $\rho:\f\to \hat{\f}$ mapping
$f\mapsto f^\rho=\{f,s(f)\}$ and the diagonal
embedding $\delta: \f\to \hat{\f}$ sending
$f\mapsto f^\de=\{f,f\}$. Accordingly,
for an arbitrary operator $A$ acting in $\f$,
$A^\rho=\{A,s(A)\}, A^\de=\{A,A\}$. The images
$f^\rho$ of $f\in\f$ are called {\em functions} (in contrast to
{\em spinors}) or {\em principle spinors} (like for adeles). 

For instance, for $\f=\mathscr{X}$,
\begin{align*}
&X^\rho:\{f_1,f_2\}\mapsto\{Xf_1,X^{-1}f_2\},\ 
&&\Gamma^\rho:\{f_1,f_2\}\mapsto\{\Gamma(f_1),
\Gamma^{-1}(f_2)\},\\
&X^\de:\{f_1,f_2\}\mapsto\{Xf_1,Xf_2\},\ 
&&\Gamma^\de:\{f_1,f_2\}\mapsto\{\Gamma(f_1),\Gamma(f_2)\},
\end{align*}
where, recall, $\Gamma(X)=q^{1/2}X$. We simply put
$$
X^\rho=\{X,X^{-1}\},\ \Ga^\rho=\{\Ga,\Ga^{-1}\},\
X^\de=\{X,X\},\ \Ga^\de=\{\Ga,\Ga\}.
$$
Obviously, $s^\rho=s=s^\de$.

If a function $f\in \f$ or an operator $A$ acting in $\f$ 
have no super-index $\de$, then
they will be treated  as $f^\rho, A^\rho$. I.e.,
by default, functions
and operators are embedded 
into $\hat{\f}$ and the algebra of spinor operators
using $\rho$.

If the operator $A$ is explicitly expressed as $\{A_1,A_2\}$,
then $A_1$ and $A_2$ must be applied to the
corresponding components of $f=\{f_1,f_2\}$. In the calculations
below, $A_i$ will be allowed to contain $s$ placed on the right,
i.e., in the form  $A_i=A_i'\cdot s$, where $A_i'$ contains no 
$s$\,. The latter can be always achieved by using the commutation 
relations between $s$ and $X, \Ga$. Then
the component
$i$ of $Af$ will be (by definition) $A_i'(f_{3-i})$. I.e., 
$s$ placed on the right {\em inside a spinor component of
the operator} will mean the switch to the 
other component ($i\mapsto 3-i$)
before applying the rest of the operator, which is $A_i'$.

For instance, $\{\Ga s,\,s-1\}(\{f_1,f_2\})=$
$\{\Ga(f_2), f_1-f_2\}$. 

We will frequently use the vertical mode for spinors: 
$$
\{f_1,f_2\}= 
\left\{\begin{array}{c}f_1 \\f_2\end{array}\right\},\
\{A_1,A_2\}= 
\left\{\begin{array}{c}A_1 \\A_2\end{array}\right\}. 
$$

\comment{
\begin{align*}
&X^\rho:\left\{\begin{array}{c}f_{1} \\f_{2}
\end{array}\right\}
\mapsto\left\{\begin{array}{c}Xf_1 \\X^{-1}f_2
\end{array}\right\}
&&\Gamma^\rho:\left\{\begin{array}{c}f_{1} \\f_{2}
\end{array}\right\}
\mapsto\left\{\begin{array}{c}\Gamma(f_1) \\\Gamma^{-1}(f_2)
\end{array}\right\}
\\
&X^\de:\left\{\begin{array}{c}f_{1} \\f_{2}
\end{array}\right\}
\mapsto\left\{\begin{array}{c}Xf_1 \\Xf_2
\end{array}\right\},\ 
&&\Gamma^\de:\left\{\begin{array}{c}f_{1} \\f_{2}
\end{array}\right\}
\mapsto\left\{\begin{array}{c}\Gamma(f_{1}) 
\\ \Gamma(f_{2})\end{array}\right\},
\end{align*}
}

\subsubsection{\sf Q-Toda via DAHA}
The $q$\~Toda {\em spinor\,} operator is the following
{\em symmetric\,} (i.e., $s$\~invariant) difference
{\em spinor\,} operator
\begin{align}\label{Lspin}
\widehat{\l}\ =\ 
\{\Gamma^{-1}+(1-X^{-2})\Gamma,\,\Gamma^{-1}+(1-X^{-2})\Gamma\}.
\end{align}
Its first component is the operator $R\!E(\mathcal{L})$
from Section \ref{sect:whitlim}; we will use the notation
and definitions
from this section.

We claim that $\widehat{\l}$
can be represented in the form $\hat{Y}+\hat{Y}^{-1}$
upon the restriction to {\em symmetric spinors}, i.e.,  
to $\{f,f\}\in \hat{\f}$. The construction of the 
{\em spinor-difference Dunkl operator}
$\hat{Y}$ goes as follows.

Let us introduce the following map on the operators
in terms of $X,\Ga$ and $s$ with values in
spinor operators: 
\begin{align}\label{aede}
&\hbox{\ae}^\de:\, X\mapsto \tilde{t}^{\ -1/2}X,\,
\Gamma \mapsto \tilde{t}^{\ -1/2}\Gamma,\, s\mapsto s 
\end{align}
for the  
{\em spinor constant\,} 
$\tilde{t}^{\ 1/2}\equal\{t^{1/2},t^{-1/2}\}$.
Spinor constants are constant diagonal matrices; they commute
with $\Ga$ and $X$ but not with $s$ 
unless they are scalar. The 
{\em spinor $R\!E$-construction\,} is 
$$
R\!E^\de:\ A\mapsto \lim_{t\to 0}\, \hbox{\ae}^\de(A).
$$
It is of course very different from the procedure
$R\!E^\rho$ from Section \ref{sect:whitlim}. The spinor-Dunkl
operators are $\hat{Y}=R\!E^\de(Y),\, \hat{Y}\,'=
R\!E^\de(Y^{-1})$.
They are inverse to each other: $\hat{Y}\hat{Y}\,'=1$.

\begin{theorem}\label{MainToda}
The map 
\begin{align*}
&Y^{\pm 1}\mapsto \hat{Y}^{\pm 1},\
\pi_-\mapsto R\!E^\de(XT),\\
&T\mapsto \hat{T}=R\!E^\de(t^{1/2}T),\
T'\mapsto \hat{T}'=R\!E^\de(t^{1/2}T^{-1})\,
\end{align*} 
can be extended to
a representation of the algebra $\,\overline{\HH}_-\,$
in the space $\,\hat{\mathscr{X}}\,$ of spinors over
$\mathscr{X}=\C[q^{\pm 1/4}][X^{\pm 1}]$. 
Correspondingly, 
\begin{align*}
&X\mapsto R\!E^\de(t^{1/2}X)=R\!E^\de(\pi_-)\circ \hat{T}',\\
&X'\mapsto R\!E^\de(t^{1/2}X^{-1})=\hat{T}
\circ R\!E^\de(\pi_-).
\end{align*}
The commutativity
of\, $T$  and\, $Y+Y^{-1}$\, in $\overline{\HH}_-$ results
in the $s$\~invariance of $\ \hat{Y}+\hat{Y}^{-1}\,$ 
and the $s$\~invariance of this operator upon its restriction
to the space of $s$\~invariant spinors, which is the
one from (\ref{Lspin}).\sq
\end{theorem}

It is clear from the construction that all 
hat-operators preserve the space of spinors for
Laurent polynomials in terms of $X^{\pm 1}$.
We will give below explicit formulas.
Upon multiplication by the Gaussian,
this $\overline{\HH}_-$\~module contains an irreducible 
submodule, the {\em spinor polynomial representation},
isomorphic to the Fourier image of the bar-polynomial 
representation times the Gaussian; see 
Section \ref{sect:QHermite}, formula (\ref{sitausi})
and Theorem \ref{SPIN-polyn} below. 
The reproducing kernel of the isomorphism between 
these two modules 
inducing $\si:\overline{\HH}_+\to$ $\overline{\HH}_-$
at the operator level is given by the {\em nonsymmetric 
$q$\~Whittaker function}; its existence was
conjectured in \cite{ChW}.

\subsubsection{\sf Spinor-Dunkl operators}
Let us calculate explicitly the operator
$\hat{Y}=R\!E^{\de}(Y)=
\lim_{t\to 0}\hbox{\ae}^\de(Y)$.

Recall that
$s$ placed on the right {\em inside a spinor component of
the operator} always mean the switch to the 
other component before applying the rest of the operator
in this component. 
\comment{
For instance (see below),
\begin{align*}
&\left\{\begin{array}{c}
\Gamma^{-1}\cdot(1-s) \\
\Gamma-\Gamma\cdot X^{-2}\cdot (1-s)\end{array}\right\}\{f_1,f_2\}
\\
&\ =\ \{\Gamma^{-1}(f_1-f_2), \Gamma(f_2)-\Gamma(X^{-2}(f_2-f_1))\}.
\end{align*}
}

Using formulas (\ref{aede}):
\begin{eqnarray*}
\hbox{\ae}^\de(Y)
&=&s\cdot(\tilde{t}^{\ -1/2}\Gamma)\cdot
\left(t^{1/2}s+\frac{t^{1/2}-t^{-1/2}}
{\tilde{t}^{-1}X^{2}-1}\cdot (s-1)\right)\\
&=&t^{1/2}\tilde{t}^{\ 1/2}\Gamma^{-1}
+\tilde{t}^{\ 1/2}\Gamma^{-1}\cdot\frac{t^{1/2}-t^{-1/2}}
{\tilde{t}X^{-2}-1}\cdot(1-s)\\
&=&\left\{\begin{array}{c}
t\Gamma^{-1}+\Gamma^{-1}\frac{t-1}{tX^{-2}-1}\cdot(1-s) \\
\Gamma+\Gamma\frac{1-t^{-1}}{t^{-1}X^{2}-1}\cdot(1-s)
\end{array}\right\}\\
\xrightarrow{t\to 0}\ \hat{Y}&=&
\left\{\begin{array}{c}
\Gamma^{-1}\cdot(1-s) \\
\Gamma-\Gamma\cdot X^{-2}\cdot (1-s)\end{array}\right\}.
\end{eqnarray*} 
Recall that $\tilde{t}^{\ 1/2}= \{t^{1/2},t^{-1/2}\}$.
A little bit more involved calculation is needed for
$\hat{Y}\,'=R\!E^{\de}(Y^{-1})$:

\begin{eqnarray*}
\hbox{\ae}^\de(Y^{-1})
&=&\left(t^{-1/2}s+\frac{t^{-1/2}-t^{1/2}}{\tilde{t}X^{-2}-1}
\cdot (s-1)\right)\cdot
(\tilde{t}^{\ 1/2}\Gamma^{-1}s)\\
&=&\left(\frac{t^{\ -1/2}\tilde{t}X^{-2}-t^{1/2}}
{\tilde{t}X^{-2}-1}\cdot s
-\frac{t^{\ -1/2}-t^{1/2}}{\tilde{t}X^{-2}-1}\right)
\cdot(\tilde{t}^{\ 1/2}\Gamma^{-1}s)\\
&=&\frac{t^{\ -1/2}\tilde{t}X^{-2}-t^{1/2}}
{\tilde{t}X^{-2}-1}\tilde{t}^{\ -1/2}\Gamma
-\frac{t^{-1/2}-t^{1/2}}{\tilde{t}X^{-2}-1}
\tilde{t}^{\ 1/2}\Gamma^{-1}s\\
&=&\left\{\begin{array}{c}
\frac{X^{-2}-1}{tX^{-2}-1}\Gamma-\frac{1-t}{tX^{-2}-1}
\Gamma^{-1}s\\
\frac{t^{-1}X^{2}-t}{t^{-1}X^{2}-1}
\Gamma^{-1}-\frac{t^{-1}-1}{t^{\,-1}X^{2}-1}\Gamma s
\end{array}\right\}\xrightarrow{t\to 0}\\
\hat{Y}\,'&=& \left\{\begin{array}{c}
(1-X^{-2})\Gamma+\Gamma^{-1}s \\
\Gamma^{-1}-\frac{1}{X^{2}}\Gamma s\end{array}\right\}\\
&=& \left\{\begin{array}{c}1-X^{-2} \\1\end{array}\right\}
\Gamma+\left\{\begin{array}{c}1 \\-X^{2}\end{array}\right\}
\Gamma^{-1}s.
\end{eqnarray*} 
Automatically, $\hat{Y}\hat{Y}\,'=1$, since these operators
were obtained by the $R\!E^{\de}$\~construction.
Now, as we claimed,
\begin{eqnarray*}
&&R\!E^{\de}(Y+Y^{-1})=\lim_{t\to 0}\hbox{\ae}^\de(Y+Y^{-1})\\
&=&\left\{\begin{array}{c}\Gamma^{-1}(1-s)+(1-X^{-2})
\Gamma+\Gamma^{-1}s
 \\\Gamma-\Gamma\frac{1}{X^{2}}(1-s)+\Gamma^{-1}
 -\frac{1}{X^{2}}\Gamma s\end{array}\right\}\\
&=&\left\{\begin{array}{c}\Gamma^{-1}+(1-X^{-2})\Gamma \\
\Gamma^{-1}+(1-X^{-2})\Gamma\end{array}\right\}
\bigl(\text{ mod }\bigl(\cdot\bigr)(s-1)\bigr).
\end{eqnarray*}

For $X$ and $X^{-1}$, we have
\begin{align}{\label{eqn:hatX}}
\hat{X}\!=\!R\!E^{\de}(t^{1/2}X)
=\lim_{t\to 0}\hbox{\ae}^\de(t^{1/2}X)
=\lim_{t\to 0}t^{1/2}\tilde{t}^{\ -1/2}X
=\left\{\begin{array}{c}X \\0\end{array}\right\},\\
\hat{X}\,'\!=\!R\!E^{\de}(t^{1/2}X^{-1})
=\lim_{t\to 0}\hbox{\ae}^\de(t^{1/2}X^{-1})
=\lim_{t\to 0}t^{1/2}\tilde{t}^{\ -1/2}X^{-1}
=\left\{\begin{array}{c}0 \\X\end{array}\right\}.\notag
\end{align}
Obviously, $\hat{X}\hat{X}\,'=0$.
Next,
\begin{align*}
&\hat{T}=R\!E^{\de}(t^{1/2}T)=
\lim_{t\to 0}\hbox{\ae}^\de(t^{1/2}T)
=\left\{\begin{array}{c}0 \\s-1\end{array}\right\},\\
&\hat{T}\,'=R\!E^{\de}(t^{1/2}T^{-1})=\lim_{t\to 0}
\hbox{\ae}^\de(t^{1/2}T^{-1})
=\left\{\begin{array}{c} 1 \\ s \end{array}\right\}.
\end{align*}

It is instructional
to check the following relations using the explicit
formulas we obtained (they of course follow from 
Theorem \ref{MainToda}): 
\begin{align}
&\hat{T}'=\hat{T}+1,\ 
\hat{T}\hat{T}\,'\,=\,0\,=\,\hat{T}\,'\hat{T},\ \, 
\hat{T}\,'\hat{X}\,'=0=\hat{X}\hat{T},
{\label{eqn:TY1}}\\
&\hat{T}\hat{Y}-\hat{Y}^{-1}\hat{T}\,=\,-\hat{Y},
\ \hat{T}\hat{Y}^{-1}-\hat{Y}\hat{T}\,=\,\hat{Y},
{\label{eqn:TY2}}\\
&\hat{T}\hat{X}-\hat{X}\,'\hat{T}=\hat{X}\,',\ 
\hat{T}\hat{X}\,'-\hat{X}\hat{T}=-\hat{X}\,', 
\ \hat{X}+\hat{X'}=X^{\de}.{\label{eqn:TX}}
\end{align}
Relations \eqref{eqn:TY2} imply that
\begin{align}\label{hatTYcom}
\hat{T}(\hat{Y}+\hat{Y}^{-1})=(\hat{Y}+\hat{Y}^{-1})\hat{T}.
\end{align}
It proves that the spinor operator $\hat{Y}+\hat{Y}^{-1}$ is
symmetric (recall that $\hat{Y}\,'=\hat{Y}^{-1}$).
Indeed, applying (\ref{hatTYcom}) to 
a symmetric spinor $\{f,f\}$, 
let $(\hat{Y}+\hat{Y}^{-1})(\{f,f\})=\{g_1,g_2\}$.
Then $\hat{T}(\{g_1,g_2\})=0$, which is possible if and only
if $g_1=g_2$.

\subsubsection{\sf Using the components}
Explicitly, the action of $\hat{Y}$ and $\hat{Y}\,'$ 
on the spinors is as follows:
\begin{align*}
&\hat{Y}(\left\{\begin{array}{c}f_{1} \\f_{2}
\end{array}\right\})
=\left\{\begin{array}{c}\Gamma^{-1}(f_{1}-f_{2}) \\
\Gamma(f_{2})-\Gamma(\frac{f_{2}-f_{1}}{X^{2}})
\end{array}\right\},
\\
&\hat{Y}\,'(\left\{\begin{array}{c}f_{1} \\f_{2}
\end{array}\right\})
=\left\{\begin{array}{c}(1-X^{-2})
\Gamma(f_{1})+\Gamma^{-1}(f_{2}) \\
\Gamma^{-1}(f_{2})-\frac{1}{X^{2}}\Gamma(f_{1})
\end{array}\right\}.
\end{align*}
It is simple but not immediate to check the relation 
$\hat{Y}\hat{Y}\,'=$\,id and other identities for 
$\hat{Y}^{\pm1}$
using the component formulas. 
The explicit formulas for $\hat{T}$ and $\hat{T}\,'$ 
are:
\begin{align}\label{Tspinors}
&\hat{T}(\left\{\begin{array}{c}f_{1} \\f_{2}
\end{array}\right\})
=\left\{\begin{array}{c}0\\f_{1}-f_{2}
\end{array}\right\},
&\hat{T}\,'(\left\{\begin{array}{c}f_{1} \\f_{2}
\end{array}\right\})
=\left\{\begin{array}{c}f_{1}\\
f_{1}\end{array}\right\}.
\end{align}
It readily gives  \eqref{eqn:TY1}, \eqref{eqn:TY2}.

Generally, there is no need to establish and check
the formulas for $\hat{X}$ and $\hat{X}\,'$ (although
they are simple). From Theorem \ref{MainToda}, 
\begin{align*}
&\hat{X}\ =\ R\!E^\de(\pi_-)\cdot \hat{T}\,',\ 
\hat{X}\,'\ =\ \hat{T}
\cdot R\!E^\de(\pi_-).
\end{align*}
Thus we need only to know $\hat{\pi}\equal
R\!E^\de(\pi_-)$, where $\pi_-=XT$.
We have
\begin{eqnarray*}
\hbox{\ae}^\de(XT)
&=&(\tilde{t}^{-1/2}X)(t^{1/2}s
+\frac{t^{1/2}-t^{-1/2}}{\tilde{t}^{-1}X^{2}-1}(s-1))\\
&=&\tilde{t}^{-1/2}t^{1/2}Xs
+\frac{X(\tilde{t}^{-1/2}t^{1/2}-\tilde{t}^{-1/2}t^{-1/2})}
{\tilde{t}^{-1}X^{2}-1}(s-1)\\
&=&\left\{\begin{array}{c}Xs \\tX^{-1}s\end{array}\right\}
+\left\{\begin{array}{c}\frac{X(1-t^{-1})}
{t^{-1}X^{2}-1}(s-1) \\
\frac{X^{-1}(t-1)}{tX^{-2}-1}(s-1)\end{array}\right\}.
\end{eqnarray*}
Taking the limit $t\to 0$, 
\begin{align*}
\hat{\pi}=\left\{\begin{array}{c}Xs \\0\end{array}\right\}
+\left\{\begin{array}{c}-X^{-1}(s-1) \\X^{-1}(s-1)
\end{array}\right\}
=\left\{\begin{array}{c}Xs-X^{-1}(s-1) \\X^{-1}(s-1)
\end{array}\right\}.
\end{align*}
Using the components,
\begin{align}\label{picompo}
\hat{\pi}:\ \left\{\begin{array}{c}f_{1} \\f_{2}
\end{array}\right\}\mapsto
\left\{\begin{array}{c}Xf_{2}+\frac{f_{1}-f_{2}}{X} \\
\frac{f_{1}-f_{2}}{X}\end{array}\right\}.
\end{align}
Check directly that  $\hat{\pi}^{2}=id$.

This formula completes the ``component presentation" of the
{\em hat-module} of $\overline{\HH}_-$
from Theorem \ref{MainToda}:
$$
T, \pi_-, Y\ \mapsto\ 
 \hat{T}, \hat{\pi}, \hat{Y}.
$$

The extension of this Theorem to arbitrary (reduced)
root systems is straightforward as well as the justification;
we will address this (and the applications) in further paper(s).
The formulas for the $\overline{Y}$\~operators are of course
getting more involved. The justifications
in the spinor $q$\~Toda theory (including global Whittaker
functions) are entirely based on DAHA theory.
We calculate and check practically everything
explicitly in this work mainly to demonstrate the practical 
aspects of the technique of spinors (and because of novelty
of this topic).

\subsubsection{\sf Spinor Whittaker function}
Let us apply the procedure $R\!E^\de$ to
the {\em global difference spherical function}
$\e_q(x,\la)$ 
from \cite{C5}, Section 5 (upon the
specialization to the case of $A_1$). We do not give here
its exact definition and do not discuss the details of the
procedure. Actually, the only point that requires comments
is using the conjugated $E$\~polynomials, $E_b^*$ in the
formula for $\e_q$ in \cite{C5}. Generally, the relation
of $\{E_b\}$ and its conjugates is via the action of
$T_{w_0}$; compare with Theorem \ref{tildee}
in the case of $A_1$. 

We arrive at the 
following spinor nonsymmetric generalization
of the function $\w_q$ from (\ref{Whitsym}) above:
\begin{align}
\Om(X,\La)= q^{x^2}q^{\la^2}\,\Bigl(1+
\sum_{m=1}^{\infty} q^{m^{2}/4}\,
\bigl(&\frac{\overline{E}_{-m}(\La)}
{\prod_{s=1}^{m} (1-q^{s})}
\left\{\begin{array}{c}X^m \\ q^mX^m\end{array}\right\}\notag\\
+\,&\frac{\overline{E}_{m}(\La)}
{\prod_{s=1}^{m-1} (1-q^{s})}
\left\{\begin{array}{c} 0 \\ X^m\end{array}\right\}
\bigr)\Bigr).\label{spinwhito}
\end{align}
Using the Pieri rules from (\ref{pienilp-}), we can present 
it as follows:
\begin{align}\label{spinwhit}
\Om= q^{x^2}q^{\la^2}\,
\sum_{m=0}^{\infty} \frac{q^{m^{2}/4}}
{\prod_{s=1}^{m} (1-q^{s})}
\left\{\begin{array}{c}X^m \overline{E}_{-m}(\La) \\ 
X^m \La^{-1}\overline{E}_{m+1}(\La)\end{array}\right\}.
\end{align}
Either of these two presentations readily gives that the
spinor- symmetrization of $\Om$ is $\{\w,\,\w\}$. We need
to apply the symmetrizer $\mathscr{P}'=T'=T+1$ to $\Om$,
equivalently, duplicate its first component;
see (\ref{Tspinors}).
Note that $\La$ is a (nonspinor) variable.

The spinor $\Om$
intertwines the bar-representation of  $\overline{\HH}_+$ and
the hat-representation of  $\overline{\HH}_-$. Namely,
\begin{align}
&\hat{Y}(\Om)=\La^{-1}(\Om),\ 
\hat{X}(\Om)=\overline{Y}\,'_\La (\Om),\ 
\hat{X}\,'(\Om)=\overline{Y}_\La (\Om),\notag\\
&\hat{\pi}(\Om)\ =\ \pi_\La(\Om),\ \, \hat{T}(\Om)\ =\ 
\overline{T}_\La(\Om), \label{Omegainter}
\end{align}
where $\overline{Y}\,'_\La,\, \overline{Y}_\La,\,\pi_\La,\,
\overline{T}_\La$\,
act on the argument $\La$; the other operators are
$X$\~operators. These (and other
related identities) follow from the general
theory for any reduced root systems (at least in
the twisted case). However, in the rank-one case
(and for $A_n$), one can use the Pieri
rules from (\ref{pienilp+}),(\ref{pienilp-})
and formulas (\ref{Ynilp+}), 
(\ref{Ynilp-}) for the direct verification.
 
Let us calculate
$\overline{Y}\,'_\La(\Om)$. 
First,
$\overline{Y}\,'_\La(\overline{E}_n(\La))=0$ for $n>0$.
Second,
$q^{\,-\la^2}\,\overline{Y}\,'_\La\,q^{\la^2}$
$=q^{-1/4}\,\overline{Y}\,'_\La\cdot\La$. 
For instance, 
$$
\overline{Y}\,'_\La(q^{\la^2})=q^{-1/4}
\overline{Y}\,'_\La(\La)\,q^{\la^2}=
\overline{Y}\,'_\La(\overline{E}_1(\La))\,q^{\la^2}=0.
$$
We see that the second spinor component 
of $\overline{Y}\,'_\La(\Om)$ vanishes, as it
is supposed to be because the second component
of  $\hat{X}(\Om)$ is obviously zero.

The first component reads as follows:
\begin{align*}
\overline{Y}\,'_\La(\Om)&= q^{x^2}q^{\la^2}\,
\sum_{m=0}^{\infty} \frac{q^{m^{2}/4-1/4}\,X^m\,
\overline{Y}\,'_\La(\La \overline{E}_{-m})
}{\prod_{s=1}^{m} (1-q^{s})}\\
&= q^{x^2}q^{\la^2}\,
\sum_{m=0}^{\infty} \frac{q^{m^{2}/4-1/4-m/2+1/2}\,X^m\,
(1-q^m)\overline{E}_{1-m}
}{\prod_{s=1}^{m} (1-q^{s})}\\
&=q^{x^2}q^{\la^2}\,
X\,\sum_{m=1}^{\infty} \frac{q^{(m-1)^{2}/4}\,X^{m-1}\,
\overline{E}_{1-m})
}{\prod_{s=1}^{m-1} (1-q^{s})},\\
\end{align*}
which coincides with the first component of $\hat{X}(\Om)$
(its second component is zero).
We have used here the nil-Pieri formula:
$$\La\overline{E}_{-n}=
(1-q^n)\overline{E}_{1-n}+\overline{E}_{n+1} \for n>0;
$$
the second term, $\overline{E}_{n+1}$,  does not
contribute to the final formula,
since $\overline{Y}\,'(\overline{E}_{n+1})=0$.

The (key) relation $\hat{Y}(\Om)=\La^{-1}(\Om)$ can be verified
directly in a similar manner. 
First,
$q^{-x^2}\,\hat{Y}\,q^{x^2}=q^{1/4}X^{-1}\hat{Y}$. Therefore
\begin{align}\label{hatYga}
q^{-x^2}\,\hat{Y}\,q^{x^2}\left\{\begin{array}{c}f_1\\
f_2\end{array}\right\}=q^{1/4}
\left\{\begin{array}{c} X^{-1}\Ga^{-1}(f_1-f_2) \\ 
X\Ga(f_2)+q^{-1}X^{-1}\Ga(f_1-f_2)\end{array}\right\}.
\end{align}
Second, $\f_m\equal \overline{E}_{-m}(\La)-
\La^{-1}\overline{E}_{m+1}(\La)=
(1-q^m)\La^{-1}E_{1-m}(\La)$ (the Pieri rules). Now,\ \,
$\La^{-1}q^{-x^2}q^{-\la^2}\,\hat{Y}(\Om)=$
\begin{align*}
&\La^{-1}\sum_{m=0}^{\infty} \frac{q^{m^{2}/4+1/4}}
{\prod_{s=1}^{m} (1-q^{s})}
\left\{\begin{array}{c} q^{-\frac{m}{2}}X^{m-1}\f_m \\ 
q^{\frac{m}{2}}X^{m+1}\La^{-1}\overline{E}_{m+1}(\La)
+q^{\frac{m}{2}-1}X^{m-1}\f_m
\end{array}\right\}\\
=&\sum_{m=0}^{\infty} \frac{q^{m^{2}/4+1/4}}
{\prod_{s=1}^{m} (1-q^{s})}
\left\{\begin{array}{c} q^{-\frac{m}{2}}(1-q^m)X^{m-1}
\overline{E}_{1-m} \\ 
q^{\frac{m}{2}}X^{m+1}\overline{E}_{m+1}
+q^{\frac{m}{2}-1}(1-q^m)X^{m-1}\overline{E}_{1-m}
\end{array}\right\}.
\end{align*}
Collecting the terms with $(1-q^m)$, we obtain that
\begin{align*}
\hat{Y}(\Om)= &\La^{-1}q^{x^2}q^{\la^2}\,
\sum_{m=1}^{\infty} \frac{q^{(m-1)^{2}/4}}
{\prod_{s=1}^{m-1} (1-q^{s})}
\left\{\begin{array}{c} X^{m-1}\overline{E}_{1-m}(\La) \\ 
q^{m-1}X^{m-1}\overline{E}_{1-m}(\La) 
\end{array}\right\}\\
+ &\La^{-1}q^{x^2}q^{\la^2}
\sum_{m=0}^{\infty} \frac{q^{(m+1)^2/4}}
{\prod_{s=1}^{m} (1-q^{s})}
\left\{\begin{array}{c} 0\\ 
X^{m+1}\overline{E}_{m+1}(\La)
\end{array}\right\},
\end{align*}
i.e., exactly the presentation from (\ref{spinwhito}) multiplied
by $\La^{-1}$.

Formulas (\ref{hatYga}), (\ref{picompo}) and (\ref{Tspinors})
result in the definition of the {\em spinor-polynomial}
representation:
$$
\mathscr{X}_{spin}=\C\oplus\ \bigl
(\oplus_{m=1}^{\ \infty} (\C\{X^m,0\}\oplus
\C\{0,X^m\})\bigr).
$$
\begin{theorem}\label{SPIN-polyn}
The space $\mathscr{X}_{spin}$ is an irreducible 
$\overline{\HH}_-$\~submodule of the space of
spinors over $\C[X^{\pm1}]$ supplied with the twisted
action:
$$
\overline{\HH}_-\ni A\mapsto q^{-x^2}\,\hat{A}\,q^{x^2}.
$$ 
More explicitly, $\mathscr{X}_{spin}$
is invariant and irreducible under the 
action of operators $\hat{T},
\hat{\pi}$  and $q^{-x^2}\,\hat{Y}\,q^{x^2}$. \sq 
\end{theorem}

The general theory of spinor nonsymmetric
Whittaker functions will be published elsewhere. Let us now
consider the technique of spinors in the differential
setting.

\setcounter{equation}{0}
\section{\sc{Differential theory}}
\subsection{{\bf The degenerate case}}
\subsubsection{\sf Degenerate DAHA}
Let us begin with the definition of 
{\em degenerate double affine 
Hecke algebra} for an arbitrary (reduced) root system $R$.
Recall that $\hat{W}=W\ltimes P^\vee$ for the 
coweight lattice $P^{\vee}$.

\begin{definition}\label{dDAHA}
The {\sf degenerate double affine Hecke algebra} $\HH'$ 
is generated by $\hat{W}$ (with the corresponding 
group relations) and pairwise commutative elements 
$y_{b}$, $b\in P$ satisfying the following relations:
\begin{align}\label{ddaha}
&s_{i}y_{b}-y_{s_{i}(b)}s_{i}\,=\,-k(b, \al_{i}^\vee) 
\text{\ \,for\ \,}i\geq 1, \\
&s_{0}y_{b}-y_{s_{0}(b)}s_{i}\,=\,k(b, \theta)
\text{\ and \ } \pi_{r}y_{b}=y_{\pi_{r}(b)}\pi_{r},\notag
\end{align}
where $y_{[b,j]}=y_b+j$, $y_{b+c}=y_b+y_c$.
\end{definition} 
Note that in contrast to the definition of DAHA
from  
(\ref{dahadef}), 
$y_b$ are labeled by $b\in P$
(not by $P^\vee$). It is convenient because
$X_a$ (to be introduced later) will be naturally labeled
by $a\in P^\vee$.

Due to the additive dependence of $y_b$ of $b$, the
exact choice ($P$ or $P^\vee$) is not too important
here; one can even take $b\in \C^n$. Similarly,
changing $(b, \al_{i}^\vee)$ to $(b, \al_{i})$ will
simply re-scale the $k$\~parameters. However, the exact
choice of the lattice is important to ensure  
the compatibility of this definition with the 
limit $q\to 1$ from $q,t$\~DAHA (see below).   
The operators $X_a$ will be (translations by) $a\in P^\vee$ 
considered as elements of $\hW\subset \HH'$. The PBW Theorem
holds for $\{X_a,y_b,W\}$. 

This algebra was introduced for
the first time as the limit $q\to 1$ of 
$q,t$\~DAHA; see \cite{C101}, Chapter 2, Section ``Degenerate
DAHA." There is another approach to its definition
via the compatibility and $\hW$\~equivariance of the
{\em affine infinite Knizhnik-Zamolodchikov equation}
from \cite{C1,Ch12}. It can be called ``elliptic AKZ"
(though no elliptic functions are used in its definition)
because this system of equations at critical level
is equivalent to the eigenvalue problem
for the elliptic deformation of the Heckman-Opdam operators.
The latter is due to 
Olshanetsky -Perelomov for $A_n$\,,\,Ochiai -Oshima -Sekiguchi
for the classical root systems, and from \cite{Ch12} for 
any (reduced) root systems.
\smallskip

Let us consider the $A_{1}$\~case.
Then $\HH'$ will be generated by $s, \pi, y$ with the 
following defining relations:
\begin{align*}
s^2=1,\ sy+ys\,=\,-k,\ \pi y\,=\,(\frac{1}{2}-y)\pi.
\end{align*}
Recall that we set $s=s_1$,\, $\om=\om_1$,\, $\pi=\om s$,\,
$y=y_{\omega}$; for instance, $\pi(\omega)=[-\om,\frac{1}{2}]$.

Letting $X=\pi s$, one has that $sXs\,=\,X^{-1}$,
$(Xs)y\,=\,(\frac{1}{2}-y)(Xs)$ and finally
\begin{align*}
\ 
X(-k-ys)\,=\,(\frac{1}{2}-y)Xs \ \Rightarrow\ 
[y, X]\,=\,\frac{1}{2}X+kXs.
\end{align*}
Similar to DAHA,  
$\HH'$ can be represented as $\lan y, s, X^{\pm1}\ran$ 
subject to the relations:
\begin{align}\label{trigA1}
sXs\,=\,X^{-1},\ sy+ys\,=\,-k,\ s^{2}\,=\,1,\  
[y, X]\,=\,\frac{1}{2}X+kXs.
\end{align}
\smallskip

This algebra can be obtained as the limit (``degeneration")
of $\HH$ from 
(\ref{dahaone}) as follows.
We set  $q=\exp(h)$, $t=q^{k}=\exp(hk)$.
Let $Y=\exp(-hy)$,\, $X=X$\, and\, $T=s+\frac{hk}{2}$.
Note that now $X$ comes from the multiplication operator
(not from translations).
The letter relation is necessary to ensure that the quadratic
relation holds modulo $\,(h^2)\,$. Indeed, then 
\begin{align*}
T^{2}\ =\ 1+hks\ =\ (t^{1/2}-t^{-1/2})T+1\mod (h^2).
\end{align*} 
Check that the coefficient
of $h$ in $TY^{-1}T=Y$ readily results in the relation
$sys+ks=-y$. 

\subsubsection{\sf Polynomial representation}
Continuing with the $A_1$\~case,
$X$ and $s$ remain the same as in the $q,t$\~case, however,
now we set $X=e^{x}$. The generator $y$ is  
mapped to the differential operator
\begin{align}\label{yformula}
y=\frac{1}{2}\frac{d}{dx}+\frac{k}{1-X^{2}}(1-s)-\frac{k}{2}\,,
\end{align}
called the trigonometric Dunkl or Cherednik-Dunkl operator.
It is simple to check directly that $sys+y=-ks$\, and that 
\begin{align*}
[y,X]\,=\,\frac{1}{2}X+\frac{k}{1-X^{-2}}(Xs-X^{-1}s)\,=\,
\frac{1}{2}X+kXs.
\end{align*}

The constant $-k/2$ in formula (\ref{yformula})
automatically results from the
limiting procedure. However, its appearance here  
can be clarified without any reference to DAHA or
degenerate DAHA.
\begin{lemma}\label{deltak}
Let $\Delta_{k}\equal (e^{x}-e^{-x})^{k}$. Then
\begin{align*}
\tilde{y}\equal\Delta_{k}\,y\,\Delta_{k}^{-1}=
\frac{1}{2}\frac{d}{dx}-\frac{k}{1-X^{-2}}s.
\end{align*}
\end{lemma}
{\em Proof.}
Indeed, we have 
\begin{eqnarray*}
\Delta_{k}\,y\,\Delta_{k}^{-1}
&=&
\frac{1}{2}\frac{d}{dx}-
\frac{k}{2}\,\frac{e^{x}+e^{-x}}{e^{x}-e^{-x}}
+\frac{k}{1-X^{-2}}(1-s)-\frac{k}{2}\\
&=&
\frac{1}{2}\frac{d}{dx}+
\frac{k}{2}\,\bigl(1-\frac{2e^{x}}{e^{x}-e^{-x}}\bigr)
+\frac{k}{1-X^{-2}}(1-s)-\frac{k}{2}\\
&=&\frac{1}{2}\frac{d}{dx}-\frac{k}{1-X^{-2}}s.
\end{eqnarray*}
\sq

Thus the constant $-k/2$\, is necessary
to make the conjugation of the
trigonometric Dunkl operator by $\Delta_k$  with
pure $s$ (but then the Laurent
polynomials will not be preserved). We mention that the trigonometric
Dunkl operators were introduced in \cite{C13} 
in terms of $(c-s)$ for an arbitrary constant $c$
(including $c=0$) and in the matrix setting. We see that
the constant $c$ can be changed using conjugations
by powers of the discriminant.
\smallskip

\rmk
For complex $k$, we need to take the function
$|e^{x}-e^{-x}|^{k}$ in the lemma (to avoid problems with
complex powers). However, 
the claim of the lemma is entirely algebraic. The best
way to proceed here algebraically   
is to conjugate by the {\em even spinor}
$$\{(e^{x}-e^{-x})^{k}, (e^{x}-e^{-x})^{k}\}$$
for any branch of $(e^{x}-e^{-x})^{k}$.
It is the first appearance of spinors in this part
of the paper.
\sq

\subsubsection{\sf The self-adjointness}
Let us first establish the connection of the trigonometric 
Dunkl operator to the $k$\~deformation of Harish-
Chandra theory of the radial parts of Laplace operators
on symmetric spaces. 
One has  
$$L\,'\equal 2y^{2}|_{\mathrm{sym}}
=\frac{1}{2}\frac{d^{2}}{dx^{2}}+
k\frac{(1+e^{-2x})}{(1-e^{-2x})}\frac{d}{dx}+
\frac{k^{2}}{2}.$$
The restriction $\,|_{\mathrm{sym}}\,$
to symmetric (even) functions simply means
that we move all $s$ to the right and then delete them.

In Harish-Chandra theory, $k$ is one-half
of the {\em root multiplicity} of the restricted root system
corresponding to the symmetric space. For instance,
$k=1$ in the so-called group case. Let us mention the
contributions of Koornwinder, Calogero, Sutherland,
Heckman, Opdam and van den Ban to developing the theory
for arbitrary $k$.
See e.g, \cite{HO1} (we do not need anything beyond
the results of this paper in this section).

Lemma \ref{deltak} readily gives that
\begin{eqnarray*}
\tilde{L}\,'\equal\Delta_{k}\,L\,'\,\Delta_{k}^{-1}
\,=\,\frac{1}{2}\frac{d^{2}}{dx^{2}}+\frac{2k(1-k)}
{(e^{x}-e^{-x})^{2}}.
\end{eqnarray*}

Now let us discuss the inner product. We set formally:
$$\lan f, g\ran\,\equal\,\int f(x)g(-x)\Delta_{k}^{2}dx.
$$
For instance, the integration here can be taken over $\R$; 
then $\Delta_{k}^{2}$ must
be understood as $|e^x-e^{-x}|^{2k}$; the functions 
$f,g$ must be chosen to ensure the convergence. 

The anti-involution
${}^+$ (formally) corresponding to the ``free" inner product 
\,$\int f(x)g(-x)dx$\, acts as follows:
\begin{align*}
x^{+}\,=\,x, \,\, (\frac{d}{dx})^{+}=\frac{d}{dx}.
\end{align*}
Then the anti-involution
$A^{\Diamond}=\Delta_{k}^{-2}\,A^{+}\,\Delta_{k}^2$
serves $\lan f,g\ran$. 

\begin{lemma}
One has
\begin{align*}
X^{\Diamond}\,=\,X^{-1},\ 
y^{\Diamond}\,=\,y,\  s^{\Diamond}\,=\,s,
\end{align*}
which implies that $(L\,')^{\Diamond}=L\,'$.
\end{lemma}
{\em Proof.}
One can check the self-adjointness of $y$ and $L'$
directly. However, the best way is via
Lemma \ref{deltak} (first, for $y$ and, second, for
$L'$). Using that $\tilde{y}^{+}=\tilde{y}$, one 
obtains that 
\begin{eqnarray*}
&y^{\Diamond}\,
=\,\Delta_{k}^{-2}\,(\Delta_{k}^{-1}\,\tilde{y}\,
\Delta_{k}^{\,})^{\,+}\,\Delta_{k}^{2}
\,=\,\Delta_{k}^{-2}\,(\Delta_{k}\,\tilde{y}\,
\Delta_{k}^{-1})\,\Delta_{k}^{2}
\,=\,\Delta_{k}^{-1}\,\tilde{y}\,\Delta_{k}=y.&
\end{eqnarray*}
\smallskip

\subsubsection{\sf The Ruijsenaars limit}
The procedure is as follows. We begin with
$\tilde{L}\,'
=\frac{1}{2}\frac{d^{2}}{dx^{2}}+\frac{2k(1-k)}
{(e^{x}-e^{-x})^{2}}$,
replace $x$ by $x+M$ and connect $M$ with $k$ by
the relation $k(1-k)=e^{2M}$. Finally, we set
$\Re M\to+\infty$. Then the resulting
operator will be  
$\frac{1}{2}\frac{d^{2}}{dx^{2}}+2e^{-x},$
the {\em Toda operator}.

Applying this method to arbitrary root systems,
one obtains a system of pairwise 
commutative Toda operators.
In contrast to $L\,'$, these
operators are {\em not} $W$\~invariant. 
The (real) {\em Whittaker function}
is their eigenfunction. Given a weight (the set
of eigenvalues), the dimension of the corresponding
space of all eigenfunctions is 
$|W|$. The ``true" Whittaker function belongs to
this space and 
can be fixed uniquely there using certain decay conditions.

Let us give a reference to paper \cite{Shim},
where this procedure was applied to the Heckman-Opdam 
functions from \cite{HO1}; their limits are, 
indeed, the {\em true} Whittaker ones. 

Note that $k$ must be arbitrary in QMBP for the 
Ruijsenaars- Etingof procedure. It is impossible
to obtain the Whittaker function directly from the
classical Harish-Chandra spherical function (which 
is for very special $k$). It is somewhat different from 
$\mathfrak{p}$\~adic theory, where the passage from
the Satake-Macdonald spherical function to the
$\mathfrak{p}$\~adic Whittaker function can be established 
via switching to the maximal unramified extension from a given
$\mathfrak{p}$\~adic field.
\smallskip

\subsection{{\bf Dunkl operator and Bessel function}}
Let $X=e^{\varepsilon x}$ with $\varepsilon>0$. 
Then the trigonometric Dunkl operator $y$
becomes 
\begin{align*}
\frac{1}{2\varepsilon}
\frac{d}{dx}+\frac{k}{2\varepsilon x}(1-s)-\frac{k}{2}+
o(\varepsilon).
\end{align*}
Letting $\varepsilon\to 0$,
\begin{align*}
\varepsilon y\to \frac{1}{2}\frac{d}{dx}+
\frac{k}{2x}(1-s).
\end{align*}
We will use the same letter $y$ on the right-hand side.
However, the {\em Dunkl operator}\, will be more convenient:
$$
\mathscr{D}\equal 2y=\frac{d}{dx}+\frac{k}{x}(1-s).
$$

This definition is due to Charles Dunkl \cite{Du}, who
introduce Dunkl (rational) operators for
arbitrary root systems and also for some groups 
generated by complex reflections. 

\subsubsection{\sf Rational DAHA}
\begin{definition}
The {\sf rational double affine Hecke algebra} 
$\HH''$ is generated by $x, y, s$ with 
the following relations:
\begin{align*}
&sxs\,=\,-x, \,\, sys\,=\,-y, \,\, s^{2}\,=\,1, \,\, 
[y,x]=\frac{1}{2}+ks.
\end{align*}
\end{definition}
It is the limit of the relations 
from (\ref{trigA1}). An abstract (and very
general) variant of this definition is actually due 
to Drinfeld \cite{Dr} (though he did not consider its
polynomial representation).   

The assignment $x\to x$, 
$y\to \mathscr{D}/2$, $s\to s$ defines the 
{\em polynomial representation} of $\HH''$ in 
$\mathbb{C}[x]$. It is an induced module from
the character of the subalgebra generated by
$y,s$ sending $y$ to $y(1)=0$ and $s$ to $s(1)=1$.
The PBW Theorem is almost immediate in the rational
setting (it also follows from the existence of the
polynomial representation). 

Upon the symmetrization of $\mathscr{D}^2$, we obtain the 
key operator
in the classical theory of Bessel functions:
\begin{align*}
L\equal\mathscr{D}^{2}|_{\mathrm{sym}}
=\frac{d^{2}}{dx^{2}}+\frac{2k}{x}\frac{d}{dx}.
\end{align*}

\begin{lemma}\label{LEMDIAM}
(i) One has 
$$
x^{k}\cdot\mathscr{D}\cdot x^{-k}=
\tilde{\mathscr{D}}\equal\frac{d}{dx}-\frac{k}{x}s,\ \  
x^{k}\cdot L\cdot x^{-k}=\tilde{L}\equal\frac{d^{2}}
{dx^{2}}+\frac{k(1-k)}{x^{2}}.
$$

(ii) Let  
$A^{\Diamond}=x^{-2k}\cdot A^{*}\cdot x^{2k}$, where the
anti-involution $*$ is as follows:
$$x^{*}\,=\,x, \quad (\frac{d}{dx})^{*}\,=\,-\frac{d}{dx};
$$ 
the anti-involution $\,\Diamond\,$ formally serves the 
bilinear symmetric form 
$\lan f, g\ran$ $=\int f(x)g(x)x^{2k}dx\,.$
One has that\, $\mathscr{D}^{\Diamond}=-\mathscr{D}$\,, and\, 
$L^{\Diamond}=L$\,.\sq
\end{lemma}

\subsubsection{\sf Bessel functions}
Assuming that $\la\neq 0$, an arbitrary solution
$\varphi_{\lambda}^{(k)}$  
of the eigenvalue problem
\begin{align}\label{Lvarphi}
L \varphi_{\lambda}^{(k)}=4\lambda^{2}\varphi_{\lambda}^{(k)}
\end{align}
analytic in a neighborhood of $x=0$ can be represented as
$$\varphi_{\lambda}^{(k)}(x)=\varphi^{(k)}(x\lambda).
$$
Here $\varphi^{(k)}$ can be readily calculated: 
\begin{eqnarray}\label{varphiform}
\varphi^{(k)}(t)=\sum_{m=0}^{\infty}
\frac{t^{2m}\Gamma(k+1/2)}{m!\Gamma(k+n+1/2)}
\end{eqnarray}
for the Gamma\~function, satisfying
$\Gamma(x+1)\,=\,x\Gamma(x)$, $\Gamma(1)=1$.
The parameter $k$ is
arbitrary here provided that $k\neq -1/2-m$ for $m\in \Z_+$.
The function $\varphi^{(k)}(t)$ is a variant of the
Bessel $J$\~function. 

See \cite{O3} (and references
therein) for the theory of multi-dimensional Bessel functions.

Notice that 
\begin{eqnarray*}
\varphi^{(k)}(t)\xrightarrow{k\to 0}
\sum_{m=0}^{\infty}\frac{(2t)^{2m}}{(2m)!}
=\frac{e^{2t}+e^{-2t}}{2},
\end{eqnarray*}
due to the relations:
\begin{align*}
\Gamma(n+1)\Gamma(n+\frac{1}{2})\,=\,2^{-2n}(2n)!\sqrt{\pi}, 
\ \, \Gamma(\frac{1}{2})\,=\,\sqrt{\pi}. 
\end{align*}

Using the passage to the Sturm-Louiville operator
$\tilde{L}$, we can control the growth of
$\varphi_{\lambda}^{(k)}$ at infinity.
   
\begin{lemma}
The differential equation 
$L \varphi =4\lambda^{2}\varphi$ has the following 
two fundamental solutions for real $x$.
If $\lambda=0$, then $1$ and $x^{1-2k}$ can be taken.
If $\lambda\neq 0$, the asymptotic behavior can be used
to fix them:
\begin{align*}
\varphi_{\lambda}^{\pm}=
x^{-k}e^{\pm 2\lambda x}(1+o(1))
\hbox{\ as\ } x\to +\infty.
\end{align*}
Any solution $\varphi$ is a linear combination of these two.
In particular, the growth of any solution as $x \to\pm\infty$
is no greater than exponential, namely, 
$O(\,x^{-\Re k} e^{\pm 2x\Re \lambda})$ 
for $\la\neq 0$. \sq
\end{lemma}
We will use this lemma only for justifying that 
the Gauss-Bessel integrals we will need below
are well defined.
The following is the classical formula; see
Introduction and Chapter 1 from \cite{C101}
for a more comprehensive exposition.

\subsubsection{\sf Hankel transform}
\begin{theorem}{\label{thm:master}}
\begin{align*}
\int^{+\infty}_{-\infty}\,
\varphi_{\lambda}^{(k)}(x)\,
\varphi_{\mu}^{(k)}(x)\,e^{-x^{2}}|x|^{2k}
dx=\Gamma(k+\frac{1}{2})\,\varphi_{\mu}^{(k)}
(\lambda)\,e^{\lambda^{2}+\mu^{2}},
\end{align*}
where $\Re k>-\frac{1}{2}$. The normalization is given by the
Euler integral:
\begin{align*}
\int^{+\infty}_{-\infty}e^{-x^{2}}|x|^{2k}
dx=\Gamma(k+\frac{1}{2}).
\end{align*}
Here one can set 
$\int^{+\infty}_{-\infty}=2\int^{+\infty}_{0}$, since all
functions are even. \sq
\end{theorem}

In order to prove Theorem \ref{thm:master}, we need  
the following definition.
\begin{definition}
The {\sf Hankel transform} for even functions 
$f$ is given by
\begin{align}\label{mastersym}
\mathbb{H}f(\lambda)=\frac{1}{\Gamma(k+\frac{1}{2})}
\int_{\mathbb{R}}f(x)\varphi_{\lambda}^{(k)}(x)|x|^{2k}dx
\end{align}
in proper functional spaces.
\end{definition}

\subsubsection{\sf Its properties}
Let us denote the operator $L$ acting in the $\la$\~space by
$L_\la$; $L$ without the suffix $\la$ will continue to be the 
operator above in terms $x$. Recall that the operator $L$ depends
on $k$; we will sometimes denote it by $L^{(k)}$. 

\begin{lemma} For any functional spaces (not only for 
even functions), provided $L$ and $\mathbb{H}$
are well defined there,
\begin{enumerate}
\item [(a)] $\mathbb{H}(L)=4\lambda^{2}$, 
$\mathbb{H}(4x^{2})=L_\la$;
\item [(b)] $e^{-x^{2}}\,L\,e^{x^{2}}=L+4x^{2}+[L,x^{2}]$.
\end{enumerate}
\end{lemma}
{\em Proof.}
Claim $(a)$ is based on the
$x\!\leftrightarrow\! \la$\~symmetry of $\varphi_{\la}^{(k)}(x)$ 
and on the self-adjointness of the operators $L$ and $x^2$
with respect to the measure we consider.

Checking $(b)$ is direct. One can also use the following
important connection with the theory of $\mathfrak{sl}(2)$.
Setting
$$
e=x^{2},\ f=-\frac{L}{4},\ h=[e,f]=x\frac{d}{dx}+\frac{1}{2}+k,
$$
we obtain a representation of this Lie algebra. Then
$e^{-x^{2}}Le^{x^{2}}$ can be interpreted and
calculated using the adjoint action of $SL_2$. It must
be {\em a priori} a linear combination of $e,f,h$; the exact 
formula is simple.
Note that the Hankel transformation becomes the group
element $s\in SL_2$ in this interpretation.
\sq

{\em Proof of theorem \ref{thm:master}}.
Let 
$\hat{\varphi}_{\mu}^{(k)}(\lambda)
\equal 
e^{-\lambda^{2}}\mathbb{H}(\varphi^{(k)}_{\mu}(x)e^{-x^{2}})$.
Due to the lemma, $\hat{\varphi}_{\mu}^{(k)}(\la)$ satisfies 
$L^{(k)}_{\lambda}\hat{\varphi}_{\mu}^{(k)}=
4\mu^{2}\hat{\varphi}_{\mu}^{(k)}$. However, this
solution is unique up to proportionality in the class 
of even analytic functions in a neighborhood
of $x=0$. Thus 
$\hat{\varphi}_{\mu}^{(k)}(\lambda)=C_\mu 
\varphi_{\mu}^{(k)}(\lambda)$. It gives 
(\ref{mastersym}) up to proportionality. 
Using the 
$\la\!\leftrightarrow\!\mu$\~symmetry on the left-hand
side of this formula and the same symmetry of
$\varphi_{\mu}^{(k)}(\lambda)$, we obtain that
$C_\mu=Ce^{\mu^2}$ for an absolute
constant $C$, which can be readily determined. 
\sq

\subsubsection{\sf Tilde-Bessel functions}
Let us try to apply the master formula to other
solutions of the eigenvalue problem
(\ref{Lvarphi}). We will manipulate
algebraically for some time, without exact 
analytic justifications.
The proof above looks very algebraic; we
even did not use that $\varphi_{\la}^{(k)}(x)$ is
even.  

For $\lambda\neq 0$, there exists another solution 
$\tilde{\varphi}_{\lambda}^{(k)}(x)=
(x\lambda)^{1-2k}\varphi_{\lambda}^{(1-k)}(x)$
of (\ref{Lvarphi}).
If $\lambda=0$, let 
$\tilde{\varphi}_{\lambda}^{(k)}(x)\equal x^{1-2k}$.
We need to assume that $\Re (k)<1/2$ to avoid
the singularity at $0$ in these solutions.

Applying the reasoning above (formally), we obtain
that 
\begin{align}\label{genmaster}
\mathbb{H}(\tilde{\varphi}_{\mu}^{k}e^{-x^{2}})
=\breve{\varphi}_{\mu}^{(k)}(\lambda)
e^{\lambda^{2}+\mu^{2}}
\end{align}
for a certain solution $\breve{\varphi}_{\mu}^{(k)}$
of the same eigenvalue problem, a linear combination
of $\varphi_{\mu}^{(k)}$ and
$\tilde{\varphi}_{\mu}^{(k)}$.

If we assume here that $0<\Re(k)<1/2$ and
set $\mu=0$, then $\tilde{\varphi}_{\mu}^{(k)}(0)=0$.
Upon obvious cancelations, we come to the following
brand new identity in the theory of Bessel functions:
$$
\int_{-\infty}^{+\infty}
\varphi_{\lambda}^{(k)}(x)|x|e^{-x^{2}}dx=e^{\lambda^{2}}.
$$
{\em Unfortunately this formula is  wrong}. Let us explain why.
\smallskip

Informally this is wrong simply because
no {\em new} identities of such a kind can be expected in the
very classical field of Bessel functions
and Hankel transform. The exact mathematical reason for this
failure is as follows. The integration by parts, necessary for the
self-adjointness claim, requires the convergence at $0$ of
the {\em first two derivatives} of the functions involved. 
The existence of the starting and the final integral can be
insufficient; one need to justify the convergence of all
intermediate integrals as well. 

The following  analytic constraints make claim $(ii)$ of 
Lemma \ref{LEMDIAM} rigorous. These conditions
are not exactly sharp, but sufficient for us.

Provided that $f,g\in C^{2}(\mathbb{R}_{+})$ and 
$f(x)|x|^{k}$, $g(x)|x|^{k}$ are absolutely integrable,
\begin{align*}
\int_{-\infty}^{+\infty}L(f)g|x|^{2k}dx=
\int_{-\infty}^{+\infty}fL(g)|x|^{2k}dx.
\end{align*}

\subsubsection{\sf Complex analytic theory}
The deduction above of (\ref{genmaster}) from the
properties of the Hankel transform is of course
formally correct; this simply gives nothing new in the
case of real integration due to the divergence at $0$
of the derivatives of the tilde-solution. 
The Laplace integration, was design exactly to 
avoid the divergences of this kind.
Let us first re-establish the usual master formula in 
the Laplace setting.
\begin{theorem}
For all $k\in \mathbb{C}$ such that 
$k\neq -\frac{1}{2}-m$, $m\in \mathbb{Z}_{+}$,
\begin{align*}
\int_{i\varepsilon +\mathbb{R}}
\varphi_{\lambda}^{(k})(x)
\varphi_{\mu}^{(k)}(x)e^{-x^{2}}(-x^{2})^{k}dx
=\frac{\pi}{\Gamma(\frac{1}{2}-k)}
\varphi_{\lambda}^{(k)}(\mu)e^{\lambda^{2}+\mu^{2}}.
\end{align*}
Here $\varepsilon>0$; the 
condition $k\neq -\frac{1}{2}-m$ is necessary for
the existence of $\varphi_{\lambda}^{(k)}(x)$.\sq
\end{theorem}

For any complex number $k$, the function
$(-x^{2})^{k}$ is defined as the function
$\exp(k\log(-x^2))$
continued along the integration path $x\in i\varepsilon+\R$ for 
the usual branch of $\log$ with the cutoff at
$\R_-$. Using $(-x^{2})^{k}$ is quite standard in classical works
on $\Gamma$ and related functions. 

Due to the Gamma-term on the right-hand side,
this integral must be zero
at $k=\frac{1}{2}+m$, $m\in \mathbb{Z}_{+}$. It is simple
to demonstrate  directly.
Indeed, 
$$(-x^{2})^{1/2}=-ix \hbox{\ \, along the path\,\ }
i\varepsilon+\R\,;
$$ 
check the point $x=i\vep$ using that 
$(\vep^2)^{1/2}=\vep$. The integrand is analytic
at zero for such $k$, so we can tend $\varepsilon\to 0$.
However the integrand is an odd function on $\R$ and,
therefore, 
\begin{align*}
\int_{i\varepsilon +\mathbb{R}}
\varphi_{\lambda}^{(k)}(x)
\varphi_{\mu}^{(k)}(x)e^{-x^{2}}(-ix)^{2m+1}dx=0.
\end{align*}

Similarly, for $\tilde{\varphi}_{\lambda}(x)
\equal(-\lambda^{2})^{1/2-k}(-x^{2})^{1/2-k}
\varphi_{\lambda}^{(1-k)}(x)$, which is the complex analytic
variant of the tilde-solution considered above,
\begin{eqnarray*}
&&\int_{i\varepsilon +\mathbb{R}}
\varphi_{\lambda}^{(k)}(x)
\tilde{\varphi}^{(k)}_{\mu}(x)(-x^{2})^{(k)}
e^{-x^{2}}dx\\
&=&\int_{i\varepsilon +\mathbb{R}}\varphi_{\lambda}^{(k)}(x)
\varphi^{(1-k)}_{\mu}(x)(-x^{2})^{1/2}dx\\
&=&\int_{\mathbb{R}}\varphi_{\lambda}^{(k)}(x)
\varphi^{(1-k)}_{\mu}(x)(-ix)dx\ =\ 0.
\end{eqnarray*}
Thus the standard solution $\varphi_{\lambda}^{(k)}(x)$
and the complex-analytic tilde-solution are orthogonal 
to each other in the master formula.
 
It is straightforward to calculate the master formula for the
tilde-solutions $\tilde{\varphi}^{(k)}_{\la}(x),
\tilde{\varphi}^{(k)}_{\mu}(x)$
coupled together in the Gauss-Bessel integral. We will provide 
the corresponding formulas below when doing the nonsymmetric master
formula.

\setcounter{equation}{0}
\section{\sc{Spinor eigenfunctions}}
We will begin with the eigenvalue problem for the
Dunkl operator. The latter is not a differential operator,
but it shares some (but not all) properties with
the first order {\em differential} operators.
\begin{lemma}\label{DUNEIG}
(i) The eigenvalue problem
\begin{align}\label{Duneigen}
\mathscr{D}\psi=2\lambda\psi, \text{ for } 
\mathscr{D}=\frac{d}{dx}+\frac{k}{x}(1-s)
\end{align}
has a unique analytic at $0$ solution $\psi=\psi_\lambda^{(k)}(x)$
satisfying $\psi(0)=1$  if
and only if $k\not \in -1/2-\Z_+.$   

(ii) Namely, it is $\psi=1$ for $\la=0$ and
$
\psi(x)=\psi^{(k)}(\lambda x)
$
for 
$$
\psi^{(k)}(t)=
\varphi^{(k)}(t)+\frac{1}{2}(\varphi^{(k)})'(t)
$$
in terms of $\varphi^{(k)}(t)$ 
from (\ref{varphiform}). 

(iii) When $\lambda=0$ and  $k=-\frac{1}{2}-m$, 
the space of analytic solutions is generated by 
$\psi=1$ and $\psi=x^{2m+1}$. When 
$\la\neq 0$ for the same $k$,
the analytic solution $\psi$ exists and
is unique up to proportionality, but vanishes at \, $0$. 
\sq
\end{lemma}
The fact that the dimension of the space of solutions
of (\ref{Duneigen}) can be $2$ (for special
values of the parameters) requires attention
and will eventually lead us to the spinor extension
of the space of functions.

\subsection{{\bf Nonsymmetric master formula}}
For $k\neq -1/2-m$, $m\in \Z_{+}$ and the function
$\psi_{\lambda}^{(k)}(x)=\psi^{(k)}(\lambda x)$ from
Lemma \ref{DUNEIG}, the following holds.
\begin{theorem}\label{NONSYMPSI}
(i) For $\Re k>-1/2$,
\begin{align*}
\int_{\mathbb{R}}
\psi_{\lambda}^{(k)}(x)\psi_{\mu}^{(k)}(x)e^{-x^{2}}|x|^{2k}dx
=\Gamma(k+\frac{1}{2})
\psi_{\lambda}(\mu)^{(k)}e^{\lambda^{2}+\mu^{2}}.
\end{align*}

(ii) Denote 
$\int_{\R}^{\varepsilon}\equal 
\frac{1}{2}(\int_{i\varepsilon +\mathbb{R}}+
\int_{-i\varepsilon +\mathbb{R}})$, then
\begin{align*}
\int^{\varepsilon}_{\mathbb{R}}
\psi_{\lambda}^{(k)}(x)\psi_{\mu}^{(k)}(x)e^{-x^{2}}
(-x^{2})^{k}dx
=\frac{\pi}{\Gamma(\frac{1}{2}-k)}
\psi_{\lambda}^{(k)}(\mu)e^{\lambda^{2}+\mu^{2}}.
\end{align*} 
\end{theorem}
{\em Proof.}
As in the symmetric theory, the formula readily results from
the basic facts concerning the {\em nonsymmetric Hankel
transform}. The (general) definition of this transform
is due to Dunkl \cite{Du2}.
Its one-dimensional version can be found in  
Hermite's works, but this was used only marginally 
in the classical theory. This transform is given by
\begin{align}\label{masternon}
\mathbb{H}_{ns}f(\lambda)=\frac{1}{\Gamma(k+\frac{1}{2})}
\int_{\mathbb{R}}f(x)\psi_{\lambda}^{(k)}(x)|x|^{2k}dx,
\end{align}
provided the existence. Its theory is actually simpler
than that of the classical symmetric Hankel transform
(at least the algebraic aspects).
We use the notation $\mathscr{D}_\la$ 
for the Dunkl operator acting in the $\la$\~space.

The following analytic conditions for the functions $f,g$
and their derivatives $f',g'$
are sufficient to ensure that
\begin{align}\label{intRD}
&\int_{\R}\,\mathscr{D}(f) g |x|^{2k}dx\ =\ 
-\int_{\R}\,f\mathscr{D}(g) |x|^{2k}dx:
\end{align}

(1) $f(x), g(x)$ are continuous and 
$f'(x), g'(x)$ exist in $\mathbb{R}\setminus 0$;

(2) the function $f(x)g(x)|x|^{2k}$ is integrable and 
continuous at $0$;

(3) $f(x)g(x)|x|^{2k-1}$, $f'(x)g(x)|x|^{2k}$, 
$f(x)g'(x)|x|^{2k}$, $f(x)g(-x)|x|^{2k}$ 

\ \ \ \ \ are integrable at zero.

For the integration $\int_{\mathbb{R}}^{\varepsilon},$
only the integrability at infinity is needed
for (\ref{intRD}). The theorem readily follows from
the following lemma.

\begin{lemma} 
For $f$ as above and provided the existence
of $\mathbb{H}_{ns}$,
\begin{enumerate}
\item [(a)] $\mathbb{H}_{ns}(\mathscr{D})=2\lambda$, 
$\mathbb{H}_{ns}(2x)=\mathscr{D}_\la$;
\item [(b)] $e^{-x^{2}}\,\mathscr{D}\,e^{x^{2}}=\mathscr{D}+2x$\,,
\end{enumerate}
where the integration in (\ref{masternon})
can be either $\int_{\mathbb{R}}$ or
$\int_{\mathbb{R}}^{\varepsilon}.$\sq
\end{lemma}

\rmk
Similar to the symmetric case, the integrals from
Theorem \ref{NONSYMPSI} in the complex case are
identically zero as $k\in 1/2+\Z_+$. It corresponds
to the vanishing condition of the inner products
associated with level-one coinvariants from 
Theorem \ref{HALLONE}. 
See also formula 
(\ref{diaform}) 
(the real case $(b)$ there). 

The affine symmetrizer $\hat{\mathscr{I}}\,$ from
(\ref{Jproj})\, 
is a $q,t$\~Jackson counterpart of the integration 
$\int_{i\vep +\R}f(x) (-x^2)^k dx$.
The zeros of the inner product
$\hat{\mathscr{I}}\,(f\, T(g))$ for $A_1$ are 
exactly in the set $1/2+\Z_+$. 

\subsubsection{\sf Using spinors}
The theory of the nonsymmetric tilde-solutions
requires the technique of spinors (already used above).
They are pairs $f=\{f_{1}, f_{2}\}$ of functions defined 
in an open set $U$
in $\R$ or $\C$.
{\em Real spinor} are defined for 
$U=\{x\in \R,\, x>0\}$; {\em complex spinors}
are defined for the set $U=\{x\in \C,\, \Im x>0\}$. The operators
act naturally on spinors; see Section \ref{sect:spinors}.
For instance, 
$$
s\{f_{1}, f_{2}\}=\{f_{2}, f_{1}\},\ x\{f_{1}, f_{2}\}=
\{xf_{1}, -xf_{2}\},\ 
\{f_{1}, f_{2}\}'=\{f_{1}', -f_{2}'\},
$$
where here and below
$f'\equal df/dx$.

The {\em super-presentation} of a spinor $f$ is 
defined to be
\begin{align*}
f=\llbracket f^{0}, f^{1}\rrbracket, \text{ where }
f^{0}=\frac{f_{1}(x)+f_{2}(x)}{2}, \,
f^{1}=\frac{f_{1}(x)-f_{2}(x)}{2}.
\end{align*}

For any two spinors, $f=\{f_{1}, f_{2}\}$, $g=\{g_{1}, g_{2}\}$, 
their product 
is given by $f\cdot g=\{f_{1}g_{1}, f_{2}g_{2}\}$. 
In the super-presentation:
\begin{align*}
f\cdot g=\llbracket f^{0}g^{0}+f^{1}g^{1}, 
f^{0}g^{1}+f^{1}g^{0}\rrbracket.
\end{align*}
It is the standard stuff about $\Z_2$\~graded algebras.

A spinor $f=\{f_{1}, f_{2}\}$  is called a 
{\em principal spinor (function)} 
if the following holds.
There must exist an open {\em connected\,} set $\tilde{U}$ 
and a function $\tilde{f}$ 
on $\tilde{U}$  such that
$U$, $U^{s}\equal s(U)\subset \tilde{U}$ and 
$f_{1}=\tilde{f}|_{U}$, $f_{2}=s(\tilde{f})|_{U}$.
 
The differentiation of spinors $\frac{d}{dx}$ 
is an odd operator defined by
$$
\frac{d}{dx}\llbracket f^{0}, f^{1}\rrbracket
=\llbracket \frac{d}{dx}f^{1}, \frac{d}{dx}f^{0}\rrbracket.
$$ 
The spinor integration is given by
$$ 
 \int_\ga \llbracket f^{0}, \ 
f^{1}\rrbracket\equal \int_{\gamma}f^{0},
$$ 
where $\gamma\subset U$ is a path in the set $U$.

\subsubsection{\sf Spinor Bessel functions}
The Dunkl spinor eigenvalue problem is
\begin{align}\label{dunspineig}
\mathscr{D}(\psi)=
\llbracket(\psi^{1})'+\frac{2k\psi^{1}}{x},(\psi^{0})'\rrbracket
=\llbracket 2\lambda \psi^{0}, 2\lambda\psi^{1}\rrbracket.
\end{align}
In the standard representation $\{\psi_1,\psi_2\}$, 
it reads as follows:
\begin{align*}
\mathscr{D}(\psi)=
\{\,\psi_1'+\frac{k(\psi_1-\psi_2)}{x},
-\psi_2'-\frac{k(\psi_2-\psi_1)}{x}\,\}=
\{\,2\lambda \psi_1, 2\lambda\psi_2\,\}.
\end{align*}

\begin{lemma}\label{DUNSOLUT}
The space of solutions of the eigenvalue problem
(\ref{dunspineig}) is always two-dimensional. 
There are three cases:
\begin{enumerate}
\item if $\lambda\neq 0$, then all the solutions are in the form
$\psi=\llbracket \varphi, \frac{\varphi'}{2\lambda}\rrbracket$
for $\varphi$ satisfying $L\varphi=4\lambda^{2}\varphi$, and  
only one of them (up to proportionality)
is a function (i.e., a principle spinor);
\item if $\lambda=0$ and $k\not\in -1/2-\Z_+$ 
then $\psi=1$ is a solution and also there is an odd 
spinor solution  $\chi_{k}$, given by 
$\chi_{k}=\llbracket 0, |x|^{-2k}\rrbracket$ in the real case and 
$\chi_{k}=\llbracket 0, (-x^{2})^{-k}\rrbracket$ in the complex case;
\item when $\lambda=0$ and $k= -1/2-m$ for $m\in \Z_+$, then
the solutions are $1$ and $x^{2m+1}$, i.e.,
both are principle spinors (functions).
 
\sq
\end{enumerate}
\end{lemma}

{\em Nonsymmetric tilde-solutions.}
For $k\notin 1/2+\mathbb{Z}_{+}$, the spinor 
\begin{align*}
\tilde{\psi}_{\lambda}^{(k)}=\chi_{k}(x)\chi_{k}(\lambda)
\psi_{\lambda}^{(-k)}(x)
\end{align*}
satisfies (\ref{dunspineig}). Actually it is a {\em bi-spinor},
in terms of $x$ and $\la$; we will skip the formal definition.

Let us incorporate the tilde-solution into
the master formula. We need to redefine the inner
product. Let 
\begin{align*}
x^{2k}\equal \left\{\begin{array}{ccc}
\llbracket\ |x|^{2k},\,0\ \rrbracket, &  
& \hbox{\ real\ case;} \\
\llbracket\, (-x^{2})^{k},\,0\ \rrbracket, &  
& \hbox{\ complex\ case.}\end{array}\right.
\end{align*}
I.e., both are even spinors (functions, if $k\in \Z$).
Note that $\chi_{k}(x)x^{2k}=\llbracket 0,1\rrbracket$ is an 
odd constant (a spinor of course).
The integration will be
\begin{eqnarray*}
\int f(x) & \equal & 
2\int_{0}^{+\infty}f^{0}(x)dx 
\text{\ \,in the real case; }\\ 
\int 
f(x) & \equal & \int_{i\varepsilon+\mathbb{R}}f^{0}(x)dx 
\text{\ \, in the complex case. }
\end{eqnarray*}

Let us check that the $\psi$\~solution and 
the $\tilde{\psi}$\~solution are orthogonal to each
other in the master formula. 
Similar to the symmetric case, we have the divergence
problem with the integration by parts, so only
the complex case will be considered.
Then the integral
\begin{align}{\label{eqn:complexint}}
\int \psi^{(k)}_{\lambda}(x)
\tilde{\psi}^{(k)}_{\mu}(x)e^{-x^{2}}x^{2k}
 \end{align}
is proportional to
\begin{align*}
I=\int_{i\varepsilon+\R}e^{-x^{2}}
(\psi_{\lambda}^{(k)}\psi_{\mu}^{(-k)}
\cdot\llbracket 0,1\rrbracket)^{0}\,dx=
\int_{i\varepsilon+\R}e^{-x^{2}}
(\psi_{\lambda}^{(k)}\psi_{\mu}^{(-k)})^1\,dx.
\end{align*} 
However, $e^{-x^2}\psi_{\lambda}^{(k)}(x)\psi_{\mu}^{(-k)}(x)$ 
is a principal spinor, i.e., a restriction
of an analytic function $F$. Therefore the component
$F^1$ is an odd function on $\R$. Letting 
$\vep\to 0$ in the integration 
path, we conclude that $I=0$.

\subsubsection{\sf Tilde master formulas}
Let us list explicitly the Gauss-Bessel integrals for the
tilde-solutions.
\begin{theorem} In the real case, 
\begin{align*}
2\int^{+\infty}_{0}
(\tilde{\psi}_{\lambda}^{(k)}
\tilde{\psi}_{\mu}^{(k)})^{0}e^{-x^{2}}|x|^{2k}dx
=\tilde{\psi}_{\lambda}^{(k)}(\mu)
e^{\lambda^{2}+\mu^{2}}\Gamma(\frac{1}{2}-k) 
\text{\ \,for\ \, }\Re k<\frac{1}{2}.
\end{align*}
In the complex case,
\begin{align*}
\int_{i\varepsilon+\mathbb{R}}
(\tilde{\psi}_{\lambda}^{(k)}\tilde{\psi}_{\mu}^{(k)})^{0}
e^{-x^{2}}(-x^{2})^{k}dx
=\frac{\pi}{\Gamma(\frac{1}{2}+k)}\tilde{\psi}_{\lambda}^{(k)}
(\mu)e^{\lambda^{2}+\mu^{2}}
\text{ as }k\notin \frac{1}{2}+\mathbb{Z}_{+};
\end{align*}
this integral is zero when $k=-1/2-m$ for $m\in \mathbb{Z}_{+}$.
\sq
\end{theorem}
We note that the spinors we integrate and those in the
right-hand side are actually {\em bi-spinors}, i.e., 
spinors in terms of $x$ and spinors in terms of $\la,\mu$. 
the formal definitions are straightforward. It suffices here
to use directly the definition:
$\tilde{\psi}_{\lambda}^{(k)}(\mu)=
\chi_{k}(\lambda)\chi_{k}(\mu)
\psi_{\lambda}^{(-k)}(\mu)$.

Let us also provide the symmetric tilde-formulas
(no spinors are needed):
\begin{align*}
2\int^{+\infty}_{0}
\tilde{\varphi}_{\lambda}^{(k)}
\tilde{\varphi}_{\mu}^{(k)}\,e^{-x^{2}}|x|^{2k}\,dx\
=\ \Gamma(\frac{3}{2}-k)\,\tilde{\varphi}_{\mu}^{(k)}
(\lambda)\,e^{\lambda^{2}+\mu^{2}},\ 
\Re k<\frac{3}{2}\,,
\end{align*}
\begin{align*}
\int_{i\varepsilon+\mathbb{R}}\tilde{\varphi}_{\lambda}^{(k)}
\tilde{\varphi}_{\mu}^{(k)}\,e^{-x^{2}}(-x^{2})^{k}\,dx
=\frac{\pi}{\Gamma(-\frac{1}{2}+k)}\,
\tilde{\varphi}_{\mu}^{(k)}(\lambda)\,
e^{\lambda^{2}+\mu^{2}},\ k\notin \frac{3}{2}+\mathbb{Z}_{+},
\end{align*}
and the latter 
integral is zero at $k=1/2-m$ for $m\in \mathbb{Z}_{+}$.

An obvious problem is in extending the nonsymmetric master formula
to all spinor solutions for arbitrary root systems. One cannot
expect the formulas to be so simple as for $A_1$, because the
Weyl groups $W$ have irreducible representations of higher 
dimensions. We do not have the general formulas at the moment. 
Similar questions can be posted for arbitrary, not necessarily 
symmetric, solutions of the $L$\~eigenvalue problems in arbitrary
ranks, when no spinors are needed.

We mention that the orthogonality relations 
for $\psi$ coupled with $\tilde{\psi}$  can be extended to
the trigonometric- differential and 
trigonometric- difference settings (any root systems), 
provided we have the $Y$\~semisimplicity.
Hopefully this can be sufficient to manage the rational case.

\subsection{{\bf Affine KZ equations}}\label{sect:AKZ}
\subsubsection{\sf Degenerate AHA and AKZ}
Let $R$ be an arbitrary (reduced) root system, $R^\vee$ its dual,
$P$ and $P^\vee$ the corresponding weight and coweight lattices. 
We set $z_{a}=(z,a)$ for $z\in \C^n$ and define the differentiation
$\partial_{b}z_{a}\equal (b,a)$ for arbitrary vectors $a,b$
(to be used mainly for $b\in P$, $a\in P^\vee$).
Let ${}^wf(z)=f(w^{-1}(z))$ for $w\in W$,
$s_{\al}$ be the reflections corresponding to the roots 
$\al$ and $\{y_{b}\}$  pairwise commutative elements
satisfying $y_{a+b}=y_a+y_b$ for $a,b\in P$.  

We will follow 
Definition \ref{dDAHA} of degenerate DAHA restricted 
to the AHA case, i.e., consider only
nonaffine reflections $s_i$; also $-k$ will be replaced by $k$. 
The relations of {\em degenerate AHA},
due to Drinfeld for $GL_n$ \cite{Dr}
 and Lusztig \cite{L}, are
\begin{align}\label{ddahanon}
&s_{i}y_{b}-y_{s_{i}(b)}s_{i}\,=\,k(b, \al_{i}^\vee), 
\text{ for }i\geq 1. 
\end{align}
The corresponding algebra will be denoted by $\h'$.

Let $\Phi\,$ be a function of $z\,$ taking its values in the abstract
algebraic span $\lan s_{\al}, y_{b}\ran$. 
The {\em affine Knizhnik-Zamolodchikov equation}, AKZ, 
is the following system of differential equations
\begin{align}{\label{eqn:AKZ}}
\partial_{b}(\Phi)
=\left(\sum_{\al\in R_{+}^\vee}\frac{k(b,\al)
s_{\al}}{e^{z_{\al}}-1}+y_{b}\right)\Phi,
\text{ where }b\in P.
\end{align}
Actually, $b$ can be arbitrary complex vectors here
and below.

\begin{theorem}
The AKZ is self-consistent and $W$\~equivariant if and only 
if the elements $s_{\al}\,$ and $y_{b}\,$ 
satisfy the relations from (\ref{ddahanon}).
The equivariance here means that if $\Phi\,$ is a solution
of AKZ, then so is $\,w({}^w\Phi(z))=w(\Phi(w^{-1}(z)))\,$.\sq
\end{theorem}

The definition of AKZ and this theorem were the starting
point of DAHA theory;
here and below see Chapter 1 of \cite{C101}.
The following construction is basically from \cite{C13},
but using the technique of spinors consistently makes 
it entirely algebraic (and essentially coinciding with
that from \cite{O2}).

In \cite{C13} and other first author's
papers, the values of AKZ were considered in
$\h'$\~modules induced from arbitrary finite-dimensional
representations of $W$ or induced from the characters of  
the polynomial algebra 
$\C[y]=\C[y_b,\, b\in P]$. In this paper we will stick
to the modules induced from $\C[y]$.

\subsubsection{\sf Spinor Dunkl operators}
The Dunkl operators will be needed here
in the following form:
\begin{align}{\label{eqn:Duntrig}}
\mathscr{D}_b^0
=\partial_b-\sum_{\al\in R_{+}^\vee}\frac{k(b,\al)
\si_{\al}}{e^{z_{\al}}-1},\where   
\si_\al(z_a)=z_{s_\al(a)}.
\end{align}
Here $\si$ stays for the action on the argument
of functions: 
$\si_u(f)(z)=f(u^{-1}z)$,\, $u\in W$.
The relation to AKZ is established via
the {\em spinor Dunkl operators\,} defined as 
a natural extension of (\ref{eqn:Duntrig}) to
the space of $W$\~spinors.

The {\em spinors} are collections
$\hat{\psi}=\{\psi_w,\, w\in W\}$ of (arbitrary) 
scalar functions with component-wise addition, 
multiplication and the differentiations by $\partial_b$. 
The action $\si_u$ for $u\in W$
is through permutations of the indices:
$$
\si_u(\hat{\psi})=\{\psi_{u^{-1}w},w\in W\}.
$$
Note the sign of $u^{-1}$, which ensures that we really have a 
representation of $W$; the spinors are actually {\em functions}
on $W\times \C^n$ so $u^{-1}$ (the dualization)
is necessary. This definition
matches the action of $W$ on functions $f$ of $z$, which 
will be considered as {\em principle spinors} under the embedding
$$
f\mapsto f^\rho\equal\{f_w={}^{w^{-1}}f,\,w\in W\}.
$$
Indeed, we have the
commutativity \,
$(\si_u(f))^\rho=\si_u(f^\rho)$. The definition of $\rho$
can be naturally extended to the operators acting on functions.

For instance, the function $z_\al$ becomes the spinor 
$\{z_{w^{-1}(\al)},w\in W\}$ under this embedding; also,  
$(\partial_b)^\rho=\{\partial_{w^{-1}(b)},w\in W\}$. 

\begin{theorem}\label{SPINDUN}
For a solution $\Phi$ of the AKZ with values
in $\h'$, let us define the
spinor $\hat{\Psi}=\{w(\Phi),w\in W\}$ for the 
action of $w\in W$ in $\h'$ by left multiplications.
Then  $\hat{\Psi}$ satisfies the following {\sf spinor
Dunkl eigenvalue problem}:
\begin{align}\label{Dspineigen}
\mathscr{D}_b^0(\hat{\Psi})=y_b\hat{\Psi},\ b\in P.
\end{align}
\end{theorem} 
{\it Proof.} The $W$\~equivariance of AKZ readily 
establishes the equivalence of this theorem with
the previous one. Explicitly, 
$\si_\al(\hat{\Psi})=\{s_\al w(\Phi),w\in W\}$ and
the relations
for the component $w=u$ of $\hat{\Psi}$
read as follows:
$$
\partial_{u^{-1}(b)}\,u(\Phi)=
\sum_{\al\in R_{+}^\vee}\frac{k(b,\al)
s_{\al}\,u(\Phi)}{\exp(z_{u^{-1}(\al)})-1}+y_b\, u(\Phi),
\ b\in P.
$$
This can be recalculated to the same AKZ system for $\Phi$
due to the $W$\~equivariance. \sq

\rmk
In \cite{C13}, an analytic variant of this construction
was used. The algebraic formalization of the 
argument from \cite{C13} can be found in Lemma 3.2 
from \cite{O2}; the proof above is very similar to that 
in \cite{O2}. This ``algebraization" can be readily extended
to the difference and elliptic theories (considered in
\cite{C101} and previous first author's works). From
the viewpoint of the applications to the isomorphism
theorems, both approaches are equivalent. 

As far as the reduction of AKZ to the Dunkl eigenvalue 
problem is concerned, arbitrary modules of $\h'$ were
considered (not only induced) in \cite{C13}. The Dunkl
operators there were given in terms of the action of
$W$ via the monodromy of AKZ (see below). 
Treating formally the corresponding $W$\~orbits as
{\em spinors}, one makes the construction entirely algebraic  
(as in Theorem \ref{SPINDUN} and in \cite{O2}). 

It is important that the monodromy can be calculated
{\em explicitly} for the asymptotically free solutions of AKZ.
For instance, these explicit formulas were
used in Theorem 4.3 from \cite{C13} to solve the
{\em real\,} (nonspinor)  Dunkl eigenvalue problem 
via AKZ in {\em functions} (not only in {\em spinors}).
The solution found in \cite{C13} using the monodromy
approach is the $G$\~function that was introduced (later)
and played the key role in paper \cite{O2}. 
\sq

\subsubsection{\sf The isomorphism theorem}
Let us apply Theorem \ref{SPINDUN} to {\em induced
representations}.
Given a one-dimensional representation $\C_{\lambda}=\C v_\circ$
of $\C[y]\,$  
defined by $y_{b}(v_\circ)=\la_b v_\circ$
for $\la_b=(\lambda, b)$, where $\la\in \C^n$,\, let 
$
\,I_{\lambda}=\Ind^{\h'}_{\C[y]}\C_\la\,
$ 
be the $\h'$\~module induced from $\C_\la$. 

We note that if  
the space of eigenvectors (pure, not generalized)
for the eigenvalue $\la$ is one-dimensional in $I_\la$, 
then there exists a rational expression in terms of 
$y_b$ serving as a projector of $I_\la$ 
onto $\C v_\circ\subset I_\la$.

Let $I_\la^*$ be $Hom(I_\la,\C)$ supplied with
the natural action of $\h'$ via the 
canonical {\em anti-involution}
of $\h'$ preserving the generators $s_i, y_b$
(reversing the order in products).
We use here that the relations in the degenerate
affine Hecke algebra are self-dual.

Next, we define the linear functional 
\,$\varpi: f\mapsto f(v_\circ)$\, on \,$I_\la^*\ni f$\,
satisfying the conditions\,
$
\varpi((y_b-\la_b)I_\la^*)=0 \for b\in P.\,
$
Assuming, that the space of $\la$\~eigenvectors 
in $I_\la$ is one-dimensional, these conditions
determine $\varpi$ uniquely up to proportionality.

The functional $\varpi$ is nonzero on
any nonzero $\h'$\~submodule $V^*\subset I_\la^*$, since 
$I_\la$ is cyclic generated by $v_\circ$.
Indeed, if $\varpi(f)=0$ for all $f\in V^*$, then
$f(\h'v_\circ)=0=f(I_\la)$ for all such $f$.  

Let $U_{0}\subset \mathbb{C}^{n}$ be a open neighborhood of 
$0$ in $\C^n$; we set $U_{0}'=\cap_{{}_{w\in W}}\,w(U_{0})$. 
We assume that $U_{0}$ satisfies the following properties
(necessary for the monodromy interpretation below): 
\begin{enumerate}
\item $U_{0}$ does not contain any zeros 
of $\prod_{\al\in R_{+}^\vee}(e^{z_{\al}}-1)$;
\item $U_{0}$ is simply connected and $U_{0}'/W$ is connected;
\end{enumerate}
$U_0^\star$ will be one of the connected components of $U_0'$
(the latter set is a disjoint union of $|W|$ connected open sets). 

By $\,Sol^{\lambda}_{AKZ}(U_{0})\,$, we denote the space of 
$I_{\lambda}^*$\~valued analytic solutions $\phi$
of the AKZ equation in $U_0$. 

Let $\,Sol_{\mathscr{D}}^{\lambda}(U^\star_{0})\,$
be the space of $W$\~{\em spinor solutions} $\hat{\psi}$
in $U^\star_{0}$ of the {\em scalar} eigenvalue problem
\begin{align}\label{Dspineig}
\mathscr{D}_b^0(\hat{\psi})=\la_b\hat{\psi},\ b\in P.
\end{align}
The spinors here are collections
$\hat{\psi}=\{\psi_w,\, w\in W\}$ of (arbitrary) 
scalar analytic functions in $U^\star_{0}$.

\begin{theorem}\label{ISOTHMFULL}
The dimension of the space 
$Sol_{\mathscr{D}}^{\lambda}(U^\star_{0})$ equals the
cardinality $|W|$ of $W$. There is an isomorphism
\begin{align}\label{solakzD}
\eta:\, Sol^{\lambda}_{AKZ}(U_{0})\ni\phi\,\mapsto\, 
\{\varpi(w(\phi))\downarrow_{U_0^\star},\, w\in W\}\,\in\,
Sol_{\mathscr{D}}^{\lambda}(U^\star_{0})
\end{align}
for the action of $w\in W$ 
on the values of $\phi$, which are from $I_\la$.
\end{theorem}
{\em Proof}.
The claim that $\eta$ is a map between the required spaces of
solutions follows from Theorem \ref{SPINDUN}.
Due to the coincidence of the dimensions of the spaces
in (\ref{solakzD}), we need only to check that $\eta$
is injective. As in \cite{C13}, this follows from the fact
that $\varpi$ is nonzero on any
$\h'$\~submodule of $I_{\la}^*$. Note that the 
construction of $\eta$ is entirely 
algebraic, so it suffices to assume that $\phi$
is defined in the same open set $U_0^\star$ as in the statement
of the theorem.

\subsubsection{\sf The monodromy interpretation}
Let $\Phi(z)$ be an invertible matrix solution
of AKZ in $U_0$ with values in $\hbox{Aut}(I_\la^*)$.
For any $w\in W$, let us define the {\em monodromy matrix} 
$\mathcal{T}_{w}$ by
\begin{align*}
w(\Phi(z))=\Phi(w(z))\mathcal{T}_{w}.
\end{align*}
Here $\Phi(w(z))$ is well defined in
$U_0\cap w^{-1}(U_0)$, so is $\mathcal{T}_{w}$. 
The matrix solution $\Phi$ is nothing but a choice of
the basis of fundamental solutions in $Sol^{\lambda}_{AKZ}(U_{0})$
(its columns). Changing the basis conjugates
all $\mathcal{T}_w$ by a constant invertible matrix. 
The matrix-valued functions $\mathcal{T}_{w}$ have the 
following properties:
\begin{enumerate}
\item[(a)] $\mathcal{T}_{w}$ are defined in $U_{0}'$ and are locally
constant;
\item[(b)] $\mathcal{T}_{uw}={}^{w^{-1}}\mathcal{T}_{u}
\mathcal{T}_{w}
=\mathcal{T}_{u}(w(z))\mathcal{T}_{w}(z)$ for $u,w\in W$.
\end{enumerate}
For each $w\in W$, let us
define its {\em $\sigma'$-action}\,: 
$$\sigma'_{w}(F)={}^{w}F\mathcal{T}_{w^{-1}}=
F(w^{-1}(z))\mathcal{T}_{w^{-1}}(z).
$$
Then $\sigma'_{1}=1$, $\sigma'_{uw}=\sigma'_{u}\sigma'_{w}$
and $\sigma'_{w}\partial_{a}=\partial_{w(a)}\sigma'_{w}$
for $u,w\in W$ and $a\in P$. We naturally set 
$\si'_{\al}=\si'_{s_\al}$ and $\si'_i=\si'_{\al_i}$.
Here $F$ can be an arbitrary function
in $U_0'$ with values in $\hbox{Aut}(I_\la^*)$.

Introducing
$$
\mathscr{D}_{b}'\equal \partial_{b}-
k\sum_{\al\in R_{+}^\vee}\frac{(\al,b)\sigma'_{\al}}
{e^{z_{\al}}-1},
$$
one readily obtains that
\begin{align}\label{abphi}
y_{b}\Phi=\left(\partial_{b}-k
\sum_{\al\in R_{+}^\vee}\frac{(\al,b)\sigma'_{\al}}
{e^{z_{\al}}-1}\right)\Phi=\mathscr{D}_{b}'\Phi.
\end{align}
We simply employ the definition of $\sigma'$ here.
The action of $\mathscr{D}_{b}'$ is given in terms
of the $W$\~action on $z$ and the {\em right} multiplications
by matrices $\t_{s_\al}$. So this action commutes with $y_b$, which
are {\em left} multiplications
by constant matrices. Therefore we can apply the
functional $\varpi$ to $\Phi$ in (\ref{abphi}), which
gives that (\ref{abphi})  holds for 
$\varpi(\Phi)$. The spinor $\hat{\Psi}$ 
from Theorem \ref{SPINDUN} is nothing but 
$\{\Psi_w=\si'_{w^{-1}}(\Phi)\downarrow U_0^\star, w\in W\}$.

\comment{
(with the natural action of $W$ on the 
indices $w$ of $\Psi_w$). 
Then we claim that 
\begin{align}\label{spinpsi}
\hat{\psi}=\{\psi_w,\,w\in W\}\equal
\varpi(\hat{\Psi})\,=\,\{\varpi(\Psi_w)\}
\end{align}
solves the eigenvalue (scalar) problem (\ref{Dspineig}) in the
space of spinors for the connected component $U_0^\star$.
If $\phi$ is the first column of $\Phi$,
then $\psi_1=\psi_{id}$ is $\varpi(\phi)$ upon the restriction
to $U_0^\star$. Thus the map claimed in the theorem is  
$$
\phi\mapsto \hat{\psi},\ \psi_1=\varpi(\phi)\downarrow_{U_0^\star}.
$$
Following \cite{C13} and using that $\varpi$ is nonzero on any
$\h'$\~submodule of $I_{\la}^*$, we prove the theorem.\sq
}

\subsubsection{\sf Connection to QMBP}
Continuing this construction, one can combine the isomorphism we 
found with the symmetrization map, which acts 
from $Sol_{\mathscr{D}}^{\lambda}(U_{0}^\star)$ to
the space of solutions  of the Heckman-Opdam system (QMBP)
in $U_0^\star$ corresponding to $\la$. 
To be exact,
the map from $Sol^{\lambda}_{AKZ}(U_{0})$ to
$Sol^{\lambda}_{QMBP}(U_{0})$
is the projection of the space of values 
onto the one-dimensional subspace of
$W$\~invariants inside $I_{\la}^*$. It gives the
Matsuo- Cherednik isomorphism theorem from \cite{Mats,C13}
(the proof follows \cite{C13}). 
The spinors do not appear in the construction
of this map and the statement of the
theorem; 
however, they provide {\em the} best way to verify it
(and dramatically reduce the proof from \cite{Mats}). 

We note that the relation of the Dunkl-spinor eigenvalue 
problem above to QMBP is actually very similar to Lemma 
\ref{DUNSOLUT}, which addresses solving the Dunkl eigenvalue
problem in {\em spinors}. Let us mention Corollary 3.4
from \cite{O2}, where a similar extension of the Dunkl
eigenvalue problem was considered.

Certain conditions on the module $I_\la$ are
necessary to ensure the {\em isomorphism} with QMBP. 
Namely, this module must be assumed {\em spherical\,},
$\h'$\~generated by $\sum_{w\in W} w(v_\circ)$,
correspondingly,
$I_\la^*$ will be {\em co-spherical\,}. See \cite{C13}
and \cite{C101}. 
\smallskip

\rmk
There are relations to the localization functor from
\cite{GGOR,VV1}. The later is, very briefly, taking the
monodromy representation of the local systems analogous
to AKZ (in more general modules). Starting with certain
rational or degenerate DAHA modules, the monodromy results in the 
representations of nonaffine (affine) $t$\~Hecke algebras. 

The monodromy is important in our approach too 
(the cocycle $\{\t_w\}$ does contain $t$). The actual
output of our approach is a complete system of eigenfunctions
of Dunkl operators in the corresponding $y$\~eigenspaces 
of the initial $\h'$\~module, the $G$\~function in 
the terminology from \cite{O2}.
Algebraically, the Dunkl operators and the operators
of multiplication by the (trigonometric) coordinates
generate the corresponding DAHA module.

The localization functor
is understood completely (so far) only in the
rational case and in the differential -trigonometric
case (corresponding to the setting of this section);
see \cite{GGOR,VV1}. Our construction and the isomorphism
theorems hold for all known families
of AKZ and Dunkl operators (including the elliptic 
theories). See \cite{C13}, \cite{Ch12}, Chapter 1
from \cite{C101} and \cite{Sto2}. The exact connection is still not
clarified. 

\medskip
  
\setcounter{equation}{0}
\section{\sc{Conclusion}}
To try to connect better the topics of this work
and to put it into perspective, we will touch upon 
the relations of DAHA, mainly the $q$\~Whittaker functions,
to the geometric quantum Langlands program, though not much is known 
in this direction. The relation of the Verlinde algebras to
the Lusztig category of the representations of quantum groups
from \cite{L1} is of key importance here; this is 
the main focus of this section. 

We will not try to review the applications (known and expected)
of the {\em ``symmetric" global $q$\~Whittaker functions}, 
including the Shintani -Casselman -Shalika formula, 
the relations to Givental-Lee theory and possible applications
in physics. See \cite{ChW} and \cite{GLO} 
for a discussion. Generally, the (coefficients
of) $q$\~Whittaker functions are expected to contain
a lot of information about quantum $K$\~theory and
$IC$\~theory of affine flag varieties. 
Givental-Lee
theory deals with quantum $K$\~theory of the flag
variety. 

We are very thankful to 
Roman Bezrukavnikov, Alexander Braverman, 
Dennis Gaitsgory, Michael Finkelberg, 
David Kazhdan, Victor Ostrik for various
discussions on quantum groups, affine Grassmannians,
quantum Langlands program and neighboring topics (though 
they do not always agree with what will follow). 

\subsection{{\bf Verlinde algebras and QG}}
The relations to DAHA are expected upon applying
$K_0$ (the Grothendieck group)
to the categories used in the quantum geometric
Langlands program and related directions.
Then these categories become
commutative rings with inner products 
and sometimes with a projective action of $PSL_2(\Z)$. 
Generally, the number of simple objects must be finite 
for the latter action to exist. 

As it was pointed out 
in Section ``Abstract Verlinde Algebras" from \cite{C101},
such rings (even if some of these structures
are missing) are very exceptional. For instance,
one can formally prove counterparts of the Macdonald
conjectures 
(the norm formulas and the evaluation-duality formulas) 
in the abstract Verlinde-type setting, establish Pieri rules 
and do more; cf. \cite{C103}. 

It is unlikely that there 
are many commutative rings with such rich structures. The major 
candidates are quotients of the 
polynomial and various similar representations of DAHA, 
including infinite-dimensional ones and the corresponding
(commutative) algebras of the $W$\~invariants.  

\subsubsection{\sf Quantum groups}
The expected connections to the Langlands program and related
projects are grouped around the following.

\newtheorem*{keyconjecture}{Conjecture}

\begin{keyconjecture} \label{CONJCH}
The commutative algebra $K_0(Rep_q\, G)$  
for the category $Rep_q\, G\,$ of finite-dimensional representations
of Lusztig's quantum group can be canonically
identified with the algebra $\mathscr{X}^W$ of $W$\~invariants 
of the polynomial representation $\mathscr{X}$ of DAHA at $t=q$, 
defined for the corresponding root system. 
It includes the roots of unity $q$.
Then sub-quotients of the $\mathscr{X}^W$ under the action 
of the subalgebra of invariants of $\HH$ (the elements 
commuting with $T_w$ for $w \in W$) correspond to
categorical sub-quotients of $Rep_q\, G\,$. Such 
sub-quotient of  $Rep_q\, G\,$
has the structure of modular category if $PSL(2,\Z)$ acts
projectively in the corresponding sub-quotient of $\mathscr{X}^W$. 
\end{keyconjecture}
\smallskip
 
For generic $q\,$, the simple
objects correspond to the classical finite-dimensional
characters, which are eigenfunctions of the $W$\~invariant 
$Y$\~operators.
The most interesting here is the case of roots of unity, 
when $\mathscr{X}$ and $\mathscr{X}^W$ become reducible.
 
For $q=e^{2\pi \imath /N}$, the algebra of
$W$\~invariants of the nonzero (canonical) irreducible
quotient of $\mathscr{X}$ can be naturally identified with
the {\em Verlinde algebra} in the special case 
$k=1\,(t=q)$; see \cite{C101}, Section 0.4. 
The projective DAHA-action of the $PSL(2,\Z)$  
leads to the Verlinde $T,S$\~operators. 

The Verlinde algebra was originally defined
in terms of integrable Kac-Moody modules with
the {\em fusion} directly related to the
{\em conformal field theory}.  
Equivalently, it is isomorphic to   
the quotient of $K_0(Rep_q\, G)$ by the modules of zero
$q$\~dimension, i.e., $K_0$ of the so-called 
{\em reduced category\,}. The equivalence of these
two approaches at roots of unity is due to Finkelberg
\cite{Fi} (\cite{KL2} apart from the 
roots of unity). It confirms the conjecture
for the {\em perfect\,} quotients of $\mathscr{X}$.

The categorical sub-quotients in the conjecture
generally cannot be expected to be tensor categories 
for Lusztig's big quantum group unless in some special cases,
including the reduced category. The first author is
grateful to Michael Finkelberg and Victor Ostrik for
clarifying discussions on these matters. 

The next case after the reduced category 
(actually the key) is the so-called {\em parallelogram
quotient} of $Rep_q\, G$. It is the category
of representations of the 
{\em small quantum group\,} \cite{AG}, which attracts a lot of
attention now. We expect that its
$K_0$ corresponds to  
the algebra of $W$\~invariants of the
{\em DAHA parallelogram module} under the same
relation $t=q$. The latter is defined for $A_1$ as 
$$
V^{-2}=\C[X,X^{-1}]/(X^{2N}+X^{-2N}-2)= 
\C[X,X^{-1}]/(X^{N}-X^{-N})^2
$$ 
in the notation from \cite{C101} Section 2.9.3;\,
its dimension is $4N$. 
Let us discuss the rank-one case in greater detail.

\subsubsection{\sf The rank-one case}
The {\em perfect} quotient of $V^{-2}$
for $q=e^{2\pi \imath/N}$ and integral $0\le k<N/2$ 
will be denoted by $V_{2N-4k}$; its dimension is $2N-4k$. 
Here one can consider half-integral $k$ too
(we will not discuss it). Let $V_{2N+4k}\,$ be the kernel
of the natural map $V^{-2}\to V_{2N-4k}$.

Both are irreducible
DAHA modules with the projective $PSL(2,\Z)$-action.
They are commutative algebras because so is $\mathscr{X}$;
$V_{2N-4k}$ is semisimple, but $V_{2N+4k}$ for $k>0$ is not.
The action of $X$ in the latter has
$\,4k\,$ Jordan $2$\~blocks ($2$\~dimensional blocks) 
with pairwise distinct eigenvalues and 
$\,2N-4k\,$ simple eigenvectors. Due to the
projective $PSL(2,\Z)$\~action (we need $\si$), 
the Jordan decomposition must be of the same 
type for $Y$ instead of $X$. 

The Jordan decomposition of $Y$ in the whole $V^{-2}$ 
is different from that of $X$.  Namely, $Y$ has 
$\,4k\,$ Jordan $2$\~blocks and the rest of it is semisimple 
(all eigenvalues are of multiplicity $2$). The 
decomposition of $X$ in $V^{-2}$ obviously consists of the
$2$\~blocks only (see its definition); their number is $2N$.
Hence, there can be no projective action of
$PSL(2,\Z)$ in $V^{-2}$ extending that in $V_{2N\pm 4k}$.

Upon taking the $W$\~invariants,\, $\dim (V^{-2})^W=2N\,,$
$$
\dim V_{2N-4k}^W\, =\, N-2k+1,\ \and
\dim V_{2N+4k}^W\, =\, N+2k-1.
$$
The latter two algebras are projective $PSL(2,\Z)$\~invariant
because the generator $T$ is fixed under this action.

Let us discuss the case $k=1$ in more detail. One has
$$
V_{2N-4}=
\C[X,X^{-1}]/(F) \for
F=\frac{X^{2N}-1}{(X^2-1)(X^2-q)}.
$$
For instance, $F=(X-q)(X+q)$ for (the minimal possible) $N=3$
and the Verlinde algebra is $\C[Z]/(Z^2-1)$ for $Z=X+X^{-1}$
in this case; $\,q=\exp(2\pi \imath/3)$.

Importantly, $Y+Y^{-1}$ acts semisimply in the invariants
of the polynomial representation for $k=1$.
It is due to the fact that the 
$(Y+Y^{-1})$\~eigenvectors in $\C[X+X^{-1}]$ do not depend 
on $q$ when $t=q$  and are proportional to 
the $SL(2)$\~Schur functions (it holds for any root systems).  
Accordingly, $(V^{-2})^W$ and $V_{2N+4}^W$ are 
$(Y+Y^{-1})$\~semisimple in this case. The spectrum
of $Y+Y^{-1}$ in $(V^{-2})^W$ is
$\{q^{i/2}+q^{-i/2}, 1\le i\le 2N\}$ for $q^{1/2}=
e^{\pi\imath/N}$; thus, $\,2,-2\,$ are simple eigenvalues and
the others are of multiplicity $2$.

The operator $X+X^{-1}$ in $(V^{-2})^W$ is {\,\em not\,}
semisimple even for $k=1$. Namely, $2,-2$ are its simple
eigenvalues, but $q^{i/2}+q^{-i/2}$ correspond to the
Jordan $2$\~blocks for  $1\le i <N$.
Since this is different from the Jordan decomposition
of $Y+Y^{-1}$ in this space, we conclude that there can be 
no projective action of $PSL(2,\Z)$ in $(V^{-2})^W$ for 
$k=1$ ($N\ge 3$).
If the conjecture above holds, then no such an action can be
expected in the parallelogram quotient of $Rep_q\, G$ 
at roots of unity extending that in the Verlinde
algebra; so it cannot be a {\em modular category}.

It is likely
that the irreducible constituents of the parallelogram
DAHA modules for integral $k$ are always projective 
$PSL(2,\Z)$\~modules, but this is known only for $A_1$;
it may be connected with \cite{Lyu}.
The parallelogram module, as the whole, has no natural
(projective) $PSL(2,\Z)$\~structure (only $\tau_-$ acts there).
\smallskip

As a related direction,
we would like to mention that Tipunin and others 
successfully calculated 
certain generalized Verlinde algebras of nonsemisimple 
type using the {\em logarithmic conformal theory}; 
see e.g., \cite{MT}. 
They obtained exactly the ones described in \cite{C101},
Proposition 2.9.6 (upon taking the $W$\~invariants). 
Technically, the (canonical)
irreducible quotient of the polynomial representation
becomes nonsemisimple for integral $N>k>N/2$ ; 
it can be identified with $V_{2N+4k}$ considered above 
for $0<k<N/2$.  

Let us also note that the limit of the minimal models as 
$c\to 1$ is important in physics applications;
the corresponding infinite-dimensional
Verlinde-type algebra is likely to be the polynomial 
DAHA representation itself. 
\smallskip

\subsection{{\bf Expected developments}}
\subsubsection{\sf Approaching the conjecture}
The most conceptual reason for the conjecture above
is a very close relation of DAHA (almost at the level of
its definition)
to $K$\~theory of affine flag varieties.
However, there are other aspects too.
\smallskip
 
{\em KZ equations.}
The affine Knizhnik-Zamolodchikov equations
and the so-called $r$\~matrix KZ (see \cite{C101},
Section 1.5) can be employed here. These KZ are 
directly connected with the {\em coinvariants} and the {\em
$\tau$\~function} for
{\em factorizable Kac-Moody algebras} associated
with $r$\~matrices (introduced in the first author's works).
Generally, the approach based on the $KZ$ equation is of key
importance in \cite{KL2}, \cite{Fi} and in \cite{Ga}, 
so this technique is certainly relevant for the conjecture.
   
\smallskip
{\em Nonsymmetric theory.}
DAHA gives the most in the nonsymmetric setting,
when we switch from the $W$\~invariant polynomials to the
whole polynomial representation. However, 
we do not know much about the geometric meaning of the
{\em nonsymmetric Macdonald polynomials}. 
There are two major general facts here. They 
are connected with the Matsumoto 
spherical functions and with the level-one 
Demazure characters; these examples are degenerate but
nevertheless important. 

Generally, 
taking the $W$\~invariants in DAHA-modules seems really
necessary to relate them to Lie-Kac-Moody theory.  
The technique of spinors, which establishes
a connection of DAHA to non-$W$\~invariant sections of 
local systems like QMBP (the Heckman-Opdam system), 
could be a bridge from the nonsymmetric theory to geometry. 
\smallskip

{\em Finite-dimensional modules.}
It is worth mentioning that the specialization $t=q$ 
used in this conjecture does not seem 
the only one related to $Rep_q\, G$. Let us
restrict ourselves to the spherical case,
which means that we will consider only the
quotients of the polynomial representation
$\mathscr{X}$. Then such modules will be 
commutative algebras and the corresponding
categories, if any,  can be expected monoidal.

Important generalizations of Verlinde
algebras can be obtained when the polynomial
representation $\mathscr{X}$ and its nonzero
irreducible quotient are considered for
the following DAHA parameters: 
 
(a)\, $t=q^k$ for {\em \, singular\,} rational
$k=-\frac{s}{d}<0\,$ and any unimodular $q$,\, 

(b)\, $t\in \C$ but $q$ is a root of unity (a variant
of the parallelogram case),  

(c)\, and when $q$ is a root of unity under the 
limits $\,t\to 0$ or $\,t\to \infty\,$,\\
although not all structures are present in these three cases. 

Only the integrality of the
structural constants of the Verlinde algebra 
will be missing in $(a)$ (since $q$ is not a root of unity);
the positivity of the Verlinde inner product will hold
for sufficiently small $\arg(q)$. 

More significantly, there will be no projective 
action of $PSL(2,\Z)$ in the cases $(b,c)$.
 The limits from $(c)$ (which are actually particular cases
of $(b)$) are very
interesting because of possible (no exact confirmations so far) 
relations to the following.
  
\subsubsection{\sf Toward Langlands program}
The (local) quantum geometric Langlands program will be
discussed here very introductory. Let $G$ be the  
simply connected Lie group over $\C$
corresponding to a given root system $R$,\, ${}^L G$ its Langlands
dual (though we mainly stick to the simply-laced $R$ in this work). 

The global ``symmetric" $q$\~Whittaker function 
can be interpreted as the Fourier transform of
$K_0(Rep_q\, {}^LG)$ for generic $q$;
we actually need $|q|<1$ here to ensure the convergence.  
The challenge is to connect it with the category
$\hbox{Whit}^c$ (see below) and the Gaitsgory-Lurie transform
$$
K_0(Rep_q\, {}^LG)\ \to\ \hbox{Whit}^c(Gr_G),
$$
a complicated functor between 
the corresponding $2$-categories. Such connection
seems almost inevitable if this transform  has something to do 
with the $q$\~Toda operators, which is exactly the key question.

The images of the simple objects of $K_0(Rep_q {}^LG)$
in  $K_0(\hbox{Whit}^c(Gr_G))$ under the Gaitsgory-Lurie
transform are of major importance; 
for many applications, knowing them is quite sufficient.
The problem is that this map cannot be fixed uniquely
at the level of $K_0$ without using involved
categorical (or other?) methods.
\smallskip

Assuming that ${}^L\!\mathscr{X}^W$, where ${}^L\!\mathscr{X}$
is the Langlands dual of $\mathscr{X}$, 
is a substitute for  $K_0(Rep_q {}^LG)$ 
(the conjecture),
its limits $t\to 0$ or $t\to\infty$ could be equally relevant
here for generic $q$ (they are connected with each other). 
Then the Fourier transform of 
$\, \lim_{t\to 0}\,{}^L\!\mathscr{X}^W$ could be, hopefully, 
a DAHA counterpart of $K_0$ of the category 
$$
\hbox{Whit}^c(Gr_G)\ =\ D\,\hbox{mod}\,^c ( G((z))/G_0)^{N((z))},
$$ 
where $N\subset G$ is the standard unipotent subgroup and 
an unramified character on $N((z))$ is needed here
to define the equivariant modules. 

Without going into detail,
it is a category of $N((z))$\~equivariant 
$c$\~twisted $D$\~modules on the affine flag 
variety $Gr_G$, which is the group $G((z))$ of 
(formal) meromorphic loops divided by 
the group $G_0=G[[z]]$ of holomorphic ones. 
The category $\hbox{Whit}^c(Gr_G)$ was proven by Gaitsgory 
in \cite{Ga} to be equivalent (for generic $q$ and 
under some technical restrictions) to $Rep_q\, {}^LG$ for
$q=e^{\pi c}$, which was conjectured by Lurie.
\smallskip

{\em $Q$\~Toda system as Hitchin system.}
The Fourier image of $\lim_{t\to 0}\,{}^L\!\mathscr{X}$ twisted by 
the Gaussian is the 
spinor polynomial representation of nil-DAHA from 
Theorems \ref{TWONIL}, called there the {\em hat-representation} 
(see also Theorem \ref{MainToda}). 
So the $W$\~invariant part of the hat-representation may be a 
candidate for $K_0(\hbox{Whit}^c(Gr_G)$.

A certain indirect confirmation is the relation of 
$\hbox{Whit}^c(Gr_G)$ to the $W$\~algebras 
and their Verlinde algebras, which, in their turn, are
connected with the DAHA-Verlinde algebras.

If one replaces the {\em Hitchin system} in the
geometric Langlands duality by the  
$q$\~Toda eigenvalue problem, then the ``symmetric"
(non-spinor) 
$q$\~Whittaker function will become the reproducing
kernel of the corresponding Fourier transform. For any
fixed set of eigenvalues, the corresponding
$q$\~Toda eigenvalue problem can be interpreted as a
$D$\~module very similar to those in the 
category $\hbox{Whit}^c(Gr_G)$ (upon the switch
from quantum groups to Kac-Moody theory). The
exact relation of this approach to the quantum Langlands
program is not established so far.
\smallskip


\subsubsection{\sf Affine flag varieties etc.}
Another source of inspiration could be Theorem 3 from
\cite{BF}, which may be more directly connected
with $q$\~Whittaker functions than the Gaitsgory-Lurie
transform. In its $K$\~theoretic variant (a conjecture),
it looks related to the Fourier duality we establish
between nil-DAHA from Theorem \ref{TWONIL}. 

If such a connection really exists, then it could
result in the $K$\~theoretical interpretation of 
the {\em spinor $q$\~Whittaker function} 
from (\ref{spinwhit}). It provides the duality
between the spinor hat-representation and the
bar-representation. The latter has a clear $K$\~theoretic
meaning; thus the former can be of geometric nature too.

Also, we expect the modular {\em translation functor\,}
and the so-called {\em wall-crossing\,} to be related to the
DAHA {\em intertwiners\,} and, more specifically,
to the analytic continuation of the asymptotic 
expansions of the global $q$\~functions from one asymptotic 
sector to another.

The mod $p\,$ methods were already used for DAHA; this is
a powerful tool. The wall-crossing is expected to be connected
with the theory of nil-DAHA; its relation to global functions
is not based on any solid evidence at the moment.

Let us outline a possible approach to
geometric theory of global functions based on their
asymptotic expansions. The definition of these functions
and the existence of their limits at infinity are from  
\cite{C5}, \cite{ChW}; let us also mention Stokman's
definition of the global functions for $C^\vee C$ and his 
recent results on the difference Harish-Chandra theory.
\smallskip

{\em Global functions geometrically.}
A complete description of the asym\~ptotic expansions of a
global function, namely, inside the asymptotic sectors, 
then at their walls, then at the walls of walls and so on, 
called the {\em  resonance conditions\,}, would fix
it uniquely as an analytic function without any reference 
to the Macdonald or $q$\~Toda operators. 

Generally, the continuation of the functions/sections from their
natural domains to the boundary requires involved tools (like 
intersection cohomology). {\em Global functions\,} are automatically
such continuations of their asymptotic expansions, so they are 
expected
to be canonical in every possible sense. In their definition,
we use that the polynomial representation multiplied
by the Gaussian is self-dual with respect to the
DAHA-Fourier transforms; the global functions are the corresponding 
reproducing kernels. It provides a conceptual explanation of their 
remarkable algebraic and analytic properties. 

The resonance theory of global 
$q,t$\~spherical and $q$\~Whittaker functions, a continuation
of the program due to Harish-Chandra, Casselman \cite{Ca} and
others, is in progress. The first development here was
the Harish-Chandra theory of asymptotic decomposition  (the first
author and Stokman), including the representation of a global
function as a weighted $W$\~summation of its asymptotic expansions.
\smallskip

{\em Associators and dilogarithm.\,} We note that DAHA can be applied 
to catch 
certain categorical structures beyond $K_0$.  Generally, changing
the asymptotic sectors of the Knizhnik-Zamolodchikov equations
gives the {\em associators\,} due to Drinfeld. In DAHA theory
we restrict ourselves to the (various) AKZ equations. 
In several examples, these associators correspond to different 
choices of maximal commutative subalgebras in AHA or DAHA
and can be calculated.

The resulting $q,t$\~pentagon-type relations in the 
limit $t\to 1$ may be connected with \cite{FG}. It is
certainly connected with the theory of asymptotic decomposition
of the global functions outlined above.
It is worth mentioning that the well-known pentagon relation 
for the quantum {\em dilogarithm}, which is nothing but
the $q$\~Gamma function, does not play any significant role in DAHA 
theory so far, though there are recent developments in this
direction in the theory of nil-DAHA. 
Adding dilogarithms to DAHA would be an important 
development.


\renewcommand\refname{\sc{References}}
\bibliographystyle{unsrt}

\end{document}

%% file: macro.tex
\renewcommand{\tilde}{\widetilde}
\renewcommand{\hat}{\widehat}

\newcommand{\Z}{{\mathbb Z}}
\newcommand{\Q}{{\mathbb Q}}
\newcommand{\N}{{\mathbb N}}
\newcommand{\C}{{\mathbb C}}
\newcommand{\R}{{\mathbb R}}

\def\HH{\mbox{${\mathcal H}$\kern-5.2pt${\mathcal H}$}}

\newtheorem{theorem}{Theorem}[section]

\newtheorem{definition}[theorem]{Definition}
\newtheorem{lemma}[theorem]{Lemma}

\newtheorem{theorem }{Theorem}[section]
\newtheorem{maintheorem }[theorem]{Main Theorem}
\newtheorem{proposition }[theorem]{Proposition}
\newtheorem{definition }[theorem]{Definition}
\newtheorem{lemma }[theorem]{Lemma}
\newtheorem{corollary }[theorem]{Corollary}
\newtheorem{notation }[theorem]{Notation}
\newtheorem{remark }[theorem]{Remark}
\newtheorem{example }[theorem]{Example}

\newtheorem{ maintheorem }[theorem]{Main Theorem}
\newtheorem{ theorem}{Theorem}[section]
\newtheorem{ proposition}[theorem]{Proposition}
\newtheorem{ definition}[theorem]{Definition}
\newtheorem{ lemma}[theorem]{Lemma}
\newtheorem{ corollary}[theorem]{Corollary}
\newtheorem{ notation}[theorem]{Notation}
\newtheorem{ remark}[theorem]{Remark}
\newtheorem{ example}[theorem]{Example}

 \newcommand{\rmk}{{\bf Comment.\ }}

\def\for{\  \hbox{ for } \ }

\def\where{\  \hbox{ where } \ }
\def\and{\  \hbox{ and } \ }
\def\and{\  \hbox{ and } \ }

\def\equal{\stackrel{\,\mathbf{def}}{= \kern-3pt =}}

\def\la{\lambda}
\def\La{\Lambda}
\def\om{\omega}
\def\Om{\Omega}

\def\th{\theta}
\def\al{\alpha}
\def\be{\beta}
\def\ga{\gamma}
\def\ep{\epsilon}

\def\de{\delta}
\def\De{\Delta}

\def\si{\sigma}
\def\Si{\Sigma}
\def\Ga{\Gamma}

\def\vph{\varphi}

\def\vep{\varepsilon}

\def\vth{{\vartheta}}

\def\tal{\tilde{\alpha}}
\def\tbe{\tilde{\beta}}

\def\tw{\widetilde w}
\def\tW{\widetilde W}

\def\tz{\tilde z}

\def\tR{\tilde R}

\def\hw{\widehat{w}}
\def\hW{\widehat{W}}
\def\hu{\hat{u}}

\def\hv{\hat{v}}

\def\I{\mathbf{I}}

\def\S{\mathbf{S}}

\def\0{\mathbf{0}}
\def\H{\mathbf{H}}

\def\f{\mathcal{F}}
\def\çF{\mathcal{F}}
\def\o{\mathcal{O}}
\def\t{\mathcal{T}}

\def\l{\mathcal{L}}
\def\m{\mathcal{M}}

\def\p{\mathcal{P}}

\def\a{\mathcal{A}}
\def\h{\mathcal{H}}
\def\c{\mathcal{C}}

\def\e{\mathcal{E}}

\def\s{\mathcal{S}}

\def\u{\mathcal{U}}
\def\w{\mathcal{W}}
\def\i{\mathcal{I}}
\def\j{\mathcal{J}}
\def\b{\mathcal{B}}

\def\lan{\langle}

\def\ran{\rangle}

 \def\dim{{\hbox{\rm dim}}_{\mathbb C}\,}

\newcommand{\Aut}{\operatorname{Aut}}
\newcommand{\Hom}{\operatorname{Hom}}

\newcommand{\Ind}{\operatorname{Ind}}

\newcommand{\ct} {\operatorname{ct}}
\newcommand{\sq}{\phantom{1}\hfill$\qed$}

\def\HH{\mathfrak{H}}

\def\HH{\hbox{${\mathcal H}$\kern-5.2pt${\mathcal H}$}}

\font\smm=msbm10 at 12pt 
\def\symbol#1{\hbox{\smm #1}}
\def\lsmash{{\symbol n}}

\def\#{\sharp}
